\documentclass[12pt]{article}

\usepackage{indentfirst}
\usepackage{amsfonts}
\usepackage{amsmath}
\usepackage{amssymb}
\usepackage{amsbsy, amsthm}
\usepackage[margin=1in]{geometry}
\usepackage{makeidx,xcolor}
\usepackage{hyperref}
\usepackage[noabbrev,capitalise]{cleveref}

\newtheorem{dfn}{Definition} [section]
\newtheorem{theorem}[dfn]{Theorem}
\newtheorem{lemma}[dfn]{Lemma}

\newtheorem{corollary}[dfn]{Corollary}
\newtheorem{conjecture}[dfn]{Conjecture}
\newenvironment{pf}{\noindent{\bf Proof.}}
{\enspace\vrule height5pt depth0pt width5pt}

\newcommand{\defn}[1]{{\emph{#1}}{\index{#1}}}

\def\X {{\mathcal X}}
\def\tw {{\rm tw}}

\def\F {{\mathcal F}}

\def\P {{\mathcal P}}

\def\Y {{\mathcal Y}}

\makeindex

\begin{document}

\title{Robertson's conjecture I. Well-quasi-ordering bounded tree-width graphs by the topological minor relation}

\maketitle

\centerline{{\bf Chun-Hung Liu}%
\footnote{Email:chliu@tamu.edu. Partially supported by NSF under awards DMS-1929851 and DMS-1954054 and CAREER award DMS-2144042.}
}
\medskip
\centerline{Department of Mathematics}
\centerline{Texas A\&M University}
\centerline{College Station, TX 77843-3368, USA}
\smallskip
\centerline{and}
\smallskip
\centerline{{\bf Robin Thomas}%
}
\medskip
\centerline{School of Mathematics}
\centerline{Georgia Institute of Technology}
\centerline{Atlanta, GA 30332-0160, USA}
\medskip

\begin{abstract}
Robertson and Seymour's celebrated Graph Minor Theorem states that graphs are well-quasi-ordered by the minor relation.
Unlike the minor relation, the topological minor relation does not well-quasi-order graphs in general.
Among all known infinite antichains with respect to the topological containment, subdivisions of a graph obtained from an arbitrarily long path by duplicating each edge can be found.
In the 1980's Robertson conjectured that this is the only obstruction.
Formally, he conjectured that for every positive integer $k$, graphs that do not contain the graph obtained from a path of length $k$ by duplicating each edge as a topological minor are well-quasi-ordered by the topological minor relation.
The case $k=1$ implies Kruskal's Tree Theorem, and the case $k=2$ implies a conjecture of V\'{a}zsonyi on subcubic graphs.

This series of papers is dedicated to a proof of Robertson's conjecture.
We prove Robertson's conjecture for graphs of bounded tree-width in this paper.
It is an essential and groundbreaking step toward the complete proof of Robertson's conjecture, and the machinery developed in this paper will be applied in future papers of the series.
This bounded tree-width case proved in this paper implies all known results about well-quasi-ordering graphs by the topological minor relation that can be proved without using the Graph Minor Theorem, and our proof in this paper is self-contained and does not use any heavily developed machinery.
\end{abstract}




\section{Introduction}

Graphs are finite and are allowed to have parallel edges and loops in this paper.
We say that a graph $G$ contains a graph $H$ as a \defn{minor} if a graph isomorphic to $H$ can be obtained from a subgraph of $G$ by contracting edges.
We say that $G$ contains $H$ as a \defn{topological minor} if a graph isomorphic to $H$ can be obtained from a subgraph of $G$ by repeatedly contracting edges incident with vertices of degree two; equivalently, some subgraph of $G$ is isomorphic to a subdivision of $H$.

A (binary) relation $\preceq$ on a set $S$ is a \defn{quasi-ordering} on $S$ if it is reflexive and transitive.
A quasi-ordering $\preceq$ on a set $S$ is a \defn{well-quasi-ordering} if for every infinite sequence $x_1,x_2,...$ over $S$, there exist $j<j'$ such that $x_j \preceq x_{j'}$.
In this case, we say that $(S,\preceq)$ is a \defn{well-quasi-ordered set}, and $S$ is \defn{well-quasi-ordered by $\preceq$}.
For simplicity, by saying that $Q$ is a quasi-order (or a well-quasi-order, respectively), we mean that $Q$ is a set equipped with a quasi-ordering (or a well-quasi-ordering, respectively) $\leq_Q$.

The history of well-quasi-ordering graphs can be traced back to a conjecture of V\'{a}zsonyi in the 1940's (see \cite{k2,k}): forests are well-quasi-ordered by the topological minor relation.
This conjecture was proved by Kruskal \cite{k} and independently by Tarkowski \cite{t_wqo_tree}; now it is usually called Kruskal's Tree Theorem.
Nash-Williams \cite{n-w_finite_tree} offered an elegant and simpler proof of this theorem a few years later.
Mader \cite{m} further generalized their results to graphs with a bounded number of disjoint cycles.

Kruskal's Tree Theorem and its labelled version are important in logic, computer science and many branches of mathematics.
For example, they are tools for showing that certain orderings have no infinite decreasing chains, and those orderings are used to prove the termination of systems of rewrite rules and the correctness of Knuth-Bendix completion procedures \cite{kb} (see \cite{g}); Higman's Lemma (Theorem \ref{Higman's lemma}) is a special case of the labelled version of Kruskal's Tree Theorem and is recently used to prove that certain categories are quasi-Gr\"{o}bner and hence the categories of their representations over a left-noetherian ring are noetherian \cite{ss}.

One of the most prominent and deepest results in graph theory is the Graph Minor Theorem \cite{rs XX}: graphs are well-quasi-ordered by the minor relation.
It was conjectured by Wagner \cite{w} and proved by Robertson and Seymour \cite{rs XX} in the 1980's.
The proof is extremely difficult and consists of around 20 papers in the Graph Minors series.
On the other hand, the tools developed in the Graph Minors series have had significant impacts in structural graph theory.
Moreover, since the topological minor relation is the same as the minor relation on subcubic graphs\footnote{A graph is \defn{subcubic} if every vertex has degree at most three. The \defn{degree} of a vertex is the number of edges incident with it, where every loop is counted twice.}, the Graph Minor Theorem confirms another conjecture of V\'{a}zsonyi: subcubic graphs are well-quasi-ordered by the topological minor relation.
Furthermore, in the (currently) last paper of the Graph Minor series, Robertson and Seymour \cite{rs XXIII} confirmed a conjecture of Nash-Williams \cite{n-w_wqo_weak}: graphs are well-quasi-ordered by the weak immersion relation.\footnote{We omit the formal definition of weak immersion as we will not discuss this relation in the rest of this paper. But we remark that if a graph $G$ contains another graph $H$ as a topological minor, then $G$ also contains $H$ as a minor and a weak immersion; the minor relation and the weak immersion relation are incomparable. See \cite{l_survey} for a survey about well-quasi-ordering on graphs with respect to various graph containments.}

Unlike the relations of minor and weak immersion, the topological minor relation does not well-quasi-order graphs in general.
For any positive integer $i$, define $G_i$ to be the graph obtained from a path of length $i$ by duplicating each edge and attaching two leaves to each end of the path.
Then it is easy to see that there exist no distinct positive integers $j$ and $j'$ such that $G_j$ contains $G_{j'}$ as a topological minor.
In fact, there are several different sequences of graphs, each containing no pair of graphs comparable by the topological minor relation.
But each such sequence contains graphs obtained from arbitrarily long paths by duplicating each edge as topological minors.
Robertson in the late 1980's conjectured that this known obstruction is the only one.

For every positive integer $k$, the \defn{Robertson chain} of length $k$ is the graph obtained from a path of length $k$ by duplicating each edge exactly once.

\begin{conjecture}[Robertson's conjecture (see \cite{d})\footnote{This conjecture is mentioned in \cite{d}. But this conjecture has been circulated in the community since the late 1980's or earlier.}] \label{Robertson conj}
For every positive integer $k$, graphs that do not contain the Robertson chain of length $k$ as a topological minor are well-quasi-ordered by the topological minor relation.
\end{conjecture}

Conjecture \ref{Robertson conj} is expected to be difficult since the case $k=2$ of Conjecture \ref{Robertson conj} implies the aforementioned conjecture of V\'{a}zsonyi about subcubic graphs, which is unknown how to be proved without using the Graph Minor Theorem. 

Conjecture \ref{Robertson conj} is much more general than known results about well-quasi-ordering graphs by the topological minor relation.
One can easily construct graphs with no long Robertson chain but far from being forests and subcubic graphs.
For example, consider a graph $G$ that is obtained from a disjoint union of arbitrarily many graphs $H_1,H_2,...,H_t$ for some integer $t$ by adding new edges between different $H_i$'s such that for each $i$, there are at most three new edges incident with vertices in $H_i$.
It is easy to see that if each $H_i$ does not contain the Robertson chain of length $r$ as a topological minor, then $G$ does not contain the Robertson chain of length $2r$ as a topological minor.
In particular, if each $H_i$ consists of one vertex, then $G$ can be an arbitrary subcubic graph; if each $H_i$ is a tree, then $r=1$ and $G$ can contain arbitrarily many disjoint cycles; if each $H_i$ is a subcubic graph, then $r=2$ and $G$ can contain arbitrarily many vertices of degree at least four.

One might notice that Conjecture \ref{Robertson conj} is not optimal since Robertson chains do not form an infinite antichain with respect to the topological minor relation.
But it becomes optimal if vertices are labelled.
Before we formally state the labelled version, we mention an equivalent definition of topological minors.

Let $H$ and $G$ be graphs.
We say that $\eta$ is a \defn{homeomorphic embedding} from $H$ to $G$ if $\eta$ is a pair of functions $(\pi_V,\pi_E)$ such that 
	\begin{itemize}
		\item $\pi_V: V(H) \rightarrow V(G)$ is an injection, and 
		\item $\pi_E$ maps each non-loop edge $xy$ of $H$ to a path in $G$ with ends $\pi_V(x)$ and $\pi_V(y)$ and maps each loop of $H$ with end $x$ to a cycle in $G$ passing through $\pi_V(x)$ such that 
			\begin{itemize}
				\item for any two different edges $e_1, e_2$ of $H$, $\pi_E(e_1) \cap \pi_E(e_2) \subseteq \bigcup_{v \in e_1 \cap e_2}\pi_V(v)$, and 
				\item for any $v \in V(H)$ and $e \in E(H)$, $\pi_V(v) \in V(\pi_E(e))$ only if $v$ is incident with $e$. 
			\end{itemize}
	\end{itemize}
In this case, we write $\eta: H \hookrightarrow G$.
Furthermore, we define $\eta(x)=\pi_V(x)$ if $x \in V(H)$, and define $\eta(x)=\pi_E(x)$ if $x \in E(H)$.
Note that $G$ contains $H$ as a topological minor if and only if a homeomorphic embedding from $H$ to $G$ exists.

The following is the labelled version of Robertson's conjecture.

\begin{conjecture} \label{Robertson conj labelled}
Let $k$ be a positive integer.
Let $G_1,G_2,...$ be graphs that do not contain a Robertson chain of length at least $k$ as a topological minor.
Let $Q$ be a set, and for each positive integer $i$, let $\phi_i: V(G_i) \rightarrow Q$ be a function.
If $Q$ is a well-quasi-order, then there exist integers $j,j'$ with $1 \leq j < j'$ and $\eta:G_j \hookrightarrow G_{j'}$ such that $\phi_j(v) \leq_Q \phi_{j'}(\eta(v))$ for every $v \in V(G_j)$.
\end{conjecture}

Clearly, Conjecture \ref{Robertson conj} is the special case of Conjecture \ref{Robertson conj labelled} when $Q$ consists of one element or $Q$ is trivial\footnote{We say that a quasi-order $Q$ is \defn{trivial} if for any elements $x,y$ of $Q$, $x \leq_Q y$ and $y \leq_Q x$.}.
And the case $k=1$ of Conjecture \ref{Robertson conj labelled} implies the labelled version of Kruskal's Tree Theorem.
In addition, Conjecture \ref{Robertson conj labelled} is optimal, since when $Q$ is non-trivial, it contains elements $x,y$ with $x \not \leq_Q y$, and subdivisions of Robertson chains of different lengths form an infinite antichain once we label the ends\footnote{An \defn{end} of a Robertson chain is an end of the original path.} by $x$ and label all other vertices by $y$.

The main objective of this series of papers is to prove Conjecture \ref{Robertson conj labelled} and provide a characterization of the topological minor ideals that are well-quasi-ordered by the topological minor relation\footnote{A statement of this characterization can be found in a survey paper of the first author \cite{l_survey}. This characterization is a strengthening of Conjecture \ref{Robertson conj} and can be viewed as the optimal form of the unlabelled version of Robertson's conjecture (Conjecture \ref{Robertson conj}).}.
We remark that a proof of Conjecture \ref{Robertson conj labelled} was provided in the PhD dissertation of the first author \cite{l_thesis}, and this series of papers is based on that.

This paper is the first paper in the series.
In this paper, we prove Conjecture \ref{Robertson conj labelled} for graphs of bounded tree-width.
(The formal definition of tree-width will be provided in Section \ref{subsec: notation}.)
This partial result of Conjecture \ref{Robertson conj labelled} is an essential step toward the complete proof of this conjecture, and the machinery developed in this paper will be crucial for other papers of the series.

The following is the main result of this paper.

\begin{theorem} \label{robertson conj bounded tree width}
Let $k$ and $w$ be positive integers.
Let $G_1,G_2,...$ be graphs of tree-width at most $w$ not containing a Robertson chain of length at least $k$ as a topological minor.
Let $Q$ be a set, and for each positive integer $i$, let $\phi_i: V(G_i) \rightarrow Q$ be a function.
If $Q$ is a well-quasi-order, then there exist integers $j,j'$ with $1 \leq j < j'$ and $\eta:G_j \hookrightarrow G_{j'}$ such that $\phi_j(v) \leq_Q \phi_{j'}(\eta(v))$ for every $v \in V(G_j)$.
\end{theorem}

The proof of Theorem \ref{robertson conj bounded tree width} is self-contained.
The only known results in the literature that we use as a black box in this paper are Higman's Lemma and \cite[Theorem 2.1]{rs IV}, where each of them has a short proof.

By the famous Grid Minor Theorem \cite{rs V}, for every planar\footnote{A graph is \defn{planar} if it can be embedded in the plane with no edge-crossing.} graph $H$, there exists an integer $w$ such that every graph that does not contain $H$ as a minor has tree-width at most $w$.
Hence we obtain the following corollary of Theorem \ref{robertson conj bounded tree width}.

\begin{corollary} \label{robertson conj excluding planar graph}
For every positive integer $k$, every planar graph $H$ and every well-quasi-order $Q$, if $G_1,G_2,...$ are graphs that do not contain $H$ as a minor and do not contain a Robertson chain of length at least $k$ as a topological minor, and for each positive integer $i$, $\phi_i: V(G_i) \rightarrow Q$ is a function, then there exist integers $j,j'$ with $1 \leq j <j'$ and $\eta:G_j \hookrightarrow G_{j'}$ such that $\phi_j(v) \leq_Q \phi_{j'}(\eta(v))$ for every $v \in V(G_j)$. 
\end{corollary}

We remark that Theorem \ref{robertson conj bounded tree width} (or Corollary \ref{robertson conj excluding planar graph}) implies all known results about well-quasi-ordering graphs by the topological minor relation that can be proved without applying the Graph Minor Theorem.
By a classical result of Erd\H{o}s and P\'{o}sa \cite{ep}, every graph that has a bounded number of pairwise disjoint cycles can be modified into a forest by deleting a bounded number of vertices, so such graphs have bounded tree-width and do not contain a long Robertson chain as a topological minor.
Therefore, Theorem \ref{robertson conj bounded tree width} implies Mader's theorem and hence the aforementioned Kruskal's Tree Theorem. 
In addition, the only known progress (as far as we are aware) on Conjecture \ref{Robertson conj} prior to the announcement of a complete proof of Conjecture \ref{Robertson conj} in the thesis of the first author \cite{l_thesis} is due to Ding \cite{d}: for every $k$, graphs that do not contain a Robertson chain of length $k$ as a minor\footnote{As we pointed out earlier, if $G$ contains $H$ as a topological minor, then $G$ contains $H$ as a minor.} are well-quasi-ordered by the topological minor relation. 
Ding's theorem is an immediate corollary of Corollary \ref{robertson conj excluding planar graph}, since every Robertson chain is planar.

The first author \cite{l_survey} showed that Conjecture \ref{Robertson conj labelled} implies a classical result of Ding \cite{d_subgraph} about subgraphs: for any positive integer $k$, graphs with no path of length $k$ are well-quasi-ordered by the subgraph relation.
We remark that the argument in \cite{l_survey} only involves graphs of bounded tree-width, so Theorem \ref{robertson conj bounded tree width} implies this Ding's result about subgraphs.

\subsection{Proof sketch and organization}

We provide a brief outline of the proof of Theorem \ref{robertson conj bounded tree width} in this subsection.
The proof involves many steps, and we include more detailed sketches for each step in later sections.

We shall use a minimal bad sequence argument to prove Theorem \ref{robertson conj bounded tree width}.
To make it work, we need a ``nice'' tree-decomposition of the graphs.

The first ``nice'' property is a ``linkedness property''.
Roughly speaking, we require that the subgraph induced by the bags in the subtree rooted at a node contains the subgraph induced by the bags in a subtree rooted at ``any" descendant of the previous node as a ``rooted'' topological minor.
The key idea to obtain this linkedness property is to convert vertex-cuts realized by the bags of a tree-decomposition into ``pseudo-edge-cuts."
Though we are not able to convert all vertex-cuts into pseudo-edge-cuts, it suffices to covert some of them such that the tree-decomposition has bounded ``elevation.''

In Section \ref{sec: construct Robertson chain}, we prove a sufficient condition that ensures the existence of a long Robertson chain topological minor. 
Due to the lack of those topological minors, this sufficient condition must be violated.
This allows us to prove the existence of pseudo-edge-cuts.
The detailed arguments and the formal definitions of pseudo-edge-cuts are included in Section \ref{sec: edge-cuts}.

Though those pseudo-edge-cuts exist, they are not necessarily realized by bags of the tree-decomposition.
Our strategy is to ``insert'' those pseudo-edge-cuts into the tree-decomposition to make them realized by bags.
One difficulty of this strategy is that inserting pseudo-edge-cuts into a tree-decomposition might covert some pseudo-edge-cuts that were realized by bags into vertex-cuts, so the insertion process possibly does not terminate.
We overcome this difficulty in Section \ref{sec: bounding elevation}.
The key is that we are manageable to repeatedly select a pseudo-edge-cut and insert it and possibly other its ``related'' cuts into the tree-decomposition such that we can ensure that those pseudo-edge-cuts that we have inserted will not be inserted again even though they are no longer realized by bags.
So we will not insert the same pseudo-edge-cut twice.
This allows us to prove the existence of a tree-decomposition of bounded width and bounded ``elevation.''
The formal definition of elevation is included in Section \ref{sec: bounding elevation}.

Section \ref{sec: tree lemma} is a preparation for proving well-quasi-ordering. 
We show that ``nicely decorated'' trees are well-quasi-ordered with respect to their ``decoration.''
The ``decoration'' can be thought as an encoding of the linkedness property.

The next goal is to achieve the ``absorption property'', which is the other ``nice property'' of a tree-decomposition we want.
Roughly speaking, the absorption property allows us to ``encode'' the subgraphs induced by the bags in the subtree rooted at the children of a given node into the bag of this given node. 
As the tree-decomposition have bounded width, each bag has bounded size, so the bags are easily well-quasi-ordered.
The main difficulty is to find and prove the correct notion of ``encoding'' that allows us to recover topological minors from the encoding and is compatible with the setting of the aforementioned linkedness property. 
We are manageable to do so.
In Section \ref{sec: new wqo bounded depth}, we formally define the notion of encoding and show that these encodings ``simulate'' the ``rooted'' topological minor relation.

In Section \ref{sec: well-behaved}, we combine the main results in previous sections to prove Theorem \ref{robertson conj bounded tree width}.
We show how to transform the linkedness property of a tree-decomposition into the ``decoration'' of trees and use the tree lemma proved in Section \ref{sec: tree lemma} to prove that graphs that admit a tree-decomposition of these two nice properties are well-quasi-ordered by the topological minor relation.
Then Theorem \ref{robertson conj bounded tree width} follows.

We remark that the remainder of the paper can be divided into parts. 
The first part consists of Sections \ref{sec: construct Robertson chain}, \ref{sec: edge-cuts} and \ref{sec: bounding elevation}.
Readers can skip Sections \ref{sec: construct Robertson chain}, \ref{sec: edge-cuts} and \ref{sec: bounding elevation} as long as they understand the statements of Lemma \ref{N-linked disjoint paths} and Theorem \ref{bounded depth} and notion used in them.
Section \ref{sec: tree lemma} forms its own part and it does not rely on any result or notion in previous sections.
Readers can skip Section \ref{sec: tree lemma} as long as they understand the statement of Theorem \ref{decorated tree lemma} and notion used in this theorem.
The rest of sections form the last part.
Section \ref{sec: new wqo bounded depth} does not rely on any result or notion in previous sections; Section \ref{sec: well-behaved} will use Lemma \ref{N-linked disjoint paths}, Theorems \ref{bounded depth} and \ref{decorated tree lemma}, and results in Section \ref{sec: new wqo bounded depth}.

We introduce common terminologies in Section \ref{subsec: notation}.
Intuition of other terminologies will be included in later sections when they are about be to formally defined.

\section{Simple terminologies} \label{subsec: notation}

\subsection{Basic notations}
Let $G$ be a graph, and let $S$ be a subset of $V(G)$.
We denote the subgraph of $G$ induced by $S$ by \defn{$G[S]$}.
We define \defn{$G-S$} to be $G[V(G)-S]$.
If $v$ is a vertex of $G$, then we define \defn{$G-v$} to be $G-\{v\}$.

Let $f$ be a function, and let $S$ be a subset of its domain.
We denote the restriction of $f$ on $S$ by \defn{$f|_S$}.
We also define \defn{$f(S)$} to be $\{f(x):x \in S\}$.
Similarly, for a sequence $\sigma=(x_1,x_2,...,x_n)$ whose entries are in the domain of $f$, we define \defn{$f(\sigma)$} to be $(f(x_1),f(x_2),...,f(x_n))$.

For a nonnegative integer $r$, we define \defn{$[r]$} to be the set $\{n \in {\mathbb N}: 1 \leq n \leq r\}$ and define \defn{$[0,r]$} to be $[r] \cup \{0\}$.

\subsection{Tree-decomposition}

We say that $(T,\X)$ is a \defn{tree-decomposition} of a graph $G$ if the following hold.
	\begin{itemize}
		\item $T$ is a tree, and $\X=(X_t: t \in V(T))$, where for every $t \in V(T)$, $X_t$ is a subset of $V(G)$ and called the \defn{bag} at $t$.
		\item $\bigcup_{t \in V(T)} X_t = V(G)$.
		\item For every edge of $G$, some bag contains all the ends of this edge.
		\item For every vertex $u \in V(G)$, the nodes of $T$ whose bags contain $u$ induce a connected subgraph of $T$.
	\end{itemize}
The \defn{adhesion} of $(T,\X)$ is $\max\{\lvert X_x \cap X_y \rvert: xy \in E(T)\}$.
The \defn{width} of $(T,\X)$ is $\max\{\lvert X_t \rvert: t \in V(T)\}-1$.
The \defn{tree-width} of a graph $G$, denoted by \defn{$\tw(G)$}, is the minimum width of a tree-decomposition of $G$.

In this paper, if $(T,\X)$ is a tree-decomposition, then we always denote the bag at $t$ by \defn{$X_t$}, and we call each vertex of $T$ a \defn{node} in order to distinguish vertices of $T$ and vertices in their bags.

\section{Looking for Robertson chains} \label{sec: construct Robertson chain}

In this section, we will prove that certain structure in a tree-decomposition of a graph can be used to construct a long Robertson chain topological minor.
It is a step toward our main structure theorem for graphs of bounded tree-width that do not contain a long Robertson chain topological minor (Theorem \ref{bounded depth}).

Roughly speaking, we will show how to construct a long Robertson chain topological minor when $G$ has a tree-decomposition that has many disjoint bags with the same size such that there are many disjoint paths in $G$ passing through all vertices in those bags and there are many extra paths that connect those paths in a certain way.
The precise form of the previous statement require certain terminologies that will be defined in later sections.
In particular, we will formally define those ``extra paths'' in Sections \ref{subsec:path_blocks} and \ref{subsec:lrjumps}.
In fact, we will prove a more general statement that allows non-disjoint bags, which leads to the notion called ``weak strips'' defined in Section \ref{subsec:weak_strip}.
Weak strips are prototypes of ``strips'', a key notion that will be defined and used in later sections to prove the main structure theorem (Theorem \ref{bounded depth}).

\subsection{Path of blocks} \label{subsec:path_blocks}

Let $G$ be a graph.
A \defn{cut-vertex} of $G$ is a vertex $v$ of $G$ such that $G-v$ has more components than $G$.
A \defn{block} of $G$ is a maximal subgraph of $G$ that does not contain any cut-vertex. 
A \defn{block tree} of $G$ is a tree $T$ with bipartition $\{A,B\}$ and a bijection $f$ between $V(T)$ and the union of the set of blocks of $G$ and the set of the cut-vertices of $G$ such that 
	\begin{itemize}
		\item $f(a)$ is a cut-vertex of $G$ for every $a \in A$,
		\item $f(b)$ is a block of $G$ for every $b \in B$, and
		\item if $a \in A$ and $b \in B$ with $ab \in E(T)$, then $f(a)$ is a cut-vertex of $G$ contained in $V(f(b))$.
	\end{itemize}
For blocks $B_1,B_2$ of $G$, we say that $H$ is a \defn{graph that is the path of blocks of $G$ from $B_1$ to $B_2$} if $H=\bigcup_t f(t)$, where the union is over all nodes $t \in B$ contained in the path in $T$ from $b_1$ to $b_2$, where $b_1,b_2$ are the nodes of $T$ with $f(b_1)=B_1$ and $f(b_2)=B_2$.

\subsection{Left jumps and right jumps} \label{subsec:lrjumps}

Let $G$ be a graph and $(T,\X)$ a tree-decomposition of $G$.
Let $k$ be a positive integer, and let $t_1,t_2$ be two nodes of $T$ with $\lvert X_{t_1} \rvert = \lvert X_{t_2} \rvert=k$.
Let $\P$ be a collection of $k$ disjoint paths in $G$ from $X_{t_1}$ to $X_{t_2}$.
Let $P$ be a member of $\P$ with $\lvert E(P) \rvert \neq \emptyset$.
For $i \in \{1,2\}$, let $v_i$ be the vertex in $X_{t_i} \cap V(P)$, and let $B_i$ be the block of $G-\bigcup_{W \in \P-\{P\}}V(W)$ containing the edge of $P$ incident with $v_{i}$.
Let $Q_{\P,P}$ be the graph that is the path of blocks of $G-\bigcup_{W \in \P-\{P\}}V(W)$ from $B_1$ to $B_2$.
If none of the blocks in $Q_{\P,P}$ is an edge, then we define $L_{\P,P}=R_{\P,P}=Q_{\P,P}$.
If some block in $Q_{\P,P}$ is an edge, then let $Q'$ be the union of the blocks in $Q_{\P,P}$ that are single edges, and define $L_{\P,P}$ (and $R_{\P,P}$, respectively) to be the component of $Q_{\P,P}-E(Q')$ containing $v_{1}$ (and $v_{2}$, respectively).
Note that $L_{\P,P}$ (and $R_{\P,P}$, respectively) consists of a vertex if $B_1$ (and $B_2$, respectively) is a single edge.
Let $G'$ be the subgraph of $G$ induced by $X_{t_1} \cup X_{t_2} \cup \bigcup_t X_t$, where the union is over all nodes $t$ in the component of $T-\{t_1,t_2\}$ containing an internal node of the $t_1$-$t_2$ path in $T$.
A \defn{right jump from $v_{1}$} (and \defn{left jump from $v_{2}$}, respectively) is a path in $G'$ from $V(L_{\P,P})$ (and $V(R_{\P,P})$, respectively) to $\bigcup_{W \in \P-\{P\}}V(W)$ internally disjoint from $V(L_{\P,P}) \cup \bigcup_{W \in \P} V(W)$ (and $V(R_{\P,P}) \cup \bigcup_{W \in \P}V(W)$, respectively). We call the graphs $Q_{\P,P}$, $L_{\P,P}$ and $R_{\P,P}$ the \defn{$(Q,\P,P)$-graph, $(L,\P,P)$-graph} and \defn{$(R,\P,P)$-graph between $X_{t_1}$ and $X_{t_2}$}, respectively.
See Figure \ref{fig_LRQjump} for an example.

\begin{figure} 
	\begin{picture}(100,230) (-5,-30)

		\thicklines
	
		\put(215,200){\circle*{5}} \put(205,200){{$z$}}
		
		\put(50,30){\circle*{5}} \put(45,20){{$v_1$}}
		\put(90,50){\circle*{5}} 
		\put(90,10){\circle*{5}}
		\put(50,30){\line(2,1){40}}
		\put(50,30){\line(2,-1){40}}
		\put(130,30){\circle*{5}}
		\put(130,30){\line(-2,1){40}}
		\put(130,30){\line(-2,-1){40}}
		\put(160,50){\circle*{5}}
		\put(160,10){\circle*{5}} \put(158,0){{$y_2$}}
		\put(130,30){\line(3,2){30}}
		\put(130,30){\line(3,-2){30}}
		\put(160,50){\line(3,0){30}}
		\put(160,50){\line(0,-4){40}}
		\put(160,50){\line(3,-4){30}}
		\put(160,10){\line(3,0){30}}
		\put(175,30){\circle*{5}} 
		\put(190,50){\circle*{5}}
		\put(190,10){\circle*{5}}
		\put(220,30){\circle*{5}}
		\put(190,50){\line(3,-2){30}}
		\put(190,10){\line(3,2){30}}

		\multiput(30,55)(5,0){40}{\line(1,0){3}}
		\multiput(230,55)(0,-5){14}{\line(0,-1){3}}
		\multiput(30,-15)(5,0){40}{\line(1,0){3}}
		\multiput(30,-15)(0,5){14}{\line(0,1){3}}
		\put(115,-10){{$L_{\P,P}$}}

		\put(240,30){\circle*{5}}
		\put(280,30){\circle*{5}}
		\put(220,30){\line(1,0){180}}
		\put(260,50){\circle*{5}}
		\put(260,50){\line(-1,-1){20}}
		\put(260,50){\line(1,-1){20}}

		\put(400,30){\circle*{5}} \put(395,20){{$v_2$}}
		\put(370,0){\circle*{5}}
		\put(340,0){\circle*{5}}
		\put(355,60){\circle*{5}} \put(350,68){{$y_3$}}
		\put(310,30){\circle*{5}}
		\put(310,30){\line(3,2){45}}
		\put(310,30){\line(1,0){90}}
		\put(310,30){\line(1,-1){30}}
		\put(310,30){\line(2,-1){60}}
		\put(340,0){\line(1,4){15}}
		\put(340,0){\line(2,1){60}}
		\put(340,0){\line(1,0){30}}
		\put(370,0){\line(-1,4){15}}
		\put(370,0){\line(1,1){30}}
		\put(400,30){\line(-3,2){45}}

		\multiput(300,75)(5,0){25}{\line(1,0){3}}
		\multiput(425,75)(0,-5){18}{\line(0,-1){3}}
		\multiput(300,-15)(5,0){25}{\line(1,0){3}}
		\multiput(300,-15)(0,5){18}{\line(0,1){3}}
		\put(305,-10){{$R_{\P,P}$}}

		\multiput(10,80)(5,0){88}{\line(1,0){3}}
		\multiput(450,80)(0,-5){22}{\line(0,-1){3}}
		\multiput(10,-30)(5,0){88}{\line(1,0){3}}
		\multiput(10,-30)(0,5){22}{\line(0,1){3}}
		\put(250,-25){{$Q_{\P,P}$}}

		\multiput(50,130)(50,0){8}{\circle*{5}}
		\put(50,130){\line(1,0){350}}
		\put(170,140){{$P_2$}}

		\put(160,10){\line(-1,2){60}}
		\multiput(100,130)(5,-10){5}{\circle*{5}}
		\put(130,90){{$Y_2$}}

		\put(70,110){\circle*{5}}
		\put(70,90){\circle*{5}}
		\put(50,30){\line(1,3){20}}
		\put(70,90){\line(0,1){20}}
		\put(70,110){\line(-1,1){20}}
		\put(75,100){{$Y_1$}}

		\put(215,200){\line(1,-1){140}}
		\multiput(215,200)(30,-30){4}{\circle*{5}}
		\put(290,100){{$Y_3$}}

		\put(315,150){\circle*{5}}
		\put(215,200){\line(2,-1){100}}
		\put(370,150){\circle*{5}}
		\put(370,150){\line(-1,0){55}}
		\put(370,150){\line(1,-4){30}}
		\put(310,160){{$Y_4$}}

		\multiput(40,10)(0,5){40}{\line(0,1){3}}
		\multiput(40,210)(5,0){38}{\line(1,0){3}}
		\multiput(230,210)(0,-5){4}{\line(0,-1){3}}
		\multiput(230,190)(-5,0){34}{\line(-1,0){3}}
		\multiput(60,190)(0,-5){36}{\line(0,-1){3}}
		\multiput(40,10)(5,0){4}{\line(1,0){3}}	
		\put(15,180){{$X_{t_1}$}}

		\multiput(390,10)(0,5){35}{\line(0,1){3}}
		\multiput(390,185)(-5,0){38}{\line(-1,0){3}}
		\multiput(200,185)(0,5){8}{\line(0,1){3}}
		\multiput(200,225)(5,0){42}{\line(1,0){3}}
		\multiput(410,225)(0,-5){43}{\line(0,-1){3}}
		\multiput(390,10)(5,0){4}{\line(1,0){3}}	
		\put(415,180){{$X_{t_2}$}}
		
	\end{picture}
	\caption{In this example, $\P$ consists of three disjoint paths between $X_{t_1}$ and $X_{t_2}$: the first path consists of the single vertex $z$, the second path is $P_2$, and the third path is a path $P$ between $v_1$ and $v_2$ contained in $Q_{\P,P}$. $Y_1$ and $Y_2$ are right jumps from $v_1$, where $Y_1$ is between $v_1$ and $V(P_2) \cap X_{t_1}$, and $Y_2$ is between $y_2$ and $V(P_2)-X_{t_1}$. $Y_3$ and $Y_4$ are left jumps from $v_2$, where $Y_3$ is between $y_3$ and $z$, and $Y_4$ is between $v_2$ and $z$.} \label{fig_LRQjump}
\end{figure}

\subsection{Weak strips} \label{subsec:weak_strip}

Let $(T,\X)$ be a tree-decomposition of a graph $G$.
For a subset $Z$ of $V(G)$ and a positive integer $s$, a \defn{weak $(Z,s)$-strip} in $(T,\X)$ is a sequence $(t_1,t_2,...,t_r)$ for some positive integer $r$ such that
	\begin{itemize}
		\item $t_1,t_2,...,t_r$ are distinct nodes of $T$ such that there exists a path in $T$ passing through them in the order listed,
		\item $Z \subseteq X_{t_i}$ and $\lvert X_{t_i}-Z \rvert = s$ for every $1 \leq i \leq r$,
		\item $X_{t_1}-Z,X_{t_2}-Z,...,X_{t_h}-Z$ are pairwise disjoint nonempty sets, and
		\item there exists a set of $\lvert X_{t_1} \rvert$ disjoint paths in $G$ from $X_{t_1}$ to $X_{t_r}$.
	\end{itemize}
See Figure \ref{fig_weak_strip} for an example.

\begin{figure} 
	\begin{picture}(100,230) (-5,-10)

		\thicklines
	
		\multiput(80,200)(60,0){6}{\circle*{5}}
		\put(80,200){\line(1,0){300}}
		\put(76,210){{$u_1$}}
		\put(136,210){{$u_2$}}
		\put(196,210){{$u_3$}}
		\put(256,210){{$u_4$}}
		\put(316,210){{$u_5$}}
		\put(376,210){{$u_6$}}
		\put(100,200){\circle*{5}}
		\put(120,200){\circle*{5}}
		\put(170,200){\circle*{5}}
		\put(220,200){\circle*{5}}
		\put(240,200){\circle*{5}}
		\put(275,200){\circle*{5}}
		\put(290,200){\circle*{5}}
		\put(305,200){\circle*{5}}
		\put(335,200){\circle*{5}}
		\put(350,200){\circle*{5}}
		\put(365,200){\circle*{5}}

		\multiput(80,150)(60,0){6}{\circle*{5}}
		\put(80,150){\line(1,0){300}}
		\put(76,160){{$v_1$}}
		\put(136,160){{$v_2$}}
		\put(196,160){{$v_3$}}
		\put(256,160){{$v_4$}}
		\put(316,160){{$v_5$}}
		\put(376,160){{$v_6$}}
		\put(110,150){\circle*{5}}
		\put(170,150){\circle*{5}}
		\put(215,150){\circle*{5}}
		\put(230,150){\circle*{5}}
		\put(245,150){\circle*{5}}
		\put(280,150){\circle*{5}}
		\put(300,150){\circle*{5}}
		\put(340,150){\circle*{5}}
		\put(360,150){\circle*{5}}

		\multiput(230,0)(0,50){3}{\circle*{5}}
		\put(215,98){{$z_1$}}
		\put(215,48){{$z_2$}}
		\put(215,-2){{$z_3$}}

	\end{picture}
	\caption{In this example, $(t_1,t_2,...,t_6)$ is a weak $(Z,2)$-strip, $Z=\{z_1,z_2,z_3\}$ and $X_{t_i} = \{u_i,v_i,z_1,z_2,z_3\}$ for every $i \in [6]$.} \label{fig_weak_strip}
\end{figure}

\begin{lemma} \label{jumps make Robertson chain}
Let $r,k$ be positive integers.
Let $(T,\X)$ be a tree-decomposition of a graph $G$.
Let $Z$ be a subset of $V(G)$.
Let $(t_1,t_2,...,t_r)$ be a weak $(Z,|X_{t_1}-Z|)$-strip in $(T,\X)$. 
Let $\P$ be a set of $|X_{t_1}|$ disjoint paths in $G$ from $X_{t_1}$ to $X_{t_r}$.
Let $P$ be a member of $\P$ disjoint from $Z$.
For every $i \in [r]$, let $v_{i}$ be the vertex in $X_{t_i} \cap V(P)$.
Assume that for every $i \in [r-1]$, either there exist two edge-disjoint paths in $G-Z$ from $v_{i}$ to $v_{i+1}$ internally disjoint from $X_{t_i} \cup X_{t_{i+1}}$, or there exist a right jump from $v_{i}$ disjoint from $Z$ and a left jump from $v_{i+1}$ disjoint from $Z$.
Assume that for every $i \in [r-1]-[1]$, if both the right jump and the left jump from $v_i$ mentioned above exist, then these two jumps can be chosen such that they intersect in at most one vertex.
If $r \geq k(k+1)\lvert X_{t_1} \rvert^{2k+2}+k+3$, then $G-Z$ contains the Robertson chain of length $k$ as a topological minor.
\end{lemma}

\begin{pf}
Without loss of generality, we may assume that $Z$ is the empty set; otherwise we delete $Z$ from $G$.
For every $i \in [r-1]$, define $Q_{i}, L_{i},R_{i}$ to be the $(Q,\P,P)$-graph, $(L,\P,P)$-graph, $(R,\P,P)$-graph between $X_{t_i}$ and $X_{t_{i+1}}$, respectively.
Denote $\P$ by $\{P_1,P_2,...,P_{\lvert X_{t_1} \rvert}\}$, and let $s$ be the index such that $P=P_s$.
For every $i \in [r-1]$, define the \defn{type} of $i$ to be $(a,b)$, where $1 \leq a \leq \lvert X_{t_1} \rvert$ and $1 \leq b \leq \lvert X_{t_1} \rvert$ are some integers such that the following hold.
	\begin{itemize}
		\item If there exist two edge-disjoint paths from $v_{i}$ to $v_{i+1}$ internally disjoint from $X_{t_{i}} \cup X_{t_{i+1}}$ in $G-Z$, then $(a,b)=(s,s)$.
		\item Otherwise, there exist a right jump $J_{R,i}$ from $v_{i}$ whose end not in $V(L_{i})$ is in $V(P_a)$, and a left jump $J_{L,i+1}$ from $v_{i+1}$ whose end not in $V(R_{i})$ is in $V(P_b)$. 
	\end{itemize}
By the assumption of this lemma, we may assume that $J_{R,i}$ and $J_{L,i}$ intersect in at most one vertex when both of them are defined.
Note that it is possible that there are more than one jump satisfying the property mentioned above, and in this case, we choose arbitrary such a jump to define the type of $i$.
Furthermore, $V(J_{R,i}) \cap V(J_{L,i+1}) \subseteq \bigcup_{j \neq s}V(P_j)$ and $\lvert V(J_{R,i}) \cap V(J_{L,i+1}) \rvert \leq 1$, since $L_{i} \neq R_{i}$ when $J_{R,i}$ and $J_{L,i+1}$ are defined.
Note that if some entry of the type of $i$ is $s$, then the corresponding two edge-disjoint paths exist since $L_{i} \neq R_{i}$, and hence both entries of the type of $i$ are $s$. 

For every $2 \leq i \leq r-1-k$, we define the \defn{$k$-type} of $i$ to be the sequence $(b_i,b_{i+1},...,b_{i+k})$, where $b_j$ is the type of $j$ for every $i \leq j \leq i+k$.
Observe that there are at most $\lvert X_{t_1} \rvert^{2k+2}$ possible $k$-types.
Since $r-k-2 \geq k(k+1)\lvert X_{t_1} \rvert^{2k+2}+1$, there exist $2 \leq i_0 < i_1< ... < i_k \leq r-k-1$ such that $i_j \equiv i_0$ (mod $k+1$) for $1 \leq j \leq k$, and the $k$-types of $i_0,i_1,...,i_k$ are the same.

We shall recursively construct a homeomorphic embedding from the Robertson chain of length $k$ to $G$ by proving the following statement.

\begin{itemize}
	\item[(*)] There exists an increasing sequence $(s_0,s_1,...,s_k)$ such that for every $m \in [0,k]$, $s_m = i_m'+m$ for some $i_m' \in \{i_0,i_1,..., i_m\}$ (so $s_m-i_0 \equiv m$ (mod $k+1$)), and $G$ contains a subgraph $S_m$ satisfying the following properties. 
	\begin{itemize}
		\item[(i)] For every $v \in V(S_m) \cap X_{t_{s_m}}$, there exists a homeomorphic embedding $\eta$ from a Robertson chain of length at least $m$ to $S_m$ such that $v=\eta(x)$ for some end $x$ of the Robertson chain.
		\item[(ii)] $S_m$ is contained in the subgraph of $G$ induced by $X_{t_{s_m}} \cup \bigcup_{t \in V(T')}X_t$, where $T'$ is the component of $T-\{t_{s_m}\}$ containing $t_1$.
		\item[(iii)] $\lvert V(S_m) \cap X_{t_{s_m}} \rvert \in \{1,2\}$.
		\item[(iv)] If $\lvert V(S_m) \cap X_{t_{s_m}} \rvert=2$, then 
			\begin{itemize}
				\item $J_{L,s_m}$ is defined, 
				\item $V(S_m) \cap X_{t_{s_m}}$ consists of $v_{s_m}$ and an end of $J_{L,s_m}$, and 
				\item if $v_{s_m}$ is not an end of $J_{L,s_m}$, then there exist a subgraph $S_m'$ of $G$ with $V(S_m') \cap X_{t_{s_m}}=V(S_m) \cap X_{t_{s_m}} - \{v_{s_m}\}$ and a homeomorphic embedding $\eta'$ from a Robertson chain of length at least $m$ to $S_m'$ such that $\eta'$ maps some end of the Robertson chain to the vertex in $V(S_m) \cap X_{t_{s_m}} - \{v_{s_m}\}$.
			\end{itemize}
		\item[(v)] If $\lvert V(S_m) \cap X_{t_{s_m}} \rvert=1$, then $V(S_m) \cap X_{t_{s_m}} = \{v_{s_m}\}$.
	\end{itemize}
\end{itemize}

Since this lemma follows from (i), it suffices to prove the statement (*).

We prove the statement (*) by induction on $m$.
When $m=0$, the statement is obviously true by choosing $s_0=i_0$ and choosing $S_0$ to be the graph consists of the single vertex $v_{i_0}$.

Now we assume that $m>0$ and $(s_0,s_1,...,s_{m-1})$ exists. 

\medskip

\noindent{\bf Claim 1:} We may assume that $s_{m-1}$ is not of type $(s,s)$.

\medskip

\noindent{\bf Proof of Claim 1:}
Assume the contrary that $s_{m-1}$ is of type $(s,s)$.
By (iv) and (v), $S_{m-1}$ contains $v_{s_{m-1}}$.
We set $s_m=s_{m-1}+1$ and $i'_m=i'_{m-1}$, and we define $S_m$ to be the graph obtained from $S_{m-1}$ by adding two edge-disjoint paths from $v_{s_{m-1}}$ to $v_{s_{m-1}+1}=v_{s_m}$ internally disjoint from $X_{t_{s_{m-1}}} \cup X_{t_{s_{m-1}+1}}$.
It is clear that $S_m$ satisfies (i)-(v).
$\Box$

\medskip

By Claim 1, we may assume that $s_{m-1}$ is of type $(a,b)$, where $a \neq s \neq b$.

We let $i'_m = i_{q+1}$, where $q$ is the number in $[0,m-1]$ such that $i_{m-1}'=i_q$, and let $s_m = i'_m+m$.
So $s_m-1=i'_m+m-1 \geq (i'_{m-1}+k+1)+m-1>s_{m-1}+1$.
In particular, $s_m>s_{m-1}$.

Since $s_m-1-i_m' \equiv s_{m-1}-i_{m-1}'$ (mod $k+1$), the type of $s_m-1$ is the same as the type of $s_{m-1}$, which equals $(a,b)$.
Since $a \neq s \neq b$, $J_{R,s_{m-1}}, J_{L,s_{m-1}+1}, J_{R,s_m-1}, J_{L,s_m}$ are defined.
Let $u_{s_{m-1}}, u_{s_{m-1}+1},u_{s_m-1},u_{s_m}$ be the end of $J_{R,s_{m-1}}, J_{L,s_{m-1}+1}, J_{R,s_m-1}, J_{L,s_m}$ contained in $L_{s_{m-1}}, R_{s_{m-1}}, L_{s_m-1},R_{s_m-1}$, respectively.
Since every block in $L_{s_{m-1}} \cup R_{s_{m-1}} \cup L_{s_m-1} \cup R_{s_m-1}$ is not a single edge, there exist connected subgraphs $Y_{m-1}, Y'_{m-1}, Y_m,Y'_m$ of $G$ disjoint from $\bigcup_{j \neq s}V(P_j)$ such that 
	\begin{itemize}
		\item $Y_{m-1}$ has an Eulerian trail\footnote{A \defn{Eulerian trial} in a graph $H$ is a walk that visits every edge of $H$ exactly once.} from $u_{s_{m-1}}$ to $u_{s_{m-1}+1}$ containing $v_{s_{m-1}}$ such that
			\begin{itemize}
				\item if $u_{s_{m-1}}$ and $v_{s_{m-1}}$ are contained in the same block in $L_{s_{m-1}}$, then $Y_{m-1}$ is a path,
				\item otherwise, $Y_{m-1}$ is the union of a path from $u_{s_{m-1}}$ to $u_{s_{m-1}+1}$ and two edge-disjoint paths from $v_{s_{m-1}}$ to a cut-vertex of $L_{s_{m-1}}$ contained in a block of $L_{s_{m-1}}$ containing $u_{s_{m-1}}$,
			\end{itemize}
		\item $Y'_{m-1}$ has an Eulerian trail from $u_{s_{m-1}}$ to $v_{s_{m-1}+1}$ containing $v_{s_{m-1}}$ such that
			\begin{itemize}
				\item if $u_{s_{m-1}}$ and $v_{s_{m-1}}$ are contained in the same block in $L_{s_{m-1}}$, then $Y'_{m-1}$ is a path,
				\item otherwise, $Y'_{m-1}$ is the union of a path from $u_{s_{m-1}}$ to $v_{s_{m-1}+1}$ and two edge-disjoint paths from $v_{s_{m-1}}$ to a cut-vertex of $L_{s_{m-1}}$ contained in a block of $L_{s_{m-1}}$ containing $u_{s_{m-1}}$,
			\end{itemize}
		\item $Y_m$ has an Eulerian trail from $u_{s_m-1}$ to $u_{s_m}$ containing $v_{s_m}$ such that
			\begin{itemize}
				\item if $u_{s_{m}}$ and $v_{s_{m}}$ are contained in the same block in $R_{s_{m}-1}$, then $Y_m$ is a path,
				\item otherwise, $Y_m$ is the union of a path from $u_{s_m-1}$ to $u_{s_m}$ and two edge-disjoint paths in $R_{s_m-1}$ from $v_{s_m}$ to a cut-vertex of $R_{s_m-1}$ contained in a block of $R_{s_m-1}$ containing $u_{s_m}$,
			\end{itemize}
		\item $Y'_m$ has an Eulerian a trail from $v_{s_m-1}$ to $u_{s_m}$ containing $v_{s_m}$ such that
			\begin{itemize}
				\item if $u_{s_{m}}$ and $v_{s_{m}}$ are contained in the same block in $R_{s_{m}-1}$, then $Y'_m$ is a path,
				\item otherwise, $Y'_m$ is the union of a path from $v_{s_m-1}$ to $u_{s_m}$ and two edge-disjoint paths in $R_{s_m-1}$ from $v_{s_m}$ to a cut-vertex of $R_{s_m-1}$ contained in a block of $R_{s_m-1}$ containing $u_{s_m}$,
			\end{itemize}
	\end{itemize}
Moreover, we define paths $Z_{m-1},Z'_{m-1},Z_m,Z_m'$ in $G$ internally disjoint from $\bigcup_{j \neq s}V(P_j)$ such that
	\begin{itemize}
		\item if $v_{s_{m-1}}$ is not an end of $J_{R,s_{m-1}}$, then $Z_{m-1}$ (and $Z'_{m-1}$, respectively) is a path from $u_{s_{m-1}}$ to $u_{s_{m-1}+1}$ (and to $v_{s_{m-1}+1}$, respectively) not containing $v_{s_{m-1}}$,  and 
		\item if $v_{s_m}$ is not an end of $J_{L,s_m}$, then $Z_{m}$ (and $Z_m'$, respectively) is a path from $u_{s_{m}-1}$ (and from $v_{s_m-1}$, respectively) to $u_{s_{m}}$ not containing $v_{s_{m}}$. 
	\end{itemize}
Finally, we define the paths $W_R,W_L,W,P'$ in $G$ such that
	\begin{itemize}
		\item $W_R$ is the subpath of $P_a$ connecting one end of $J_{R,s_{m-1}}$ and one end of $J_{R,s_m-1}$, and $W_L$ is the subpath of $P_b$ connecting one end of $J_{L,s_{m-1}+1}$ and one end of $J_{L,s_m}$, 
		\item if $a=b$, then $W$ is the subpath of $P_a$ connecting one end of $J_{R,s_{m-1}}$ and one end of $J_{L,s_m}$, and
		\item $P'$ is the subpath of $P$ connecting $v_{s_{m-1}+1}$ and $v_{s_m-1}$.
	\end{itemize}

Now, we are ready to construct $S_m$ (and $S_m'$ if required).

\medskip

\noindent{\bf Claim 2:} We may assume that $\lvert V(S_{m-1}) \cap X_{t_{s_{m-1}}} \rvert=2$ and $V(J_{R,s_{m-1}}) \cap X_{t_{s_{m-1}}} \cap V(S_{m-1})-\{v_{s_{m-1}}\} \neq \emptyset$.

\medskip

\noindent{\bf Proof of Claim 2:}
We assume the contrary that $\lvert V(S_{m-1}) \cap X_{t_{s_{m-1}}} \rvert=1$ or $V(J_{R,s_{m-1}}) \cap X_{t_{s_{m-1}}} \cap V(S_{m-1})-\{v_{s_{m-1}}\} = \emptyset$.
By (iv) and (v), $v_{s_{m-1}} \in V(S_{m-1}) \cap X_{t_{s_{m-1}}}$.
If $a=b$, then define $S_m$ to be the graph obtained from $S_{m-1}$ by adding $J_{R,s_{m-1}} \cup Y'_{m-1} \cup P' \cup Y'_m \cup J_{L,s_m} \cup W$.
If $a \neq b$, then define $S_m$ to be the graph obtained from $S_{m-1}$ by adding $J_{R,s_{m-1}} \cup Y_{m-1} \cup J_{L,s_{m-1}+1} \cup W_L \cup J_{L,s_m} \cup Y_m \cup J_{R,s_m-1} \cup W_R$.
Note that $V(S_m) \cap X_{t_{s_m}}$ contains $v_{s_m}$, and if $V(S_m) \cap X_{t_{s_m}} - \{v_{s_m}\} \neq \emptyset$, then one end of $J_{L,s_m}$ is in $X_{t_{s_m}}-\{v_{s_m}\}$, and $V(S_m) \cap X_{t_{s_m}}-\{v_{s_m}\}$ consists of this vertex.
Hence $S_m$ satisfies (i)-(iii) and (v), and to prove (iv), it suffices to show that $S_m'$ exists when $\lvert V(S_m) \cap X_{t_{s_m}} \rvert =2$ and $v_{s_m}$ is not an end of $J_{L,s_m}$.
In this case, if $a=b$, then we define $S_m'$ to be the graph obtained from $S_{m-1}$ by adding $J_{R,s_{m-1}} \cup Y'_{m-1} \cup P' \cup Z'_m \cup J_{L,s_m} \cup W$; if $a \neq b$, then we define $S'_m$ to be the graph obtained from $S_{m-1}$ by adding $J_{R,s_{m-1}} \cup Y_{m-1} \cup J_{L,s_{m-1}+1} \cup W_L \cup J_{L,s_m} \cup Z_m \cup J_{R,s_m-1} \cup W_R$.
Then $S_m$ and $S_m'$ satisfy (iv).
$\Box$

\medskip

By Claim 2, we may assume that $\lvert V(S_{m-1}) \cap X_{t_{s_{m-1}}} \rvert=2$ and $V(J_{R,s_{m-1}}) \cap X_{t_{s_{m-1}}} \cap V(S_{m-1})-\{v_{s_{m-1}}\} \neq \emptyset$.

By (iv), $V(S_{m-1}) \cap X_{t_{s_{m-1}}}$ consists of $v_{s_{m-1}}$ and an end of $J_{L,s_{m-1}}$.
So the vertex in $V(S_{m-1}) \cap X_{t_{s_{m-1}}} - \{v_{s_{m-1}}\}$ is a common vertex of $J_{L,s_{m-1}}$ and $J_{R,s_{m-1}}$.
Hence $v_{s_{m-1}}$ is not a common vertex of $J_{L,s_{m-1}}$ and $J_{R,s_{m-1}}$ by our assumption.

\medskip

\noindent{\bf Claim 3:} We may assume that $v_{s_{m-1}}$ is not an end of $J_{L,s_{m-1}}$.

\medskip

\noindent{\bf Proof of Claim 3:}
We assume the contrary that $v_{s_{m-1}}$ is an end of $J_{L,s_{m-1}}$.
Hence the ends of $J_{L,s_{m-1}}$ are the two vertices in $V(S_{m-1}) \cap X_{t_{s_{m-1}}}$, and $u_{s_{m-1}} \neq v_{s_{m-1}}$.
If $a=b$, then define $S_m$ to be the graph obtained from $S_{m-1}$ by adding $J_{R,s_{m-1}} \cup Z'_{m-1} \cup P' \cup Y'_m \cup J_{L,s_m} \cup W$.
If $a \neq b$, then define $S_m$ to be the graph obtained from $S_{m-1}$ by adding $J_{R,s_{m-1}} \cup Z_{m-1} \cup J_{L,s_{m-1}+1} \cup W_L \cup J_{L,s_m} \cup Y_m \cup J_{R,s_m-1} \cup W_R$.
Note that $V(S_m) \cap X_{t_{s_m}}$ contains $v_{s_m}$, and if $V(S_m) \cap X_{t_{s_m}} - \{v_{s_m}\} \neq \emptyset$, then one end of $J_{L,s_m}$ is in $X_{t_{s_m}}-\{v_{s_m}\}$, and $V(S_m) \cap X_{t_{s_m}}-\{v_{s_m}\}$ consists of this vertex.
Hence $S_m$ satisfies (i)-(iii) and (v), and to prove (iv), it suffices to show that $S_m'$ exists when $\lvert V(S_m) \cap X_{t_{s_m}} \rvert =2$ and $v_{s_m}$ is not an end of $J_{L,s_m}$.
In this case, if $a=b$, then we define $S_m'$ to be the graph obtained from $S_{m-1}$ by adding $J_{R,s_{m-1}} \cup Z'_{m-1} \cup P' \cup Z'_m \cup J_{L,s_m} \cup W$; if $a \neq b$, then we define $S'_m$ to be the graph obtained from $S_{m-1}$ by adding $J_{R,s_{m-1}} \cup Z_{m-1} \cup J_{L,s_{m-1}+1} \cup W_L \cup J_{L,s_m} \cup Z_m \cup J_{R,s_m-1} \cup W_R$.
Then $S_m$ and $S_m'$ satisfy (iv).
$\Box$

\medskip

So we may assume that $v_{s_{m-1}}$ is not an end of $J_{L,s_{m-1}}$ by Claim 3.
Hence $S_{m-1}'$ is defined by (iv).
If $a=b$, then define $S_m$ to be the graph obtained from $S'_{m-1}$ by adding $J_{R,s_{m-1}} \cup Y'_{m-1} \cup P' \cup Y'_m \cup J_{L,s_m} \cup W$.
If $a \neq b$, then define $S_m$ to be the graph obtained from $S'_{m-1}$ by adding $J_{R,s_{m-1}} \cup Y_{m-1} \cup J_{L,s_{m-1}+1} \cup W_L \cup J_{L,s_m} \cup Y_m \cup J_{R,s_m-1} \cup W_R$.
Hence $S_m$ satisfies (i)-(iii) and (v), and to prove (iv), it suffices to show that $S_m'$ exists when $\lvert V(S_m) \cap X_{t_{s_m}} \rvert =2$ and $v_{s_m}$ is not an end of $J_{L,s_m}$.
In this case, if $a=b$, then we define $S_m'$ to be the graph obtained from $S'_{m-1}$ by adding $J_{R,s_{m-1}} \cup Y'_{m-1} \cup P' \cup Z'_m \cup J_{L,s_m} \cup W$; if $a \neq b$, then we define $S'_m$ to be the graph obtained from $S'_{m-1}$ by adding $J_{R,s_{m-1}} \cup Y_{m-1} \cup J_{L,s_{m-1}+1} \cup W_L \cup J_{L,s_m} \cup Z_m \cup J_{R,s_m-1} \cup W_R$.
Then $S_m$ and $S_m'$ satisfy (iv).

This proves the statement (*) and hence this lemma. 
\end{pf}

\section{Looking for edge-cuts}
\label{sec: edge-cuts}

The goal of this section is to prove the existence of ``pseudo-edge-cuts'', which is a step toward our main structure theorem for graphs with no long Robertson chain topological minor (Theorem \ref{bounded depth}) that will be used in proving well-quasi-ordering results in later sections.
Roughly speaking, a pseudo-edge-cut is a separation $(A,B)$ such that some vertices in $A \cap B$ are incident with most one edge whose other end is in $A-B$; so deleting those edges and other vertices in $A \cap B$ can separate the graph, and the number of deleted edges and vertices is at most $|A \cap B|$.
The key idea for constructing a pseudo-edge-cut is based on the lack of certain kind of jumps due to the absence of a long Robertson chain topological minor proved in Section \ref{sec: construct Robertson chain}.

We will define necessary notions in Sections \ref{subsec:rooted_tree_decomposition}-\ref{subsec:pseudo_ec_strips} and then prove the main result of this section in Section \ref{subsec:breaking_strips}.

\subsection{Rooted tree-decomposition} \label{subsec:rooted_tree_decomposition}

A \defn{rooted tree} is a directed graph whose underlying graph is a tree such that all but one node has in-degree one.
The vertex in a rooted tree with in-degree not one is called the \defn{root}.
It is easy to see that the root has in-degree zero.
For every non-root node $v$, the tail $u$ of the edge with head $v$ is the \defn{parent} of $v$, and we say that $v$ is a \defn{child} of $u$ in this case.
If there exists a directed path from a node $x$ to another node $y$, then $x$ is an \defn{ancestor} of $y$, and $y$ is a \defn{descendant} of $x$.
Note that every node is an ancestor and a descendant of itself.
We say that $x$ is a \defn{proper ancestor} (and \defn{proper descendant}, respectively) of $y$ if $x \neq y$ and $x$ is an ancestor (and descendant, respectively).

If $x,y$ are nodes of a rooted tree $T$, then we denote the set of all nodes of the path in the underlying graph of $T$ from $x$ to $y$ by \defn{$xTy$}.
When there is no danger for creating confusion, we also denote the directed path in $T$ between $x$ and $y$ by $xTy$.

We say that a tree-decomposition $(T,\X)$ is a \defn{rooted tree-decomposition} if $T$ is a rooted tree.
In this case, for every node $t$, we define \defn{$(T,\X)\uparrow t$} to be $\bigcup_s X_s$, where the union is taken over all descendants $s$ of $t$; we define \defn{$(T,\X)\downarrow t$} to be $\bigcup_s X_s$, where the union is taken over all nodes $s$ where either $s=t$ or $s$ is a non-descendant of $t$.
When the rooted tree-decomposition is clear in the context, we simply denote $(T,\X) \uparrow t$ and $(T,\X) \downarrow t$ by \defn{$\uparrow t$} and \defn{$\downarrow t$}, respectively.\footnote{It might be helpful to absorb these notations if the readers imagine that the tree $T$ is drawn in the plane such that the root is drawn at the bottom.}
Note that $\downarrow t \cap \uparrow t = X_t$.

Given a rooted tree-decomposition $(T,\X)$, we say that a node $t_1$ of $T$ is a \defn{precursor} of a node $t_2$ of $T$ if $t_1$ is a proper ancestor of $t_2$ with $\lvert X_{t_1} \rvert = \lvert X_{t_2} \rvert$ and $\lvert X_t \rvert \geq \lvert X_{t_1} \rvert$ for all $t \in t_1Tt_2$.

\subsection{Separations}

A \defn{separation} $(A,B)$ of a graph $G$ is an ordered pair of subsets of $V(G)$ with $A \cup B = V(G)$ such that no edge of $G$ has one end in $A-B$ and one end in $B-A$.
The \defn{order} of $(A,B)$ is $\lvert A \cap B \rvert$.
We remark that in the literature, a separation is often an ordered pair of edge-disjoint subgraphs whose union is the entire graph.
But our arguments in this paper about separations only depend on the vertex-sets of the two edge-disjoint subgraphs, so we define a separation to be an order pair of subsets of vertices in this paper.

When $(T,\X)$ is a rooted tree-decomposition and $t$ is a node of $T$, the \defn{separation given by $t$ in $(T,\X)$}, denoted by \defn{$(A_t,B_t)$}, is the separation $(\downarrow t, \uparrow t)$.

We say that a separation $(A,B)$ \defn{separates} two subsets $X,Y$ of $V(G)$ if $X$ and $Y$ are not subsets of $A \cap B$, and either $X \subseteq A$ and $Y \subseteq B$, or $X \subseteq B$ and $Y \subseteq A$.
When $X$ consists of one vertex, we say $(A,B)$ separates $Y$ and the vertex in $X$.

\subsection{Pseudo-edge-cuts and strips} \label{subsec:pseudo_ec_strips}

Let $(A,B)$ be a separation of a graph $G$.
We say that a vertex $v \in A \cap B$ is \defn{pointed for $(A,B)$} if $v$ is incident with at most one edge of $G$ whose other end is in $A-B$.
For each subset $Z$ of $V(G)$, we say that $(A,B)$ is a \defn{pseudo-edge-cut modulo $Z$} if every vertex in $A \cap B-Z$ is pointed for $(A,B)$.

For any positive integer $s$ and subset $Z$ of $V(G)$, a \defn{$(Z,s)$-strip} in a rooted tree-decomposition $(T,\X)$ of $G$ is a sequence of nodes $(t_1,t_2,...,t_h)$ of $T$ such that the following hold.
	\begin{itemize}
		\item $t_i$ is a precursor of $t_{i+1}$ for every $1 \leq i \leq h-1$.
		\item $(t_1,t_2,...,t_h)$ is a weak $(Z,s)$-strip.
		\item There exists no node $t$ in $t_1Tt_h$ such that $\lvert X_t \rvert = \lvert X_{t_1} \rvert$ and the separation $(A_t,B_t)$ given by $t$ in $(T,\X)$ is a pseudo-edge-cut modulo $Z$. 
	\end{itemize}
The \defn{length} of this $(Z,s)$-strip is $h$.
The \defn{$(Z,s)$-depth} of $(T,\X)$ is the maximum length of a $(Z,s)$-strip in $(T,\X)$.

Let $\alpha$ be a positive integer. 
We say that a separation $(A,B)$ of $G$ \defn{$\alpha$-breaks} a $(Z,s)$-strip $(t_1,...,t_h)$ in $(T,\X)$ if there exist distinct indices $1 \leq i_1<i_2<...<i_{\alpha}<j_1<j_2<...<j_{\alpha} \leq h$ such that $\downarrow t_{i_\alpha} \subseteq A$ and $\uparrow t_{j_1} \subseteq B$. 

\subsection{Breaking strips} \label{subsec:breaking_strips}

Let $G$ be a graph.
Let $(T,\X)$ be a rooted tree-decomposition of $G$.
For every node $t$ of $T$ that has a precursor, and for every precursor $t'$ of $t$, we say that a set $\P$ of $\lvert X_t \rvert$ disjoint paths $P_1,...,P_{\lvert X_t \rvert}$ in $G$ from $X_{t'}$ to $X_t$ of size $\lvert X_t \rvert$ is a \defn{set of foundation paths between $t'$ and $t$}.
When a set of foundation paths $\P$ between $t'$ and $t$ and a member $P$ of $\P$ with $V(P) \cap X_t \neq V(P) \cap X_{t'}$ are given, we define $L_{\P,P}, R_{\P,P},Q_{\P,P}$ to be the $(L,\P,P)$-graph, $(R,\P,P)$-graph, $(Q,\P,P)$-graph between $X_{t'}$ and $X_t$, respectively.
(Recall the definitions in Section \ref{subsec:lrjumps}.)
A \defn{parent-side jump at $t$ with respect to $t'$ at the vertex in $X_t \cap V(P)$} is a path in $G[\uparrow t' \cap \downarrow t]$ from $V(R_{\P,P})$ to $\bigcup_{W \in \P-\{P\}} V(W) - (X_{t'} \cap X_t)$ internally disjoint from $V(R_{\P,P}) \cup \bigcup_{W \in \P}V(W)$.
A \defn{child-side jump at $t'$ with respect to $t$ at the vertex in $X_{t'} \cap V(P)$} is a path in $G[\uparrow t' \cap \downarrow t]$ from $V(L_{\P,P})$ to $\bigcup_{W \in \P-\{P\}} V(W) - (X_{t'} \cap X_t)$ internally disjoint from $V(L_{\P,P}) \cup \bigcup_{W \in \P}V(W)$.
A parent-side (or child-side) jump at $t$ is \defn{ambiguous} if both its ends are in $X_t$; otherwise it is \defn{unambiguous}.
See Figure \ref{fig_pcjump} for an example.

\begin{figure} 
	\begin{picture}(100,230) (-5,-30)

		\thicklines
	
		\put(215,200){\circle*{5}} \put(205,200){{$z$}}
		
		\put(50,30){\circle*{5}} \put(45,20){{$v_1$}}
		\put(90,50){\circle*{5}} 
		\put(90,10){\circle*{5}}
		\put(50,30){\line(2,1){40}}
		\put(50,30){\line(2,-1){40}}
		\put(130,30){\circle*{5}}
		\put(130,30){\line(-2,1){40}}
		\put(130,30){\line(-2,-1){40}}
		\put(160,50){\circle*{5}}
		\put(160,10){\circle*{5}} \put(158,0){{$y_2$}}
		\put(130,30){\line(3,2){30}}
		\put(130,30){\line(3,-2){30}}
		\put(160,50){\line(3,0){30}}
		\put(160,50){\line(0,-4){40}}
		\put(160,50){\line(3,-4){30}}
		\put(160,10){\line(3,0){30}}
		\put(175,30){\circle*{5}} 
		\put(190,50){\circle*{5}}
		\put(190,10){\circle*{5}}
		\put(220,30){\circle*{5}}
		\put(190,50){\line(3,-2){30}}
		\put(190,10){\line(3,2){30}}

		\multiput(30,55)(5,0){40}{\line(1,0){3}}
		\multiput(230,55)(0,-5){14}{\line(0,-1){3}}
		\multiput(30,-15)(5,0){40}{\line(1,0){3}}
		\multiput(30,-15)(0,5){14}{\line(0,1){3}}
		\put(115,-10){{$L_{\P,P}$}}

		\put(240,30){\circle*{5}}
		\put(280,30){\circle*{5}}
		\put(220,30){\line(1,0){180}}
		\put(260,50){\circle*{5}}
		\put(260,50){\line(-1,-1){20}}
		\put(260,50){\line(1,-1){20}}

		\put(400,30){\circle*{5}} \put(395,20){{$v_2$}}
		\put(370,0){\circle*{5}}
		\put(340,0){\circle*{5}}
		\put(355,60){\circle*{5}} \put(350,68){{$y_4$}}
		\put(310,30){\circle*{5}}
		\put(310,30){\line(3,2){45}}
		\put(310,30){\line(1,0){90}}
		\put(310,30){\line(1,-1){30}}
		\put(310,30){\line(2,-1){60}}
		\put(340,0){\line(1,4){15}}
		\put(340,0){\line(2,1){60}}
		\put(340,0){\line(1,0){30}}
		\put(370,0){\line(-1,4){15}}
		\put(370,0){\line(1,1){30}}
		\put(400,30){\line(-3,2){45}}

		\multiput(300,75)(5,0){25}{\line(1,0){3}}
		\multiput(425,75)(0,-5){18}{\line(0,-1){3}}
		\multiput(300,-15)(5,0){25}{\line(1,0){3}}
		\multiput(300,-15)(0,5){18}{\line(0,1){3}}
		\put(305,-10){{$R_{\P,P}$}}

		\multiput(10,80)(5,0){88}{\line(1,0){3}}
		\multiput(450,80)(0,-5){22}{\line(0,-1){3}}
		\multiput(10,-30)(5,0){88}{\line(1,0){3}}
		\multiput(10,-30)(0,5){22}{\line(0,1){3}}
		\put(250,-25){{$Q_{\P,P}$}}

		\multiput(50,130)(50,0){8}{\circle*{5}}
		\put(50,130){\line(1,0){350}}
		\put(170,140){{$P_2$}}

		\put(160,10){\line(-1,2){60}}
		\multiput(100,130)(5,-10){5}{\circle*{5}}
		\put(130,90){{$Y_2$}}

		\put(70,110){\circle*{5}}
		\put(70,90){\circle*{5}}
		\put(50,30){\line(1,3){20}}
		\put(70,90){\line(0,1){20}}
		\put(70,110){\line(-1,1){20}}
		\put(75,100){{$Y_1$}}

		\put(215,115){\circle*{5}}
		\put(130,30){\line(1,1){85}}
		\put(215,115){\line(0,1){85}}
		\put(180,95){{$Y_3$}}
		\put(125,20){{$y_3$}}

		\put(215,200){\line(1,-1){140}}
		\multiput(215,200)(30,-30){4}{\circle*{5}}
		\put(290,100){{$Y_4$}}

		\put(370,90){\circle*{5}}
		\put(370,90){\line(1,-2){30}}
		\put(370,90){\line(3,4){30}}
		\put(360,100){{$Y_5$}}

		\multiput(40,10)(0,5){40}{\line(0,1){3}}
		\multiput(40,210)(5,0){38}{\line(1,0){3}}
		\multiput(230,210)(0,-5){4}{\line(0,-1){3}}
		\multiput(230,190)(-5,0){34}{\line(-1,0){3}}
		\multiput(60,190)(0,-5){36}{\line(0,-1){3}}
		\multiput(40,10)(5,0){4}{\line(1,0){3}}	
		\put(15,180){{$X_{t'}$}}

		\multiput(390,10)(0,5){35}{\line(0,1){3}}
		\multiput(390,185)(-5,0){38}{\line(-1,0){3}}
		\multiput(200,185)(0,5){8}{\line(0,1){3}}
		\multiput(200,225)(5,0){42}{\line(1,0){3}}
		\multiput(410,225)(0,-5){43}{\line(0,-1){3}}
		\multiput(390,10)(5,0){4}{\line(1,0){3}}	
		\put(415,180){{$X_{t}$}}
		
	\end{picture}
	\caption{In this example, $t'$ is a precursor of $t$ with $|X_{t'}|=|X_t|=3$, and $\P$ is a set of foundation paths between $t'$ and $t$ consisting of $|X_t|=3$ disjoint paths: the first path consists of the single vertex $z$, the second path is $P_2$, and the third path is a path $P$ between $v_1 \in X_{t'}$ and $v_2 \in X_t$ contained in $Q_{\P,P}$. And $Y_1$ is a path between $v_1$ and $V(P_2) \cap X_{t'}$ internally disjoint from $V(L_{\P,P}) \cup \{z\} \cup V(P_2) \cup V(P)$, and $Y_2$ is a path between $y_2$ and $V(P_2)-X_{t'}$ internally disjoint from $V(L_{\P,P}) \cup \{z\} \cup V(P_2) \cup V(P)$. So $Y_1$ and $Y_2$ are child-side jumps at $t'$ with respect to $t$ at $v_1$, where $Y_1$ is ambiguous and $Y_2$ is unambiguous. Moreover, $Y_3$ is a path between $y_3$ and $z \in X_{t'} \cap X_{t}$ internally disjoint from $V(L_{\P,P}) \cup \{z\} \cup V(P_2) \cup V(P)$. Since $z \in X_{t'} \cap X_{t}$, $Y_3$ is not a child-side jump at $t'$ with respect to $t$ from $v_1$. And $Y_4$ is a path between $y_4$ and $z \in X_t \cap X_{t'}$ internally disjoint from $V(R_{\P,P}) \cup \{z\} \cup V(P_2) \cup V(P)$. Since $z \in X_t \cap X_{t'}$, $Y_4$ is not a parent-side jump at $t$ with respect to $t'$. In addition, $Y_5$ is a path between $v_2$ and a vertex in $V(P_2) \cap X_t$ internally disjoint from $V(R_{\P,P}) \cup \{z\} \cup V(P_2) \cup V(P)$. So $Y_5$ is an ambiguous parent-side jump at $t$ with respect to $t'$ from $v_2$.} \label{fig_pcjump}
\end{figure}

\subsubsection{Limiting jumps}

\begin{lemma} \label{breaking_weak_0}
For any positive integers $k,w,\alpha$, there exists an integer $f = f(k,w,\alpha)$ such that the following holds.
Let $G$ be a graph that does not contain the Robertson chain of length $k$ as a topological minor.
Let $(T,\X)$ be a rooted tree-decomposition of $G$ of width at most $w$.
Let $Z \subseteq V(G)$, $s$ be a positive integer, and $(t_1,t_2,...,t_{f+1})$ be a $(Z,s)$-strip in $(T,\X)$.
Let $\P = \{P_1,P_2,...,P_{\lvert Z \rvert+s}\}$ be a set of foundation paths between $t_1$ and $t_{f+1}$, where $P_i$ is a one-vertex path with $V(P_i) \subseteq Z$ for every $s+1 \leq i \leq s+\lvert Z \rvert$.
Then there exists a subsequence $(t_1',t_2',...,t'_{\alpha+1})$ of $(t_1,t_2,...,t_{f+1})$ such that for every $j' \in [s]$, the following statements hold. 
	\begin{enumerate}
		\item For every $\ell \in [\alpha]$, there exist no two edge-disjoint paths in $G-\bigcup_{W \in \P-\{P_{j'}\}}V(W)$ from the vertex in $X_{t'_{\ell}} \cap V(P_{j'})$ to the vertex in $X_{t'_{\ell+1}} \cap V(P_{j'})$.
		\item Either
			\begin{itemize}
				\item for every $j \in [\alpha+1]-[1]$, there exists no parent-side jump at $t'_{j}$ with respect to $t'_{j-1}$ at the vertex in $X_{t'_{j}} \cap V(P_{j'})$ disjoint from $Z$, or 
				\item for every $j \in [\alpha]$, there exists no unambiguous child-side jump at $t'_{j}$ with respect to $t'_{j+1}$ at the vertex in $X_{t'_{j}} \cap V(P_{j'})$ disjoint from $Z$.
			\end{itemize}
	\end{enumerate}
\end{lemma}

\begin{pf}
Define $f=2^{w+2}(\alpha+1) (k(k+1)(w+1)^{2k+2}+k+3)(w+1)k$.
We shall show that $f$ satisfies the conclusion of the lemma.

Let $G$, $(T,\X)$, $Z$, $s$, $(t_1,t_2,...,t_{f+1})$ and $\P$ be the ones mentioned in the lemma.
Note that subpaths of the members of $\P$ form a set of foundation paths between $t_i$ and $t_{i+1}$ for each $i \in [f]$.
In addition, $\lvert Z \rvert+s =\lvert X_{t_1} \rvert \leq w+1$ since the width of $(T,\X)$ is at most $w$.

If there exist $i \in [f-k]$ and $P \in \{P_1,P_2,...,P_s\}$ such that for each $j \in [0,k-1]$, there exist two edge-disjoint paths in $G-\bigcup_{W \in \P-\{P\}}V(W)$ from the vertex in $X_{t_{i+j}} \cap V(P)$ to the vertex in $X_{t_{i+j+1}} \cap V(P)$ internally disjoint from $X_{t_{i+j}} \cup X_{t_{i+j+1}}$, then there exists a homeomorphic embedding from a Robertson chain of length at least $k$ to $G$, a contradiction.
So for each $P \in \{P_1,P_2,...,P_s\}$ and $i \in [f-k]$, there exists $j \in [0,k-1]$ such that there do not exist two edge-disjoint paths in $G-\bigcup_{W \in \P-\{P\}}V(W)$ from the vertex in $X_{t_{i+j}} \cap V(P)$ to the vertex in $X_{t_{i+j+1}} \cap V(P)$ internally disjoint from $X_{t_{i+j}} \cup X_{t_{i+j+1}}$.

Therefore, for each $P \in \{P_1,P_2,...,P_s\}$ and $i \in [f-k]$, there do not exist two edge-disjoint paths in $G-\bigcup_{W \in \P-\{P\}}V(W)$ from the vertex in $X_{t_i} \cap V(P)$ to the vertex in $X_{t_{i+k}} \cap V(P)$ internally disjoint from $X_{t_i} \cup X_{t_{i+k}}$.

To simplify the notation, by taking a subsequence of $(t_1,t_2,...t_{f+1})$ of length $f/k$, we can call $t_{ki+1}$ as $t_i$ for each $i \in [f/k]$, and assume that for each $i \in [f/k]$ and $P \in \{P_1,P_2,...,P_s\}$, there do not exist two edge-disjoint paths in $G-\bigcup_{W \in \P-\{P\}}V(W)$ from the vertex in $X_{t_i} \cap V(P)$ to the vertex in $X_{t_{i+1}} \cap V(P)$.
Note that it implies that for any $i,j$ with $1 \leq i <j \leq f/k$, there do not exist two edge-disjoint paths in $G-\bigcup_{W \in \P-\{P\}}V(W)$ from the vertex in $X_{t_i} \cap V(P)$ to the vertex in $X_{t_{j}} \cap V(P)$.

Let $\beta= (k(k+1)\lvert X_{t_1} \rvert^{2k+2}+k+3)s$.
Suppose that there exists a sequence $(i_0,i_1,i_2,...,i_{\beta+1})$ with $i_0=1 < 2 \leq i_1<i_2<...<i_{\beta} \leq f/k-1 <f/k=i_{\beta+1}$ such that for each $j \in [\beta]$, there exists a vertex $v_j $ in $X_{t_{i_j}}-Z$ such that there exist a parent-side jump at $t_{i_j}$ with respect to $t_{i_{j-1}}$ at $v_j$ disjoint from $Z$ and an unambiguous child-side jump at $t_{i_j}$ with respect to $t_{i_{j+1}}$ at $v_j$ disjoint from $Z$.
By the pigeonhole principle, there exists a sequence $(i_0',i_1',...,i_{\beta/s+1}')$ with $i_0'=1<2 \leq i_1'<i_2'<...<i_{\beta/s}'<i'_{\beta/s+1}$ such that there exists a member $P$ of $\{P_1,P_2,...,P_s\}$ such that for each $j \in [\beta/s]$, there exists a parent-side jump at $t_{i_j'}$ with respect to $t_{i'_{j-1}}$ at the vertex in $(X_{t_{i_j'}}-Z) \cap V(P)$ disjoint from $Z$ and an unambiguous child-side jump at $t_{i_j'}$ with respect to $t_{i'_{j+1}}$ at the vertex in $(X_{t_{i_j'}}-Z) \cap V(P)$ disjoint from $Z$.
By Lemma \ref{jumps make Robertson chain}, $G$ contains the Robertson chain of length $k$ as a topological minor, a contradiction.

So there exist at most $\beta-1$ nodes $t_i$ in $\{t_2,...,t_{f/k-1}\}$ such that each $X_{t_{i}}-Z$ contains a vertex $v$ such that there exist a parent-side jump at $t_i$ with respect to $t_{i-1}$ at $v$ disjoint from $Z$ and an unambiguous child-side jump at $t_i$ with respect to $t_{i+1}$ at $v$ disjoint from $Z$.
Let $f_1 = \lceil (f/k-2-\beta+1)/\beta \rceil$.
By the pigeonhole principle, there exists $i$ with $2 \leq i \leq f/k-f_1-1$ such that for each $j \in [f_1-1]$, $X_{t_{i+j}}-Z$ contains no vertex $v$ such that both a parent-side jump at $t_{i+j}$ with respect to $t_{i+j-1}$ at $v$ disjoint from $Z$ and an unambiguous child-side jump at $t_{i+j}$ with respect to $t_{i+j+1}$ at $v$ disjoint from $Z$ exist.

Let $f_2 = \lceil f_1/2^s \rceil$.
By the pigeonhole principle, there exists a subsequence $(t_1',t_2',...,t'_{f_2})$ of $(t_1,t_2,...,t_{f+1})$ such that for each $j' \in [s]$, there exist no two edge-disjoint paths in $G-\bigcup_{W \in \P-\{P_{j'}\}}V(W)$ from the vertex in $X_{t'_{\ell}} \cap V(P_{j'})$ to the vertex in $X_{t'_{\ell+1}} \cap V(P_{j'})$ for all $\ell \leq [f_2-1]$, and either 
	\begin{itemize}
		\item for every $j \in [f_2]-[1]$, there exist no parent-side jumps at $t'_{j}$ with respect to $t'_{j-1}$ at the vertex in $X_{t'_{j}} \cap V(P_{j'})$ disjoint from $Z$, or 
		\item for every $j \in [f_2-1]$, there exist no unambiguous child-side jumps at $t'_{j}$ with respect to $t'_{j+1}$ at the vertex in $X_{t'_{j}} \cap V(P_{j'})$ disjoint from $Z$.
	\end{itemize}
Since $f_2 \geq \alpha+1$, this proves the lemma.
\end{pf}

\bigskip

Now we further strengthen Lemma \ref{breaking_weak_0}.

\begin{lemma} \label{breaking_weak}
For any positive integers $k,w,\alpha$, there exists an integer $f = f(k,w,\alpha)$ such that the following holds.
Let $G$ be a graph that does not contain the Robertson chain of length $k$ as a topological minor.
Let $(T,\X)$ be a rooted tree-decomposition of $G$ of width at most $w$.
Let $Z \subseteq V(G)$, $s$ be a positive integer, and $(t_1,t_2,...,t_{f+1})$ be a $(Z,s)$-strip in $(T,\X)$.
Let $\P = \{P_1,P_2,...,P_{\lvert Z \rvert+s}\}$ be a set of foundation paths between $t_1$ and $t_{f+1}$, where $P_i$ is a one-vertex path with $V(P_i) \subseteq Z$ for every $i \in [s+\lvert Z \rvert]-[s]$. 
Then there exists a subsequence $(t_1',t_2',...,t'_{\alpha+1})$ of $(t_1,t_2,...,t_{f+1})$ such that for every $j' \in [s]$, the following statements hold. 
	\begin{enumerate}
		\item For every $\ell \in [\alpha]$, there exist no two edge-disjoint paths in $G-\bigcup_{W \in \P-\{P_{j'}\}}V(W)$ from the vertex in $X_{t'_{\ell}} \cap V(P_{j'})$ to the vertex in $X_{t'_{\ell+1}} \cap V(P_{j'})$.
		\item Either
			\begin{enumerate}
				\item for every $j \in [\alpha+1]-[1]$, there exists no parent-side jump at $t'_{j}$ with respect to $t'_{j-1}$ at the vertex in $X_{t'_{j}} \cap V(P_{j'})$ disjoint from $Z$, or 
				\item for every $j \in [\alpha]$, 
					\begin{enumerate}
						\item there exists no unambiguous child-side jump at $t'_{j}$ with respect to $t'_{j+1}$ at the vertex in $X_{t'_{j}} \cap V(P_{j'})$ disjoint from $Z$, and 
						\item if $j \geq 2$, then there exists $v \in X_{t'_j}-(Z \cup V(P_{j'}))$ such that $v$ is an end of every parent-side jump at $t'_{j}$ with respect to $t'_{j-1}$ at the vertex in $X_{t'_{j}} \cap V(P_{j'})$ disjoint from $Z$ and every child-side jump at $t'_{j}$ with respect to $t'_{j+1}$ at the vertex in $X_{t'_{j}} \cap V(P_{j'})$ disjoint from $Z$.
					\end{enumerate}
			\end{enumerate}
	\end{enumerate}
\end{lemma}

\begin{pf}
Let $\alpha_0=k(k+1)\lvert X_{t_1} \rvert^{2k+2}+k+3$.
Let $\alpha_1 = ((\alpha_0-1)s+1)(\alpha+1)$.
Define $f = f_0(k,w,\alpha_1)$, where $f_0$ is the number $f$ mentioned in Lemma \ref{breaking_weak_0}.

By Lemma \ref{breaking_weak_0}, there exists a subsequence $(t_1',t_2',...,t'_{\alpha_1+1})$ of $(t_1,t_2,...,t_{f+1})$ such that for every $j' \in [s]$, the following statements hold. 
	\begin{itemize}
		\item For every $\ell \in [\alpha_1]$, there exist no two edge-disjoint paths in $G-\bigcup_{W \in \P-\{P_{j'}\}}V(W)$ from the vertex in $X_{t'_{\ell}} \cap V(P_{j'})$ to the vertex in $X_{t'_{\ell+1}} \cap V(P_{j'})$.
		\item Either
			\begin{itemize}
				\item for every $j \in [\alpha_1+1]-[1]$, there exists no parent-side jump at $t'_{j}$ with respect to $t'_{j-1}$ at the vertex in $X_{t'_{j}} \cap V(P_{j'})$ disjoint from $Z$, or 
				\item for every $j \in [\alpha_1]$, there exists no unambiguous child-side jump at $t'_{j}$ with respect to $t'_{j+1}$ at the vertex in $X_{t'_{j}} \cap V(P_{j'})$ disjoint from $Z$.
			\end{itemize}
	\end{itemize}

For any $j' \in [s]$ and $j \in [\alpha_1]-[1]$, we say that $(j,j')$ is \defn{bad} if 
	\begin{itemize}
		\item there exist a parent-side jump at $t'_{j}$ with respect to $t'_{j-1}$ at the vertex in $X_{t'_{j}} \cap V(P_{j'})$ disjoint from $Z$ and a child-side jump at $t'_{j}$ with respect to $t'_{j+1}$ at the vertex in $X_{t'_{j}} \cap V(P_{j'})$ disjoint from $Z$, and
		\item for every $v \in X_{t'_j}-(Z \cup V(P_{j'}))$, either there exists a parent-side parent-side jump at $t'_{j}$ with respect to $t'_{j-1}$ at the vertex in $X_{t'_{j}} \cap V(P_{j'})$ disjoint from $Z \cup \{v\}$, or there exists a child-side jump at $t'_{j}$ with respect to $t'_{j+1}$ at the vertex in $X_{t'_{j}} \cap V(P_{j'})$ disjoint from $Z \cup \{v\}$.
	\end{itemize}
Note that for any bad $(j,j')$, there exist a parent-side jump at $t'_{j}$ with respect to $t'_{j-1}$ at the vertex in $X_{t'_{j}} \cap V(P_{j'})$ disjoint from $Z$ and a child-side jump at $t'_{j}$ with respect to $t'_{j+1}$ at the vertex in $X_{t'_{j}} \cap V(P_{j'})$ disjoint from $Z$ such that these two jumps only intersect in at most one vertex.
	If there exist $j' \in [s]$ and $\alpha_0$ elements in $[\alpha_1]-[1]$ such that $(j,j')$ is bad, then Lemma \ref{jumps make Robertson chain} implies that $G$ contains the Robertson chain of length $k$ as a topological minor, a contradiction.

Hence for every $j' \in [s]$, there exist at most $\alpha_0-1$ elements in $[\alpha_1]-[1]$ such that $(j,j')$ is bad.
So there are at most $s(\alpha_0-1)$ bad pairs $(j',j)$.
Since $\alpha_1 \geq ((\alpha_0-1)s+1)\alpha+(\alpha_0-1)s$, by the pigeonhole principle, there exist $i^* \in [\alpha_1+1-\alpha]$ such that the subsequence $(t'_{i^*}, t'_{i^*+1},...,t'_{i^*+\alpha-1})$ satisfies the conclusion of this lemma.
\end{pf}

\subsubsection{Parent-side and child-side progress}

Let $G$ be a graph.
Let $(T,\X)$ be a rooted tree-decomposition of $G$.
Let $t$ be a node of $T$, and let $t'$ be a precursor of $t$.
Let $\P$ be a set of foundation paths between $t'$ and $t$.

For a member $P$ of $\P$ with $V(P) \cap X_t \neq V(P) \cap X_{t'}$, we say that $P$ is \defn{parent-side static with respect to $t'$ and $t$} if 
	\begin{itemize}
		\item there exist no two edge-disjoint paths in $G-\bigcup_{W \in \P-\{P\}}V(W)$ from the vertex in $V(P) \cap X_{t'}$ to the vertex in $V(P) \cap X_t$ internally disjoint from $X_t \cup X_{t'}$, and 
		\item there exists no parent-side jump at $t$ with respect to $t'$ at the vertex in $V(P) \cap X_t$.
	\end{itemize}

Assume that $P$ is parent-side static with respect to $t'$ and $t$.
Let $u_P$ be the vertex in $V(P) \cap X_{t'}$ and $v_P$ the vertex in $V(P) \cap X_t$.
Since there exists no parent-side jump at $t$ with respect to $t'$ at $v_P$ and there exist no two edge-disjoint paths in $G-\bigcup_{W \in \P-\{P\}}V(W)$ from $u_P$ to $v_P$ internally disjoint from $X_t \cup X_{t'}$, there exists a unique vertex $w$ in the $(R,\P,P)$-graph $R_{\P,P}$ such that $w$ is contained in a single-edge block with edge-set $\{e\}$ of the $(Q,\P,P)$-graph $Q_{\P,P}$, and we define $R_{\P,P}'$ to be the component of $(G-(X_t-\{v_P\}))-e$ containing $R_{\P,P}$. 
Since $P$ is parent-side static with respect to $t'$ and $t$, there exists a separation $(L_P,M_P)$ of $G$ such that $L_P = \downarrow t-V(R'_{\P,P}-\{w\})$ and $M_P = \uparrow t \cup V(R'_{\P,P})$.
Note that $L_P \cap M_P = (X_t-\{v_P\}) \cup \{w\}$, and $e$ is the unique edge incident with $w$ whose other end is in $L_P-M_P$. 
See Figure \ref{fig_pstatic} for an example.

\begin{figure} 
	\begin{picture}(100,280) (-5,-70)

		\thicklines
	
		\put(215,200){\circle*{5}} \put(205,200){{$z$}}
		
		\put(50,30){\circle*{5}} \put(45,20){{$u_P$}}
		\put(90,50){\circle*{5}} 
		\put(90,10){\circle*{5}}
		\put(50,30){\line(2,1){40}}
		\put(50,30){\line(2,-1){40}}
		\put(130,30){\circle*{5}}
		\put(130,30){\line(-2,1){40}}
		\put(130,30){\line(-2,-1){40}}
		\put(160,50){\circle*{5}}
		\put(160,10){\circle*{5}} 
		\put(130,30){\line(3,2){30}}
		\put(130,30){\line(3,-2){30}}
		\put(160,50){\line(3,0){30}}
		\put(160,50){\line(0,-4){40}}
		\put(160,50){\line(3,-4){30}}
		\put(160,10){\line(3,0){30}}
		\put(175,30){\circle*{5}} 
		\put(190,50){\circle*{5}}
		\put(190,10){\circle*{5}}
		\put(220,30){\circle*{5}}
		\put(190,50){\line(3,-2){30}}
		\put(190,10){\line(3,2){30}}

		\put(240,30){\circle*{5}}
		\put(280,30){\circle*{5}}
		\put(220,30){\line(1,0){180}}
		\put(260,50){\circle*{5}}
		\put(260,50){\line(-1,-1){20}}
		\put(260,50){\line(1,-1){20}}

		\put(400,30){\circle*{5}} \put(395,20){{$v_P$}}
		\put(370,0){\circle*{5}}
		\put(340,0){\circle*{5}}
		\put(355,60){\circle*{5}} 
		\put(310,30){\circle*{5}} \put(305,38){{$w$}} \put(288,33){{$e$}}
		\put(310,30){\line(3,2){45}}
		\put(310,30){\line(1,0){90}}
		\put(310,30){\line(1,-1){30}}
		\put(310,30){\line(2,-1){60}}
		\put(340,0){\line(1,4){15}}
		\put(340,0){\line(2,1){60}}
		\put(340,0){\line(1,0){30}}
		\put(370,0){\line(-1,4){15}}
		\put(370,0){\line(1,1){30}}
		\put(400,30){\line(-3,2){45}}

		\multiput(300,75)(5,0){25}{\line(1,0){3}}
		\multiput(425,75)(0,-5){18}{\line(0,-1){3}}
		\multiput(300,-15)(5,0){25}{\line(1,0){3}}
		\multiput(300,-15)(0,5){18}{\line(0,1){3}}
		\put(355,-10){{$R_{\P,P}$}}

		\multiput(10,80)(5,0){88}{\line(1,0){3}}
		\multiput(450,80)(0,-5){22}{\line(0,-1){3}}
		\multiput(10,-30)(5,0){88}{\line(1,0){3}}
		\multiput(10,-30)(0,5){22}{\line(0,1){3}}
		\put(150,-25){{$Q_{\P,P}$}}

		\multiput(50,130)(50,0){8}{\circle*{5}}
		\put(50,130){\line(1,0){350}}
		\put(170,140){{$P_2$}}

		\multiput(40,10)(0,5){40}{\line(0,1){3}}
		\multiput(40,210)(5,0){38}{\line(1,0){3}}
		\multiput(230,210)(0,-5){4}{\line(0,-1){3}}
		\multiput(230,190)(-5,0){34}{\line(-1,0){3}}
		\multiput(60,190)(0,-5){36}{\line(0,-1){3}}
		\multiput(40,10)(5,0){4}{\line(1,0){3}}	
		\put(42,180){{$X_{t'}$}}

		\multiput(390,10)(0,5){35}{\line(0,1){3}}
		\multiput(390,185)(-5,0){38}{\line(-1,0){3}}
		\multiput(200,185)(0,5){8}{\line(0,1){3}}
		\multiput(200,225)(5,0){42}{\line(1,0){3}}
		\multiput(410,225)(0,-5){43}{\line(0,-1){3}}
		\multiput(390,10)(5,0){4}{\line(1,0){3}}	
		\put(393,180){{$X_{t}$}}

		\multiput(187,250)(0,-5){16}{\line(0,1){3}}
		\multiput(190,175)(5,0){39}{\line(-1,0){3}}
		\multiput(380,173)(0,-5){18}{\line(0,1){3}}
		\multiput(380,88)(-5,0){17}{\line(-1,0){3}}
		\multiput(297,88)(0,-5){26}{\line(0,-1){3}}
		\put(190,240){{$M_{P}$}}

		\multiput(430,250)(0,-5){34}{\line(0,1){3}}
		\multiput(430,85)(-5,0){23}{\line(-1,0){3}}
		\multiput(317,85)(0,-5){32}{\line(0,-1){3}}
		\put(295,-65){{$L_{P}$}}
		
	\end{picture}
	\caption{In this example, $t'$ is a precursor of $t$ with $|X_{t'}|=|X_t|=3$, and $\P$ is a set of foundation paths between $t'$ and $t$ consisting of $|X_t|=3$ disjoint paths: the first path consists of the single vertex $z$, the second path is $P_2$, and the third path is a path $P$ between $u_P \in X_{t'}$ and $v_P \in X_t$ contained in $Q_{\P,P}$. Moreover, $R_{\P,P}=R'_{\P,P}$ in this example. And $L_P = \downarrow t-V(R'_{\P,P}-\{w\})$ and $M_P = \uparrow t \cup V(R'_{\P,P})$.} \label{fig_pstatic}
\end{figure}

The \defn{parent-side progress of $t$ with respect to $t'$} is the separation $$(\bigcap_{P'} L_{P'}, \bigcup_{P'} M_{P'}),$$ where the intersection and the union are taken over all parent-side static members $P'$ of $\P$ with respect to $t'$ and $t$. 
Observe that the order of $(\bigcap_{P'} L_{P'}, \bigcup_{P'} M_{P'})$ is $\lvert X_t \rvert$, and every vertex in $\bigcap_{P'} L_{P'} \cap (\bigcup_{P'} M_{P'})-X_t$ is pointed for $(\bigcap_{P'} L_{P'}, \bigcup_{P'} M_{P'})$.

For a member $P$ with $V(P) \cap X_t \neq V(P) \cap X_{t'}$, we say that $P$ is \defn{child-side static with respect to $t'$ and $t$} if 
	\begin{itemize}
		\item there exist no two edge-disjoint paths in $G-\bigcup_{W \in \P-\{P\}}V(W)$ from the vertex in $V(P) \cap X_{t'}$ to the vertex in $V(P) \cap X_t$ internally disjoint from $X_t \cup X_{t'}$, and 
		\item there exists no unambiguous child-side jump at $t'$ with respect to $t$ at the vertex in $V(P) \cap X_{t'}$.
	\end{itemize}
Notice that unlike the parent-side static case, it is possible that there exists an ambiguous child-side jump in the child-side static case.

Moreover, for a member $P$ with $V(P) \cap X_t \neq V(P) \cap X_{t'}$ and a precursor $t''$ of $t'$, we say that $P$ is \defn{strongly child-side static with respect to $t'',t',t$} if
	\begin{itemize}
		\item $P$ is child-side static with respect to $t'$ and $t$, and
		\item there exists $v \in X_{t'}-V(P)$ such that $v$ is an end of every parent-side jump at $t'$ with respect to $t''$ at the vertex in $X_{t'} \cap V(P)$ and every child-side jump at $t'$ with respect to $t$ at the vertex in $X_{t'} \cap V(P)$.
	\end{itemize}

Assume that $P$ is a member of $\P$ with $V(P) \cap X_t \neq V(P) \cap X_{t'}$ such that $P$ is child-side static with respect to $t'$ and $t$.
Let $u_P$ be the vertex in $V(P) \cap X_{t'}$ and $v_P$ the vertex in $V(P) \cap X_t$.
Since there exists no unambiguous child-side jump at $t'$ with respect to $t$ at $u_P$ and there exist no two edge-disjoint paths in $G-\bigcup_{W \in \P-\{P\}}V(W)$ from $u_P$ to $v_P$ internally disjoint from $X_t \cup X_{t'}$, we know that there exists a unique vertex $w$ in the $(L,\P,P)$-graph $L_{\P,P}$ such that $w$ is contained in a single-edge block with edge-set $\{e\}$ of the $(Q,\P,P)$-graph $Q_{\P,P}$, and we define $L_{\P,P}'$ to be the component of $(G-(X_{t'}-\{u_P\}))-e$ containing $L_{\P,P}$. 
In this case, there exists a separation $(L_P',M_P')$ of $G$ such that $L'_P = \downarrow t \cup V(L'_{\P,P}) \cup \{w'\}$ and $M'_P = \uparrow t - V(L'_{\P,P})$, where $w'$ is the end of $e$ other than $w$.
Note that $L_P' \cap M_P' = (X_t-\{u_P\}) \cup \{w'\}$, and $e$ is the unique edge incident with $w'$ whose other end is in $L_P'-M_P'$. 

For a precursor $t''$ of $t'$, the \defn{child-side progress of $t'$ with respect to $t,t''$} is the separation $$(\bigcup_{P'} L_{P'}', \bigcap_{P'} M_{P'}'),$$ where the union and the intersection are taken over all members $P'$ of $\P$ that are strongly child-side static with respect to $t'',t',t$ but not parent-side static with respect to $t''$ and $t'$. 
Observe that the order of $(\bigcup_{P'} L_{P'}', \bigcap_{P'} M_{P'}')$ is $\lvert X_t \rvert$, and every vertex in $(\bigcup_{P'} L_{P'}') \cap \bigcap_{P'} M_{P'}'-X_t$ is pointed for $(\bigcup_{P'} L_{P'}', \bigcap_{P'} M_{P'}')$.

\subsubsection{From progress to pseudo-edge-cuts}

Now we show how to use parent-side progress and child-side progress to obtain a pseudo-edge-cut.

\begin{lemma} \label{progress_shift}
Let $G$ be a graph.
Let $(T,\X)$ be a rooted tree-decomposition of $G$.
Let $t_1,t_2,t_3$ be nodes of $T$ such that $t_i$ is a precursor of $t_{i+1}$ for each $i \in [2]$.
Assume that the bags $X_{t_1},X_{t_2},X_{t_3}$ are pairwise disjoint sets of size $s$, for some positive integer $s$.
Let $\P$ be a set of foundation paths $\{P_1,P_2,...,P_s\}$ between $X_{t_1}$ and $X_{t_3}$.
If for every $i \in [s]$, $P_i$ is parent-side static with respect to $t_1$ and $t_2$ or is strongly child-side static with respect to $t_1,t_2,t_3$, then there exists a pseudo-edge-cut $(L^*,M^*)$ of $G$ modulo $\emptyset$ of order $s$ such that $\downarrow t_1 \subseteq L^*$ and $\uparrow t_3 \subseteq M^*$.
\end{lemma}

\begin{pf}
Since for every $i \in [s]$, $P_i$ is parent-side static with respect to $t_1$ and $t_2$ or is strongly child-side static with respect to $t_1,t_2,t_3$, we can reindex the members of $\P$ to assume that there exists $r \in [0,s]$ such that for every $i \in [r]$, $P_i$ is parent-side static with respect to $t_1$ and $t_2$, and for every $j \in [s]-[r]$, $P_j$ is not parent-side static with respect to $t_1$ and $t_2$ but is strongly child-side static with respect to $t_1,t_2,t_3$.

Let $(L,M)$ be the parent-side progress of $t_2$ with respect to $t_1$.
Let $(L',M')$ be the child-side progress of $t_2$ with respect to $t_3,t_1$.
Let $$(L^*,M^*)=(L \cup (\uparrow t_2 \cap L'-(X_{t_2}\cap (\bigcup_{\ell=1}^r V(P_\ell)))), M' \cup (\downarrow t_2 \cap M-(X_{t_2}\cap (\bigcup_{\ell=r+1}^s V(P_\ell))))).$$
We shall prove that $(L^*,M^*)$ is a pseudo-edge-cut $(L^*,M^*)$ of $G$ modulo $\emptyset$ of order $s$ such that $\downarrow t_1 \subseteq L^*$ and $\uparrow t_3 \subseteq M^*$.

\medskip

\noindent{\bf Claim 1:} $L^* \cup M^* = V(G)$.

\noindent{\bf Proof of Claim 1:}
Suppose to the contrary that there exists $v \in V(G)-(L^* \cup M^*)$.
Then $v \in (M-L) \cap (L'-M')$.
So if $v \in X_{t_2} \cap (\bigcup_{\ell=r+1}^sV(P_\ell))$, then $v \in \uparrow t_2 \cap L'-(X_{t_2}\cap (\bigcup_{\ell=1}^r V(P_\ell)) \subseteq L^*$, a contradiction; if $v \in X_{t_2} \cap (\bigcup_{\ell=1}^rV(P_\ell))$, then $v \in \downarrow t_2 \cap M-(X_{t_2}\cap (\bigcup_{\ell=r+1}^s V(P_\ell))) \subseteq M^*$, a contradiction.
Hence $v \not \in X_{t_2}$.
So if $v \in \downarrow t_2$, then $v \in \downarrow t_2 \cap (M-L) \subseteq M^*$, a contradiction; if $v \in \uparrow t_2$, then $v \in \uparrow t_2 \cap (L'-M') \subseteq L^*$, a contradiction.
$\Box$

\medskip

\noindent{\bf Claim 2:} There exists no edge $xy$ of $G$ such that $x \in L^*-(M^* \cup L \cup X_{t_2})$ and $y \in  X_{t_2} \cap (\bigcup_{\ell=1}^rV(P_\ell))-L^*$.

\noindent{\bf Proof of Claim 2:}
Suppose to the contrary that such an edge $xy$ of $G$ exists.
Since $x \in L^*-(L \cup X_{t_2})$, $x \in \uparrow t_2 \cap L'-X_{t_2}$.
Since $x \not \in M' \cup X_{t_2}$ and $y \in X_{t_2} \cap (\bigcup_{\ell=1}^rV(P_\ell))$, there exists a child-side jump $J$ at $t_2$ with respect to $t_3$ at the vertex in $X_{t_2} \cap V(P_{j'})$ for some $j' \in [s]-[r]$ such that $J$ contains $xy$.
Since $P_{j'}$ is child-side static with respect to $t_2$ and $t_3$, $J$ is an ambiguous child-side jump at $t_2$ with respect to $t_3$ at the vertex in $X_{t_2} \cap V(P_{j'})$.
So the ends of $J$ and are the vertex in $X_{t_2} \cap V(P_{j'})$ and $y$.
Since $P_{j'}$ is strongly child-side static with respect to $t_1,t_2,t_3$ but not parent-side static with respect to $t_1$ and $t_2$, there exists a parent-side jump $J'$ at $t_2$ with respect to $t_1$ at the vertex in $X_{t_2} \cap V(P_{j'})$ such that $y$ is an end of $J'$.
So $J'$ is a parent-side jump at $t_2$ with respect to $t_1$ at $y \in X_{t_2} \cap (\bigcup_{\ell=1}^rV(P_\ell))$, contradicting that $P_\ell$ is parent-side static with respect to $t_1$ and $t_2$ for every $\ell \in [r]$. 
$\Box$

\medskip

\noindent{\bf Claim 3:} $(L^*,M^*)$ is a separation of $G$.

\noindent{\bf Proof of Claim 3:}
Suppose to the contrary that $(L^*,M^*)$ is not a separation of $G$.
Since $L^* \cup M^* = V(G)$ by Claim 1, there exists an edge $xy$ of $G$ such that $x \in L^*-M^*$ and $y \in M^*-L^*$.

We first suppose that $y \in \downarrow t_2 \cap M-(X_{t_2}\cap (\bigcup_{\ell=r+1}^s V(P_\ell)))$.
If $x \in X_{t_2}$, then since $x \not \in M'$, we have $x \in X_{t_2} \cap \bigcup_{\ell=r+1}^sV(P_\ell)$, so some path containing $xy$ is a parent-side jump at $t_2$ with respect to $t_1$ at the vertex in $X_{t_2} \cap V(P_{j'})$ for some $j' \in [r]$, a contradiction.
If $x \in L-X_{t_2}$, then since $y \not \in L$ and $(L,M)$ is a separation, we know $x \in L \cap M-X_{t_2} \subseteq \downarrow t_2 \cap M - X_{t_2} \subseteq M^*$, a contradiction.
So $x \not \in L \cup X_{t_2}$.
Since $x \in L^*-M^*$, $x \in \uparrow t_2 \cap L'-X_{t_2}$.
Since $xy \in E(G)$, $y \in \uparrow t_2$, so $y \in \uparrow t_2 \cap \downarrow t_2 \cap M-(X_{t_2}\cap (\bigcup_{\ell=r+1}^s V(P_\ell))) \subseteq X_{t_2} \cap (\bigcup_{\ell=1}^rV(P_\ell))$, contradicting Claim 2.  

Hence $y \in M' - (\downarrow t_2 \cap M-(X_{t_2}\cap (\bigcup_{\ell=r+1}^s V(P_\ell))))$.

Suppose $x \in \uparrow t_2 \cap L'-(X_{t_2}\cap (\bigcup_{\ell=1}^r V(P_\ell)))$.
Since $x \not \in M^*$, $x \in L'-M'$.
Since $xy \in E(G)$, $y \in L' \cap M'$.
Since $y \not \in L^*$, $y \in L' \cap M' \cap X_{t_2} \cap \bigcup_{\ell=1}^rV(P_\ell)$.
If $x \in X_{t_2}$, then $x \in X_{t_2} \cap \bigcup_{\ell=r+1}^sV(P_\ell)$, so $xy$ forms a parent-side jump at $t_2$ with respect to $t_1$ at the vertex in $X_{t_2} \cap V(P_{j'})$ for some $j' \in [s]$, a contradiction.
So $x \not \in X_{t_2}$.
Since $L \subseteq \downarrow t_2$ and $x \in \uparrow t_2 - (M^* \cup X_{t_2})$, we have $x \in L^*-(M^* \cup L \cup X_{t_2})$, contradicting Claim 2.

Hence $x \in L-(\uparrow t_2 \cap L'-(X_{t_2}\cap (\bigcup_{\ell=1}^r V(P_\ell))))$.
Since $L \subseteq \downarrow t_2$ and $x \not \in M^*$, $x \not \in M - (X_{t_2} \cap \bigcup_{\ell=r+1}^sV(P_\ell))$.
Since $X_{t_2} \cap \bigcup_{\ell=r+1}^sV(P_\ell) \subseteq \uparrow t_2 \cap L'-(X_{t_2}\cap (\bigcup_{\ell=1}^r V(P_\ell)))$, $x \not \in M$.

If $y \in X_{t_2}$, then $y \in M' \cap X_{t_2} \subseteq X_{t_2} \cap (\bigcup_{\ell=1}^r V(P_\ell)))$, a contradiction.
So $y \in M'-X_{t_2} \subseteq \uparrow t_2 - X_{t_2}$.
Since $xy \in E(G)$, $x \in \uparrow t_2 \subseteq M$, a contradiction.
$\Box$

\medskip

\noindent{\bf Claim 4:} $L^* \cap M^* \subseteq (L \cap M \cap \bigcup_{\ell=1}^rV(P_\ell)) \cup (L' \cap M' \cap \bigcup_{\ell=r+1}^sV(P_\ell))$.

\noindent{\bf Proof of Claim 4:}
By the definition of $L^*$ and $M^*$, we have $L^* \cap M^* \subseteq (L \cap M') \cup (L \cap M-(X_{t_2} \cap \bigcup_{\ell=r+1}^sV(P_\ell))) \cup (L' \cap M' - (X_{t_2} \cap \bigcup_{\ell=1}^rV(P_\ell)))$.

Note that $L \cap M' \subseteq \downarrow t_2 \cap \uparrow t_2 = X_{t_2}$. 
So if $v$ is a vertex in $L \cap M'$, then there exists $j \in [s]$ such that $v$ is the vertex in $X_{t_2} \cap V(P_j)$ and $V(P_j) \cap L \cap M = \{v\} = V(P_j) \cap L' \cap M'$.
Hence $L \cap M' \subseteq (L \cap M-(X_{t_2} \cap \bigcup_{\ell=r+1}^sV(P_\ell))) \cup (L' \cap M' - (X_{t_2} \cap \bigcup_{\ell=1}^rV(P_\ell)))$. 

Since $L \cap M \cap \bigcup_{\ell=r+1}^sV(P_\ell) \subseteq X_{t_2}$, $L \cap M-(X_{t_2} \cap \bigcup_{\ell=r+1}^sV(P_\ell)) = L \cap M \cap \bigcup_{\ell=1}^rV(P_\ell)$.
Similarly, $L' \cap M' - (X_{t_2} \cap \bigcup_{\ell=1}^rV(P_\ell)) = L' \cap M' \cap \bigcup_{\ell=r+1}^sV(P_\ell)$.
Therefore, $L^* \cap M^* \subseteq (L \cap M \cap \bigcup_{\ell=1}^rV(P_\ell)) \cup (L' \cap M' \cap \bigcup_{\ell=r+1}^sV(P_\ell))$.
$\Box$

\medskip

\noindent{\bf Claim 5:} $(L^*,M^*)$ is a pseudo-edge-cut modulo $\emptyset$.

\noindent{\bf Proof of Claim 5:}
Suppose to the contrary that $(L^*,M^*)$ is not a pseudo-edge-cut modulo $\emptyset$.
By Claim 3, $(L^*,M^*)$ is a separation of $G$, so there exists $v \in L^* \cap M^*$ such that $v$ is not pointed for $(L^*,M^*)$.
By Claim 4, $v \in (L \cap M \cap \bigcup_{\ell=1}^rV(P_\ell)) \cup (L' \cap M' \cap \bigcup_{\ell=r+1}^sV(P_\ell))$.

We first suppose that $v \in L \cap M \cap \bigcup_{\ell=1}^rV(P_\ell)$.
For every $j \in [r]$, since $P_j$ is parent-side static with respect to $t_1$ and $t_2$, the vertex in $L \cap M \cap V(P_j)$ is pointed for $(L,M)$.
So $v$ is pointed for $(L,M)$.
Since $v$ is not pointed for $(L^*,M^*)$, $v$ has a neighbor $u$ in $(L^*-M^*)-(L-M)$.
Note that 
	\begin{align*}
		(L^*-M^*)-(L-M) \subseteq & ((L-M^*)-(L-M)) \cup (\uparrow t_2 \cap L' - ((X_{t_2} \cap \bigcup_{\ell=1}^rV(P_\ell)) \cup M^*))\\
		\subseteq & (L \cap M - M^*) \cup ((\uparrow t_2 \cap L' - (X_{t_2} \cup M'))) \cup (X_{t_2} \cap \bigcup_{\ell=r+1}^sV(P_\ell)) \\
		\subseteq & (X_{t_2} \cap \bigcup_{\ell=r+1}^sV(P_\ell)) \cup (\uparrow t_2 \cap L' - (X_{t_2} \cup M')).
	\end{align*}
If $u \in X_{t_2} \cap \bigcup_{\ell=r+1}^sV(P_\ell)$, then there exists a parent-side jump at $t_2$ with respect to $t_1$ at the vertex in $X_{t_2} \cap V(P_{j'})$ for some $j' \in [r]$, a contradiction.
So $u \in \uparrow t_2 \cap L'-(X_{t_2} \cup M')$.
Hence $v \in L \cap M \cap \bigcup_{\ell=1}^rV(P_\ell) \cap X_{t_2}$.
So the path consists of $xy$ is an unambiguous child-side jump at $t_2$ with respect to $t_3$ at the vertex in $X_{t_2} \cap V(P_{j'})$ for some $j' \in [s]-[r]$, a contradiction.

So $v \in L' \cap M' \cap \bigcup_{\ell=r+1}^sV(P_\ell)$.
Since $v$ is pointed for $(L',M')$, $v$ has a neighbor in $(L^*-M^*)-(L'-M')$.
But $L^* \subseteq L'$ and $M^* \supseteq M'$, so $(L^*-M^*)-(L'-M') = \emptyset$, a contradiction.
$\Box$

\medskip

By Claim 5 and the existence of $\{P_1,...,P_s\}$, $(L^*,M^*)$ is a pseudo-edge-cut modulo $\emptyset$ of size $s$.
Note that $\downarrow t_1 \subseteq L \subseteq L^*$ and $\uparrow t_3 \subseteq M' \subseteq M^*$.
This proves the lemma.
\end{pf}

\subsubsection{Breaking}

\begin{lemma} \label{breaking a long strip}
For any positive integers $k,w,\alpha$, there exists an integer $f=f(k,w,\alpha)$ such that the following holds.
Let $G$ be a graph that does not contain the Robertson chain of length $k$ as a topological minor.
Let $(T,\X)$ be a rooted tree-decomposition of $G$ of width at most $w$.
If there exist $Z \subseteq V(G)$, a positive integer $s$ and a $(Z,s)$-strip $(t_1,t_2,...,t_{f+1})$ in $(T,\X)$, then there exists a pseudo-edge-cut $(A,B)$ modulo $Z$ of order $\lvert X_{t_1} \rvert$ $\alpha$-breaking $(t_1,...,t_{f+1})$.
\end{lemma}

\begin{pf}
Define $f = f_{\ref{breaking_weak}}(k,w,2\alpha+4)$, where $f_{\ref{breaking_weak}}$ is the function $f$ mentioned in Lemma \ref{breaking_weak}.
We shall show that $f$ satisfies the conclusion of the lemma.

Let $G$, $(T,\X)$, $Z$, $s$ and $(t_1,t_2,...,t_{f+1})$ be the ones as mentioned in the lemma.
Since $(t_1,t_2,...,t_{f+1})$ is a $(Z,s)$-strip, there exists a set $\P=\{P_1,P_2,...,P_{\lvert Z \rvert +s}\}$ of foundation paths between $t_1$ and $t_{f+1}$ such that $P_i$ is a one-vertex path with $V(P_i) \subseteq Z$ for every $i\in [s+\lvert Z \rvert]-[s]$.

By Lemma \ref{breaking_weak}, there exist a subsequence $(t_1',t_2',...,t_{2\alpha+5}')$ of $(t_1,t_2,...,t_{f+1})$ such that for each $j'\in [s]$, there exist no two edge-disjoint paths in $G-\bigcup_{W \in \P-\{P_{j'}\}}V(W)$ from the vertex in $X_{t'_{\ell}} \cap V(P_{j'})$ to the vertex in $X_{t'_{{\ell+1}}} \cap V(P_{j'})$ for any $\ell \in [2\alpha+4]$, and either 
	\begin{itemize}
		\item for every $j \in [2\alpha+5]-[1]$, there exist no parent-side jumps at $t'_{j}$ with respect to $t'_{j-1}$ at the vertex in $X_{t'_{j}} \cap V(P_{j'})$ disjoint from $Z$, or 
		\item for every $j \in [2\alpha+4]$, 
			\begin{itemize}
				\item there exist no unambiguous child-side jumps at $t'_{j}$ with respect to $t'_{j+1}$ at the vertex in $X_{t'_{j}} \cap V(P_{j'})$ disjoint from $Z$, and
				\item if $j \geq 2$, then there exists $v \in X_{t'_j}-(Z \cup V(P_{j'}))$ such that $v$ is an end of every parent-side jump at $t'_{j}$ with respect to $t'_{j-1}$ at the vertex in $X_{t'_{j}} \cap V(P_{j'})$ disjoint from $Z$ and every child-side jump at $t'_{j}$ with respect to $t'_{j+1}$ at the vertex in $X_{t'_{j}} \cap V(P_{j'})$ disjoint from $Z$.
			\end{itemize}
	\end{itemize}

For every $t \in V(T)$, let $X'_t = X_t-Z$.
Let $\X' = (X'_t: t \in V(T))$.
Then $(T,\X')$ is a rooted tree-decomposition of $G-Z$ of width at most $w$.
Note that $\{P_1,P_2,...,P_s\}$ is a set of foundation paths between $X'_{t_1'}$ and $X'_{t'_{2\alpha+5}}$.
Then for any $j \in [2\alpha+4]-[1]$ and $j' \in [s]$, either $P_{j'}$ is parent-side static with respect to $t'_{j-1}$ and $t'_j$, or $P_{j'}$ is strongly child-side static with respect to $t'_{j-1},t'_j,t'_{j+1}$.

Let $i^*=\alpha+2$.
By Lemma \ref{progress_shift}, there exists a pseudo-edge-cut $(L^*_0,M^*_0)$ modulo $\emptyset$ of order $s$ in $G-Z$ such that $(T,\X')\downarrow t'_{i^*-1} \subseteq L^*_0$ and $(T,\X') \uparrow t'_{i^*+1} \subseteq M_0^*$.

Define $(L^*,M^*) = (L^*_0 \cup Z, M^*_0 \cup Z)$.
Then $(L^*,M^*)$ is a pseudo-edge-cut modulo $Z$ of order $s+\lvert Z \rvert$ in $G$ such that $(T,\X) \downarrow t'_{\alpha+1} \subseteq L^*$ and $(T,\X) \uparrow t'_{\alpha+3} \subseteq M^*$.
Therefore, $(L^*,M^*)$ $\alpha$-breaks $(t_1',t_2',...,t'_{2\alpha+5})$ and hence $\alpha$-breaks $(t_1,t_2,...,t_{f+1})$.
This proves the lemma.
\end{pf}

\subsubsection{Anti-pointed vertices and breaking coherently}

Let $(A,B)$ be a separation of a graph $G$.
We say that a vertex $v \in A \cap B$ is \defn{anti-pointed for $(A,B)$} if $v$ is incident with at most one edge whose other end is in $B-A$.

Let $(T,\X)$ be a rooted tree-decomposition of a graph $G$.
Let $t_1,t_2$ be nodes of $T$, where $t_1$ is an ancestor of $t_2$.
We say that a vertex $v \in X_{t_1} \cap X_{t_2}$ is \defn{coherent} for $t_1,t_2$, if the following two statements hold.
	\begin{itemize}
		\item Either $v$ is not pointed for $(A_{t_1},B_{t_1})$, or there exists $i \in \{0,1\}$ such that both the number of edges between $v$ and $A_{t_1}-B_{t_1}$ and the number of edges between $v$ and $A_{t_2}-B_{t_2}$ are equal to $i$.
		\item Either $v$ is not anti-pointed for $(A_{t_2},B_{t_2})$, or there exists $i \in \{0,1\}$ such that both the number of edges between $v$ and $B_{t_1}-A_{t_1}$ and the number of edges between $v$ and $B_{t_2}-A_{t_2}$ are equal to $i$.
	\end{itemize}

\begin{lemma} \label{breaking a long strip better}
For any positive integers $k,w,\alpha$, there exist integers $f(k,w,\alpha)$ and $g(k,w,\alpha) = 9^w (f(k,w,\alpha)+1)$ such that the following hold.
Let $G$ be a graph that does not contain the Robertson chain of length $k$ as a topological minor.
Let $(T,\X)$ be a rooted tree-decomposition of $G$ of width at most $w$.
If there exist $Z \subseteq V(G)$, a positive integer $s$ and a $(Z,s)$-strip $(t_1,t_2,...,t_{g(k,w,\alpha)})$, then there exist a $(Z,s)$-strip $R$ and a pseudo-edge-cut $(A,B)$ modulo $Z$ of order $\lvert X_{t_1} \rvert$ such that the following hold.
	\begin{enumerate}
		\item $R$ is a subsequence of $(t_1,t_2,...,t_{g(k,w,\alpha)})$ of length $f(k,w,\alpha)+1$.
		\item For any pairs of nodes $t,t'$ in $R$, every vertex in $Z$ is coherent for $t$ and $t'$. 
		\item $(A,B)$ $\alpha$-breaks $R$.
	\end{enumerate}
\end{lemma}

\begin{pf}
Define $f(k,w,\alpha)$ to be the number $f(k,w,\alpha)$ mentioned in Lemma \ref{breaking a long strip}.

For each $t_i$ and each $v \in Z$, 
	\begin{itemize}
		\item either 
			\begin{itemize}
				\item $v$ is not pointed for $(A_{t_i},B_{t_i})$, or 
				\item $v$ is incident with exactly $x$ edges whose other end is in $A_{t_i}-B_{t_i}$ for some $x \in \{0,1\}$, 
			\end{itemize}
			and 
		\item either 
			\begin{itemize}
				\item $v$ is not anti-pointed for $(A_{t_i},B_{t_i})$, or 
				\item $v$ is incident with exactly $y$ edges whose other end is in $B_{t_i}-A_{t_i}$ for some $y \in \{0,1\}$.
			\end{itemize}
	\end{itemize}
Hence for each $t_i$ and $v \in Z$, there are nine possibilities mentioned above, so we can use nine colors to color each pair $(t_i,v)$.
So for each $t_i$, there are at most $9^{\lvert Z \rvert}$ different colors for the pairs in $\{(t_i,v): v \in Z\}$.
Since $\lvert Z \rvert \leq w$, we can use at most $9^w$ different colors to color each $t_i$ according to the colors of $(t_i,v)$ for $v \in Z$.
Therefore, there are at least $g(k,w,\alpha)/9^w \geq f(k,w,\alpha)+1$ nodes in $\{t_i: 1 \leq i \leq g(k,w,\alpha)\}$ having the same color.
Let $R$ be a subsequence of $(t_1,t_2,...,t_{g(k,w,\alpha)})$ with $f(k,w,\alpha)+1$ entries with the same color.
Then $R$ is a $(Z,s)$-strip such that every vertex in $v \in Z$ is coherent for any pair of nodes in $R$.
Finally, by Lemma \ref{breaking a long strip}, there exists a pseudo-edge-cut modulo $Z$ with order $\lvert X_{t_1} \rvert$ $\alpha$-breaking $R$.
This proves the lemma.
\end{pf}

\section{Bounding the elevation}
\label{sec: bounding elevation}

In this section, we prove the first key ingredient (Theorem \ref{bounded depth}) in this paper, which is a structure theorem for graphs of bounded tree-width that does not contain a long Robertson chain as a topological minor.
Roughly speaking, the theorem states that every such a graph has a tree-decomposition of small width such that its bags provide sufficiently many pseudo-edge-cuts to break long strips, and there are disjoint paths between the bags with the same size unless the bags are separated by separations of smaller order given by other bags.
The existence of such a tree-decomposition is crucial for future sections that address well-quasi-ordering.

\subsection{Intuition}

We need a number of terminologies in order to develop the necessary machinery for proving Theorem \ref{bounded depth}.
We explain the intuition behind those terminologies and the lemmas that we will prove.

Lemma \ref{breaking a long strip better}, proved in the previous section, shows that the absence of a long Robertson chain implies the existence of pseudo-edge-cuts that break long strips.
However, those pseudo-edge-cuts possibly are not separations given by bags of the tree-decomposition.
So we have to ``insert'' pseudo-edge-cuts into the current tree-decomposition to make them given by bags.
The main technical issue is that inserting a pseudo-edge-cut into a tree-decomposition can destroy other pseudo-edge-cuts given by bags in the previous tree-decomposition so that the new tree-decomposition is not better than the old one.
The key solution for overcoming this technical issue is to insert the ``correct'' pseudo-edge-cut and its ``partner'' and design a way to evaluate how good a tree-decomposition is.

The formal definition of the ``partner'' is the reflection that we define in Section \ref{subsec:reflection}.
Roughly speaking, a reflection of a separation $(A,B)$ is another separation obtained by replacing some vertices in $A \cap B$ pointed for $(A,B)$ by their unique neighbors in $A-B$.
Inserting both a separation and its reflection can reduce the possibilities for destroying existing pseudo-edge-cuts given by bags, comparing to only inserting a separation.
As we want the reflection $(A',B')$ of $(A,B)$ also separates two given sets $X$ and $Y$ that are separated by $(A,B)$, we have to assume that vertices in $X$ are not in $A \cap B - (A' \cap B')$, which is the motivation of the notion for a separation strongly separating two sets introduced in Section \ref{subsec:reflection}.
In Section \ref{subsec:reflection}, we formally define those notions and show that we can transform a pair of separations $(A_1,B_1)$ and $(A_2,B_2)$, where one of them is a reflection of the other, into another pair of separations that have extra properties so that we have better control for the new tree-decomposition after we insert both of them into a tree-decomposition.

Now we consider how to evaluate how good a tree-decomposition is.
Intuitively we would like to say that a tree-decomposition is better if its bags ``provide'' more pseudo-edge-cuts.
But the naive way that counts the number of pseudo-edge-cuts given by bags does not work: it is possible that such counts is significantly decreased after inserting a nicely chosen pseudo-edge-cut and its reflection because many pseudo-edge-cuts given by bags in the old tree-decomposition are no longer given by bags in the new one.
The key solution is to find a way to collect the information in the new tree-decomposition about what pseudo-edge-cuts in the old tree-decomposition are destroyed.
The formal form of this information is the ``incorporation'' defined in Section \ref{subsec:incorporation}.
It leads us to define the ``signature'' of a tree-decomposition in Section \ref{subsec:incorporation} to evaluation how good a tree-decomposition is.
Roughly speaking, a tree-decomposition has higher signature has more incorporated pseudo-edge-cuts, where we prefer pseudo-edge-cuts with smaller order and more pointed vertices.
It is the motivation of the ``breadth'' defined in Section \ref{subsec:incorporation}.

By considering the tree-decomposition with the highest signature, we have a \linebreak tree-decomposition that incorporates many separations (or pseudo-edge-cuts).
But incorporated separations are not necessarily given by bags, and our goal is to obtain a tree-decomposition whose bags give pseudo-edge-cuts that break strips.
So the next step is to show that the tree-decomposition with the highest signature actually contains the bags that give desired pseudo-edge-cuts.
The formal form of ``containing the bags that give desired pseudo-edge-cuts'' is the ``integration'' defined in Section \ref{subsec:integration}.
We will prove this statement in Section \ref{subsec:integration}, which is the most technical part of this section and relies on the preparation developed in other subsections of this section.
A proof sketch will be provided there.

It remains to consider the other desired property of our tree-decomposition: we want disjoint paths between bags with the same size unless some bag gives a separation of smaller order that witness the non-existence of the paths.
Such a linkedness property is a standard one used in the literature.
However, we are not able to prove that a tree-decomposition can have this linkedness property as well as having the aforementioned property about pseudo-edge-cuts.
We can only achieve a weaker version of the linkedness property while keeping the property about pseudo-edge-cuts.
This weaker version of the linkedness property is described in Section \ref{subsec:linkedness} and is sufficient for us to prove results on well-quasi-ordering in later sections.
In Section \ref{subsec:incorporation}, we will prove that the tree-decomposition with the highest signature have this weak linkedness property.

Finally, we will combine everything together in Section \ref{subsec:elevation} to prove the main structure theorem that will use be for proving well-quasi-ordering results in future sections.

\subsection{Linkedness} \label{subsec:linkedness}

\begin{lemma} \label{separating bags to subtrees applied}
Let $(T,\X)$ be a rooted tree-decomposition of a graph $G$.
Let $t_1,t_2 \in V(T)$ be such that $t_1$ is an ancestor of $t_2$.
Let $(A, B)$ be a separation of $G$ such that $X_{t_1}\subseteq A$ and  $X_{t_2}\subseteq B$.
Then there exists a separation $(A^*,B^*)$ of $G$ such that the following hold.
	\begin{enumerate}
		\item[{\rm(i)}] $A^* \cap B^* = A \cap B$.
		\item[{\rm(ii)}] $\downarrow t_1 \subseteq A^*$ and $\uparrow t_2 \subseteq B^*$.
		\item[{\rm(iii)}] For every node $s$ of $T$  that is a descendant of $t_1$ and not an ancestor or a descendant of $t_2$, 
if $X_s \subseteq A^*$ or $X_s \subseteq B^*$, then either $\uparrow s \subseteq A^*$, or $\uparrow s \subseteq B^*$.
	\end{enumerate}
\end{lemma}

\begin{pf}
Let us say that a node $s\in V(T)$ is a {\em side node} if $s$ is a descendant of $t_1$ and not an ancestor or a descendant of $t_2$.
Let $(A,B)$ be a separation as in the statement of this lemma.
We say that a side node $s\in V(T)$ is {\em bad (for $(A,B)$)} if $X_s \subseteq A$ or $X_s \subseteq B$, but $\uparrow s \not\subseteq A$ and $\uparrow s \not\subseteq B$. 
We say that $t_1$ is {\em bad} if $\downarrow t_1 \not\subseteq A$ and
we say that $t_2$ is {\em bad} if $\uparrow t_2 \not\subseteq B$.
We proceed by induction on the number of bad nodes.
If there are no bad nodes, then $(A,B)$ satisfies the conclusion of the lemma.

If $t_1$ is bad, then every vertex $v\in\downarrow t_1-A$ belongs to a component $C$ of $G-(A\cap B)$
that is disjoint from $A$, and hence from $X_{t_1}$, and hence it is disjoint from $\uparrow t_2$ and $\uparrow s$
for every side node $s$.
Let $A'=A\cup\bigcup V(C)$ and $B'=B-\bigcup V(C)$, where the union is over all components $C$  of $G-(A\cap B)$
contained in $\downarrow t_1-A$. Then $(A',B')$ is a separation of $G$ with  $X_{t_1}\subseteq A'$ and  $X_{t_2}\subseteq B'$
that satisfies $A' \cap B' = A \cap B$.
Furthermore, every node that is bad for $(A',B')$ is bad for $(A,B)$, but $t_1$ is no longer bad for $(A',B')$.
Thus the conclusion of the lemma follows by induction applied to the separation  $(A',B')$.

An analogous argument applies when $t_2$ is bad. We may therefore assume that there exists a bad side node $s$.
Thus $X_s \subseteq A$ or $X_s \subseteq B$, and $\uparrow s \not\subseteq A$ and $\uparrow s \not\subseteq B$. 
We assume that $X_s \subseteq A$, because the case $X_s \subseteq B$ is analogous.
Since $\uparrow s \not\subseteq A$ and $X_s \subseteq A$, there exists a component $C$ of $G-(A\cap B)$ contained in $\uparrow s-A$.
Every such component is disjoint from $\downarrow t_1$, $\uparrow t_2$, and $\uparrow s'$ for every bad side node $s'$ that is not a descendant or an ancestor of $s$.
Let $A'=A\cup\bigcup V(C)$ and $B'=B-\bigcup V(C)$, where the union is over all components $C$  of $G-(A\cap B)$
contained in $\uparrow s-A$. Then $(A',B')$ is a separation of $G$ with  $X_{t_1}\subseteq A'$ and  $X_{t_2}\subseteq B'$
that satisfies $A' \cap B' = A \cap B$.

We claim that every node that is bad for $(A',B')$ is bad for $(A,B)$.
To see that let $s'$ be bad for $(A',B')$.
Then $s'$ is not a descendant of $s$, because $\uparrow s\subseteq A'$. 
If it is an ancestor of $s$, then $X_{s'}\cap A=X_{s'}\cap A'$ and $X_{s'}\cap B=X_{s'}\cap B'$. 
Then $\uparrow s' \not\subseteq A$ and $\uparrow s' \not\subseteq B$, because the same holds for $s$.
 Thus $s'$ is bad for $(A,B)$.
If $s'$ is neither a descendant nor an ancestor of $s$, then  $X_{s'}\cap A=X_{s'}\cap A'$ and $X_{s'}\cap B=X_{s'}\cap B'$,
and $\uparrow s'\cap A=\uparrow s'\cap A'$ and $\uparrow s'\cap B=\uparrow s'\cap B'$.
Thus, again, $s'$ is bad for $(A,B)$.

But $s$ is no longer bad for $(A',B')$.
Thus the conclusion of the lemma follows by induction applied to the separation  $(A',B')$.
\end{pf}

\bigskip

Let $N$ be a positive integer.
We say that a rooted tree-decomposition $(T,\X)$ of a graph $G$ is \defn{$N$-linked} if the following holds.
	\begin{itemize}
		\item If $t_1,t_2 \in V(T)$, where $t_1$ is a precursor of $t_2$, such that $\uparrow t_2$ contains at least $N$ vertices of $G$ each of which cannot be separated from $\downarrow t_1$ by a separation of order less than $\lvert X_{t_1} \rvert$ given by a node of $T$, then there does not exist a separation $(A,B)$ of $G$ of order less than $\lvert X_{t_1} \rvert$ with $\downarrow t_1 \subseteq A$ and $\uparrow t_2 \subseteq B$.
	\end{itemize}
We say that a rooted tree-decomposition $(T,\X)$ of a graph $G$ is \defn{weakly $N$-linked} if the following holds.
	\begin{itemize}
		\item If $t_1,t_2,...,t_{N+1}$ are nodes of $T$ such that $t_i$ is a precursor of $t_{i+1}$ for all $i \in [N]$, and the sets $X_{t_i}$ are distinct for $i \in [N+1]$, then there exist $\lvert X_{t_1} \rvert$ disjoint paths in $G$ from $X_{t_1}$ to $X_{t_2}$.
	\end{itemize}

The following lemma shows that every $N$-linked rooted tree-decomposition is weakly $N$-linked.

\begin{lemma} \label{N-linked disjoint paths}
Let $N$ be a positive integer and let $(T,\X)$ be an $N$-linked rooted tree-decomposition of a graph $G$.
Let $t_1,t_2,...,t_{N+1}$ be nodes of $T$ such that $t_i$ is a precursor of $t_{i+1}$ for all $i \in [N]$.
If the sets $X_{t_i}$ are distinct for $i \in [N+1]$, then there exist $\lvert X_{t_1} \rvert$ disjoint paths in $G$ from $X_{t_1}$ to $X_{t_2}$.
\end{lemma}

\begin{pf}
Suppose that there do not exist $\lvert X_{t_1} \rvert$ disjoint paths in $G$ from $X_{t_1}$ to $X_{t_2}$.
So $X_{t_1} \neq X_{t_2}$ and $\lvert X_{t_1} \rvert \geq 1$.
Since $X_{t_1}, X_{t_2}, ..., X_{t_{N+1}}$ are $N+1$ distinct sets with the same size, $X_{t_{i+1}}-\bigcup_{j=1}^iX_{t_j} \neq \emptyset$ for every $i \in [N]$ by the definition of a tree-decomposition.
Since $t_i$ is a precursor of $t_{i+1}$ for every $i \in [N]$, we know $\lvert X_t \rvert \geq \lvert X_{t_1} \rvert$ for all $t \in t_1Tt_{N+1}$, so $\uparrow t_2$ contains at least $\lvert X_{t_2} \rvert + N-1 \geq N$ vertices each of which cannot be separated from $\downarrow t_1$ by a separation of order less than $\lvert X_{t_1} \rvert$ given by a node of $T$.

Let $(A,B)$ be a separation of $G$ with minimum order such that $X_{t_1} \subseteq A$ and $X_{t_2} \subseteq B$.
So the order of $(A,B)$ is less than $\lvert X_{t_1} \rvert$ by Menger's theorem.
By Lemma \ref{separating bags to subtrees applied}, there exists a separation $(A^*,B^*)$ with $\lvert A^* \cap B^* \rvert = \lvert A \cap B \rvert < \lvert X_{t_1} \rvert$ such that $\downarrow t_1 \subseteq A^*$ and $\uparrow t_2 \subseteq B^*$.
So $(T,\X)$ is not $N$-linked, a contradiction.
\end{pf}

\subsection{Incorporation} \label{subsec:incorporation}

We say that a separation $(A,B)$ of a graph $G$ \defn{weakly separates} two subsets $Y,Z$ of $V(G)$ if either 
	\begin{itemize}
		\item $Y \subseteq A$ and $Z \subseteq B$, or 
		\item $Y \subseteq B$ and $Z \subseteq A$.
	\end{itemize}

Let $G$ be a graph and let $(A,B)$ be a separation of $G$.
Recall that we say that a vertex $v \in A \cap B$ is \defn{pointed for $(A,B)$} if $v$ is incident with at most one edge whose other end is in $A-B$; we say that a vertex $v \in A \cap B$ is \defn{anti-pointed for $(A,B)$} if $v$ is incident with at most one edge whose other end is in $B-A$.
We further define the following:
	\begin{itemize}
		\item We say that $v \in A \cap B$ is \defn{doubly pointed for $(A,B)$} if it is pointed and anti-pointed for $(A,B)$.
		\item The \defn{thickness} of $(A,B)$ is the number of vertices in $A \cap B$ not pointed for $(A,B)$.
		\item The \defn{breadth} of $(A,B)$ is the sequence $(\lvert A \cap B \rvert, j)$, where $j$ is the thickness of $(A,B)$.
	\end{itemize}
In this paper, sequences are compared by the lexicographic order.
So the breadth of a separation $(A,B)$ is smaller than the breadth of a separation $(C,D)$ if and only if either the order of $(A,B)$ is less than the order of $(C,D)$, or they have the same order but the thickness of $(A,B)$ is smaller than the thickness of $(C,D)$.

Let $(T,\X)$ be a rooted tree-decomposition of a graph $G$.
Recall that the separation given by a node $t$ of $T$ is the separation $(\downarrow t, \uparrow t)$ and is denoted by $(A_t,B_t)$.
Let $(A,B)$ be a separation of $G$.
We say that $(A,B)$ is \defn{incorporated} in $(T,\X)$ if there exists $S \subseteq V(T)$ such that the following hold.
\begin{itemize}
	\item[(INC1)] For every $t \in S$, the breadth of the separation $(A_t,B_t)$ given by $t$ in $(T,\X)$ is at most the breadth of $(A,B)$.
	\item[(INC2)] $B = \bigcup_{t \in S} B_t$.
	\item[(INC3)] $\sum_{i,j} a_{i,j}2^{2i^2+j} \leq 2^{2\lvert A \cap B \rvert^2 + \ell}$, where $\ell$ is the thickness of $(A,B)$, and $a_{i,j}$ is the number of nodes $t$ in $S$ such that $(A_t,B_t)$ has breadth $(i,j)$.
\end{itemize}
We say that $S$ is a \defn{witness set} of the incorporation.
In fact, (INC1) follows from (INC3), but (INC1) is included for better clarity.

\begin{lemma} \label{incorporation basic}
If $(A,B)$ is a separation that is incorporated in a rooted tree-decomposition $(T,\X)$ of a graph $G$ with witness set $S$, then the following hold.
	\begin{enumerate}
		\item $\lvert S \rvert \leq 2^{2\lvert A \cap B \rvert^2+\lvert A \cap B \rvert}$.
		\item If $(A,B)$ is a separation weakly separating two sets $Y,Z$ with minimum order, and $(C,D)$ is a separation given by a node in $S$ of breadth at least the breadth of $(A,B)$ weakly separating $Y,Z$, then $(A,B)=(C,D)$.
		\item There exists $S' \subseteq S$ satisfying (INC1)-(INC3) such that no node in $S'$ is a proper ancestor of another node in $S'$.
	\end{enumerate}
\end{lemma}

\begin{pf}
It is clear that $\lvert S \rvert \leq 2^{2\lvert A \cap B \rvert^2+\lvert A \cap B \rvert}$ by (INC3) since $(A,B)$ has thickness at most $\lvert A \cap B \rvert$.

Let $s \in S$ be such that $(A_s,B_s)$ is a separation of breadth at least the breadth of $(A,B)$.
By (INC1), the breadth of $(A_s,B_s)$ equals the breadth of $(A,B)$.
By (INC3), $s$ is the unique member of $S$.
By (INC2), $B_s=B$.
Since $(A,B)$ is a separation of minimum order weakly separating $Y,Z$, and $(A_s,B_s)$ weakly separates $Y,Z$ and the order of $(A,B)$ is equal to the order of $(A_s,B_s)$, $A_s \cap B_s = A \cap B$.
So $(A_s,B_s) = (A,B)$.

If $s_1,s_2$ are two nodes in $S$ such that $s_1$ is a proper ancestor of $s_2$, then $S-\{s_2\}$ also satisfies (INC1)-(INC3).
So we may repeat this process to remove nodes in $S$ that are proper descendants of other nodes in $S$, until no node in $S$ is a proper ancestor of another node in $S$.
\end{pf}

\begin{lemma} \label{not incorporated}
Let $(T,\X)$ be a rooted tree-decomposition of a graph $G$.
Let $w$ be a nonnegative integer and let $N = (w+1) \cdot 2^{2(w+1)(w+2)}+1$.
Let $t_1,t_2$ be two nodes of $T$, where $t_1$ is an ancestor of $t_2$.
Let $(A,B)$ be a separation of $G$ with minimum order such that $\downarrow t_1 \subseteq A$ and $\uparrow t_2 \subseteq B$. 
If $|A \cap B| \leq w+1$, and $\uparrow t_2$ contains $N$ vertices each of which cannot be separated from $\downarrow t_1$ by a separation given by a node of $T$ whose breadth is strictly less than the breadth of $(A,B)$, then either $(A,B)$ is given by a node of $T$, or $(A,B)$ is not incorporated.
\end{lemma}

\begin{pf}
Suppose that $(A,B)$ is not given by a node of $T$, but $(A,B)$ is incorporated.
Let $S$ be a subset of $V(T)$ satisfying (INC1)-(INC3) witnessing that $(A,B)$ is incorporated.
Since $(A,B)$ is a separation weakly separating $\downarrow t_1$ and $\uparrow t_2$ with minimum order, Statement 2 of Lemma \ref{incorporation basic} implies that every separation given by a node in $S$ has breadth strictly less than the breadth of $(A,B)$.

Note that $\uparrow t_2 \subseteq \uparrow t_1$ and $(A,B)$ is a separation of $G$ with minimum order such that $\downarrow t_1 \subseteq A$ and $\uparrow t_2 \subseteq B$.
So the order of $(\downarrow t_1,\uparrow t_1)$ is at least the order of $(A,B)$.
That is, $\lvert X_{t_1} \rvert \geq \lvert A \cap B \rvert$.

\medskip

\noindent{\bf Claim 1:} $\downarrow t_1 \not \subseteq X_t$ for every $t \in S$.

\noindent{\bf Proof of Claim 1:}
Suppose to the contrary that there exists a node $t \in S$ such that $\downarrow t_1 \subseteq  X_t$.
Then $\lvert X_{t_1} \rvert \leq \lvert \downarrow t_1 \rvert \leq \lvert X_t \rvert \leq \lvert A \cap B \rvert \leq \lvert X_{t_1} \rvert$.
So $\downarrow t_1 = X_t$, which has size $\lvert A \cap B \rvert$.
Since $\downarrow t_1 \subseteq A$ and $X_t \subseteq B$ (by (INC2)), $X_t=\downarrow t_1 = A \cap B$.
Furthermore, since $B_t \subseteq B$ and $A \cap B=X_t=A_t \cap B_t$, we know $A \subseteq A_t$.
So for each vertex in $A \cap B$ not pointed for $(A,B)$, it is in $A_t \cap B_t$ and not pointed for $(A_t,B_t)$.
Hence the breadth of $(A_t,B_t)$ is at least the breadth of $(A,B)$, a contradiction.
$\Box$

\medskip

\noindent{\bf Claim 2:} No node in $S$ is an ancestor of $t_1$.

\noindent{\bf Proof of Claim 2:}
Suppose to the contrary that some node $t \in S$ is an ancestor of $t_1$.
So $\uparrow t_1 \subseteq B_t$.
By (INC2), $B_t \subseteq B$.
So $X_{t_1} \subseteq \downarrow t_1 \cap B_t \subseteq A \cap B$.
But $\lvert X_{t_1} \rvert \geq \lvert A \cap B \rvert$.
This implies that $X_{t_1}=A \cap B$.
Hence if $A=\downarrow t_1$, then $B=\uparrow t_1$ and $(A,B)$ is given by $t_1$, a contradiction.
So $A - (\downarrow t_1) \neq \emptyset$.
But $A-(\downarrow t_1) \subseteq \uparrow t_1 \subseteq B_t \subseteq B$.
Therefore, $(A \cap B) - (\downarrow t_1) \neq \emptyset$.
But $A \cap B = X_{t_1} \subseteq \downarrow t_1$, a contradiction.
$\Box$

\medskip

Since $\uparrow t_2\subseteq B = \bigcup_{t \in S}B_t$, for each vertex $v \in \uparrow t_2$, there exists $t_v \in S$ such that $v \in B_{t_v}$.
Assume now that $v \in \uparrow t_2$ cannot be separated from $\downarrow t_1$ by a separation given by a node of $T$ of breadth strictly less than the breadth of $(A,B)$.
Since the breadth of each $(A_{t_v},B_{t_v})$ is less than the breadth of $(A,B)$, $(A_{t_v},B_{t_v})$ does not separate $\downarrow t_1$ and $v$.
By Claim 1, $\downarrow t_1 \not \subseteq X_{t_v}$.
So either $v \in X_{t_v}$, or $\downarrow t_1 \not \subseteq A_{t_v}$.
If $\downarrow t_1 \not \subseteq A_{t_v}$, then $t_v$ is not $t_1$ or a descendant of $t_1$, so $t_v$ is not an ancestor of $t_1$ nor a descendant of $t_1$ by Claim 2, and hence $v \in B_{t_v} \cap \uparrow t_2 \subseteq X_{t_v}$.
That is, $v \in X_{t_v} \subseteq \bigcup_{s \in S} X_s$.

Therefore, there are at least $N$ vertices contained in $\bigcup_{s \in S} X_s$.
For each $s \in S$, $\lvert X_s \rvert \leq \lvert A \cap B \rvert \leq w+1$ by (INC1).
So $\bigcup_{s \in S} X_s$ contains at most $(w+1)\lvert S \rvert \leq (w+1)2^{2(w+1)(w+2)}<N$ vertices by Statement 1 of Lemma \ref{incorporation basic}, a contradiction.
\end{pf}

\bigskip

Let $(T,\X)$ be a rooted tree-decomposition of a graph $G$.
For all nonnegative integers $i,j$ with $j \leq i$, let $b_{i,j}$ be the number of separations of $G$ of breadth $(i,j)$ incorporated in $(T,\X)$.
For each nonnegative integer $k$, let $b_k$ be the sequence $(b_{k,0},b_{k,1},...,b_{k,k})$.
The \defn{signature} of $(T,\X)$ is the sequence $b=(b_0,b_1,...,b_{\lvert V(G) \rvert})$.
If $b'=(b'_0,b'_1,...,b'_{\lvert V(G) \rvert})$ is the signature of another rooted tree-decomposition of $G$, where $b'_k=(b'_{k,0},b'_{k,1},...,b'_{k,k})$,
then we say that $b'$ is {\em greater} than $b$ if there exist integers $i,j\in\{0,1,\ldots,\lvert V(G) \rvert\}$ such that $j\le i$,
$b_0=b'_0$, $b_1=b'_1,\ldots,b_{i-1}=b'_{i-1}$, $b_{i,0}=b'_{i,0}$, $b_{i,1}=b'_{i,1},\ldots,b_{i,j-1}=b'_{i,j-1}$, and $b_{i,j}<b'_{i,j}$. 

\begin{lemma} \label{N-linked}
Let $w$ be a nonnegative integer and let $N=(w+1) \cdot 2^{2(w+1)(w+2)}+1$.
Let $(T,\X)$ be a rooted tree-decomposition of a graph $G$ of width at most $w$.
If $(T,\X)$ is not $N$-linked, then there exists a rooted tree-decomposition $(T^*,\X^*)$ of $G$ of width no more than the 
width of $(T,\X)$ but with signature  greater than the signature of $(T,\X)$.
\end{lemma}

\begin{pf}
Since $(T,\X)$ is not $N$-linked, there exist $t_1,t_2 \in V(T)$, where $t_1$ is a precursor of $t_2$, such that $\uparrow t_2$ contains at least $N$ vertices each of which cannot be separated from $\downarrow t_1$ by a separation of order less than $\lvert X_{t_1} \rvert$ given by a node of $T$, but there exists a separation $(A,B)$ of order less than $\lvert X_{t_1} \rvert$ with $\downarrow t_1 \subseteq A$ and $\uparrow t_2 \subseteq B$.
We may assume that the order of $(A,B)$ is as small as possible, so there exist $\lvert A \cap B \rvert$ disjoint paths in $G$ from $\downarrow t_1$ to $\uparrow t_2$ and hence from $X_{t_1}$ to $X_{t_2}$.
Subject to the minimality of the order of $(A,B)$, we further assume that $\sum_{v \in A \cap B} d_v$ is as small as possible, where $d_v$ is the minimum distance from a node whose bag contains $v$ to the path $t_1Tt_2$.

Since the order of $(A,B)$ is less than $\lvert X_{t_1} \rvert$, we know $X_{t_1} \subseteq A$, $X_{t_1} \not \subseteq B$, $X_{t_2} \subseteq B$, and $X_{t_2} \not \subseteq A$.

\medskip

\noindent{\bf Claim 1:} $(A,B)$ is not given by a node of $T$.

\noindent{\bf Proof of Claim 1:}
Suppose to the contrary that $(A,B)$ is given by a node $t$ of $T$.
Then $t \not\in t_1Tt_2$, since  $\lvert X_t \rvert \geq \lvert X_{t_1} \rvert$  as $t_1$ is a precursor of $t_2$.
If $t$ is an ancestor  of $t_1$, then $X_{t_1}\subseteq B_t=B$,
a contradiction.
So $t$ is not an ancestor  of $t_1$. Then $X_{t_2}\subseteq A_t=A$, a contradiction.
$\Box$

\medskip

Let us say that a node $t\in V(T)$ is a {\em side node} if $t$
is a descendant of $t_1$ but not an ancestor or a descendant of $t_2$.
By Lemma \ref{separating bags to subtrees applied}, we may further assume, by replacing the separation $(A,B)$, that for every side node $t$, 
if $X_t \subseteq A$ or $X_t \subseteq B$, then either $\uparrow t \subseteq A$ or $\uparrow t \subseteq B$.
Note that every separation of breadth less than the breadth of $(A,B)$ has order less than $\lvert X_{t_1} \rvert$, so by Lemma \ref{not incorporated} and Claim 1, $(A,B)$ is not incorporated.

Now we construct a new tree-decomposition.
For each vertex $v$ of $G$, let $t_v$ be a node of $T$ such that $v \in X_{t_v}$.
Let $T'$ be a copy of $T$ and let $T''$ be a copy of the maximal subtree of $T$ rooted at $t_1$.
For each node $t$ of $T$, we denote the copy of $t$ in $T'$ by $t'$; for each node $t$ that is a descendant of $t_1$, we denote the copy of $t$ in $T''$ by $t''$.
Define $T^*$ to be the rooted tree obtained from $T' \cup T''$ by adding a new node $t^*$ and directed edges $t_2't^*$ and $t^*t_1''$.
Define $X^*_{t^*} = A \cap B$.
For each $t' \in V(T')$, define $X^*_{t'} = (X_t \cap A) \cup \{v \in A \cap B: t \in t_vTt_2\}$; for each $t'' \in V(T'')$, define $X^*_{t''} = (X_t \cap B) \cup \{v \in A \cap B: t \in t_1Tt_v\}$.

\medskip

\noindent{\bf Claim 2:} $(T^*,\X^*)$ is a rooted tree-decomposition of $G$ of width at most the width of $(T,\X)$ such that for every side node $t$ of $T$, $((T^*,\X^*)\downarrow t', (T^*,\X^*)\uparrow t') = (A_t \cup B, B_t \cap A)$ and $((T^*,\X^*)\downarrow t'', (T^*,\X^*)\uparrow t'') = (A_t \cup A, B_t \cap B)$.

\noindent{\bf Proof of Claim 2:}
It is straightforward to check that $(T^*,\X^*)$ is a rooted tree-decomposition.
(Alternatively, this can be proved by the arguments in \cite[Page 542]{bd}, since $(T,\X) \downarrow t_1 \subseteq A$.)
Since there exist $\lvert A \cap B \rvert$ disjoint paths in $G$ from $X_{t_1}$ to $X_{t_2}$, it is not hard to show that the width of $(T^*,\X^*)$ is at most the width of $(T,\X)$.
It is clear that for every side node $t$ of $T$, $((T^*,\X^*)\downarrow t', (T^*,\X^*)\uparrow t') = (A_t \cup B, B_t \cap A)$ and $((T^*,\X^*)\downarrow t'', (T^*,\X^*)\uparrow t'') = (A_t \cup A, B_t \cap B)$.
$\Box$

\medskip

Note that $(A,B)$ is incorporated in $(T^*,\X^*)$ since $t^*$ gives the separation $(A,B)$.
Hence, to prove this lemma, by Claim 2, it suffices to show that every separation of $G$ of breadth no more than the breadth of $(A,B)$ incorporated in $(T,\X)$ is incorporated in $(T^*,\X^*)$.

Let $(C,D)$ be a separation of $G$ incorporated in $(T,\X)$ of breadth no more than the breadth of $(A,B)$.
It suffices to show that $(C,D)$ is incorporated in $(T^*,\X^*)$.

Let $S$ be a subset of $V(T)$ witnessing the incorporation of $(C,D)$. 
By Lemma \ref{incorporation basic}, we may assume that no node in $S$ is a proper ancestor of another node in $S$.
Since the order of $(C,D)$ is no more than $(A,B)$, every node in $S$ has bag size less than $\lvert X_{t_1} \rvert$.
Since $t_1$ is a precursor of $t_2$, no node in $S$ is in $t_1Tt_2$.

Define the following sets:
	\begin{itemize}
		\item $S_1 = \{s \in S: \hbox{$s$ belongs to the component of $T-t_1$ containing the root}\}$ and $S^*_1=\{s':s\in S_1\}$.
             	\item $S_2= \{s\in S$: $s$ is a side node and $\uparrow s \subseteq A\}$ and $S^*_2=\{s':s\in S_2\}$.
		\item $S_3 = \{s \in S$: $s$ is a side node,  $\uparrow s \not\subseteq A$ and $\uparrow s \subseteq B\}$ and $S^*_3=\{s'':s\in S_3\}$. 
		\item $S_4= \{s \in S$: $s$ is a descendant of $t_2\}$ and $S^*_4=\{s'':s\in S_4\}$. 
             	\item $S_5 = \{s \in S$: $s$ is a side node,  $\uparrow s \not\subseteq A$ and $\uparrow s\not \subseteq B\}$ and $S^*_5=\{s',s'':s\in S_5\}$. 
		\item $S^*=S_1^* \cup S_2^* \cup S_3^* \cup S_4^*\cup S_5^*$.
	\end{itemize}
Since $\downarrow t_1 \subseteq A$ and $\uparrow t_2 \subseteq B$, the following statements hold by Claim 2:
	\begin{itemize}
            \item If $s\in S_1\cup S_2$, then $(T^*,\X^*)\uparrow{s'}=(T,\X) \uparrow s$.
            \item If $s\in S_3\cup S_4$, then $(T^*,\X^*)\uparrow{s''}=(T,\X) \uparrow s$.
            \item If $s\in S_5$, then $(T^*,\X^*)\uparrow{s'}\cup(T^*,\X^*)\uparrow{s''}=(T,\X) \uparrow s$.
	    \item For $s\in S_5$ the order of $(B_s \cap A, A_s \cup B)$ equals $\lvert X^*_{s'} \rvert$.
	    \item For $s\in S_5$ the order of $(B_s \cap B, A_s \cup A)$ equals $\lvert X^*_{s''} \rvert$.
	\end{itemize}
Since no node in $S$ is in $t_1Tt_2$, we have $S=S_1 \cup S_2 \cup S_3 \cup S_4\cup S_5$, so $S^*$ satisfies (INC2), and (INC1) immediately follows from Claim 3 below.

\medskip

\noindent{\bf Claim 3:} For every side node $s \in S_5$, $\lvert X^*_{s'} \rvert < \lvert X_s \rvert$ and $\lvert X^*_{s''} \rvert < \lvert X_s \rvert$.

\noindent{\bf Proof of Claim 3:}
Since $(A,B)$ satisfies the conclusion of Lemma \ref{separating bags to subtrees applied}, $A\not\supseteq X_s\not\subseteq B$.
Suppose to the contrary that $\lvert X^*_{s'} \rvert \geq \lvert X_s \rvert$.
By  submodularity, the order of $(B_s \cup A,A_s \cap B)$ is at most the order of $(A,B)$, since the order of $(B_s \cap A, A_s \cup B)$ equals $\lvert X^*_{s'} \rvert$ and the order of $(A_s,B_s)$ equals $\lvert X_s \rvert$.
Furthermore, $\downarrow t_1 \subseteq A \subseteq B_s \cup A$ and $\uparrow t_2 \subseteq A_s \cap B$.
But $\sum_{v \in (B_s \cup A) \cap (A_s \cap B)} d_v < \sum_{v \in A \cap B} d_v$ unless $A \cap B \subseteq A_s$.
So $A \cap B \subseteq A_s$ by the minimality of $(A,B)$.
But it implies that $X^*_{s'} = X_s \cap A \subset X_s$, so $\lvert X^*_{s'} \rvert < \lvert X_s \rvert$, a contradiction.
This proves that $\lvert X^*_{s'} \rvert < \lvert X_s \rvert$.
Similarly, $\lvert X^*_{s''} \rvert < \lvert X_s \rvert$.
$\Box$

\medskip

Therefore $S^*$ satisfies (INC1) and (INC2). 
It remains to prove (INC3).
Note that for each $s \in S$, if both $s',s''$ are contained in $S^*$, then $s\in S_5$, so $\lvert X^*_{s'} \rvert < \lvert X_s \rvert$ and $\lvert X^*_{s''} \rvert < \lvert X_s \rvert$, and hence $2^{2|X^*_{s'}|^2+|X^*_{s'}|}+2^{2|X^*_{s''}|^2+|X^*_{s''}|} \leq 2^{2|X_s|^2+j}$, where $j$ is the integer such that the breadth of the separation given by $s$ is $(|X_s|,j)$.
So $S^*$ satisfies (INC3).
This proves that $(C,D)$ is incorporated in $(T^*,\X^*)$ and completes the proof of this lemma.
\end{pf}

\subsection{Reflection} \label{subsec:reflection}

Let $(A,B)$ be a separation of $G$, and let $U,V$ be two subsets of $V(G)$.
We say that $(A,B)$ \defn{strongly separates} $U,V$ if 
	\begin{itemize}
		\item $U \subseteq A$, $V \subseteq B$, and 
		\item every vertex in $A \cap B-(U \cap V)$ is pointed for $(A,B)$ and is not in $U$.
	\end{itemize}
Note that $(A,B)$ strongly separates $U,V$ does not imply that $(B,A)$ strongly separates $V,U$.

Let $G$ be a graph and let $Z \subseteq V(G)$.
Let $(A_2,B_2)$ be a separation of $G$ with $Z \subseteq A_2 \cap B_2$ such that every vertex in $A_2 \cap B_2-Z$ is pointed for $(A_2,B_2)$.
Let $W$ be a subset of $A_2 \cap B_2$ such that every vertex in $W$ is doubly pointed for $(A_2,B_2)$ and no vertex in $W$ is adjacent to some vertex in $A_2 \cap B_2-(Z \cup W)$.
We say that the separation $(A_1,B_1)$ of a graph $G$ is the \defn{reflection of $(A_2,B_2)$ with respect to $Z,W$} if 
	\begin{itemize}
		\item $A_1=A_2-(A_2 \cap B_2-(W \cup Z))$ and 
		\item $B_1=B_2 \cup \{u \in A_2-B_2: uv \in E(G)$ for some $v \in A_2 \cap B_2-(W \cup Z)\}$.
	\end{itemize}
(See Figure \ref{fig_reflection} for an example.)
Notice that $A_1 \cap B_1 \cap A_2 \cap B_2=W \cup Z$.
We say that a separation is a \defn{reflection of $(A_2,B_2)$ with respect to $Z$} if it is the reflection with respect to $Z,W'$ for some subset $W'$ of $A_2 \cap B_2$ such that every vertex in $W'$ is doubly pointed for $(A_2,B_2)$ and $W'$ is not adjacent to any vertex in $A_2 \cap B_2-(Z \cup W')$.

\begin{figure} 
	\begin{picture}(100,230) (-35,-30)

		\thicklines
	
		\put(280,150){\circle*{5}} \put(285,150){{$v_2$}}
		\put(280,150){\line(-1,0){200}}
		\put(280,180){\circle*{5}} \put(285,180){{$v_1$}}
		\put(280,180){\line(-1,0){200}}
		\multiput(265,220)(0,-5){49}{\line(0,-1){3}}
		\multiput(265,220)(5,0){10}{\line(1,0){3}}
		\multiput(313,220)(0,-5){49}{\line(0,-1){3}}
		\multiput(265,-23)(5,0){10}{\line(1,0){3}}
		\put(268,210){{$A_2 \cap B_2$}}
		\put(215,210){{$A_2-B_2$}}
		\put(320,210){{$B_2-A_2$}}

		\put(80,150){\circle*{5}} \put(65,150){{$u_2$}}
		\put(80,180){\circle*{5}} \put(65,180){{$u_1$}}
		\multiput(60,220)(5,0){12}{\line(1,0){3}}
		\multiput(60,130)(0,5){18}{\line(0,1){3}}
		\multiput(120,220)(0,-5){16}{\line(0,-1){3}}
		\multiput(120,140)(5,0){37}{\line(1,0){3}}
		\multiput(303,140)(0,-5){35}{\line(0,-1){3}}
		\multiput(60,130)(5,0){40}{\line(1,0){3}}
		\multiput(260,130)(0,-5){33}{\line(0,-1){3}}
		\multiput(260,-33)(5,0){9}{\line(1,0){3}}
		\put(68,210){{$A_1 \cap B_1$}}
		\put(10,210){{$A_1-B_1$}}
		\put(125,210){{$B_1-A_1$}}

		\multiput(280,0)(0,30){5}{\circle*{5}}
		\put(275,110){{$w_1$}}
		\put(230,120){\line(1,0){100}}
		\put(283,80){{$w_2$}}
		\put(230,90){\line(1,0){100}}
		\put(280,90){\line(0,-1){30}}
		\put(275,50){{$z_1$}}
		\put(280,60){\line(1,0){50}}
		\put(280,60){\line(-5,1){50}}
		\put(280,60){\line(-5,-1){50}}
		\put(275,20){{$z_2$}}
		\put(280,30){\line(-5,1){50}}
		\put(280,30){\line(-5,-1){50}}
		\put(280,30){\line(5,1){50}}
		\put(280,30){\line(5,-1){50}}
		\put(275,-10){{$z_3$}}
		\put(230,0){\line(1,0){100}}

	\end{picture}
	\caption{An example of the reflection of $(A_2,B_2)$ with respect to $Z,W$. Here $W=\{w_1,w_2\}$, $Z=\{z_1,z_2,z_3\}$, $A_2 \cap B_2 = \{v_1,v_2\} \cup W \cup Z$ and $A_1 \cap B_1 = \{u_1,u_2\} \cup W \cup Z$.} \label{fig_reflection}
\end{figure}

\begin{lemma} \label{reflection_basic}
Let $G$ be a graph and $Z$ a subset of $V(G)$.
Let $(A_2,B_2)$ be a separation of $G$ with $Z \subseteq A_2 \cap B_2$ such that every vertex in $A_2 \cap B_2-Z$ is pointed for $(A_2,B_2)$.
Let $W$ be a subset of $A_2 \cap B_2$ such that every vertex in $W$ is doubly pointed for $(A_2,B_2)$ and no vertex in $W$ is adjacent to some vertex in $A_2 \cap B_2-(Z \cup W)$.
Let $(A_1,B_1)$ be the reflection of $(A_2,B_2)$ with respect to $Z,W$.
If there exist $X,Y \subseteq V(G)$ with $X \cap Y=Z$ such that $(A_2,B_2)$ strongly separates $X,Y$, and there exist $\lvert A_2 \cap B_2 \rvert$ disjoint paths in $G$ from $X$ to $Y$, then every vertex in $A_1 \cap B_1-Z$ is anti-pointed for $(A_1,B_1)$, and every vertex in $W$ is doubly-pointed for both $(A_1,B_1)$ and $(A_2,B_2)$.
\end{lemma}

\begin{pf}
Since $(A_2,B_2)$ strongly separates $X,Y$, every vertex in $A_2 \cap B_2 - (X \cap Y)=A_2 \cap B_2-Z$ is not in $X$.
Since every vertex in $A_2 \cap B_2 - Z$ is pointed for $(A_2,B_2)$ and there exist $\lvert A_2 \cap B_2 \rvert$ disjoint paths in $G$ from $X$ to $Y$, the edges with one end in $A_2 \cap B_2-Z$ and with one end in $A_2-B_2$ form a matching.
So the order of $(A_1,B_1)$ equals the order of $(A_2,B_2)$, every vertex in $A_1 \cap B_1 - (W \cup Z)$ is anti-pointed for $(A_1,B_1)$.
Since every vertex in $W$ is not adjacent to any vertex in $A_2 \cap B_2-(Z \cup W)$, the set of edges with one end in $W$ and one end in $B_1-A_1$ equals the set of edges with one end in $W$ and one end in $B_2-A_2$. 
So every vertex in $W$ is doubly pointed for both $(A_1,B_1)$ and $(A_2,B_2)$.
\end{pf}

\bigskip

Notice that as long as $(A_2,B_2)$ strongly separates $X,Y$, there exists a reflection of $(A_2,B_2)$ with respect to $X \cap Y$, as we can take $W=\emptyset$.
Observe that if there exist $\lvert X \rvert$ disjoint paths from $X$ to $Y$, and $(C,D)$ is a separation of order $\lvert X \rvert$ with $X \subseteq C$, $Y \subseteq D$ and $Y-X \subseteq D-C$ such that every vertex in $C \cap D-(X \cap Y)$ is anti-pointed for $(C,D)$, then $(C,D)$ is a reflection of some pseudo-edge-cut modulo $X \cap Y$ with respect to $X \cap Y$.

\begin{lemma} \label{flipping components}
Let $(T,\X)$ be a rooted tree-decomposition of a graph $G$.
Let $t_0,t_1,t_2,t_3$ be nodes of $T$ such that $t_i$ is a precursor of $t_{i+1}$ for $i \in \{0,1,2\}$, and $X_{t_0} \cap X_{t_1}=X_{t_1} \cap X_{t_2}=X_{t_2} \cap X_{t_3} = \bigcap_{j=0}^3 X_j$. 
Assume that there exist $\lvert X_{t_0} \rvert$ disjoint paths $P_1,...,P_{\lvert X_{t_0} \rvert}$ in $G$ from $X_{t_0}$ to $X_{t_3}$, and every vertex in $X_{t_0} \cap X_{t_3}$ is coherent for $t_0,t_3$.
Let $(A_2,B_2)$ be a separation of $G$ of order $\lvert X_{t_0} \rvert$ strongly separating $\downarrow t_0$ and $\uparrow t_3$.
Let $W$ be a subset of $A_2 \cap B_2$ such that every vertex in $W$ is doubly pointed for $(A_2,B_2)$, and $W$ is not adjacent to any vertex in $A_2 \cap B_2 - (W \cup (X_{t_1} \cap X_{t_2}))$.
Let $(A_1,B_1)$ be the reflection of $(A_2,B_2)$ with respect to $X_{t_0} \cap X_{t_3}, W$.

If $A_2 \supseteq \downarrow t_0 \cup (\downarrow t_1 \cap (\bigcup_{i=1}^{\lvert X_{t_0} \rvert} V(P_i)))$ and $B_1 \supseteq \uparrow t_3 \cup (\uparrow t_2 \cap (\bigcup_{i=1}^{\lvert X_{t_0} \rvert} V(P_i)))$, then there exist separations $(A_1',B_1'),(A_2',B_2')$ such that the following hold.
	\begin{enumerate}
		\item $A_1' \cap B_1' = A_1 \cap B_1$ and $A_2' \cap B_2' = A_2 \cap B_2$.
		\item Every vertex in $A_2' \cap B_2'$ pointed (and doubly pointed, respectively) for $(A_2,B_2)$ is pointed (and doubly pointed, respectively) for $(A_2',B_2')$.  
			In particular, every vertex in $A_2' \cap B_2'-(X_{t_1} \cap X_{t_2})$ is pointed for $(A_2',B_2')$, and every vertex in $W$ is doubly pointed for $(A_2',B_2')$.
		\item $(A_1',B_1')$ is the reflection of $(A_2',B_2')$ with respect to $X_{t_1} \cap X_{t_2},W$. 
		\item $A_2' \supseteq \downarrow t_0 \cup (\downarrow t_1 \cap (\bigcup_{i=1}^{\lvert X_{t_0} \rvert} V(P_i)))$ and $B_1' \supseteq \uparrow t_3 \cup (\uparrow t_2 \cap (\bigcup_{i=1}^{\lvert X_{t_0} \rvert} V(P_i)))$.
		\item For every node $t$ that is a descendant of $t_0$ but not an ancestor of $t_3$, if $X_t \subseteq A_1'$ or $X_t \subseteq B_2'$, then $\uparrow t \subseteq A_1'$ or $\uparrow t \subseteq B_2'$.	
		\item Every vertex in $(A_1-A_1') \cup (A_1'-A_1) \cup (B_2-B_2') \cup (B_2'-B_2)$ is contained in some component of $G-(W \cup (X_{t_0} \cap X_{t_3}))$ disjoint from $X_{t_0} \cup X_{t_3}$.
		\item $(A_2',B_2')$ strongly separates $\downarrow t_0$ and $\uparrow t_3$, and every vertex in $A_1' \cap B_1' - (X_{t_1} \cap X_{t_2})$ is anti-pointed for $(A_1',B_1')$.
	\end{enumerate}
\end{lemma}

\begin{pf}
We first show that Statement 7 follows from other statements.

\medskip

\noindent{\bf Claim 1:} Statements 1-4 imply Statement 7.

\noindent{\bf Proof of Claim 1:} 
Statements 3 and 4 imply that $\downarrow t_0 \subseteq A_2'$ and $\uparrow t_3 \subseteq B_2'$.
Since $(A_2,B_2)$ strongly separates $\downarrow t_0$ and $\uparrow t_3$, Statements 1 and 2 imply that every vertex in $A_2' \cap B_2'-(\downarrow t_0 \cap \uparrow t_3) = A_2 \cap B_2 - (X_{t_0} \cap X_{t_3}) = A_2 \cap B_2 - (X_{t_1} \cap X_{t_2})$ is pointed for $(A_2',B_2')$ and is not in $\downarrow t_0$.
So $(A_2',B_2')$ strongly separates $\downarrow t_0$ and $\uparrow t_3$.
Since $\lvert A_2' \cap B_2' \rvert = \lvert A_2 \cap B_2 \rvert = \lvert X_{t_0} \rvert$ and there exist $\lvert X_{t_0} \rvert$ disjoint paths from $\downarrow t_0$ to $\uparrow t_3$, Lemma \ref{reflection_basic} implies that every vertex in $A_1' \cap B_1' - (X_{t_1} \cap X_{t_2})$ is anti-pointed for $(A_1',B_1')$.
So Statement 7 holds.
$\Box$

\medskip

We say a node of $T$ is a \defn{side node} if it is a descendant of $t_0$ but not an ancestor or a descendant of $t_3$.
We define a weakening of Statement 5.
	\begin{itemize}
		\item[5'.] For every side node $t$, if $X_t \subseteq A_1'$ or $X_t \subseteq B_2'$, then $\uparrow t \subseteq A_1'$ or $\uparrow t \subseteq B_2'$.	
	\end{itemize}

\medskip

\noindent{\bf Claim 2:} Statements 3, 4 and 5' imply Statement 5.

\noindent{\bf Proof of Claim 2:} 
Statements 3 and 4 imply that $\uparrow t_3 \subseteq B_2'$.
So for every descendant $t$ of $t_3$, $\uparrow t \subseteq \uparrow t_3 \subseteq B_2'$.
So Statement 5' implies Statement 5.
$\Box$

\medskip

Hence to prove this lemma, it suffices to prove Statements 1-4, 5' and 6.

We say that a side node $t$ is \defn{bad} for a pair of separations $(C_1,D_1)$ and $(C_2,D_2)$, where $(C_1,D_1)$ is a reflection of $(C_2,D_2)$ with respect to $X_{t_0} \cap X_{t_3}$, if $X_t \subseteq C_1$ or $X_t \subseteq D_2$, but $\uparrow t \not \subseteq C_1$ and $\uparrow t \not\subseteq D_2$.

Let $k$ be the number of bad side nodes for $(A_1,B_1),(A_2,B_2)$.
We shall prove this lemma by induction on $k$.

Clearly, when we choose $(A_i',B_i')$ to be $(A_i,B_i)$ for $i \in \{1,2\}$, Statements 1-4 and 6 hold.
And Statement 5' holds when $k=0$.

So we may assume that $k \geq 1$ and this lemma holds for all smaller $k$.

Let $Z=W \cup (X_{t_0} \cap X_{t_3})$, and let $M$ be the set consisting of the edges with one end in $A_1 \cap B_1-Z$ and one end in $A_2 \cap B_2-Z$. 
Since $W$ is not adjacent to any vertex in $A_2 \cap B_2 - Z$, and there exist $\lvert X_{t_1} \rvert$ disjoint paths in $G$ from $X_{t_0}$ to $X_{t_3}$, and every vertex in $A_2 \cap B_2-Z$ is pointed for $(A_2,B_2)$, $M$ is a matching.

Let $G_A$ (and $G_B$, respectively) be the union of the components of $G-(M \cup Z)$ intersecting $\downarrow t_0 \cup (\downarrow t_1 \cap (\bigcup_{i=1}^{\lvert X_{t_0} \rvert} V(P_i)))$ but disjoint from $\uparrow t_3$ (and intersecting $\uparrow t_3 \cup (\uparrow t_2 \cap (\bigcup_{i=1}^{\lvert X_{t_0} \rvert} V(P_i)))$ but disjoint from $\downarrow t_0$, respectively).
Note that $V(G_A) \subseteq A_1$ and $V(G_B) \subseteq B_2$ by the existence of $P_1,...,P_{\lvert X_{t_0} \rvert}$.
Since $P_1,P_2,...,P_{\lvert X_{t_0} \rvert}$ are disjoint paths in $G$ from $X_{t_0}$ to $X_{t_3}$, $A_1 \cap B_1-Z \subseteq G_A$ and $A_2 \cap B_2-Z \subseteq G_B$.
So every component of $G-(M \cup Z)$ disjoint from $G_A \cup G_B$ is a component of $G-Z$.

In addition, the set of nodes $t$ of $T$ with $X_t \cap V(G_A) \neq\emptyset$ (and $X_t \cap V(G_B) \neq \emptyset$, respectively) induces a connected subgraph of $T$.
Furthermore, since $V(G_A) \cap Z=V(G_B) \cap Z=\emptyset$, we have $V(G_A) \cap B_2=\emptyset$ and $V(G_B) \cap A_1=\emptyset$.

\medskip

\noindent{\bf Claim 3:} Let $C$ be a component of $G-(M \cup Z)$ disjoint from $G_A \cup G_B$.
	\begin{itemize}
		\item If $w$ is a vertex in $A_2 \cap B_2$ pointed for $(A_2,B_2)$, then $w$ is pointed for $(A_2 \cup V(C), B_2-V(C))$ and for $(A_2 - V(C), B_2 \cup V(C))$.
		\item If $w$ is a vertex in $A_1 \cap B_1$ anti-pointed for $(A_1,B_1)$, then $w$ is anti-pointed for $(A_1 \cup V(C), B_1-V(C))$ and for $(A_1 - V(C), B_1 \cup V(C))$.
	\end{itemize}

\noindent{\bf Proof of Claim 3:}
We assume that $w$ is pointed for $(A_2,B_2)$.
The case that $w$ is anti-pointed for $(A_1,B_1)$ can be proved analogously and we omit the proof.

Suppose that $w$ is pointed for $(A_2,B_2)$ but not pointed for $(A_2 \cup V(C), B_2-V(C))$.
Then $w$ is adjacent to some vertex in $C$.
Since $C$ is disjoint from $G_A \cup G_B$, it is a component of $G-Z$, so $w \in Z$.
If $w \not \in X_{t_0} \cap X_{t_3}$, then $w \in W$ is doubly pointed for $(A_2,B_2)$, and the neighbors of $w$ contained in $(A_2-B_2) \cup (B_2-A_2)$ belong to $G_A \cup G_B$ since there are $\lvert X_{t_0} \rvert$ disjoint paths between $X_{t_0}$ and $X_{t_3}$, so $w$ is not adjacent to vertices in $C$, a contradiction. 
So $w \in X_{t_0} \cap X_{t_3}$.

Since $C$ is disjoint from $G_B$, $C$ is disjoint from $\uparrow t_3$.
So $V(C) \subseteq A_{t_3}$.
Since $A_2 \cup V(C) \subseteq A_{t_3}$ and $w$ is not pointed for $(A_2 \cup V(C), B_2 - V(C))$, $w$ is incident with at least two edges whose other ends are in $(A_2 \cup V(C))-(B_2-V(C)) \subseteq A_{t_3}-B_{t_3}$.
So $w$ is not pointed for $(A_{t_3},B_{t_3})$.
Since $w$ is pointed for $(A_2,B_2)$, $w$ is pointed for $(A_{t_0},B_{t_0})$.
Since every vertex in $X_{t_0} \cap X_{t_3}$ is coherent for $t_0,t_3$, $w$ is pointed for $(A_{t_3},B_{t_3})$, a contradiction.

Therefore, if $w$ is pointed for $(A_2,B_2)$, then it is pointed for $(A_2 \cup V(C), B_2-V(C))$.
In addition, since $A_2-V(C) \subseteq A_2$, $w$ is pointed for $(A_2-V(C),B_2 \cup V(C))$ if $w$ is pointed for $(A_2,B_2)$.
This proves the claim.
$\Box$

\medskip

Let $s$ be a side node bad for $(A_1,B_1), (A_2,B_2)$.
We assume that $X_s \subseteq A_1$ and $\uparrow s-A_1 \neq \emptyset$.
(The case that $X_s \subseteq B_2$ and $\uparrow s-B_2 \neq \emptyset$ can be proved analogously, so we omit the proof of that case.)

Since $V(G_A) \subseteq A_1$, $\uparrow s-A_1$ is clearly disjoint from $G_A$.
Since $V(G_B) \cap A_1 = \emptyset$ and $X_s \subseteq A_1$, $X_s$ is disjoint from $G_B$.
Since the set of nodes of $T$ whose bags intersect $G_B$ induces a connected subgraph of $T$, and $t_3$ is in that set, we know that $\uparrow s-A_1$ is disjoint from $G_B$. 

Define $(A^*_1,B^*_1)$ (and $(A^*_2,B^*_2)$, respectively) to be the separation obtained from $(A_1,B_1)$ (and $(A_2,B_2)$, respectively) by removing all components of $G-(M \cup Z)$ intersecting $\uparrow s-A_1$ from $B_1$ (and $B_2$, respectively) and adding them into $A_1$ (and $A_2$, respectively).
Since $\uparrow s \subseteq A_1^*$, $s$ is not a bad side node for $(A^*_1,B^*_1),(A^*_2,B^*_2)$.

Let $C$ be the union of the components that we moved.

\medskip

\noindent{\bf Claim 4:} Statements 1-4, 6 and 7 hold if $(A_1',B_1')$ and $(A_2',B_2')$ are replaced by $(A^*_1,B^*_1)$ and $(A^*_2,B^*_2)$, respectively.

\noindent{\bf Proof of Claim 4:}
Clearly, Statement 1 holds if $(A_1',B_1')$ and $(A_2',B_2')$ are replaced by $(A^*_1,B^*_1)$ and $(A^*_2,B^*_2)$, respectively.

Since $\uparrow s-A_1$ is disjoint from $G_A \cup G_B$, so is $C$.
Hence by Claim 3, Statement 2 holds if $(A_1',B_1')$ and $(A_2',B_2')$ are replaced by $(A^*_1,B^*_1)$ and $(A^*_2,B^*_2)$, respectively.
In particular, every vertex in $A_2 \cap B_2-(X_{t_1} \cap X_{t_2})$ is pointed for $(A_2 \cup V(C), B_2-V(C))$, and every vertex in $A_1 \cap B_1-(X_{t_1} \cap X_{t_2})$ is anti-pointed for $(A_1 \cup V(C), B_1-V(C))$.
So $(A_1^*,B_1^*)$ is the reflection of $(A_2^*,B_2^*)$ with respect to $X_{t_0} \cap X_{t_3}, W$.
Hence Statement 3 holds if $(A_1',B_1')$ and $(A_2',B_2')$ are replaced by $(A^*_1,B^*_1)$ and $(A^*_2,B^*_2)$, respectively.

Furthermore, $A_2^* \supseteq A_2 \supseteq \downarrow t_0 \cup (\downarrow t_1 \cap (\bigcup_{i=1}^{\lvert X_{t_0} \rvert} V(P_i)))$.
Since $C$ is disjoint from $G_A \cup G_B$, $C$ is disjoint from $\uparrow t_3 \cup (\uparrow t_2 \cap (\bigcup_{i=1}^{\lvert X_{t_0} \rvert} V(P_i)))$. 
Therefore, $B_1^* \supseteq \uparrow t_3 \cup (\uparrow t_2 \cap (\bigcup_{i=1}^{\lvert X_{t_0} \rvert} V(P_i)))$.
Hence Statement 4 holds if $(A_1',B_1')$ and $(A_2',B_2')$ are replaced by $(A^*_1,B^*_1)$ and $(A^*_2,B^*_2)$, respectively.

By Claim 1, Statement 7 holds if $(A_1',B_1')$ and $(A_2',B_2')$ are replaced by $(A^*_1,B^*_1)$ and $(A^*_2,B^*_2)$, respectively.

Since $C$ is a union of components of $G-(M \cup Z)$ disjoint from $G_A \cup G_B$, $C$ is a union of components of $G-Z$ disjoint from $X_{t_0} \cup X_{t_3}$.
Since $(A_1-A_1^*) \cup (A_1^*-A_1) \cup (B_2-B_2^*) \cup (B_2^*-B_2) \subseteq V(C)$, every vertex in $(A_1-A_1^*) \cup (A_1^*-A_1) \cup (B_2-B_2^*) \cup (B_2^*-B_2)$ is contained in some component of $G-Z$ disjoint from $X_{t_0} \cup X_{t_3}$.
Hence Statement 6 holds if $(A_1',B_1')$ and $(A_2',B_2')$ are replaced by $(A^*_1,B^*_1)$ and $(A^*_2,B^*_2)$, respectively.
$\Box$

\medskip

\noindent{\bf Claim 5:} The number of bad side nodes for $(A^*_1,B^*_1),(A^*_2,B^*_2)$ is less than the number of bad side nodes for $(A_1,B_1),(A_2,B_2)$.

\noindent{\bf Proof of Claim 5:}
As $s$ is bad for $(A_1,B_1),(A_2,B_2)$ but not for $(A^*_1,B^*_1),(A^*_2,B^*_2)$, it suffices to prove that every bad side node for $(A^*_1,B^*_1),(A^*_2,B^*_2)$ is a bad side node for $(A_1,B_1),(A_2,B_2)$.

Suppose to the contrary that there exists a side node $s'$ that is bad for $(A^*_1,B^*_1),(A^*_2,B^*_2)$ but is not bad for $(A_1,B_1),(A_2,B_2)$.
Since $\uparrow s \subseteq A_1^*$, $s'$ is not a descendant of $s$. 

Suppose that $s'$ is an ancestor of $s$.
Since $s'$ is bad for $(A_1^*,B_1^*),(A_2^*,B_2^*)$, either $X_{s'} \subseteq A_1^*$ or $X_{s'} \subseteq B_2^*$.
Since $V(C) \subseteq \uparrow s-X_s$, either $X_{s'} \subseteq A_1$ or $X_{s'} \subseteq B_2$.
Since $s'$ is not bad for $(A_1,B_1),(A_2,B_2)$, either $\uparrow s' \subseteq A_1$ or $\uparrow s' \subseteq B_2$.
But $\uparrow s \not \subseteq A_1$, so $\uparrow s' \not \subseteq A_1$, and hence $\uparrow s' \subseteq B_2$.
Then $\uparrow s \subseteq B_2$, so $s$ is not bad for $(A_1,B_1),(A_2,B_2)$, a contradiction.

So $s'$ is not an ancestor or a descendant of $s$.
Since $V(C)$ is a subset of $\uparrow s$ and is disjoint from $X_s$, $V(C)$ is disjoint from $\uparrow s'$.
Therefore, $s'$ is bad for $(A^*_1,B^*_1),(A^*_2,B^*_2)$ if and only if it is bad for $(A_1,B_1),(A_2,B_2)$, a contradiction.
$\Box$

\medskip

By Claim 4, $(A_1^*,B_1^*)$ and $(A_2^*,B_2^*)$ satisfy the same condition for this lemma as $(A_1,B_1)$ and $(A_2,B_2)$.
By Claim 5, we can apply induction to obtain separations $(A_1',B_1')$ and $(A_2',B_2')$ satisfying Statements 1-7 (where all $(A_1,B_1)$ and $(A_2,B_2)$ are replaced by $(A_1^*,B_1^*)$ and $(A_2^*,B_2^*)$, respectively).
Note that Statements 3-5 and 7 are irrelevant with $(A_1,B_1)$ and $(A_2,B_2)$.
Since Statements 1,2 and 6 hold if $(A_1',B_1')$ and $(A_2',B_2')$ are replaced by $(A^*_1,B^*_1)$ and $(A^*_2,B^*_2)$, respectively, we know that those statements also hold (without replacing $(A_1,B_1)$ and $(A_2,B_2)$ by $(A_1^*,B_1^*)$ and $(A_2^*,B_2^*)$, respectively).
This proves the lemma.
\end{pf}

\subsection{Integration} \label{subsec:integration}

Recall that for a rooted tree-decomposition $(T,\X)$ of a graph $G$, and for every node $t$ of $T$, the separation of $G$ given by $t$ is the separation $(\downarrow t,\uparrow t)$ and denoted by $(A_t,B_t)$.

Let $N$ be a positive integer.
We say that a rooted tree-decomposition $(T,\X)$ of a graph $G$ is \defn{$N$-integrated} if for any nodes $t_0,t_1,t_2,t_3$ of $T$ and separation $(A,B)$ of $G$ satisfying statements (TE1)-(TE7) defined below, there exists a separation given by a node in $t_0Tt_3$ of breadth equal to the breadth of $(A,B)$:
	\begin{itemize}
		\item[(TE1)] $t_i$ is an ancestor of $t_{i+1}$ for $i \in \{0,1,2\}$.
		\item[(TE2)] $\lvert X_{t_i} \rvert = \lvert X_{t_0} \rvert$ for $i \in [3]$.
		\item[(TE3)] There exist $\lvert X_{t_0} \rvert$ disjoint paths in $G$ from $X_{t_0}$ to $X_{t_3}$. 
		\item[(TE4)] $X_{t_i} \cap X_{t_j} = \bigcap_{\ell=0}^3 X_{t_\ell}$ for all $i,j$ with $0 \leq i < j \leq 3$. 
		\item[(TE5)] Every vertex in $\bigcap_{i=0}^3 X_{t_i}$ is coherent for $t_0,t_3$. 
		\item[(TE6)] $(A,B)$ strongly separates $\downarrow t_1$ and $\uparrow t_2$ and has breadth $(\lvert X_{t_1} \rvert, k)$, where $k$ is the number of vertices in $X_{t_0} \cap X_{t_3}$ non-pointed for $(A_{t_3},B_{t_3})$.
		\item[(TE7)] $\uparrow t_3$ contains at least $N$ vertices of $G$ each of which cannot be separated from $\downarrow t_0$ by a separation of breadth less than the breadth of $(A,B)$ given by a node in $T$.
	\end{itemize}

The main result in this subsection is the following lemma (Lemma \ref{N-integrated}), which states that if a rooted tree-decomposition is not $N$-integrated for some large $N$, then there exists another rooted tree-decomposition of width no more than the previous one but having greater signature.
We sketch its proof here.
We first prove that there is a pair of separation $(A_1,B_1)$ and its reflection $(A_2,B_2)$ satisfying certain nice properties, and we choose an ``optimal'' such pair based on certain minimization conditions (i.e.\ conditions (a)-(h)).
We show that this pair of separations satisfy other properties (Claims 2 and 3), and then we insert this pair of separations into the tree-decomposition to obtain a new tree-decomposition.
The goal is to show that the new tree-decomposition is desired.
It is easy to show that the width of the new tree-decomposition is not larger than the old one.
We can show that some separation is incorporated in the new tree-decomposition but not in the old one.
So it suffices to fix a separation that is incorporated in the old tree-decomposition and show that it is also incorporated in the new tree-decomposition.
By suitably modifying the set witnessing the incorporation in the old tree-decomposition, we reduce the problem so that it suffices to show that for any fixed special kind of tree node $s$ in the old tree-decomposition, the separation given by $s$ has larger breadth than the separation given by the ``image tree node of $s$'' in the new tree-decomposition.
Suppose to the contrary that it is not true.
We can obtain some information about $s$ (Claim 5) and show that there is another pair of a separation and its reflection with nice properties (Claims 6-10).
By the minimality of $(A_1,B_1)$ and $(A_2,B_2)$, we obtain extra information about the separation given by $s$ (Claims 11 and 12).
Using this extra information, we can construct another pair of a separation and its reflection with nice properties (Claims 13-16).
Again, the minimality of $(A_1,B_1)$ and $(A_2,B_2)$ implies further extra information about the separation given by $s$ (Claim 17), allowing us to construct another pair of a separation and its reflection with nice properties (Claims 18-22) that contradict the minimality of $(A_1,B_1)$ and $(A_2,B_2)$ to complete the proof.

\begin{lemma} \label{N-integrated}
Let $w$ be a nonnegative integer and let $N=(w+1) \cdot 2^{2(w+1)(w+2)}+1$.
Let $(T,\X)$ be a rooted tree-decomposition of a graph $G$ of width at most $w$.
If $(T,\X)$ is not $N$-integrated, then there exists a rooted tree-decomposition $(T^*,\X^*)$ of $G$ of width no more than the width of $(T,\X)$ but the signature is greater than $(T,\X)$.
\end{lemma}

\begin{pf}
Let $t_0,t_1,t_2,t_3$ be nodes of $T$ and $(A,B)$ a separation witnessing that $(T,\X)$ is not $N$-integrated. 
That is, they satisfy (TE1)-(TE7), but no separation given by a node in $t_0Tt_3$ has breadth equal to the breadth of $(A,B)$.

Note that the combination of (TE3), (TE5) and (TE6) implies that the breadth of $(A,B)$ is minimum among all separations strongly separating $\downarrow t_0$ and $\uparrow t_3$.
Hence, every separation given by a node in $t_0Tt_3$ has breadth larger than the breadth of $(A,B)$.
In particular, $(A_{t_i},B_{t_i})$ have breadth greater than $(A,B)$ for all $i \in \{0,1,2,3\}$, and $\lvert A \cap B \rvert \leq \lvert X_{t_0} \rvert$.
Since there are $\lvert X_{t_0} \rvert$ disjoint paths in $G$ from $X_{t_0}$ to $X_{t_3}$ by (TE3), $\lvert A \cap B \rvert = \lvert X_{t_0} \rvert$.

If $X_{t_i}=X_{t_j}$ for some $0 \leq i < j \leq 3$, then $X_{t_0}=X_{t_1}=X_{t_2}=X_{t_3}$ by (TE2) and (TE4), so $(A_{t_1},B_{t_1})$ and $(A,B)$ have the same breadth by (TE3) and (TE6), a contradiction. 
Hence $X_{t_0},X_{t_1},X_{t_2},X_{t_3}$ are pairwise distinct.
So for each $i \in [0,2]$, $t_{i+1}$ is a precursor of $t_i$.

We say that a node of $T$ is a \defn{side node} if it is a descendant of $t_0$ but not an ancestor or a descendant of $t_3$.

Let $P_1,P_2,...,P_{\lvert X_{t_0} \rvert}$ be disjoint paths in $G$ from $X_{t_0}$ to $X_{t_3}$.
Define $G_A = \downarrow t_0 \cup (\downarrow t_1 \cap (\bigcup_{i=1}^{\lvert X_{t_0} \rvert} V(P_i)))$ and define $G_B = \uparrow t_3 \cup (\uparrow t_2 \cap (\bigcup_{i=1}^{\lvert X_{t_0} \rvert} V(P_i)))$.

\medskip

\noindent{\bf Claim 1:} There exist separations $(A_1,B_1)$ and $(A_2,B_2)$ satisfying the following. 
	\begin{itemize}
		\item[(a)] The breadth of $(A_2,B_2)$ equals the breadth of $(A,B)$.
		\item[(b)] $(A_1,B_1)$ is a reflection of $(A_2,B_2)$ with respect to $X_{t_0} \cap X_{t_3}$, and every vertex $v$ in $A_1 \cap B_1-(X_{t_0} \cap X_{t_3})$ is anti-pointed for $(A_1,B_1)$. 
		\item[(c)] $\downarrow t_0 \subseteq A_1$ and $\uparrow t_3 \subseteq B_2$.
		\item[(d)] $A_2 \supseteq G_A$ and $B_1 \supseteq G_B$. 
		\item[(e)] For every $t \in V(T)$ that is a descendant of $t_0$ but not an ancestor of $t_3$, if $X_t \subseteq A_1$ or $X_t \subseteq B_2$, then $\uparrow t \subseteq A_1$ or $\uparrow t \subseteq B_2$.
	\end{itemize}

\noindent{\bf Proof of Claim 1:}
If we take $(A_2,B_2)=(A,B)$ and take $(A_1,B_1)$ to be a reflection of $(A,B)$ with respect to $X_{t_0} \cap X_{t_3}$, then (a)-(d) hold by (TE4)-(TE6).
If we further apply Lemma \ref{flipping components} to the chosen $(A_1,B_1)$ and $(A_2,B_2)$, then we obtain a pair of separations $(A_1',B_1')$ and $(A_2',B_2')$ satisfying (b)-(e) such that the breadth of $(A_2',B_2')$ is not strictly greater than the breadth of $(A,B)$.
But $(A,B)$ is the separation with minimum breadth strongly separating $\downarrow t_0$ and $\uparrow t_3$, so the breadth of $(A_2',B_2')$ equals the breadth of $(A,B)$ and hence (a) is satisfied.
This shows the existence of the desired $(A_1,B_1),(A_2,B_2)$. 
$\Box$

\medskip

By Claim 1, there exist separations $(A_1,B_1)$ and $(A_2,B_2)$ satisfying (a)-(e) and the following.

	\begin{itemize}
		\item[(f)] Subject to (a)-(e), the number of side nodes $t$ such that either $A_1 \cap B_1 \cap A_2 \cap B_2-A_t \neq \emptyset$, or there exists an edge $uv$ of $G$ such that $u \in A_1 \cap B_1-((A_2 \cap B_2) \cup A_t)$ and $v \in A_2 \cap B_2 -((A_1 \cap B_1) \cup A_t)$ is as small as possible. 
	\end{itemize}

We say that a side node $t$ is \defn{(f)-bad} for separations $(C_1,D_1)$ and $(C_2,D_2)$ if $(C_1,D_1)$ and $(C_2,D_2)$ satisfy (a)-(e), but either $C_1 \cap D_1 \cap C_2 \cap D_2-A_t \neq \emptyset$ or there exists an edge $uv$ of $G$ such that $u \in C_1 \cap D_1-((C_2 \cap D_2) \cup A_t)$ and $v \in C_2 \cap D_2-((C_1 \cap D_1) \cup A_t)$.

That is, $(A_1,B_1)$ and $(A_2,B_2)$ are separations satisfying (a)-(e) such that the number of (f)-bad side nodes is minimum.

We say that a side node $t$ is \defn{(g)-bad} for separations $(C_1,D_1)$ and $(C_2,D_2)$ if $(C_1,D_1)$ and $(C_2,D_2)$ satisfy (a)-(f), $t$ is not (f)-bad for $(C_1,D_1)$ and $(C_2,D_2)$, and either 
			\begin{itemize}
				\item some vertex in $C_2 \cap D_2 \cap X_t - (C_1 \cap D_1)$ is adjacent to a vertex in $C_1 \cap D_1-((C_2 \cap D_2) \cup A_t)$ and a vertex $D_2-(C_2 \cup A_t \cup \bigcup_{i=1}^{\lvert X_{t_0} \rvert}V(P_i))$, or
				\item some vertex in $C_1 \cap D_1 \cap X_t - (C_2 \cap D_2)$ is adjacent to a vertex in $C_2 \cap D_2-((C_1 \cap D_1) \cup A_t)$ and a vertex in $C_1-(D_1 \cup A_t \cup \bigcup_{i=1}^{\lvert X_{t_0} \rvert}V(P_i))$.
			\end{itemize}

We further assume that $(A_1,B_1)$ and $(A_2,B_2)$ satisfy (a)-(f) and the following.
	\begin{itemize}
		\item[(g)] Subject to (a)-(f), $\sum_q (\lvert V(T) \rvert+1)^{d_q}$ is minimum, where the sum is over all (g)-bad side nodes $q$ for $(A_1,B_1)$ and $(A_2,B_2)$, and $d_q$ is the distance in $T$ from $q$ to $t_0Tt_3$.
	\end{itemize}

We say that a side node $t$ is \defn{(h)-bad} for separations $(C_1,D_1)$ and $(C_2,D_2)$ if $(C_1,D_1)$ and $(C_2,D_2)$ satisfy (a)-(g), $t$ is not (f)-bad for $(C_1,D_1)$ and $(C_2,D_2)$, and either $C_1 \cap D_1-A_t \neq \emptyset$ or $C_2 \cap D_2 -A_t \neq \emptyset$.

We further assume that $(A_1,B_1)$ and $(A_2,B_2)$ satisfy (a)-(g) and the following.
	\begin{itemize}
		\item[(h)] Subject to (a)-(g), the number of (h)-bad side nodes for $(A_1,B_1)$ and $(A_2,B_2)$ is as small as possible. 
	\end{itemize}

That is, $(A_1,B_1)$ and $(A_2,B_2)$ are separations satisfying (a)-(g), and subject to those, the number of (h)-bad side nodes is minimum.

\medskip

\noindent{\bf Claim 2:} $(A_2,B_2)$ is not given by a node of $t$.

\noindent{\bf Proof of Claim 2:}
Suppose that $(A_2,B_2)$ is given by a node $t$ of $T$.
Since the breath of the separation given by any node in $t_0Tt_3$ is greater than $(A,B)$, $(A_2,B_2)$ is not given by a node in $t_0Tt_3$ by (a).
So $t \not \in t_0Tt_3$.

Suppose that $t$ is a descendant of $t_3$, then $\downarrow t_3 \subseteq A_t=A_2$.
By (c), $\uparrow t_3 \subseteq B_2$, so $A_{t_3} \subseteq A_2$ and $B_{t_3} \subseteq B_2$.
Since the breadth of $(A_2,B_2)$ is at most the breadth of $(A_{t_3},B_{t_3})$, $(A_2,B_2)=(A_{t_3},B_{t_3})$, a contradiction.

So $t$ is not a descendant of $t_3$.
Suppose that $t$ is an ancestor of $t_0$.
So $X_{t_0} \cup X_{t_1} \subseteq B_2$.
By (d), $A_2 \supseteq G_A \supseteq X_{t_0} \cup X_{t_1}$.
So $A_2 \cap B_2 \supseteq X_{t_0} \cup X_{t_1}$.
But $X_{t_0} \neq X_{t_1}$, so $\lvert X_{t_0} \rvert = \lvert A_2 \cap B_2 \rvert \geq \lvert X_{t_0} \cup X_{t_1} \rvert > \lvert X_{t_0} \rvert$, a contradiction.
So $t$ is not an ancestor of $t_0$.
Since $t \not \in t_0Tt_3$ and $t$ is not a descendant of $t_3$ or an ancestor of $t_0$, $\uparrow t_3 \subseteq A_t=A_2$.
By (c), $\uparrow t_3 \subseteq A_2 \cap B_2$.
However, $\uparrow t_3$ contains at least $N$ vertices by (TE7), and $N>w+1$, a contradiction.
Therefore, $(A_2,B_2)$ is not given by a node of $t$.
$\Box$

\medskip

\noindent{\bf Claim 3:} $(A_2,B_2)$ is not incorporated.

\noindent{\bf Proof of Claim 3:}
Since there are $\lvert X_{t_0} \rvert = \lvert A_2 \cap B_2 \rvert$ disjoint paths from $X_{t_0}$ to $X_{t_3}$, $(A_2,B_2)$ is a separation of $G$ with minimum order such that $\downarrow t_0 \subseteq A_2$ and $\uparrow t_3 \subseteq B_2$ (by (c) and (d)).
Furthermore, the order of $(A_2,B_2)$ is $\lvert A \cap B \rvert=\lvert X_{t_1} \rvert \leq w+1$, and $\uparrow t_3$ contains at least $N$ vertices of $G$ each of which cannot be separated from $\downarrow t_0$ by a separation of breadth less than the breadth of $(A_2,B_2)$ given by a node of $T$ (by (TE7) and (a)).
By Lemma \ref{not incorporated} and Claim 2, $(A_2,B_2)$ is not incorporated.
$\Box$

\medskip

For each vertex $v$ of $G$, define $t_v$ to be a node of $T$ with $v \in X_{t_v}$.

Let $T'$ be a copy of $T$ and let $T''$ be a copy of the maximal subtree of $T$ rooted at $t_0$.
For each node $t$ of $T$, let $t'$ be the copy of $t$ in $T'$; for each node $t$ that is a descendant of $t_0$, let $t''$ be the copy of $t$ in $T''$.

Define $T^*$ to be the rooted tree obtained from $T' \cup T''$ by adding a path $q_0q_1q_2...q_{k+1}$ and new edges $t_3'q_0,q_{k+1}t_0''$, where $k=\lvert A_1 \cap B_1 - (A_2 \cap B_2) \rvert$.
We define the following.
	\begin{itemize}
		\item For each node $t'$ of $T'$, define $X^*_{t'} = (X_t \cap A_1) \cup \{v \in A_1 \cap B_1: t \in t_vTt_3\}$.
		\item For each node $t''$ of $T''$, define $X^*_{t''}=(X_t \cap B_2) \cup \{v \in A_2 \cap B_2: t \in t_vTt_0\}$.
		\item Define $X^*_{q_0} = A_1 \cap B_1$ and $X^*_{q_{k+1}}=A_2 \cap B_2$.
		\item Let $u_1,u_2,...,u_k$ be the vertices in $A_1 \cap B_1-(A_2 \cap B_2)$, and let $v_1,v_2,...,v_k$ be the vertices in $A_2 \cap B_2 - (A_1 \cap B_1)$ such that $u_i$ is adjacent to $v_i$ for each $i \in [k]$.
			For each $i \in [k]$, define $X^*_{q_i}=\{v_1,v_2,...,v_i,u_i,u_{i+1},...,u_k\} \cup (A_1 \cap B_1 \cap A_2 \cap B_2)$. 
	\end{itemize}
Then $(T^*,\X^*)$ is a rooted tree-decomposition of $G$, which can be proved straightforwardly as in the proof of Lemma \ref{N-linked}.
Note that every node in $t_0Tt_3$ has bag size at least $\lvert A \cap B \rvert$ in $(T,\X)$, and some node in $t_0Tt_3$ has bag size at least $\lvert A \cap B \rvert+1$ in $(T,\X)$ since there exist $\lvert X_{t_0} \rvert=\lvert X_{t_3} \rvert$ disjoint paths in $G$ between distinct sets $X_{t_0}$ and $X_{t_3}$.
Since $|A_1 \cap B_1|=|A_2 \cap B_2|=|A \cap B|$ and there are $|A \cap B|$ disjoint paths between $X_{t_0}$ and $X_{t_3}$, it is straightforward to show that the width of $(T^*,\X^*)$ is at most the width of $(T,\X)$.

Observe that $(A_2,B_2)$ is incorporated in $(T^*,\X^*)$ with witness set $\{q_{k+1}\}$.
Recall that $(A_2,B_2)$ is not incorporated in $(T,\X)$ by Claim 3.

Hence, to prove this lemma, it suffices to prove that every separation of $G$ of breadth at most $(A_2,B_2)$ incorporated in $(T,\X)$ is incorporated in $(T^*,\X^*)$.

Let $(C,D)$ be a separation of breadth at most the breadth of $(A_2,B_2)$ incorporated in $(T,\X)$, and let $S$ be a witness set for the incorporation of $(C,D)$ in $(T,\X)$.
By Lemma \ref{incorporation basic}, we may assume that no node in $S$ is an ancestor of another node in $S$.
By (INC1), no node in $S$ is in the path $t_0Tt_3$.

Define the following sets.
	\begin{itemize}
		\item $S_1 = \{s \in S: \hbox{$s$ belongs to the component of $T-t_0$ containing the root}\}$ and $S^*_1=\{s':s\in S_1\}$.
             	\item $S_2= \{s\in S$: $s$ is a side node and $\uparrow s \subseteq A_1\}$ and $S^*_2=\{s':s\in S_2\}$.
		\item $S_3 = \{s \in S$: $s$ is a side node, $\uparrow s \not\subseteq A_1$ and  $\uparrow s \subseteq B_2\}$ and $S^*_3=\{s'':s\in S_3\}$. 
		\item $S_4= \{s \in S$: $s$ is a descendant of $t_3\}$ and $S^*_4=\{s'':s\in S_4\}$. 
             	\item $S_5 = \{s \in S$: $s$ is a side node, $\uparrow s \not\subseteq A_1$ and $\uparrow s\not \subseteq B_2\}$ and $S^*_5=\{s',s'':s\in S_5\}$. 	
		\item $S^*=S_1^* \cup S_2^* \cup S_3^* \cup S_4^*\cup S_5^*$.
	\end{itemize}
To prove this lemma, it suffices to prove that $(C,D)$ is incorporated in $(T^*,\X^*)$, and $S^*$ is a witness set for $(C,D)$ being incorporated in $(T^*,\X^*)$.

By (c), (e) and the fact that $A_1 \cup B_2 = V(G)$, we have the following.
	\begin{itemize}
            \item If $s\in S_1\cup S_2$, then $(T^*,\X^*)\uparrow{s'}=(T,\X) \uparrow s$.
            \item If $s\in S_3\cup S_4$, then $(T^*,\X^*)\uparrow{s''}=(T,\X) \uparrow s$.
            \item If $s\in S_5$, then $(T^*,\X^*)\uparrow{s'}\cup(T^*,\X^*)\uparrow{s''}=(T,\X) \uparrow s$.
	    \item If $s \in S_5$, then $(A_{s'},B_{s'})=(A_s \cup B_1, B_s \cap A_1)$ in $(T^*,\X^*)$, and $(A_{s''},B_{s''})=(A_2 \cup A_s, B_2 \cap B_s)$ in $(T^*,\X^*)$.
	\end{itemize}
Since no node in $S$ is in $t_0Tt_3$, it follows that $S=S_1 \cup S_2 \cup S_3 \cup S_4\cup S_5$, and hence  $S^*$ satisfies (INC2).

Now we prove that $S^*$ satisfies (INC1) and (INC3).

Let $s$ be an element of $S$.
If $s' \in S_1^* \cup S_2^*$ (or $s'' \in S_3^* \cup S_4^*$, respectively), then it is clear that the separation given by $s'$ (or $s''$, respectively) has breadth at most the breadth of the separation given by $s$.

Assume that $s \in S$ with $s',s'' \in S^*$.
So $s \in S_5$ and $s',s'' \in S_5^*$.
Hence $s$ is a side node.

To prove this lemma, it suffices to show that the breadth of $(A_{s'},B_{s'})$ and the breadth of $(A_{s''},B_{s''})$ are less than the breadth of $(A_s,B_s)$.
We will only prove that the breadth of $(A_{s'},B_{s'})$ is less than the breadth of $(A_s,B_s)$ since the case for $(A_{s''},B_{s''})$ can be proved analogously.

Since there exist $\lvert X_{t_0} \rvert$ disjoint paths in $G$ from $X_{t_0}$ to $X_{t_3}$, $\max\{\lvert X^*_{s'} \rvert,\lvert X^*_{s''} \rvert\} \leq \lvert X_s \rvert$.

\medskip

\noindent{\bf Claim 4:} The breadth of $(A_{s'},B_{s'})$ is at most the breadth of $(A_s,B_s)$.

\noindent{\bf Proof of Claim 4:}
Suppose to the contrary that the breadth of $(A_{s'},B_{s'})$ is greater than the breadth of $(A_s,B_s)$.
In particular, $\lvert X^*_{s'} \rvert = \lvert X_s \rvert$. 
By the submodularity, the order of $(A_s \cap B_1, B_s \cup A_1)$ is at most the order of $(A_1,B_1)$.
By (c), $(A_s \cap B_1, B_s \cup A_1)$ separates $\downarrow t_0$ and $\uparrow t_3$.
Since there exist $\lvert X_{t_0} \rvert$ disjoint paths from $X_{t_0}$ to $X_{t_3}$, the order of $(A_s \cap B_1, B_s \cup A_1)$ equals the order of $(A_1,B_1)$.
So every vertex in $X_s \cap A_1 \cap B_1$ is adjacent to a vertex in $(A_s \cap B_1) - (B_s \cup A_1)$; otherwise $X_{t_0}$ and $X_{t_3}$ can be separated by a separation of order less than the order of $(A,B)$.

Let $v \in X^*_{s'}$.
If $v \in X^*_{s'}-B_1$, then $v \in X_s-B_1$.
If $v \in X^*_{s'}-B_1$ and $v$ is pointed for $(A_s,B_s)$, then it is pointed for $(A_s \cup B_1, B_s \cap A_1)$.
If $v \in X^*_{s'}-A_s$, then $v \in (A_1 \cap B_1) -(X_{t_0} \cap X_{t_3})$, so $v$ is anti-pointed for $(A_1,B_1)$ by (b) and is pointed for $(A_s \cup B_1, B_s \cap A_1)$.

Since the breadth of $(A_{s'},B_{s'})$ is greater than the breadth of $(A_s,B_s)$, we may assume that $v \in X^*_{s'} \cap A_s \cap B_1$ and $v$ is pointed for $(A_s,B_s)$ but not pointed for $(A_{s'},B_{s'})$.
Then $v$ is adjacent to a vertex in $A_s \cap B_1 - (B_s \cup A_1)$.
So $v$ has no neighbor in $A_s \cap A_1-B_s$.

If $v \not \in X_{t_0} \cap X_{t_3}$, then $v$ has no neighbor in $B_1 \cap B_s$ since $v$ is anti-pointed for $(A_1,B_1)$ by (b), so $v$ is pointed for $(A_{s'},B_{s'})$, a contradiction.
So $v \in X_{t_0} \cap X_{t_3}$.
Since $v$ is pointed for $(A_s,B_s)$, it is anti-pointed for $(A_{t_3},B_{t_3})$.
So it is anti-pointed for $(A_{t_0},B_{t_0})$ as $v$ is coherent for $t_0,t_3$ by (TE5).
Hence $v$ is pointed for $(A_{s'},B_{s'})$, a contradiction. 
$\Box$

\medskip

Suppose toward a contradiction that the breadth of $(A_{s'},B_{s'})$ is at least the breadth of $(A_s,B_s)$.
By Claim 4, $(A_s,B_s)$ and $(A_{s'}, B_{s'})$ have the same breadth.

Then every vertex in $X_s-A_1$ is pointed for $(A_s,B_s)$, since every vertex in $X^*_{s'}-X_s$ is pointed for $(A_{s'},B_{s'})$.

\medskip

\noindent{\bf Claim 5:} $s$ is not a descendant of $t_2$.

\noindent{\bf Proof of Claim 5:}
Suppose to the contrary that $s$ is a descendant of $t_2$.
Since $s \in S_5$, $X_s-A_1 \neq \emptyset$ by (e).
Since $(A_s,B_s)$ and $(A_{s'},B_{s'})$ have the same breadth, $\lvert X_s \rvert = \lvert X^*_{s'} \rvert$.
So one of $P_1,...,P_{\lvert X_{t_0} \rvert}$, say $Q$, passes through a vertex $x$ in $X_{t_0}-X_{t_3}$, a vertex $y$ in $A_1 \cap B_1 \cap (\uparrow s-X_s)$, a vertex $z$ in $X_s-A_1$, and a vertex $b \in X_{t_3}-X_{t_0}$ in the order listed. 
Since $s$ is a descendant of $t_0$ and $y \in \uparrow s-X_s$, the subpath of $Q$ between $x$ and $y$ contains a vertex $a$ in $X_s$.
Since $s$ is a descendant of $t_2$, the subpath of $Q$ between $x$ and $a$ contains a vertex $a'$ in $X_{t_2}$.
Note that $a'$ is the unique vertex belonging to $V(Q) \cap X_{t_2}$, since $P_1,P_2,...,P_{\lvert X_{t_0} \rvert}$ are $\lvert X_{t_2} \rvert$ disjoint paths intersecting $X_{t_2}$.
Note that $y \in \uparrow s-X_s$, so $y \not \in X_{t_2}$ and $y \neq a'$.

Since $t_0 \in t_0Tt_3$, the breadth of $(A,B)$ is smaller than the breadth of $(A_{t_0},B_{t_0})$.
So the order of $(A,B)$ is at most $\lvert X_{t_0} \rvert$.
By (a) and (b), the order of $(A_1,B_1)$ equals $\lvert X_{t_0} \rvert$.
So $y$ is the unique vertex belonging to $V(Q) \cap A_1 \cap B_1$.
Since $B_1 \supseteq G_B \supseteq X_{t_2}$ by (d), $a' \in B_1$.
Since $y \in A_1 \cap B_1-X_{t_2}$, $Q$ passes through $x,y,a'$ in the order listed, a contradiction. 
$\Box$

\medskip

Define $(A_1',B_1') = (A_1 \cup B_s, B_1 \cap A_s)$.

\medskip

\noindent{\bf Claim 6:} Every vertex in $A_1' \cap B_1'-(X_{t_0} \cup X_{t_3})$ is anti-pointed for $(A_1',B_1')$.

\noindent{\bf Proof of Claim 6:}
Let $v$ be a vertex is in $A_1' \cap B_1'-(X_{t_0} \cap X_{t_3})$.
If $v \in A_1$, then $v \in A_1 \cap B_1 \cap A_s-(X_{t_0} \cap X_{t_3})$, so $v$ is anti-pointed for $(A_1,B_1)$ by (b), and hence is anti-pointed by $(A_1',B_1')$ since $B_1'-A_1' \subseteq B_1-A_1$.
If $v \in B_1-A_1$, then $v \in X_s-A_1$, so $v$ is pointed for $(A_s,B_s)$, and hence $v$ is anti-pointed for $(A_1',B_1')$ since $B_1'-A_1' \subseteq A_s-B_s$.
Hence every vertex in $A_1' \cap B_1'-(X_{t_0} \cup X_{t_3})$ is anti-pointed for $(A_1',B_1')$.
$\Box$

\medskip

\noindent{\bf Claim 7:} For every vertex $u \in A_1' \cap B_1'-(X_{t_0} \cup X_{t_3})$, there exists a unique neighbor $u'$ of $u$ in $B_1'-A_1'$.
Furthermore, some member of $\{P_1,P_2,...,P_{\lvert X_{t_0} \rvert}\}$ contains both $u$ and $u'$, and $u' \in \downarrow t_3$.
In addition, if $u_1,u_2$ are distinct vertices in $A_1' \cap B_1'-(X_{t_0} \cup X_{t_3})$, then $u_1' \neq u_2'$.

\noindent{\bf Proof of Claim 7:}
Since $\lvert X^*_{s'} \rvert = \lvert X_s \rvert$, $\lvert A_1' \cap B_1' \rvert \leq \lvert A_1 \cap B_1 \rvert$ by the submodularity.
Since $X_{t_0} \subseteq A_1 \subseteq A_1'$ and $X_{t_3} \subseteq B_2 \cap A_s \subseteq B_1 \cap A_s = B_1'$, the existence of $P_1,P_2,...,P_{\lvert X_{t_0} \rvert}$ implies that $\lvert A_1' \cap B_1' \rvert \geq \lvert X_{t_0} \rvert \geq \lvert A \cap B \rvert = \lvert A_1 \cap B_1 \rvert$.
So $\lvert A_1' \cap B_1' \rvert = \lvert A \cap B \rvert = \lvert X_{t_0} \rvert$.
Hence for every vertex $u \in A_1' \cap B_1'$, there exists a unique member $P_u$ of $\{P_1,P_2,...,P_{\lvert X_{t_0} \rvert}\}$ containing $u$.
Since $B_1 \supseteq G_B$ (by (d)) and $\uparrow t_2 \subseteq A_s$ (by Claim 5), we have $\uparrow t_2 \cap (\bigcup_{i=1}^{\lvert X_{t_0} \rvert}V(P_i)) \subseteq B_1 \cap A_s=B_1'$.
So for each $u \in A_1' \cap B_1'-(X_{t_0} \cap X_{t_3})$, the subpath of $P_u$ between $u$ and $X_{t_3} \cap V(P_u)$ contains at least two vertices.
So for every $u \in A_1' \cap B_1' - (X_{t_0} \cap X_{t_3})$, there exists a neighbor $u'$ of $u$ contained in the subpath of $P_u$ between $u$ and $X_{t_3} \cap V(P_u)$.
Note that $u' \in \downarrow t_3$.
Since $\lvert X_{t_0} \rvert = \lvert A_1' \cap B_1' \rvert$, $u' \in B_1'-A_1'$.
By Claim 6, $u'$ is the unique neighbor of $u$ in $B_1'-A_1'$.
Since $P_1,P_2,...,P_{\lvert X_{t_0} \rvert}$ are disjoint, for distinct $u_1,u_2 \in A_1' \cap B_1'-(X_{t_0} \cap X_{t_3})$, $u_1' \neq u_2'$.
$\Box$

\medskip

\noindent{\bf Claim 8:} There exists a separation $(A_2',B_2')$ such that $(A_1',B_1')$ is the reflection of $(A_2',B_2')$ with respect to $X_{t_0} \cap X_{t_3}, \emptyset$.

\noindent{\bf Proof of Claim 8:}
For each $u \in A_1' \cap B_1'-(X_{t_0} \cap X_{t_3})$, by Claim 7, there exists a unique neighbor $u'$ of $u$ in $B_1'-A_1'$.
Let $R = A_1' \cap B_1'-(X_{t_0} \cap X_{t_3})$.
Let $R' =\{u': u \in R\}$.
Define $A_2' = A_1' \cup R'$ and define $B_2' = B_1' - R$.
By Claim 7, every vertex in $A_2' \cap B_2'-(X_{t_0} \cap X_{t_3})=R'$ is pointed for $(A_2',B_2')$.
So $(A_1',B_1')$ is the reflection of $(A_2',B_2')$ with respect to $X_{t_0} \cap X_{t_3},\emptyset$.
$\Box$

\medskip

Define $(A_2',B_2')$ to be the separation such that $(A_1',B_1')$ is the reflection of $(A_2',B_2')$ with respect to $X_{t_0} \cap X_{t_3}, \emptyset$. 
The existence of $(A_2',B_2')$ follows from Claim 8. 

\medskip

\noindent{\bf Claim 9:} There exist separations $(A_1^*,B_1^*)$ and $(A^*_2,B^*_2)$ satisfying (a)-(e) such that $A_1^* \cap B_1^* = A_1' \cap B_1'$, $A_2^* \cap B_2^* = A_2' \cap B_2'$, and every vertex in $(A_1^*-A_1') \cup (A_1'-A_1^*) \cup (B_2^*-B_2') \cup (B_2'-B_2^*)$ is contained in some component of $G-(X_{t_0} \cap X_{t_3})$ disjoint from $X_{t_0} \cup X_{t_3} \cup (A_1' \cap B_1') \cup (A_2' \cap B_2')$.

\noindent{\bf Proof of Claim 9:}
Since $s \not \in t_0Tt_3$ and $s$ is not a descendant of $t_2$ by Claim 5, we know $A_s \supseteq \uparrow t_2 \supseteq G_B$.
So $B_1' \supseteq G_B$ by (d). 
Hence $B_2' \supseteq \uparrow t_3$.
Since $A_1' \supseteq A_1$ and $(A_1',B_1')$ is the reflection of $(A_2',B_2')$ with respect to $X_{t_0} \cap X_{t_3}, \emptyset$, we have $A_2' \supseteq A_2$.
Since $A_2 \supseteq G_A$ by (d), $A_2' \supseteq G_A$. 
So $(A_1',B_1')$ and $(A_2',B_2')$ satisfy (b)-(d).

Since $(A_1',B_1')$ and $(A_2',B_2')$ satisfy (d), $(A_2',B_2')$ strongly separates $\downarrow t_0$ and $\uparrow t_3$.
Recall that the breadth of $(A,B)$ is minimum among all separations strongly separating $\downarrow t_0$ and $\uparrow t_3$.
So the breadth of $(A_2',B_2')$ is at least the breadth of $(A,B)$.
Since every vertex in $A_2' \cap B_2'-(X_{t_0} \cap X_{t_3})$ is pointed for $(A_2',B_2')$, and every vertex in $X_{t_0} \cap X_{t_3}$ is coherent for $t_0,t_3$, the breadth of $(A_2',B_2')$ is at most the breadth of $(A,B)$.
Hence $(A_1',B_1')$ and $(A_2',B_2')$ satisfy (a)-(d).

By Lemma \ref{flipping components}, there exist separations $(A_1^*,B_1^*)$ and $(A_2^*,B_2^*)$ satisfying (a)-(e) such that $A_1^* \cap B_1^* = A_1' \cap B_1'$, $A_2^* \cap B_2^* = A_2' \cap B_2'$, and every vertex in $(A_1^*-A_1') \cup (A_1'-A_1^*) \cup (B_2^*-B_2') \cup (B_2'-B_2^*)$ is contained in some component of $G-(X_{t_0} \cap X_{t_3})$ disjoint from $X_{t_0} \cup X_{t_3}$ since $(A_1',B_1')$ is the reflection of $(A_2',B_2')$ with respect to $X_{t_0} \cap X_{t_3},\emptyset$. 
By the existence of $P_1,P_2,...,P_{\lvert X_{t_0} \rvert}$, every component of $G-(X_{t_0} \cap X_{t_3})$ disjoint from $X_{t_0} \cup X_{t_3}$ is disjoint from $(A_1' \cap B_1') \cup (A_2' \cap B_2')$.
$\Box$

\medskip

Since $\lvert X^*_{s'} \rvert = \lvert X_s \rvert$ and there exist $\lvert A_1 \cap B_1 \rvert$ disjoint paths in $G$ from $X_{t_0}$ to $X_{t_3}$, if $A_1 \cap B_1-A_s = \emptyset$, then $X_s \subseteq A_1$, so $s \not \in S_5$ by (e), a contradiction.
So $A_1 \cap B_1-A_s \neq \emptyset$.

\medskip

\noindent{\bf Claim 10:} Every side node that is (f)-bad for $(A_1^*,B_1^*)$ and $(A_2^*,B_2^*)$ is (f)-bad for $(A_1,B_1)$ and $(A_2,B_2)$.

\noindent{\bf Proof of Claim 10:}
Suppose to the contrary that $t$ is a side node that is (f)-bad for $(A_1^*,B_1^*)$ and $(A_2^*,B_2^*)$ but not (f)-bad for $(A_1,B_1)$ and $(A_2,B_2)$.
So $A_1 \cap B_1 \cap A_2 \cap B_2-A_t = \emptyset$ and for every edge $xy \in E(G)$ with $x \in A_1 \cap B_1-(A_2 \cap B_2)$ and $y \in A_2 \cap B_2-(A_1 \cap B_1)$, either $x \in A_t$ or $y \in A_t$.

Since $(A_1^*,B_1^*)$ is the reflection of $(A_2^*,B_2^*)$ with respect to $X_{t_0} \cap X_{t_3}, \emptyset$, we know $A_1^* \cap B_1^* \cap A_2^* \cap B_2^* -A_t \subseteq X_{t_0} \cap X_{t_3}-A_t = \emptyset$.
So there exists an edge $u'v'$ of $G$ such that $u' \in A_1^* \cap B_1^* -((A_2^* \cap B_2^*) \cup A_t)$ and $v' \in A_2^* \cap B_2^*-((A_1^* \cap B_1^*) \cup A_t)$.

Suppose that $u' \in A_1 \cap B_1$.
Since $u' \not \in A_t$ and $A_1 \cap B_1 \cap A_2 \cap B_2-A_t =\emptyset$, $u' \not \in A_2 \cap B_2$.
So $u'$ is anti-pointed for $(A_1,B_1)$.
Since $v' \in A_2^* \cap B_2^*-(A_1^* \cap B_1^*) = A_2' \cap B_2'-(A_1' \cap B_1') \subseteq B_1'-A_1' \subseteq B_1-A_1$, $v'$ is the unique neighbor of $u'$ in $B_1-A_1$.
So $v' \in A_2 \cap B_2-(A_1 \cap B_1)$.
Since $u' \in A_1 \cap B_1 - ((A_2 \cap B_2) \cup A_t)$ and $u'v' \in E(G)$, we know $v' \not \in A_2 \cap B_2 - ((A_1 \cap B_1) \cup A_t)$, so $v' \in A_t$, contradicting $v' \in A_2^* \cap B_2^*-((A_1^* \cap B_1^*) \cup A_t)$.

So $u' \not \in A_1 \cap B_1$.
Since $A_1^* \cap B_1^* - (A_1 \cap B_1) = A_1' \cap B_1' - (A_1 \cap B_1) \subseteq X_s$, $u' \in X_s-(A_t \cup A_1)$. 
In particular, $X_s-A_t \neq \emptyset$.
So $t$ is an ancestor of $s$.
That is, $A_t \subseteq A_s$.

Since $\lvert X^*_{s'} \rvert = \lvert X_s \rvert$, there exist $x \in A_1 \cap B_1-A_s$ and $i \in [\lvert X_{t_0} \rvert]$ such that $P_i$ contains both $x$ and $u'$, and the subpath $Q_i$ of $P_i$ between $x$ and $u'$ is contained in $G[B_s]$ and is internally disjoint from $X_s$.
Let $y$ be the neighbor of $x$ in $Q_i$.
Since $\lvert A_1 \cap B_1 \rvert = \lvert A \cap B \rvert = \lvert X_{t_0} \rvert$, $x$ is the only vertex in $V(P_i) \cap A_1 \cap B_1$, so $y \in B_1-A_1$.
Since $A_1 \cap B_1 \cap A_2 \cap B_2-A_t = \emptyset$ and $A_t \subseteq A_s$, $x \not \in A_2 \cap B_2$.
Hence $x \in A_1 \cap B_1 - ((A_2 \cap B_2) \cup A_t)$, and $y \in A_2 \cap B_2 - (A_1 \cap B_1)$ by (b).
Since $xy \in E(G)$ and $t$ is not (f)-bad for $(A_1, B_1)$ and $(A_2,B_2)$, $y \in A_t \subseteq A_s$.
Since $x \in B_s-A_s$ and $xy \in E(G)$, $y \in X_s$.
So $y=u' \not \in A_t$, a contradiction.
$\Box$

\medskip

\noindent{\bf Claim 11:} $A_2 \cap B_2 - A_s = \emptyset$.

\noindent{\bf Proof of Claim 11:}
Suppose to the contrary that $A_2 \cap B_2 -A_s \neq \emptyset$.
Let $v \in A_2 \cap B_2 -A_s$.
Since $A_1^* \cap B_1^*-A_s = A_1' \cap B_1'-A_s = \emptyset$, $s$ is not (f)-bad for $(A_1^*,B_1^*)$ and $(A_2^*,B_2^*)$.
Since $(A_1,B_1)$ and $(A_2,B_2)$ satisfy (f), by Claim 10, $s$ is not (f)-bad for $(A_1,B_1)$ and $(A_2,B_2)$.
So $v \not \in A_1 \cap B_1 \cap A_2 \cap B_2-A_s$.
Hence $v \in A_2 \cap B_2 - ((A_1 \cap B_1) \cup A_s)$.
Let $u$ be the neighbor of $v$ in $A_2-B_2$.
So $u \in A_1 \cap B_1 -(A_2 \cap B_2)$.
Since $s$ is not (f)-bad for $(A_1,B_1)$ and $(A_2,B_2)$, $u \in A_s$.
Since $v \not \in A_s$ and $uv \in E(G)$, $u \in B_s$.
That is, $u \in A_1 \cap B_1 \cap X_s -(A_2 \cap B_2) \subseteq A_1' \cap B_1' - (X_{t_0} \cap X_{t_3})$.
Since $u$ is anti-pointed for $(A_1,B_1)$, $v$ is the unique neighbor of $u$ in $B_1-A_1$.
By Claim 7, $u$ has a neighbor $u'$ in $B_1'-A_1' = (B_1 \cap A_s)-(A_1 \cup B_s)$.
However, $u'$ and $v$ are distinct neighbors of $u$ in $B_1-A_1$, a contradiction.
$\Box$

\medskip

\noindent{\bf Claim 12:} $X_s-A_1 \subseteq A_2 \cap B_2$.

\noindent{\bf Proof of Claim 12:}
Let $v \in X_s-A_1$.
Since $\lvert X^*_{s'} \rvert = \lvert X_s \rvert$, there exists $i \in [\lvert X_{t_0} \rvert]$ such that $P_i$ contains $v$ and a vertex $x$ in $A_1 \cap B_1-A_s$, and the subpath $P_i'$ of $P_i$ between $v$ and $x$ is contained in $G[B_s \cap B_1]$ and is internally disjoint from $X_s$.
By Claim 11, $x \not \in A_2 \cap B_2$.
So the neighbor $y$ of $x$ in $P_i'$ is in $A_2 \cap B_2$.
By Claim 11, $y \in A_s$.
Since $x \not \in A_s$ and $xy \in E(G)$, $y \in X_s$.
So $v=y \in A_2 \cap B_2$.
$\Box$

\medskip

Since $(A_1,B_1)$ is a reflection of $(A_2,B_2)$ with respect to $X_{t_0} \cap X_{t_3}$, there exists $F \subseteq A_1 \cap B_1 \cap A_2 \cap B_2$ with $F \cup (X_{t_0} \cap X_{t_3})=A_1 \cap B_1 \cap A_2 \cap B_2$ such that $(A_1,B_1)$ is the reflection of $(A_2,B_2)$ with respect to $X_{t_0} \cap X_{t_3}, F$.
Since $F \subseteq A_2 \cap B_2$, $F \subseteq A_s$ by Claim 11.

Let $(A_1'',B_1'') = (A_1',B_1')$.

\medskip

\noindent{\bf Claim 13:} There exists a separation $(A_2'',B_2'')$ such that $(A_1'',B_1'')$ is the reflection of $(A_2'',B_2'')$ with respect to $X_{t_0} \cap X_{t_3}, F$, and $(A''_1,B''_1)$ and $(A_2'',B_2'')$ satisfy (b)-(d).

\noindent{\bf Proof of Claim 13:}
For every $u \in A_1' \cap B_1'-(X_{t_0} \cap X_{t_3})$, let $u'$ be the unique neighbor of $u$ in $B_1'-A_1'=B_1''-A_1''$ mentioned in Claim 7.

Let $R = (A_1 \cap B_1 \cap A_s-(A_2 \cap B_2)) \cup (X_s-A_1)$.
Note that $R \subseteq A_1' \cap B_1'-(X_{t_0} \cap X_{t_3})$.
Let $R' = \{u': u \in R\}$.
Let $A_2'' = A_1'' \cup R'$ and $B_2'' = B_1''-R$.
By Claim 11, $A_1'' \cap B_1''-R = A_1 \cap B_1 \cap A_2 \cap B_2 = F \cup (X_{t_0} \cap X_{t_3})$.

We first show that $(A_2'',B_2'')$ is a separation.
Since $R \subseteq A_1''$, $A_2'' \cup B_2''=V(G)$.
Suppose there exists $ab \in E(G)$ such that $a \in A_2''-B_2''$ and $b \in B_2''-A_2''$.
Then $b \in B_1''-(R \cup A_1'' \cup R')$.
Since $(A_1'',B_1'')$ is a separation, $a \in B_1''$.
Since $a \not \in B_2''=B_1''-R$, $a \in R$.
Since $b \in B_1''-A_1''=B_1'-A_1'$ and $a \in R \subseteq A_1' \cap B_1'-(X_{t_0} \cap X_{t_3})$, $b=a' \in R'$ by the definition of $a'$, a contradiction.

Hence $(A_2'',B_2'')$ is a separation of $G$.
Note that $A_2'' \cap B_2'' = (A_1'' \cap B_1''-R) \cup R' = (A_1 \cap B_1 \cap A_2 \cap B_2) \cup R'$.

Now we show that $(A_1'',B_1'')$ and $(A_2'',B_2'')$ satisfy (c) and (d).
Note that as shown in the proof of Claim 9, $(A_1',B_1')$ and $(A_2',B_2')$ satisfy (c) and (d).
Since $A_1'' = A_1'$ and $B_2'' \supseteq B_2'$, $(A_1'',B_1'')$ and $(A_2'',B_2'')$ satisfy (c).
By the definition of $R$ and $A_1'$, we know $A_2'' \supseteq A_2$.
Since $A_2'' \supseteq A_2$ and $B_1''=B_1'$, $(A_1'',B_1'')$ and $(A_2'',B_2'')$ satisfy (d).

Now we show that $(A_1'',B_1'')$ is the reflection of $(A_2'',B_2'')$ with respect to $X_{t_0} \cap X_{t_3}, F$.

Since $(A_1'',B_1'')$ and $(A_2'',B_2'')$ satisfy (c) and (d), they separate $\downarrow t_0$ and $\uparrow t_3$.
Since every vertex in $X_{t_0} \cap X_{t_3}$ is coherent for $t_0,t_3$ and every vertex in $F$ is doubly pointed for $(A_1,B_1)$ and $(A_2,B_2)$, every vertex in $F \cap X_{t_0} \cap X_{t_3}$ is doubly pointed for $(A_1'',B_1'')$ and $(A_2'',B_2'')$ and is not adjacent to any vertex in $A_2'' \cap B_2'' - (F \cup (X_{t_0} \cap X_{t_3}))$.

Since $B_1''-A_1'' \subseteq B_1-A_1$ and every vertex in $F-(X_{t_0} \cap X_{t_3})$ is anti-pointed for $(A_1,B_1)$, every vertex in $F-(X_{t_0} \cap X_{t_3})$ is anti-pointed for $(A_1'',B_1'')$ and $(A_2'',B_2'')$.
Since $A_2'' \cap B_2'' - (F \cup (X_{t_0} \cap X_{t_3})) \subseteq R'$, every vertex in $F-(X_{t_0} \cap X_{t_3})$ is not adjacent to any vertex in $A_2'' \cap B_2'' - (F \cup (X_{t_0} \cap X_{t_3}))$.

Hence every vertex in $F$ is not adjacent to any vertex in $A_2'' \cap B_2'' - (F \cup (X_{t_0} \cap X_{t_3}))$.
So to show that $(A_1'',B_1'')$ is the reflection of $(A_2'',B_2'')$ with respect to $X_{t_0} \cap X_{t_3}, F$, it suffices to show that every vertex in $A_2'' \cap B_2''-(F \cup (X_{t_0} \cap X_{t_3})) \subseteq R'$ is pointed for $(A_2'',B_2'')$, and every vertex in $F-(X_{t_0} \cap X_{t_3})$ is doubly pointed for $(A_2'',B_2'')$.

Any neighbor of some vertex in $R'$ in $A_2''-B_2''$ is in $A_1'' \cap B_1''-B_2'' \subseteq (R \cup (A_1 \cap B_1 \cap A_2 \cap B_2))-(B_1''-R) \subseteq R$ by Claim 11.
Since the set of edges between $R$ and $R'$ is a matching, every vertex in $R'$ is pointed for $(A_2'',B_2'')$.

Let $c \in F-(X_{t_0} \cap X_{t_3})$.
So $c$ is doubly pointed for $(A_2,B_2)$ and anti-pointed for $(A_1,B_1)$ by (b). 
Since $F \subseteq A_1' \cap B_1'$ (by Claim 11), $c \in R$ and hence $c'$ is defined.
Since $c'$ is a neighbor of $c$ in $B_1'-A_1' \subseteq B_1-A_1$, $c$ has no neighbor in $(B_1-A_1) \cap A_1'$. 
Note that $A_2''-B_2'' = (A_1' \cup R') - (B_1'-R) = (A_1 \cup B_s \cup R') - (B_1 \cap A_s-R)$.

Let $c_0$ be a neighbor of $c$ in $A_2''-B_2''$.
Note that $A_2''-B_2'' = (A_1'' \cup R')-B_2'' = A_1''-B_2''$.
So $c_0 \in A_1''-B_2''$ and hence $c_0 \not \in R'$.
If $c_0 \not \in A_1$, then $c_0 \in (B_1-A_1) \cap A_1'' = (B_1-A_1) \cap A_1'$, so $c_0$ is a neighbor of $c$ in $B_1-A_1$ other than $c' \in B_1-A_1$, contradicting that $c$ is anti-pointed for $(A_1,B_1)$. 
So $c_0 \in A_1 = (A_2-B_2) \cup (A_1 \cap B_1 \cap A_2 \cap B_2)$.
Since $c_0 \not \in B_2''=B_1'-R \supseteq A_1 \cap B_1 \cap A_2 \cap B_2$, $c_0 \in A_2-B_2$.

This implies that every neighbor of $c$ in $A_2''-B_2''$ is in $A_2-B_2$.
Since $c$ is pointed for $(A_2,B_2)$, $c$ is pointed for $(A_2'',B_2'')$.
Since $c$ is anti-pointed for $(A_1,B_1)$ and $B_2''-A_2'' \subseteq B_1''-A_1'' \subseteq B_1-A_1$, $c$ is anti-pointed for $(A_2'',B_2'')$.
Therefore, every vertex in $F-(X_{t_0} \cap X_{t_3})$ is doubly pointed for $(A_2'',B_2'')$.

Hence $(A_1'',B_1'')$ is the reflection of $(A_2'',B_2'')$ with respect to $X_{t_0} \cap X_{t_3}, F$.
Therefore, $(A_1'',B_1'')$ and $(A_2'',B_2'')$ satisfy (b).
$\Box$

\medskip

\noindent{\bf Claim 14:} There exist separations $(A_1^{**},B_1^{**})$ and $(A_2^{**},B_2^{**})$ satisfying (a)-(e) such that $A_1^{**} \cap B_1^{**} = A_1'' \cap B_1''$, $A_2^{**} \cap B_2^{**} = A_2'' \cap B_2''$, and every vertex in $(A_1^{**}-A_1'') \cup (A_1''-A_1^{**}) \cup (B_2^{**}-B_2'') \cup (B_2''-B_2^{**})$ is contained in some component of $G-(F \cup (X_{t_0} \cap X_{t_3}))$ disjoint from $X_{t_0} \cup X_{t_3} \cup (A_1'' \cap B_1'') \cup (A_2'' \cap B_2'')$.

\noindent{\bf Proof of Claim 14:}
By Claim 13, $(A_1'',B_1'')$ and $(A_2'',B_2'')$ satisfy (b)-(d).
Since $(A_1'',B_1'')$ and $(A_2'',B_2'')$ satisfy (c) and (d), $(A_2'',B_2'')$ strongly separates $\downarrow t_0$ and $\uparrow t_3$.
Recall that the breadth of $(A,B)$ is minimum among all separations strongly separating $\downarrow t_0$ and $\uparrow t_3$.
So the breadth of $(A_2'',B_2'')$ is at least the breadth of $(A,B)$.
Since $|X^*_{s}|=|X_s|$, we know $|A_2'' \cap B_2''| = |A_1'' \cap B_1''| = |A_1' \cap B_1'| \leq |A \cap B|$ by the submodularity.
Since every vertex in $A_2'' \cap B_2''-(X_{t_0} \cap X_{t_3})$ is pointed for $(A_2'',B_2'')$, and every vertex in $X_{t_0} \cap X_{t_3}$ is coherent for $t_0,t_3$, the breadth of $(A_2'',B_2'')$ is at most the breadth of $(A,B)$.
Hence $(A_1'',B_1'')$ and $(A_2'',B_2'')$ satisfy (a)-(d).
Then this claim follows from Lemma \ref{flipping components}.
$\Box$

\medskip

\noindent{\bf Claim 15:} The set of (f)-bad side nodes for $(A_1^{**},B_1^{**})$ and $(A_2^{**}, B_2^{**})$ equals the set of (f)-bad side nodes for $(A_1,B_1)$ and $(A_2,B_2)$.

\noindent{\bf Proof of Claim 15:}
Since $(A_1,B_1)$ and $(A_2,B_2)$ satisfy (f), it suffices to show that every (f)-bad side node for $(A_1^{**},B_1^{**})$ and $(A_2^{**},B_2^{**})$ is (f)-bad for $(A_1,B_1)$ and $(A_2,B_2)$.
Suppose to the contrary that there exists a side node $t$ that is (f)-bad for $(A_1^{**},B_1^{**})$ and $(A_2^{**},B_2^{**})$ but not (f)-bad for $(A_1,B_1)$ and $(A_2,B_2)$.

Since $A_1^{**} \cap B_1^{**} \cap A_2^{**} \cap B_2^{**} = F \cup (X_{t_0} \cap X_{t_3}) = A_1 \cap B_1 \cap A_2 \cap B_2$ and $t$ is not (f)-bad for $(A_1,B_1)$ and $(A_2,B_2)$, there exist $x \in A_1^{**} \cap B_1^{**} - ((A_2^{**} \cap B_2^{**}) \cup A_t)$ and $y \in A_2^{**} \cap B_2^{**} - ((A_1^{**} \cap B_1^{**}) \cup A_t)$ such that $xy \in E(G)$.
Note that $A_1^{**} \cap B_1^{**} - (A_2^{**} \cap B_2^{**}) = A_1' \cap B_1' - (A_2'' \cap B_2'') \subseteq (A_1 \cap B_1 \cap A_s-(A_2 \cap B_2)) \cup (X_s-A_1)$. 

Suppose that $x \in A_1 \cap B_1 \cap A_s-(A_2 \cap B_2)$.
Then $x$ is anti-pointed for $(A_1,B_1)$.
Since $y \in A_2^{**} \cap B_2^{**}-(A_1^{**} \cap B_1^{**}) = A_2'' \cap B_2''-(A_1'' \cap B_1'') \subseteq B_1''-A_1'' \subseteq B_1-A_1$, $y$ is the unique neighbor of $x$ in $B_1-A_1$.
So $y \in A_2 \cap B_2-(A_1 \cap B_1)$.
Hence $y \in (A_2 \cap B_2) - ((A_1 \cap B_1) \cup A_t)$.
Since $x \in A_1 \cap B_1-((A_2 \cap B_2) \cup A_t)$, $t$ is (f)-bad, a contradiction.

Hence $x \in X_s-A_1$.
Since $x \not \in A_t$, $t$ is an ancestor of $s$.
That is, $A_t \subseteq A_s$.

Since $\lvert X^*_{s'} \rvert = \lvert X_s \rvert$, there exist $a \in A_1 \cap B_1-A_s$ and $i \in [\lvert X_{t_0} \rvert]$ such that $P_i$ contains both $a$ and $x$, and the subpath $Q_i$ of $P_i$ between $a$ and $x$ is contained in $G[B_s]$ and is internally disjoint from $X_s$.
Let $b$ be the neighbor of $a$ in $Q_i$.
Since $\lvert A_1 \cap B_1 \rvert = \lvert A \cap B \rvert = \lvert X_{t_0} \rvert$, $a$ is the only vertex in $V(P_i) \cap A_1 \cap B_1$, so $b \in B_1-A_1$.
Since $t$ is not (f)-bad for $(A_1,B_1)$ and $(A_2,B_2)$, $A_1 \cap B_1 \cap A_2 \cap B_2-A_t = \emptyset$.
Since $A_t \subseteq A_s$, $a \not \in A_2 \cap B_2$.
Hence $a \in A_1 \cap B_1 - ((A_2 \cap B_2) \cup A_t)$.
So $b \in A_2 \cap B_2 - (A_1 \cap B_1)$ by (b).
Since $ab \in E(G)$ and $t$ is not (f)-bad for $(A_1,B_1)$ and $(A_2,B_2)$, $b \in A_t \subseteq A_s$.
Since $a \in B_s-A_s$, $b \in X_s$.
So $b=x \not \in A_t$, a contradiction.
$\Box$

\medskip

Claim 15 implies that $(A_1^{**},B_1^{**})$ and $(A_2^{**},B_2^{**})$ satisfy (a)-(f).

\medskip

\noindent{\bf Claim 16:} Let $t$ be a side node that is (g)-bad for $(A_1^{**},B_1^{**})$ and $(A_2^{**},B_2^{**})$.
Then $t$ is not a descendant of $s$, and if $t$ is not an ancestor of $s$, then $t$ is (g)-bad for $(A_1,B_1)$ and $(A_2,B_2)$.

\noindent{\bf Proof of Claim 16:}
Since $A_1^{**} \cap B_1^{**} = A_1' \cap B_1'$, $A_1^{**} \cap B_1^{**} \subseteq A_s$.
By Claims 7, 13 and 14, we know $A_2^{**} \cap B_2^{**} \subseteq A_s$.
Since $t$ is (g)-bad for $(A_1^{**},B_1^{**})$ and $(A_2^{**},B_2^{**})$, $((A_1^{**} \cap B_1^{**}) \cup (A_2^{**} \cap B_2^{**}))-A_t \neq \emptyset$, so $t$ is not a descendant of $s$.

Suppose to the contrary that $t$ is not an ancestor of $s$, and $t$ is not (g)-bad for $(A_1,B_1)$ and $(A_2,B_2)$.

By Claim 15, if $t$ is (f)-bad for $(A_1,B_1)$ and $(A_2,B_2)$, then $t$ is (f)-bad for $(A_1^{**},B_1^{**})$ and $(A_2^{**},B_2^{**})$, so $t$ is not (g)-bad for $(A_1^{**},B_1^{**})$ and $(A_2^{**},B_2^{**})$, a contradiction.
Hence $t$ is not (f)-bad for $(A_1,B_1)$ and $(A_2,B_2)$.

Since $t$ is (g)-bad for $(A_1^{**},B_1^{**})$ and $(A_2^{**},B_2^{**})$, either
	\begin{itemize}
		\item some vertex in $A_2^{**} \cap B_2^{**} \cap X_t - (A_1^{**} \cap B_1^{**})$ is adjacent to a vertex in $A_1^{**} \cap B_1^{**}-((A_2^{**} \cap B_2^{**}) \cup A_t)$ and a vertex in $B_2^{**}-(A_2^{**} \cup A_t \cup \bigcup_{i=1}^{\lvert X_{t_0} \rvert}V(P_i))$, or
		\item some vertex in $A_1^{**} \cap B_1^{**} \cap X_t - (A_2^{**} \cap B_2^{**})$ is adjacent to a vertex in $A^{**}_2 \cap B^{**}_2-((A^{**}_1 \cap B^{**}_1) \cup A_t)$ and a vertex in $A_1^{**}-(B^{**}_1 \cup A_t \cup \bigcup_{i=1}^{\lvert X_{t_0} \rvert}V(P_i))$.
	\end{itemize}

We first suppose that the former holds.
That is, there exist $v \in A_2^{**} \cap B_2^{**} \cap X_t - (A_1^{**} \cap B_1^{**})$ such that $v$ is adjacent to a vertex $u \in A_1^{**} \cap B_1^{**}-((A_2^{**} \cap B_2^{**}) \cup A_t)$ and a vertex $a \in B_2^{**}-(A_2^{**} \cup A_t \cup \bigcup_{i=1}^{\lvert X_{t_0} \rvert}V(P_i))$.

Since $A_1^{**} \cap B_1^{**} = A_1' \cap B_1'$ and $u \not \in A_t \supseteq B_s$, $u \in A_1 \cap B_1-A_t$.
If $u \in A_2 \cap B_2$, then $u \in A_1 \cap B_1 \cap A_2 \cap B_2-A_t$, so $t$ is (f)-bad for $(A_1,B_1)$ and $(A_2,B_2)$, a contradiction.
So $u \in A_1 \cap B_1 - ((A_2 \cap B_2) \cup A_t)$.

Since $u \in A_1^{**} \cap B_1^{**}-(A_2^{**} \cap B_2^{**})$ and $v \in A_2^{**} \cap B_2^{**} - (A_1^{**} \cap B_1^{**})$, $v$ is the unique neighbor of $u$ in $B_1^{**}-A_1^{**}$.
If $v \not \in B_1''-A_1''$, then by Claim 14, $v$ is contained in some component of $G-(F \cup (X_{t_0} \cup X_{t_3})) = G-(A_1 \cap B_1 \cap A_2 \cap B_2)$ disjoint from $A_1'' \cap B_1''$ containing $u \in A_1^{**} \cap B_1^{**}=A_1'' \cap B_1''$, a contradiction.
So $v$ is in $B_1''-A_1'' \subseteq B_1-A_1$.

Since $u \in A_1 \cap B_1 - ((A_2 \cap B_2) \cup A_t)$ and $v \in B_1-A_1$, we know that $v$ is the unique neighbor of $u$ in $B_1-A_1$.
Hence $v \in A_2 \cap B_2 - (A_1 \cap B_1)$.
Since $uv \in E(G)$ and $u \not \in A_t$, $v \in B_t$.
If $v \not \in X_t$, then $v \not \in A_t$, so $t$ is (f)-bad for $(A_1,B_1)$ and $(A_2,B_2)$, a contradiction.
Hence $v \in A_2 \cap B_2 \cap X_t - (A_1 \cap B_1)$. 
Since $v \in A_2 \cap B_2 \cap X_t - (A_1 \cap B_1)$ and $u \in A_1 \cap B_1 - (A_2 \cap B_2)$, we know that $u$ is the unique neighbor of $v$ in $A_2-B_2$.
Recall that $v$ is adjacent to $u \in A_1^{**} \cap B_1^{**}$ and $a \in B_2^{**}-A_2^{**}$.
So $a \neq v$, and hence $a \in B_2-A_t$.
Since $A_1'' \supseteq A_1$ and $A_1'' \cap B_1'' \cap A_2'' \cap B_2'' = F \cup (X_{t_0} \cap X_{t_3}) \subseteq A_1 \cap B_1 \cap A_2 \cap B_2$, we know $A_2'' \supseteq A_2$.

Since $v \in A_2 \cap B_2 \cap X_t - (A_1 \cap B_1)$ is adjacent to $u \in A_1 \cap B_1 - ((A_2 \cap B_2) \cup A_t)$ and $t$ is not (g)-bad for $(A_1,B_1)$ and $(A_2,B_2)$, we know $a \not \in B_2-(A_2 \cup A_t \cup \bigcup_{i=1}^{|X_{t_0}|}V(P_i))$.
Since $a \in B_2-(A_t \cup \bigcup_{i=1}^{\lvert X_{t_0} \rvert}V(P_i))$, we know $a \in A_2 \cap B_2-(A_t \cup \bigcup_{i=1}^{\lvert X_{t_0} \rvert}V(P_i)) \subseteq A_2'' \cap B_2-(A_t \cup \bigcup_{i=1}^{\lvert X_{t_0} \rvert}V(P_i))$.
Since $a \not \in A_2^{**} \cap B_2^{**}=A_2'' \cap B_2''$, $a \in B_2-B_2''$.
So $a \in A_1'' = A_1 \cup B_s$.
Since $a \not \in A_t$ and $t$ is not an ancestor or descendant of $s$, we know $a \not \in B_s$, so $a \in A_1$.
Since $a \in A_2 \cap B_2 \cap A_1$, $a \in A_1 \cap B_1 \cap A_2 \cap B_2$.
But $a \not \in A_t$, so $t$ is (f)-bad for $(A_1,B_1)$ and $(A_2,B_2)$, a contradiction.

Therefore, the latter condition for $t$ being (g)-bad for $(A_1^{**},B_1^{**})$ and $(A_2^{**},B_2^{**})$ holds.
That is, there exist $x \in A_1^{**} \cap B_1^{**} \cap X_t - (A_2^{**} \cap B_2^{**})$ adjacent to a vertex $y \in A^{**}_2 \cap B^{**}_2-((A^{**}_1 \cap B^{**}_1) \cup A_t)$ and a vertex $z \in A_1^{**}-(B^{**}_1 \cup A_t \cup \bigcup_{i=1}^{\lvert X_{t_0} \rvert}V(P_i))$.

Suppose $z \not \in A_1''$.
By Claim 14, $z$ is contained in some component $C$ of $G-(F \cup (X_{t_0} \cap X_{t_3}))$ disjoint from $A_1'' \cap B_1'' = A_1^{**} \cap B_1^{**}$.
But since $xz \in E(G)$ and $x \not \in A_1^{**} \cap B_1^{**} \cap A_2^{**} \cap B_2^{**}= F \cup (X_{t_0} \cap X_{t_3})$, $C$ contains $x \in A_1^{**} \cap B_1^{**}$, a contradiction.

So $z \in A_1''$.
Since $z \in A_1'-A_t$ and $t$ is not an ancestor or descendant of $s$, $z \in A_1-B_s$.
Since $z \not \in A_1^{**} \cap B_1^{**} = A_1' \cap B_1'$ and $z \in A_1-B_s$, $z \in A_1-(B_1 \cup B_s)$.
Since $xz \in E(G)$, $x \in A_1 \cap A_s$.
Since $x \in A_1^{**} \cap B_1^{**} = A_1' \cap B_1'$, $x \in A_1 \cap B_1 \cap A_s$.
Since $x \not \in X_{t_0} \cap X_{t_3}$, $x$ has a unique neighbor in $B_1-A_1$.
Since $y \in A_2^{**} \cap B_2^{**}-(A_1^{**} \cup B_1^{**}) = A_2'' \cap B_2'' - (A_1'' \cap B_1'')$, $y \not \in A_1'' \supseteq A_1$.
Hence $y$ is the unique neighbor of $x$ in $B_1-A_1$.
Since $x \not \in A_1^{**} \cap B_1^{**} \cap A_2^{**} \cap B_2^{**} = F \cup (X_{t_0} \cap X_{t_3}) = A_1 \cap B_1 \cap A_2 \cap B_2$, $x \in A_1 \cap B_1 - (A_2 \cap B_2)$, and hence $y \in A_2 \cap B_2 - (A_1 \cap B_1)$.

Therefore, $x \in A_1 \cap B_1 \cap X_t - (A_2 \cap B_2)$ is adjacent to $y \in A_2 \cap B_2 - ((A_1 \cap B_1) \cup A_t)$ and $z \in A_1-(B_1 \cup A_t \cup \bigcup_{i=1}^{\lvert X_{t_0} \rvert}V(P_i))$.
So $t$ is (g)-bad for $(A_1,B_1)$ and $(A_2,B_2)$, a contradiction.
$\Box$

\medskip

\noindent{\bf Claim 17:} Every vertex in $X_s-A_1$ is doubly pointed for $(A_2,B_2)$.

\noindent{\bf Proof of Claim 17:}
Suppose to the contrary that there exists $u \in X_s-A_1$ not doubly pointed for $(A_2,B_2)$.
By Claim 12, $u \in A_2 \cap B_2$.
Recall that every vertex in $X_s-A_1$ is pointed for $(A_s,B_s)$ (as shown right before Claim 5).
So $u$ is pointed for $(A_s,B_s)$.

If $u \in X_{t_3}$, then since $s$ is not a descendant of $t_2$ by Claim 5, $u \in X_{t_2} \cap X_{t_3} = X_{t_0} \cap X_{t_3} \subseteq A_1$, a contradiction.
So $u \not \in X_{t_3}$.
In particular, $u \not \in X_{t_0} \cap X_{t_3}$.
Hence $u$ is pointed for $(A_2,B_2)$.
Since $u$ is not doubly pointed for $(A_2,B_2)$, $u$ is not anti-pointed for $(A_2,B_2)$.
Since $u$ is pointed for $(A_s,B_s)$ but not anti-pointed for $(A_2,B_2)$, there exists a neighbor $v$ of $u$ in $(B_2-A_2) \cap B_s$.
Since $v \not \in A_2$, $v \not \in A_1$.
So if $v \in X_s$, then $v \in X_s-A_1 \subseteq A_2 \cap B_2$ by Claim 12, a contradiction.
Hence $v \in B_2 - (A_2 \cup A_s)$.

Suppose that there exists $P \in \{P_1,P_2,...,P_{\lvert X_{t_0} \rvert}\}$ such that $v \in V(P)$.
Let $P'$ be the subpath of $P$ between $v$ and the vertex in $V(P) \cap X_{t_3}$.
Since $v \not \in A_s$, $V(P') \cap X_s \neq \emptyset$.
Since $\lvert V(P) \cap A_2 \cap B_2 \rvert=1$ and $(V(P) \cap X_{t_3}) \cup \{v\} \subseteq B_2-A_2$, we know $V(P') \subseteq B_2-A_2$.
Hence $V(P') \cap X_s -A_2 \neq \emptyset$.
But $X_s-A_2 = X_s-(A_1 \cup A_2) = (X_s-A_1)-A_2 \subseteq A_2 \cap B_2-A_2 = \emptyset$ by Claim 12, a contradiction.

So $v \in B_2-(A_2 \cup A_s \cup \bigcup_{i=1}^{\lvert X_{t_0} \rvert}V(P_i))$.
Since $\lvert X^*_{s'} \rvert = \lvert X_s \rvert$, $u$ is adjacent to some vertex $z$ in $((A_2-B_2) \cap \bigcup_{i=1}^{\lvert X_{t_0} \rvert}V(P_i))-A_s \subseteq A_2-(B_2 \cup A_s)$. 
Since $u \in X_s-A_1$, $u \in A_2 \cap B_2-(A_1 \cap B_1)$.
Since $z \in A_2-B_2$, $z \in A_1 \cap B_1-(A_2 \cap B_2)$.
So $z \in A_1 \cap B_1- ((A_2 \cap B_2) \cup A_s)$.
Since $u$ is adjacent to $z$ and $v$, $s$ is (g)-bad for $(A_1,B_1)$ and $(A_2,B_2)$.

Since $B_s \subseteq A_1''$, $A_2^{**} \cap B_2^{**}-A_s = A_2'' \cap B_2''-A_s = \emptyset$.
Also, $A_1^{**} \cap B_1^{**}-A_s = A_1'' \cap B_1''-A_s = \emptyset$.
So $s$ is not (g)-bad for $(A_1^{**},B_1^{**})$ and $(A_2^{**},B_2^{**})$.

Let $Q$ be the set consisting of the side nodes $q$ such that $q$ is (g)-bad for $(A_1^{**},B_1^{**})$ and $(A_2^{**},B_2^{**})$ but not (g)-bad for $(A_1,B_1)$ and $(A_2,B_2)$.
Since $s$ is (g)-bad for $(A_1,B_1)$ and $(A_2,B_2)$ but not (g)-bad for $(A_1^{**},B_1^{**})$ and $(A_2^{**},B_2^{**})$, we know that $Q \neq \emptyset$ by (g).
By Claim 16, every element of $Q$ is an ancestor of $s$ but not a descendant of $s$.
So $s \not \in Q$ and every element of $Q$ is a proper ancestor of $s$.
Since $\lvert Q \rvert \leq \lvert V(T) \rvert$, $(\lvert V(T) \rvert+1)^{d_s} > \sum_{q \in Q} (\lvert V(T) \rvert+1)^{d_q}$, where $d_s$ and $d_q$ are the distance in $T$ from $s$ and $q$, respectively, to $t_0Tt_3$.
Hence $(A_1,B_1)$ and $(A_2,B_2)$ do not satisfy (g), a contradiction.
$\Box$

\medskip

Let $(C_1,D_1)=(A_1',B_1')$.

\medskip

\noindent{\bf Claim 18:} There exists a separation $(C_2,D_2)$ such that $(C_1,D_1)$ is the reflection of $(C_2,D_2)$ with respect to $X_{t_0} \cap X_{t_3}, F \cup (X_s-A_1)$, and $(C_1,D_1)$ and $(C_2,D_2)$ satisfy (b)-(d).

\noindent{\bf Proof of Claim 18:}
For every $u \in A_1' \cap B_1'-(X_{t_0} \cap X_{t_3})$, let $u'$ be the unique neighbor of $u$ in $B_1'-A_1'=D_1-C_1$ mentioned in Claim 7.

Let $R = A_1 \cap B_1 \cap A_s-(A_2 \cap B_2)$.
Note that $R \subseteq A_1' \cap B_1'-(X_{t_0} \cap X_{t_3})$.
Let $R' = \{u': u \in R\}$.
Let $C_2 = C_1 \cup R'$ and $D_2 = D_1-R$.

We first show that $(C_2,D_2)$ is a separation.
Since $R \subseteq C_2$, $C_2 \cup D_2=V(G)$.
Suppose that there exists $ab \in E(G)$ such that $a \in C_2-D_2$ and $b \in D_2-C_2$.
Then $b \in D_1-(R \cup C_1 \cup R')$.
Since $(C_1,D_1)$ is a separation, $a \in D_1$.
Since $a \not \in D_2=D_1-R$, $a \in R$.
Since $b \in D_2-C_2 \subseteq D_1-C_1=B_1'-A_1'$ and $a \in R \subseteq A_1' \cap B_1'-(X_{t_0} \cap X_{t_3})$, we know $b=a' \in R'$, a contradiction.

Hence $(C_2,D_2)$ is a separation of $G$.
Note that $C_2 \cap D_2 = (C_1 \cap D_1-R) \cup R' = (A_1 \cap B_1 \cap A_2 \cap B_2) \cup (X_s-A_1) \cup R'$.

Now we show that $(C_1,D_1)$ and $(C_2,D_2)$ satisfy (c) and (d).
Note that as shown in the proof of Claim 9, $(A_1',B_1')$ and $(A_2',B_2')$ satisfy (c) and (d).
Since $C_1 = A_1'$ and $(A_1',B_1')$ is the reflection of $(A_2',B_2')$ with respect to $X_{t_0} \cap X_{t_3},\emptyset$, we know $D_2 \supseteq B_2'$, so $(C_1,D_1)$ and $(C_2,D_2)$ satisfy (c).
Suppose that $C_2 \not \supseteq G_A$.
Since $C_2 \supseteq C_1 \supseteq A_1$ and $A_2 \supseteq G_A$ by (d), $G_A \cap A_2-A_1 \neq \emptyset$.
Note that $A_2 = A_1 \cup (A_2 \cap B_2)$, so $G_A-C_2 \subseteq A_2-C_2 \subseteq A_2-C_1 \subseteq A_2 \cap B_2-(A_1 \cup B_s)$.
Hence each vertex in $G_A-C_2$ is a neighbor of some vertex in $A_1 \cap B_1 \cap A_s-(A_2 \cap B_2)=R$ contained in $A_2 \cap B_2-(A_1 \cup B_s) \subseteq B_1'-A_1'$.
So $G_A-C_2 \subseteq R' \subseteq C_2$, a contradiction.
Hence $C_2 \supseteq G_A$. 
Since $C_2 \supseteq G_A$ and $D_1=B_1'$, $(C_1,D_1)$ and $(C_2,D_2)$ satisfy (d).

Now we show that $(C_1,D_1)$ is the reflection of $(C_2,D_2)$ with respect to $X_{t_0} \cap X_{t_3}, F \cup (X_s-A_1)$.

Since $(C_1,D_1)$ and $(C_2,D_2)$ satisfy (c) and (d), they separate $\downarrow t_0$ and $\uparrow t_3$.
Since every vertex in $X_{t_0} \cap X_{t_3}$ is coherent for $t_0,t_3$ and every vertex in $F$ is doubly pointed for $(A_1,B_1)$ and $(A_2,B_2)$, every vertex in $F \cap X_{t_0} \cap X_{t_3}$ is doubly pointed for $(C_1,D_1)$ and $(C_2,D_2)$ and is not adjacent to any vertex in $C_2 \cap D_2 - (F \cup (X_{t_0} \cap X_{t_3}))$.

Recall $C_2 \cap D_2 = (C_1 \cap D_1-R) \cup R' = (A_1 \cap B_1 \cap A_2 \cap B_2) \cup (X_s-A_1) \cup R'$.
So $F-(X_{t_0} \cap X_{t_3}) \subseteq A_1 \cap B_1 \cap A_2 \cap B_2 \subseteq C_1 \cap D_1 \cap C_2 \cap D_2$.

Since $D_1-C_1 \subseteq B_1-A_1$ and every vertex in $F-(X_{t_0} \cap X_{t_3})$ is anti-pointed for $(A_1,B_1)$, every vertex in $F-(X_{t_0} \cap X_{t_3})$ is anti-pointed for $(C_1,D_1)$ and $(C_2,D_2)$.
Since every vertex in $X_s-A_1$ is pointed for $(A_s,B_s)$ (as mentioned right before Claim 5), every vertex in $X_s-A_1$ is anti-pointed for $(A_1',B_1')=(C_1,D_1)$ and hence is anti-pointed for $(C_2,D_2)$.
Since $C_2 \cap D_2 - ((F \cup (X_{t_0} \cap X_{t_3})) \cup (X_s-A_1)) \subseteq R'$, every vertex in $(F-(X_{t_0} \cap X_{t_3})) \cup (X_s-A_1)$ is not adjacent to any vertex in $C_2 \cap D_2 - (F \cup (X_{t_0} \cap X_{t_3}) \cup (X_s-A_1))$. 

Hence every vertex in $F \cup (X_s-A_1)$ is not adjacent to any vertex in $C_2 \cap D_2 - (F \cup (X_{t_0} \cap X_{t_3}) \cup (X_s-A_1))$, and $(F-(X_{t_0} \cap X_{t_3})) \cup (X_s-A_1)$ is anti-pointed for $(C_2,D_2)$, and every vertex in $F \cap X_{t_0} \cap X_{t_3}$ is doubly pointed for $(C_2,D_2)$.
So to show that $(C_1,D_1)$ is the reflection of $(C_2,D_2)$ with respect to $X_{t_0} \cap X_{t_3}, F \cup (X_s-A_1)$, it suffices to show that every vertex in $C_2 \cap D_2-(F \cup (X_{t_0} \cap X_{t_3}) \cup (X_s-A_1)) \subseteq R'$ is pointed for $(C_2,D_2)$, and every vertex in $(F-(X_{t_0} \cap X_{t_3})) \cup (X_s-A_1)$ is pointed for $(C_2,D_2)$.

Any neighbor of some vertex in $R'$ in $C_2-D_2$ is in $C_1 \cap D_1-D_2 \subseteq R$.
Since the set of edges between $R$ and $R'$ is a matching by Claim 7, every vertex in $R'$ is pointed for $(C_2,D_2)$.

So it remains to show that every vertex in $(F-(X_{t_0} \cap X_{t_3})) \cup (X_s-A_1)$ is pointed for $(C_2,D_2)$.

Since $(A_1'',B_1'')$ is the reflection of $(A_2'',B_2'')$ with respect to $X_{t_0} \cap X_{t_3},F$, we know that $C_2-D_2 \subseteq A_2''-B_2''$ and every vertex in $F-(X_{t_0} \cap X_{t_3})$ is pointed for $(A_2'',B_2'')$, so every vertex in $F-(X_{t_0} \cap X_{t_3})$ is pointed for $(C_2,D_2)$.

Suppose there exists a vertex $c \in X_s-A_1$ not pointed for $(C_2,D_2)$.
Note that $C_2-D_2 = (C_1 \cup R')-(D_1-R) = (C_1-D_1) \cup R =(A_1'-B_1') \cup R$. 
Since $c \not \in A_1$, every neighbor of $c$ is in $B_1$.
So every neighbor of $c$ in $C_2-D_2$ is in $((A_1'-B_1') \cup R) \cap B_1 \subseteq (B_1-A_s) \cup R$.
Since every vertex in $R = A_1 \cap B_1 \cap A_s - (A_2 \cap B_2)$ is anti-pointed for $(A_1,B_1)$ and there exist $|X_{t_0}|$ disjoint paths from $X_{t_0}$ to $X_{t_3}$ and $|X_s|=|X^*_{s'}|$, the neighbors of $R$ in $B_1-A_1$ are contained in $R'$, so $c$ is not adjacent to $R$.
Hence every neighbor of $c$ in $C_2-D_2$ is in $B_1-A_s$.
By Claim 17, $c \in X_s-A_1$ is doubly pointed for $(A_2,B_2)$.
By Claim 7, $c'$ is a neighbor of $c$ in $(B_1'-A_1') \cap (\bigcup_{i=1}^{\lvert X_{t_0} \rvert}V(P_i))$, so $c'$ is the unique neighbor of $c$ in $B_2-A_2$.
Hence $c$ has no neighbor in $(B_2-A_2) \cap A_1' \supseteq (B_2-A_2) \cap B_s$.
Recall that every neighbor of $c$ in $C_2-D_2$ is in $B_1-A_s$, so every neighbor of $c$ in $C_2-D_2$ is in $(B_1-(B_2 \cup A_s)) \cup (B_1 \cap A_2-A_s) = B_1 \cap A_2-A_s$.
Since $A_2 \cap B_2 \subseteq A_s$ by Claim 11, $c$ has no neighbor in $A_2 \cap B_2$.
So every neighbor of $c$ in $C_2-D_2$ is in $A_2-B_2$.
Since $c$ is pointed for $(A_2,B_2)$, $c$ is pointed for $(C_2,D_2)$, a contradiction.

Hence every vertex in $(F-(X_{t_0} \cap X_{t_3})) \cup (X_s-A_1)$ is pointed for $(C_2,D_2)$.
This shows that $(C_1,D_1)$ is the reflection of $(C_2,D_2)$ with respect to $X_{t_0} \cap X_{t_3}, F \cup (X_s-A_1)$.
Therefore, $(C_1,D_1)$ and $(C_2,D_2)$ satisfy (b).
$\Box$

\medskip

\noindent{\bf Claim 19:} There exist separations $(C_1^*,D_1^*)$ and $(C_2^*,D_2^*)$ satisfying (a)-(e) such that $C_1^* \cap D_1^* = C_1 \cap D_1$, $C_2^* \cap D_2^* = C_2 \cap D_2$, and every vertex in $(C_1^*-C_1) \cup (C_1-C_1^*) \cup (D_2^*-D_2) \cup (D_2-D_2^*)$ is contained in some component of $G-(F \cup (X_{t_0} \cap X_{t_3}) \cup (X_s-A_1))$ disjoint from $X_{t_0} \cup X_{t_3} \cup (C_1 \cap D_1) \cup (C_2 \cap D_2)$.

\noindent{\bf Proof of Claim 19:}
By Claim 18, $(C_1,D_1)$ and $(C_2,D_2)$ satisfy (b)-(d).
Since $(C_1,D_1)$ and $(C_2,D_2)$ satisfy (c) and (d), $(C_2,D_2)$ strongly separates $\downarrow t_0$ and $\uparrow t_3$.
Recall that the breadth of $(A,B)$ is minimum among all separations strongly separating $\downarrow t_0$ and $\uparrow t_3$.
So the breadth of $(C_2,D_2)$ is at least the breadth of $(A,B)$.
Since every vertex in $C_2 \cap D_2-(X_{t_0} \cap X_{t_3})$ is pointed for $(C_2,D_2)$, and every vertex in $X_{t_0} \cap X_{t_3}$ is coherent for $t_0,t_3$, the breadth of $(C_2,D_2)$ is at most the breadth of $(A,B)$.
Hence $(C_1,D_1)$ and $(C_2,D_2)$ satisfy (a)-(d).
Then this claim follows from Lemma \ref{flipping components}.
$\Box$

\medskip

\noindent{\bf Claim 20:} The set of (f)-bad side nodes for $(C^*_1,D^*_1)$ and $(C^*_2, D^*_2)$ equals the set of (f)-bad side nodes for $(A_1,B_1)$ and $(A_2,B_2)$.

\noindent{\bf Proof of Claim 20:}
Since $(A_1,B_1)$ and $(A_2,B_2)$ satisfy (f), it suffices to show that every (f)-bad side node for $(C^*_1,D^*_1)$ and $(C^*_2,D^*_2)$ is (f)-bad for $(A_1,B_1)$ and $(A_2,B_2)$.
Suppose to the contrary that there exists a side node $t$ that is (f)-bad for $(C_1^*,D_1^*)$ and $(C_2^*,D_2^*)$ but not (f)-bad for $(A_1,B_1)$ and $(A_2,B_2)$.

Suppose that $C_1^* \cap D_1^* \cap C_2^* \cap D_2^*-A_t \neq \emptyset$.
Note that $C_1^* \cap D_1^* \cap C_2^* \cap D_2^* = F \cup (X_{t_0} \cap X_{t_3}) \cup (X_s-A_1) = (A_1 \cap B_1 \cap A_2 \cap B_2) \cup (X_s-A_1)$.
Since $t$ is not (f)-bad for $(A_1,B_1)$ and $(A_2,B_2)$, $(X_s-A_1)-A_t \neq \emptyset$.
Hence $t$ is an ancestor of $s$.
Let $a \in (X_s-A_1)-A_t$.
Note that $a \in X_s \cap A_2 \cap B_2-A_1$ by Claim 12.
Since $\lvert X^*_{s'} \rvert = \lvert X_s \rvert$, the neighbor $b$ of $a$ in $A_1 \cap B_1-(A_2 \cap B_2)$ is contained in $B_s-A_s \subseteq B_t-A_t$.
Hence $ab \in E(G)$, $b \in A_1 \cap B_1-((A_2 \cap B_2) \cup A_t)$ and $a \in A_2 \cap B_2-((A_1 \cap B_1) \cup A_t)$.
So $t$ is (f)-bad for $(A_1,B_1)$ and $(A_2,B_2)$, a contradiction.

Hence there exist $x \in C_1^* \cap D_1^* - ((C_2^* \cap D_2^*) \cup A_t)$ and $y \in C_2^* \cap D_2^* - ((C_1^* \cap D_1^*) \cup A_t)$ such that $xy \in E(G)$.
Recall that $A_1^{**} \cap B_1^{**} = A_1'' \cap B_1'' = C_1^* \cap D_1^*$ and $(A_1^{**},B_1^{**})$ is the reflection of $(A_2^{**},B_2^{**})$ with respect to $X_{t_0} \cap X_{t_3}, F$.
Since $(C_1^*,D_1^*)$ is the reflection of $(C_2^*,D_2^*)$ with respect to $X_{t_0} \cap X_{t_3}, F'$ for some $F' \supseteq F$, we know that $x \in C_1^* \cap D_1^* - ((C_2^* \cap D_2^*) \cup A_t) \subseteq A_1^{**} \cap B_1^{**} - ((A_2^{**} \cap B_2^{**}) \cup A_t)$ and $y \in C_2^* \cap D_2^* - ((C_1^* \cap D_1^*) \cup A_t) \subseteq A_2^{**} \cap B_2^{**} - ((A_1^{**} \cap B_1^{**}) \cup A_t)$.
So $t$ is (f)-bad for $(A_1^{**},B_1^{**})$ and $(A_2^{**},B_2^{**})$.
Hence by Claim 15, $t$ is (f)-bad for $(A_1,B_1)$ and $(A_2,B_2)$, a contradiction.
$\Box$

\medskip

Claims 19 and 20 imply that $(C_1^*,D_1^*)$ and $(C_2^*,D_2^*)$ satisfy (a)-(f).

\medskip

\noindent{\bf Claim 21:} The set of (g)-bad side nodes for $(C^*_1,D^*_1)$ and $(C^*_2, D^*_2)$ equals the set of (g)-bad side nodes for $(A_1,B_1)$ and $(A_2,B_2)$.

\noindent{\bf Proof of Claim 21:}
Since $(A_1,B_1)$ and $(A_2,B_2)$ satisfy (g), it suffices to show that every (g)-bad side node for $(C_1^*,D_1^*)$ and $(C_2^*,D_2^*)$ is (g)-bad for $(A_1,B_1)$ and $(A_2,B_2)$.
Suppose to the contrary that there exists a side node $t$ that is (g)-bad for $(C_1^*,D_1^*)$ and $(C_2^*,D_2^*)$ but not (g)-bad for $(A_1,B_1)$ and $(A_2,B_2)$.

By Claim 20, if $t$ is (f)-bad for $(A_1,B_1)$ and $(A_2,B_2)$, then $t$ is (f)-bad for $(C_1^*,D_1^*)$ and $(C_2^*,D_2^*)$, so $t$ is not (g)-bad for $(C_1^*,D_1^*)$ and $(C_2^*,D_2^*)$, a contradiction.
Hence $t$ is not (f)-bad for $(A_1,B_1)$ and $(A_2,B_2)$.

We first suppose that there exist $v \in C_2^* \cap D_2^* \cap X_t - (C_1^* \cap D_1^*)$ such that $v$ is adjacent to a vertex $u \in C_1^* \cap D_1^*-((C_2^* \cap D_2^*) \cup A_t)$ and a vertex $a \in D_2^*-(C_2^* \cup A_t \cup \bigcup_{i=1}^{\lvert X_{t_0} \rvert}V(P_i))$.
Since $C_2^* \cap D_2^* - (C_1^* \cap D_1^*)$ is contained in the intersection of the neighborhood of vertices in $C_1^* \cap D_1^*-(C_1^* \cap D_1^* \cap C_2^* \cap D_2^*) \subseteq A_1' \cap B_1'-((A_1 \cap B_1 \cap A_2 \cap B_2) \cup (X_s-A_1)) \subseteq A_1 \cap B_1 \cap A_s$ and $B_1'-A_1' \subseteq B_1-A_1$, we know $C_2^* \cap D_2^* - (C_1^* \cap D_1^*) \subseteq A_2 \cap B_2-(A_1 \cap B_1)$.
Since $v \in C_2^* \cap D_2^* \cap X_t - (C_1^* \cap D_1^*)$, we know $v \in A_2 \cap B_2 \cap X_t-(A_1 \cap B_1)$.
By the existence of $P_1,P_2,...,P_{\lvert X_{t_0} \rvert}$, $u \in A_1 \cap B_1 - (A_2 \cap B_2)$, so $u \in A_1 \cap B_1 - ((A_2 \cap B_2)\cup A_t)$.
Note that $u$ is the unique neighbor of $v$ in $A_2-B_2$.
So $a \in B_2-(C_2^* \cup A_t \cup \bigcup_{i=1}^{\lvert X_{t_0} \rvert}V(P_i))$.
Since $t$ is not (g)-bad for $(A_1,B_1)$ and $(A_2,B_2)$, $a \in A_2 \cap B_2-(C_2^* \cup A_t \cup \bigcup_{i=1}^{\lvert X_{t_0} \rvert}V(P_i))$.
Since $X_s-A_1 \subseteq A_2 \cap B_2$ by Claim 12 and $X_s-A_1 \subseteq C_1^* \cap D_1^* \cap C_2^* \cap D_2^* \subseteq C_2^*$, we know $a \in A_2 \cap B_2-(X_s-A_1)$.
If $a \in A_1$, then since $va \in E(G)$ and $v \in A_2 \cap B_2-(A_1 \cap B_1)$ and $u$ is the unique neighbor of $v$ in $A_2-B_2$, we know $a \in A_1 \cap B_1 \cap A_2 \cap B_2 \subseteq C_2^*$, a contradiction.
So $a \in A_2 \cap B_2-(A_1 \cup X_s)$.
Since $A_2 \cap B_2 \subseteq A_s$ by Claim 11, $a \in A_2 \cap B_2-(A_1 \cup B_s)$. 
So $a$ is the unique neighbor of some vertex $b$ in $A_1 \cap B_1 \cap A_s-(A_2 \cap B_2) \subseteq A_1' \cap B_1'-(A_2 \cap B_2)$ in $B_1-A_1$.
Since $a \not \in B_s$, $a$ is the unique neighbor of $b \in (A_1' \cap B_1')-(A_2 \cap B_2)$ in $B_1'-A_1'$.
Since $X_s-A_1 \subseteq A_2 \cap B_2$ by Claim 12, $b \in C_1^* \cap D_1^*-(F \cup (X_{t_0} \cup X_{t_3}) \cup (X_s-A_1))$.
Hence $a \in C_2^*$, a contradiction.

Therefore, since $t$ is (g)-bad for $(C_1^*,D_1^*)$ and $(C_2^*,D_2^*)$, we know that there exist $x \in C_1^* \cap D_1^* \cap X_t - (C_2^* \cap D_2^*)$ adjacent to a vertex $y \in C_2^* \cap D_2^*-((C_1^* \cap D_1^*) \cup A_t)$ and a vertex $z \in C_1^*-(D_1^* \cup A_t \cup \bigcup_{i=1}^{\lvert X_{t_0} \rvert}V(P_i))$.

So $x \in A_1 \cap B_1 \cap X_t - (A_2 \cap B_2)$.
By the existence of $P_1,P_2,...,P_{\lvert X_{t_0} \rvert}$, we have $y \in A_2 \cap B_2 - ((A_1 \cap B_1) \cup A_t)$, and $y$ is the unique neighbor of $x$ in $B_1-A_1$.
So $z \in A_1$.
Since $z \not \in \bigcup_{i=1}^{\lvert X_{t_0} \rvert}V(P_i)$, $z \not \in A_1 \cap B_1$.
So $z \in A_1-(B_1 \cup A_t \cup \bigcup_{i=1}^{\lvert X_{t_0} \rvert}V(P_i))$.
Hence $t$ is (g)-bad for $(A_1,B_1)$ and $(A_2,B_2)$, a contradiction.
$\Box$

\medskip

By Claim 21, $(C_1^*,D_1^*)$ and $(C_2^*,D_2^*)$ satisfy (a)-(g).

\medskip

\noindent{\bf Claim 22:} Every side node that is (h)-bad for $(C_1^*,D_1^*)$ and $(C_2^*,D_2^*)$ is (h)-bad for $(A_1,B_1)$ and $(A_2,B_2)$.

\noindent{\bf Proof of Claim 22:}
Suppose to the contrary that there exists an (h)-bad side node $t$ for $(C_1^*,D_1^*)$ and $(C_2^*,D_2^*)$ but not (h)-bad for $(A_1,B_1)$ and $(A_2,B_2)$.
Since $t$ is (h)-bad for $(C_1^*,D_1^*)$ and $(C_2^*,D_2^*)$, $t$ is not (f)-bad for $(C_1^*,D_1^*)$ and $(C_2^*,D_2^*)$, so $C_1^* \cap D_1^* \cap C_2^* \cap D_2^*-A_t = \emptyset$.
Note that $t$ is not (f)-bad for $(A_1,B_1)$ and $(A_2,B_2)$ by Claim 20.

Suppose that there exists $v \in C_1^* \cap D_1^*-A_t$.
Since $C_1^* \cap D_1^* \cap C_2^* \cap D_2^*-A_t = \emptyset$, $v \in C_1^* \cap D_1^*-((C_2^* \cap D_2^*) \cup A_t) \subseteq A_1 \cap B_1 - ((A_2 \cap B_2) \cup A_t)$ by Claims 18 and 19.
So $t$ is (h)-bad for $(A_1,B_1)$ and $(A_2,B_2)$, a contradiction.

Hence there exists $u \in C_2^* \cap D_2^*-A_t$.
Since $C_1^* \cap D_1^* \cap C_2^* \cap D_2^*-A_t = \emptyset$, $u \in C_2^* \cap D_2^*-((C_1^* \cap D_1^*) \cup A_t) \subseteq A_2 \cap B_2 - ((A_1 \cap B_1) \cup A_t)$ by Claims 18 and 19.
So $t$ is (h)-bad for $(A_1,B_1)$ and $(A_2,B_2)$, a contradiction.
$\Box$

\medskip

Note that $C_1^* \cap D_1^* = C_1 \cap D_1 \subseteq A_s$ and $C_2^* \cap D_2^* = C_2 \cap D_2 \subseteq A_s$.
Hence $s$ is not (h)-bad for $(C_1^*,D_1^*)$ and $(C_2^*,D_2^*)$.
By Claim 11, $s$ is not (f)-bad.
Since $A_1 \cap B_1 -A_s \neq \emptyset$ (as mentioned right above Claim 10), $s$ is (h)-bad for $(A_1,B_1)$ and $(A_2,B_2)$.
Therefore, $(A_1,B_1)$ and $(A_2,B_2)$ do not satisfy (h) by Claim 22, a contradiction.

This proves that the breadth of $(A_{s'},B_{s'})$ is less than the breadth of $(A_s,B_s)$ and completes the proof of this lemma.
\end{pf}

\subsection{Elevation} \label{subsec:elevation}

Let $(T,\X)$ be a rooted tree-decomposition of a graph $G$.
The \defn{elevation} of $(T,\X)$ is the maximum $(Z,s)$-depth among all $Z \subseteq V(G)$ and positive integers $s$.
(Recall that the $(Z,s)$-depth is the maximum length of a $(Z,s)$-strip in $(T,\X)$ defined in Section \ref{subsec:pseudo_ec_strips}.)

Now we are ready to prove the main theorem of this section.

\begin{theorem} \label{bounded depth}
For any positive integers $k,w$, there exist integers $N=(w+1)2^{2(w+1)(w+2)}+1$ and $f=f(k,w)$ such that if $G$ is a graph of tree-width at most $w$ not containing the Robertson chain of length $k$ as a topological minor, then there exists an $N$-linked and $N$-integrated rooted tree-decomposition $(T,\X)$ of $G$ of width $\tw(G)$ and of elevation at most $f(k,w)$ such that for every edge $xy$ of $T$, either $X_x \subseteq X_y$ or $X_y \subseteq X_x$.
\end{theorem}

\begin{pf}
Define $f(k,w)=g_{\ref{breaking a long strip better}}(k,w,N+2)$, where $g_{\ref{breaking a long strip better}}$ is the function $g$ mentioned in Lemma \ref{breaking a long strip better}. 
Let $(T,\X)$ be a rooted tree-decomposition of $G$ of width $\tw(G)$, and subject to this, the signature of $(T,\X)$ is as large as possible, and subject to these, the number of edges $xy$ of $T$ such that $X_x \not \subseteq X_y$ and $X_y \not \subseteq X_x$ is as small as possible.
By Lemma \ref{N-linked}, $(T,\X)$ is $N$-linked.
By Lemma \ref{N-integrated}, $(T,\X)$ is $N$-integrated.

If there exists an edge $xy$ of $T$ such that $X_x \not \subseteq X_y$ and $X_y \not \subseteq X_x$, then subdividing $xy$ and defining the bag of the new node to be $X_x \cap X_y$ will result in a rooted tree-decomposition of $G$ of width $\tw(G)$ and signature at least $(T,\X)$, but smaller number of ``bad'' edges, a contradiction.
So for every edge $xy$ of $T$, either $X_x \subseteq X_y$ or $X_y \subseteq X_x$.

To prove this theorem, it suffices to show that the elevation of $(T,\X)$ is at most $f(k,w)$.

Suppose to the contrary that the elevation of $(T,\X)$ is greater than $f(k,w)$.
So there exists a $(Z,s)$-strip in $(T,\X)$ with length at least $f(k,w)$, for some $Z \subseteq V(G)$ and integer $s \in [w+1]$.

Suppose $s=w+1$.
Then $w=\tw(G)$ and all nodes contained in the shortest path in $T$ passing through all nodes in this $(Z,s)$-strip have bag size exactly $w+1$.
But it implies that all of those bags are identical since $X_x \subseteq X_y$ or $X_y \subseteq X_x$ for every $xy \in E(T)$, contradicting the definition of a $(Z,s)$-strip.

So $s \in [w]$.
By Lemma \ref{breaking a long strip better}, there exist a $(Z,s)$-strip $R$ in $(T,\X)$ and a pseudo-edge-cut $(A,B)$ modulo $Z$ of order $\lvert Z \rvert+s$ such that every vertex in $Z$ is coherent for all pairs of nodes in $R$, and $(A,B)$ is a separation $(N+2)$-breaking $R$.
So there exist $t_1,t_2,...,t_{2N+4} \in R$, where $t_i$ is a precursor of $t_{i+1}$ for $i \in [2N+3]$, such that $\downarrow t_{N+2} \subseteq A$ and $\uparrow t_{N+3} \subseteq B$.
Note that $X_{t_1},X_{t_{N+4}},X_{t_{N+5}},...,X_{t_{2N+4}}$ are distinct.
By Lemma \ref{N-linked disjoint paths}, since $(T,\X)$ is $N$-linked, there exist $\lvert X_{t_1} \rvert$ disjoint paths in $G$ from $X_{t_1}$ to $X_{t_{N+4}}$.
Hence, $(A,B)$ is a separation of $G$ strongly separating $\downarrow t_{N+1}$ and $\uparrow t_{N+3}$ of breadth $(\lvert X_{t_1} \rvert, r)$, where $r$ is the number of vertices in $Z$ non-pointed for $(A_{t_{N+4}},B_{t_{N+4}})$, since every vertex in $Z$ is coherent for all pairs of nodes in $R$.
Note that since $s>0$ and the sets $X_{t_i}-Z$ are nonempty and pairwise disjoint for all $i \in [2N+4]$, we know that for each integer $j$ with $j \geq N+5$, there exists a vertex $v_j \in X_{t_j}-X_{t_{j-1}}$.
So $\{v_j: N+5 \leq j \leq 2N+4\}$ contains $N$ vertices that cannot be separated from $\downarrow t_N$ by a separation of breadth less than the breadth of $(A,B)$ given by a node of $T$.
Hence $t_N,t_{N+1},t_{N+3},t_{N+4}$ satisfy (TE1)-(TE7).
However, since $(T,\X)$ is $N$-integrated, there exists a node $t \in t_{N}Tt_{N+4}$ such that the breadth of $(A_t,B_t)$ equals the breadth of $(A,B)$.
This implies that $(A_t,B_t)$ is a pseudo-edge-cut of order $|X_{t_1}|$ modulo $Z$, which is impossible since $R$ is a $(Z,s)$-strip. 
Therefore, $(T,\X)$ has elevation at most $f(k,w)$.
\end{pf}

\section{A tree lemma} \label{sec: tree lemma}

The goal of this section is to prove Theorem \ref{decorated tree lemma}, which is a strengthening of \cite[Theorem 2.2]{rs IV} and a form of the minimal bad sequence argument for proving results on well-quasi-ordering.
We will use it in later sections.

\subsection{An old tree lemma}

A \defn{stable set} in a finite or infinite graph is a subset of pairwise nonadjacent vertices.
We say that a subset $I$ of vertices of an infinite graph $G$ is \defn{rich} in $G$ if no infinite subset of $I$ is a stable set.

To strengthen \cite[Theorem 2.2]{rs IV}, we shall use the following weaker version of it.

\begin{theorem}[{\cite[Theorem 2.1]{rs IV}}] \label{rs weak tree lemma}
Let $T_1,T_2,...$ be a countable sequence of disjoint rooted trees. 
Let $M$ be an infinite graph with $V(M) = V(T_1 \cup T_2 \cup ...)$ such that for $i'>i \geq 1$, if $u \in V(T_i)$ is adjacent to $w \in V(T_{i'})$ in $M$, then $u$ is adjacent in $M$ to all non-root ancestors of $w$.
If the set of the roots of $T_1, T_2,...$ is a stable set of $M$, then there is an infinite stable set $X$ of $M$ such that $\lvert X \cap V(T_i) \rvert \leq 1$ for each $i \geq 1$ and such that the set of heads of all edges of $T_1 \cup T_2 \cup ...$ with tails in $X$ in rich in $M$. 
\end{theorem}

\subsection{Decorated trees} \label{subsec:decorated_trees}

We need terminologies to state the main result of this section (Theorem \ref{decorated tree lemma}).

For every nonnegative integer $n$ and any set $X$, we define $[X]^{\leq n}$ to be the collection of subsets of $X$ with size at most $n$.

Let $X$ be a set and $n$ be a nonnegative integer.
Let $T$ be a rooted tree, and let $\phi,\tau$ be functions from $E(T)$ to $[X]^{\leq n}$ such that $\tau(e) \subseteq \phi(e)$ for all $e \in E(T)$.
Let $N$ be a nonnegative integer and $\mu: E(T) \rightarrow \{0,1,2,...,N\}$. 
For $v,w \in V(T)$, we say that $v$ \defn{precedes} $w$ in $T$ with respect to $(\phi,\tau,\mu)$ if the following hold.
	\begin{itemize}
		\item $v$ is not the root of $T$.
		\item $v$ is an ancestor of $w$.
		\item Let $e,f$ be the edges of $T$ with heads $v,w$, respectively.
Then $\lvert \phi(e) \rvert = \lvert \phi(f) \rvert$, $\tau(e)=\tau(f)$ and $\mu(e)=\mu(f)$.
		\item $\lvert \phi(g) \rvert \geq \lvert \phi(f) \rvert$ for all edges $g$ in $vTw$.
		\item $\mu(g) \geq \mu(e)$ for all edges $g$ in $vTw$ with $\lvert \phi(g) \rvert = \lvert \phi(e) \rvert$.
	\end{itemize}
Let $m$ be a positive integer.
We say that $(T,\phi,\tau,\mu)$ is \defn{$(n,m,N)$-decorated} if the following hold.
	\begin{itemize}
		\item If $e,e',e'' \in E(T)$ and they appear on a directed path in $T$ in the order listed, then $\phi(e) \cap \phi(e'') \subseteq \phi(e')$.
		\item $t \leq m$ whenever $t$ is a positive integer such that the following hold.
			\begin{itemize}
				\item $P$ is a directed path in $T$.
				\item $e_1,e_2,...,e_t \in E(P)$ are distinct directed edges appearing in $P$ in the order listed with $\lvert \phi(e_1) \rvert = \lvert \phi(e_2) \rvert = ... = \lvert \phi(e_t) \rvert$. 
				\item There exists $Z$ such that $\phi(e_i) \cap \phi(e_j) = Z$ for all $1 \leq i < j \leq t$.
				\item $\lvert \phi(e) \rvert \geq \lvert \phi(e_1) \rvert$ for every edge $e$ of $P$.
				\item $\tau(e) \not \subseteq Z$ and $\mu(e) = \mu(e_1)$ for every edge $e$ of $P$ with $\lvert \phi(e) \rvert = \lvert \phi(e_1) \rvert$.
			\end{itemize}
	\end{itemize}

\subsection{Lemmas for the main tree lemma}
The main goal of Section \ref{sec: tree lemma} is to prove Theorem \ref{decorated tree lemma}, which will be proved in Section \ref{subsec:main_tree_lemma}.
In this subsection, we prove lemmas towards the proof of Theorem \ref{decorated tree lemma}.

We say that $(D,(T_i,\phi_i,\tau_i,\mu_i)_{i \in {\mathbb N}},N,h,d,L)$ is a \defn{standard tuple for Theorem \ref{decorated tree lemma}} if the following hold.
	\begin{itemize}
		\item $T_1,T_2,...$ is an infinite sequence of pairwise disjoint rooted trees.
		\item $N,h,d$ are integers with $N,h \geq 0$ and $d > 0$, and $L$ is a set.
		\item For each $i \geq 1$, $\mu_i:E(T_i) \rightarrow \{0,1,2,...,N\}$ is a function, and $\phi_i$ and $\tau_i$ are functions from $E(T)$ to $[L]^{\leq h}$ such that $(T_i,\phi_i,\tau_i,\mu_i)$ is $(h,d,N)$-decorated.
		\item $D$ is an infinite graph with $V(D) = V(T_1 \cup T_2 \cup ...)$ such that for $i'>i \geq 1$, if $u \in V(T_i)$ is adjacent to $w \in V(T_{i'})$ in $D$, and $v \in V(T_{i'})$ precedes $w$ in $T_{i'}$ with respect to $(\phi_{i'},\tau_{i'},\mu_{i'})$, then $u$ is adjacent to $v$ in $D$.
		\item The roots of $T_1,T_2,...$ form a stable set in $D$.
	\end{itemize}

\begin{lemma} \label{tree lemma bounded height}
Let $(D,(T_i,\phi_i,\tau_i,\mu_i)_{i \in {\mathbb N}},N,h,d,L)$ be a standard tuple for Theorem \ref{decorated tree lemma}.
If there exists a positive integer $\ell$ such that each $T_i$ has no directed path of length $\ell$, then there exists an infinite stable set $S$ of $D$ such that $\lvert S \cap V(T_i) \rvert \leq 1$ for each $i \geq 1$ and such that the set of heads of all edges of $T_1 \cup T_2 \cup ...$ with tails in $S$ is rich in $D$.
\end{lemma}

\begin{pf}
Let $R$ be the set of the roots of $T_1,T_2, \cdots$.
We shall prove this lemma by induction on $\ell$.
When $\ell=1$, each $T_i$ contains only one node, so the lemma holds by choosing $S=R$ since the empty set is rich.
So we may assume that $\ell \geq 2$ and this lemma holds for all smaller $\ell$.

Let $C$ be the set of the children of all nodes in $R$.
If $C$ is rich, then we are done by choosing $S=R$.
So we may assume that $C$ contains an infinite stable set $W$ in $D$.
Since each $T_i$ is finite, $W$ intersects $V(T_i)$ for infinitely many integers $i$.
So we can take an infinite subset $W'$ of $W$ such that $\lvert W' \cap T_i \rvert \leq 1$ for all $i$.
Let $i_1<i_2<...$ be the indices $j$ such that $W' \cap V(T_j) \neq \emptyset$.
For each $j \geq 1$, define $T_j'$ to be the maximal subtree of $T_{i_j}$ rooted at the vertex in $W' \cap V(T_{i_j})$.
So each $T_j'$ has no directed path of length $\ell-1$.
Furthermore, the roots of $T_j'$ form the set $W'$, which is a stable set in $D$.
Let $D'=D[V(T_1' \cup T_2' \cup ...)]$.
By the induction hypothesis, there exists a stable set $S$ in $D'$ such that $\lvert S \cap T_j' \rvert \leq 1$ for all $j$ and the set of all heads of all edges of $T_1' \cup T_2' \cup ...$ with tails in $S$ is rich in $D'$.
Note that the set of all heads of all edges of $T_1' \cup T_2' \cup ...$ with tails in $S$ in $D'$ is the same as the set of all heads of all edges of $T_1 \cup T_2 \cup ...$ with tails in $S$ in $D$.
This proves the lemma.
\end{pf}

\bigskip

Let $(D,(T_i,\phi_i,\tau_i,\mu_i)_{i \in {\mathbb N}},N,h,d,L)$ be a standard tuple for Theorem \ref{decorated tree lemma}.
For $i \geq 1$ and $F_i \subseteq E(T_i)$, we define the \defn{$F_i$-contraction of $(T_i,\phi_i,\tau_i,\mu_i)$} to be the tuple $(T_i',\phi_i',\tau_i',\mu_i')$ as follows.
	\begin{itemize}
		\item Define $T_i'$ to be the rooted tree obtained from $T_i$ by contracting each component of $T_i-F_i$ into a node such that the root of $T_i'$ is the node obtained by contracting the component of $T_i-F_i$ containing the root of $T_i$.
		\item Define $\phi_i'=\phi|_{F_i}$, $\tau_i'=\tau_i|_{F_i}$, and $\mu_i'=\mu_i|_{F_i}$.
	\end{itemize}
And we define the \defn{$(\bigcup_{i \geq 1}F_i)$-contraction of $D$}, denoted by $D'$, to be the infinite graph with $V(D')=V(T_1' \cup T_2' \cup ...)$ such that for any positive integers $i$ and $j$ and for any $x \in V(T_i')$ and $y \in V(T_j')$, $x$ and $y$ are adjacent in $D'$ if and only if $i \neq j$ and the root of the component of $T_i-F_i$ contracted into $x$ and the root of the component of $T_j-F_j$ contracted into $y$ are adjacent in $D$.

\begin{lemma} \label{still rich}
Let $(D,(T_i,\phi_i,\tau_i,\mu_i)_{i \in {\mathbb N}},N,h,d,L)$ be a standard tuple for Theorem \ref{decorated tree lemma}.
For each $i \geq 1$, let $F_i \subseteq E(T_i)$, and denote the $F_i$-contraction of $(T_i,\phi_i,\tau_i,\mu_i)$ by $(T_i',\phi',\tau_i',\mu_i')$.
Denote the $(\bigcup_{i \geq 1} F_i)$-contraction of $D$ by $D'$.
Assume that there exists an infinite stable set $S'$ of $D'$ such that $\lvert S' \cap V(T_i') \rvert \leq 1$ for each $i \geq 1$ and such that the set of heads of all edges of $T_1' \cup T_2' \cup ...$ with tails in $S'$ is rich in $D'$.
Let $i_1<i_2<...$ be the indices $j$ such that $S' \cap V(T_j') \neq \emptyset$.
For each $j \geq 1$, define $T_j''$ to be the component of $T_{i_j}-F_{i_j}$ contracted into the node in $S' \cap V(T_{i_j}')$.
Define $D''=D[V(T_1'' \cup T_2'' \cup ...)]$.

If there exists an infinite stable set $S''$ of $D''$ such that $\lvert S'' \cap V(T_i'') \rvert \leq 1$ for each $i \geq 1$ and such that the set of heads of all edges of $T_1'' \cup T_2'' \cup ...$ with tails in $S''$ is rich in $D''$, then there exists an infinite stable set $S$ of $D$ such that $\lvert S \cap V(T_i) \rvert \leq 1$ for each $i \geq 1$ and such that the set of heads of all edges of $T_1 \cup T_2 \cup ...$ with tails in $S$ is rich in $D$.
\end{lemma}

\begin{pf}
We claim that we can choose $S$ to be $S''$.
Clearly, $\lvert S'' \cap V(T_i) \rvert = \lvert S'' \cap V(T_i'') \rvert \leq 1$ for each $i \geq 1$.
Suppose that the set of heads of all edges of $T_1 \cup T_2 \cup ...$ with tails in $S$ contains an infinite stable set $R$ in $D$.
Since $R \cap V(T_1'' \cup T_2'' \cup ...)$ is finite by assumption, the set $R-V(T_1'' \cup T_2'' \cup ...)$, denoted by $R'$, is an infinite stable set in $D$.
Since each node in $R'$ is the root of some component of $T_i-F_i$ contracted into a child in $T_i'$ of a node in $S'$, we obtain an infinite subset $R''$ of heads of the set of all edges of $T_1' \cup T_2' \cup ...$ with tails in $S'$ such that $R''$ is stable in $D'$ by the definition of $D'$, a contradiction.
This proves the lemma.
\end{pf}

\begin{lemma} \label{tree lemma same mu simple}
Let $(D,(T_i,\phi_i,\tau_i,\mu_i)_{i \in {\mathbb N}},N,h,d,L)$ be a standard tuple for Theorem \ref{decorated tree lemma}.
For each $i \geq 1$, let $W_i$ be a subset of $\{e \in E(T_i): \lvert \phi_i(e) \rvert= \min_{e' \in E(T_i)} \lvert \phi_i(e') \rvert\}$ and let $(T_i^W,\phi_i^W,\tau_i^W,\mu_i^W)$ to be the $W_i$-contraction of $(T_i,\phi_i, \allowbreak \tau_i,\mu_i)$.
If for every $i \geq 1$, there exists a nonnegative integer $p_i$ such that $\mu_i(e)=p_i$ for all edges $e$ of $T_i$, then for any nodes $v$ and $w$ of $T_i^W$, $v$ precedes $w$ in $T_i^W$ with respect to $(\phi_i^W,\tau_i^W,\mu_i^W)$ if and only if the root of the component of $T_i-W_i$ contracted to $v$ precedes the root of the component of $T_i-W_i$ contracted to $w$ in $T_i$ with respect to $(\phi_i,\tau_i,\mu_i)$. 
\end{lemma}

\begin{pf}
This lemma immediately follows from the assumption that for each $i \geq 1$, $W_i$ is a subset of $\{e \in E(T_i): \lvert \phi_i(e) \rvert= \min_{e' \in E(T_i)} \lvert \phi_i(e') \rvert\}$ and $\mu_i(e)=p_i$ for all edges $e$ of $T_i$.
\end{pf}

\begin{lemma} \label{tree lemma disjoint labels same size same mu}
Let $(D,(T_i,\phi_i,\tau_i,\mu_i)_{i \in {\mathbb N}},N,h,d,L)$ be a standard tuple for Theorem \ref{decorated tree lemma}.
If for each $i \geq 1$, $\phi_i(e') \cap \phi_i(e'') = \emptyset$ for every pair of distinct edges $e',e''$ of $T_i$, and there exist nonnegative integers $h_i,p_i$ with $\lvert \phi_i(e) \rvert=h_i$ and $\mu_i(e)=p_i$ for all $e \in E(T_i)$, then there exists an infinite stable set $S$ of $D$ such that $\lvert S \cap V(T_i) \rvert \leq 1$ for each $i \geq 1$ and such that the set of heads of all edges of $T_1 \cup T_2 \cup ...$ with tails in $S$ is rich in $D$.
\end{lemma}

\begin{pf}
For each $i \geq 1$, we define $F_i=\{e \in E(T_i): \tau_i(e)=\emptyset\}$ and define $(T_i',\phi_i',\tau_i',\mu_i')$ to be the $F_i$-contraction of $(T_i,\phi_i,\tau_i,\mu_i)$.
Define $D'$ to be the $(\bigcup_{i \geq 1}F_i)$-contraction of $D$. 

For each $i$, since $\lvert \phi_i'(e) \rvert=h_i$ and $\mu_i'(e)=p_i$ for all $e \in E(T_i')$, by Lemma \ref{tree lemma same mu simple}, if $x \in V(T_i'), y \in V(T_j')$ with $i<j$ and $x$ is adjacent to $y$ in $D'$, then $x$ is adjacent to all non-root ancestors of $y$ in $T_{j}'$.
By Theorem \ref{rs weak tree lemma}, there exists an infinite stable set $S' \subseteq V(D')$ in $D'$ such that $\lvert S' \cap V(T_i') \rvert \leq 1$ for each $i \geq 1$ and such that the set of the children of the members of $S'$ is rich in $D'$.

Let $i_1<i_2<...$ be the indices $j$ such that $S' \cap V(T_j') \neq \emptyset$.
For each $j \geq 1$, define $T_j''$ to be the component of $T_{i_j}-F_{i_j}$ contracted into the node in $S' \cap V(T_{i_j}')$.
Define $D''=D[V(T_1'' \cup T_2'' \cup ...)]$.

Note that for each $i \geq 1$ and each edge $e$ in a component of $T_i-F_i$, $\tau_i(e) \neq \emptyset$. 
Since $(T_i,\phi_i,\tau_i,\mu_i)$ is $(h,d,N)$-decorated and the sets $\phi_i(e)$ (for $e \in E(T_i)$) are pairwise disjoint sets with the same size and $\mu_i(e)$ (for $e \in E(T_i)$) is a constant, each component of $T_i-F_i$ has no directed path with length $d+1$.
By Lemma \ref{tree lemma bounded height}, there exists an infinite stable set $S$ in $D''$ such that $\lvert S \cap V(T_i'') \rvert \leq 1$ for all $i \geq 1$ and the set of heads of all edges of $T_1'' \cup T_2'' \cup ...$ with tails in $S$ is rich in $D''$.
Then the lemma follows from Lemma \ref{still rich}.
\end{pf}

\begin{lemma} \label{tree lemma same mu}
Let $(D,(T_i,\phi_i,\tau_i,\mu_i)_{i \in {\mathbb N}},N,h,d,L)$ be a standard tuple for Theorem \ref{decorated tree lemma}.
If for every $i \geq 1$, there exists a nonnegative integer $p_i$ such that $\mu_i(e)=p_i$ for all edges $e$ of $T_i$, then there exists an infinite stable set $S$ of $D$ such that $\lvert S \cap V(T_i) \rvert \leq 1$ for each $i \geq 1$ and such that the set of heads of all edges of $T_1 \cup T_2 \cup ...$ with tails in $S$ is rich in $D$.
\end{lemma}

\begin{pf}
Let $b=\min\{\lvert \phi_i(e) \rvert: i \geq 1, e \in E(T_i)\}$.
Note that $b \geq 0$ exists.
We shall prove this lemma by induction on the lexicographic order of $(h,h-b)$.
When $h=0$, $\phi_i(e)=\emptyset$ for all $i \geq 1$ and $e \in E(T_i)$, so every non-root node precedes all its descendants, so this lemma follows from Theorem \ref{rs weak tree lemma}.
Hence we may assume that $h \geq 1$ and this lemma holds for all pairs lexicographically smaller $(h,h-b)$.

\medskip

\noindent{\bf Claim 1:} We may assume that $\lvert \phi_i(e) \rvert$ is a constant for all $i \geq 1$ and $e \in E(T_i)$.

\noindent{\bf Proof of Claim 1:}
For each $i \geq 1$, we define $W_i=\{e \in E(T_i): \lvert \phi_i(e) \rvert=b\}$ and define $(T_i^W,\phi_i^W,\tau_i^W,\mu_i^W)$ to be the $W_i$-contraction of $(T_i,\phi_i,\tau_i,\mu_i)$.
Define $D_W$ to be the $(\bigcup_{i \geq 1}W_i)$-contraction of $D$. 
 
For each $i$, since $\lvert \phi_i^W(e') \rvert=b$ for all $e' \in E(T_i^W)$ and $\mu_i(e)=p_i$ for all $e \in E(T_i)$, by Lemma \ref{tree lemma same mu simple}, if $x \in V(T_i^W), y \in V(T_j^W)$ with $i<j$ and $x$ is adjacent to $y$ in $D_W$, then $x$ is adjacent to all nodes of $T_j^W$ preceding $y$ in $D_W$.

Note that for every $i \geq 1$ and for every edge $e$ contained in a component of $T_i-W_i$, $\lvert \phi_i(e) \rvert \geq b+1$.
So the induction hypothesis and Lemma \ref{still rich} imply that it suffices to prove that the lemma holds for the standard tuple $(D_W, (T_i^W,\phi_i^W,\tau_i^W,\mu_i^W)_{i \in {\mathbb N}},N,h,d,L)$ for Theorem \ref{decorated tree lemma}.

Since $\lvert \phi_i(e) \rvert=b$ for all $i \geq 1$ and $e \in E(T_i^W)$, it suffices to prove this lemma for the standard tuple $(D,(T_i,\phi_i,\tau_i,\mu_i)_{i \in {\mathbb N}}, \allowbreak N,h,d,L)$ for Theorem \ref{decorated tree lemma} with the extra assumption that $\lvert \phi_i(e) \rvert$ is a constant for all $i \geq 1$ and $e \in E(T_i)$.
$\Box$

\medskip

\noindent{\bf Claim 2:} We may assume that for each $i \geq 1$, there exists a nonempty set $Q_i \in [L]^{\leq h}$ such that 
	\begin{itemize}
		\item[(i)] $\lvert \phi_i(e) \rvert=\lvert Q_i \rvert$ is a constant and $\phi_i(e) \cap Q_i \neq \emptyset$ for all $e \in E(T_i)$, and
		\item[(ii)] if $e,e'$ are distinct edges appearing in a directed path in $T_i$ in the order listed, then $\phi_i(e') \cap Q_i \subseteq \phi_i(e) \cap Q_i$.
	\end{itemize}

\noindent{\bf Proof of Claim 2:}
For each $i \geq 1$, we define $F_i$ to be a maximal subset of $E(T_i)$ with the following properties.
	\begin{itemize}
		\item Every edge incident with the root of $T_i$ belongs to $F_i$.
		\item If $e_1,e_2$ are distinct elements in $F_i$, and $P$ is a directed path in $T_i$ with $E(P) \cap F_i = \{e_1,e_2\}$ such that $e_1$ is incident with the source of $P$ and $e_2$ is incident with the sink of $P$, then $\phi_i(e_1) \cap \phi_i(e_2)=\emptyset$, and for every $e \in E(P)-\{e_2\}$, $\phi_i(e) \cap \phi_i(e_1) \neq \emptyset$. 
	\end{itemize}
So for each $i \geq 1$, the sets $\phi_i(e)$ (for $e \in F_i$) are pairwise disjoint.
Since $\lvert \phi_i(e) \rvert$ is a constant for all $i$ and $e$, by Lemmas \ref{tree lemma same mu simple} and \ref{tree lemma disjoint labels same size same mu}, there is an infinite stable set of the $(\bigcup_{i \geq 1}F_i)$-contraction of $D$ intersecting each $F_i$-contraction of $T_i$ in at most one node such that the set of heads of all tree edges of the $F_i$-contractions of $T_i$ (for all $i \geq 1$) with tails in this stable set is rich in the $(\bigcup_{i \geq 1} F_i)$-contraction of $D$.
By Lemma \ref{still rich}, to prove this lemma, we may restrict the problem to the components of $T_i-F_i$.
Note that for each $i \geq 1$, the component of $T_i-F_i$ containing the root of $T_i$ has only one node; for each component $C$ of $T_i-F_i$ not containing the root of $T_i$, we know that $\phi_i(e) \cap \phi_i(e_C) \neq \emptyset$, where $e_C$ is the edge in $F_i \subseteq E(T_i)$ having the root of $C$ as its head, and if $e,e'$ are distinct edges appearing in a directed path in $C$ in the order listed, then $e_C,e,e'$ are distinct edges appearing in a directed path in $T_i$ in the order listed, so $\phi_i(e') \cap \phi_i(e_C) \subseteq \phi_i(e) \cap \phi_i(e_C)$ since $(T_i,\phi_i,\tau_i,\mu_i)$ is $(h,d,N)$-decorated.

In other words, to prove this lemma, it suffices to prove that this lemma holds for the standard tuple $(D,(T_i,\phi_i,\tau_i,\mu_i)_{i \in {\mathbb N}}, \allowbreak N,h,d,L)$ for Theorem \ref{decorated tree lemma} with the claimed extra assumption.
$\Box$

\medskip

\noindent{\bf Claim 3:} We may assume that for each $i \geq 1$, there exist $h_i \in [h]$, a nonempty set $Q_i \in [L]^{\leq h}$ and a nonempty set $X_i$ such that $\lvert \phi_i(e) \rvert = h_i$ and $\phi_i(e) \cap Q_i=X_i$ for all $e \in E(T_i)$.

\noindent{\bf Proof of Claim 3:}
For each $i \geq 1$, we define $F_i'$ to be a maximal subset of $E(T_i)$ with the following properties.
	\begin{itemize}
		\item Every edge of $T_i$ incident with the root of $T_i$ belongs to $F_i'$. 
		\item If $e_1$ and $e_2$ are distinct elements of $F_i'$, and $P$ is a directed path in $T_i$ with $E(P) \cap F_i'=\{e_1,e_2\}$ such that $e_1$ is incident with the source of $P$ and $e_2$ is incident with the sink of $P$, then $\phi_i(e_1) \cap Q_i \neq \phi_i(e_2) \cap Q_i$, and for every $e \in E(P)-\{e_2\}$, $\phi_i(e) \cap Q_i = \phi_i(e_1) \cap Q_i$.  
	\end{itemize}
By (ii) in Claim 2, for each $i \geq 1$ and for each directed path $P$ in $T_i$, the sets $\phi_i(e) \cap Q_i$ (for $e \in F_i' \cap E(P)$) are pairwise distinct.

For each $i \geq 1$, define $(T_i',\phi_i',\tau_i',\mu_i')$ to be the $F_i'$-contraction of $(T_i,\phi_i,\tau_i,\mu_i)$.
Since there are at most $2^h$ different subsets of $Q_i$, there exists no directed path in $T_i'$ with length $2^h+1$.
So by Lemma \ref{tree lemma bounded height}, there exists an infinite stable set of the $(\bigcup_{i \geq 1}F_i')$-contraction of $D$ intersecting each $F_i'$-contraction of $T_i$ in at most one node such that the set of heads of all tree edges of the $F_i'$-contractions of $T_i$ (for all $i \geq 1$) with tails in this stable set is rich in the $(\bigcup_{i \geq 1} F_i')$-contraction of $D$.
By Lemma \ref{still rich}, to prove this lemma, we may restrict the problem to the components of $T_i-F_i'$.

In other words, to prove this lemma, it suffices to show that this lemma holds for the tuple $(D,(T_i,\phi_i,\tau_i,\mu_i)_{i \in {\mathbb N}}, \allowbreak N,h,d,L)$ with the claimed extra assumption.
$\Box$

\medskip

\noindent{\bf Claim 4:} We may assume that for each $i \geq 1$, there exist $h_i \in [h]$ and an element $x_i$ such that 
	\begin{itemize}
		\item $\lvert \phi_i(e) \rvert = h_i$ and $x_i \in \phi_i(e)$ for every $e \in E(T_i)$, and 
		\item either $x_i \in \tau_i(e)$ for all $e \in E(T_i)$, or $x_i \not \in \tau_i(e)$ for all $e \in E(T_i)$.
	\end{itemize}

\noindent{\bf Proof of Claim 4:}
For each $i \geq 1$, let $z_i$ be an element of $X_i$, define $F_i''=\{e \in E(T_i): z_i \in \tau_i(e)\}$, and define the $F_i''$-contraction of $(T_i,\phi_i,\tau_i,\mu_i)$ to be $(T_i'',\phi_i'',\tau_i'',\mu_i'')$.
Since for each $i \geq 1$, $\lvert \phi_i(e) \rvert=h_i$ for all $e \in E(T_i)$ by Claim 3, a node $v$ of $T_i''$ precedes a node $w$ of $T_i''$ with respect to $(\phi_i'',\tau_i'',\mu_i'')$ if and only if the root of the component of $T_i-F_i''$ contracted into $v$ precedes the root of the component of $T_i-F_i''$ contracted into $w$ in $T_i$ with respect to $(\phi_i,\tau_i,\mu_i)$ by Lemma \ref{tree lemma same mu simple}.

Hence by Lemma \ref{still rich}, we may reduce the problem to the one for the standard tuple given by the $F_i''$-contraction and to the one for the standard tuple given by the components of $T_i-F_i''$.
So to prove this lemma, it suffices to show that this lemma holds for the tuple $(D,(T_i,\phi_i,\tau_i,\mu_i)_{i \in {\mathbb N}}, \allowbreak N,h,d,L)$ with the claimed extra assumption.
$\Box$

\medskip

For each $i \geq 1$ and each edge $e \in E(T_i)$, define $\phi_e^*(e)=\phi_i(e)-\{x_i\}$, $\tau_i^*(e)=\tau_i(e)-\{x_i\}$.
It is straightforward to verify that for each $i \geq 1$, $(T_i,\phi_i^*,\tau_i^*,\mu_i)$ is $(h-1,d,N)$-decorated, and for distinct nodes $v,w \in T_i$, $v$ precedes $w$ in $T_i$ with respect to $(\phi_i,\tau_i,\mu_i)$ if and only if $v$ precedes $w$ in $T_i$ with respect to $(\phi_i^*,\tau_i^*,\mu_i)$.
Then this lemma immediately follows from the induction hypothesis.
\end{pf}

\subsection{Main tree lemma} \label{subsec:main_tree_lemma}

The following is the main theorem of this section.

\begin{theorem} \label{decorated tree lemma}
Let $T_1,T_2,...$ be an infinite sequence of pairwise disjoint rooted trees.
Let $N,h,d$ be integers with $N,h \geq 0$ and $d > 0$, and let $L$ be a set.
For each $i \geq 1$, let $\phi_i,\tau_i$ be functions from $E(T)$ to $[L]^{\leq h}$ and $\mu_i$ a function from $E(T)$ to $\{0,1,...,N\}$ such that $(T_i,\phi_i,\tau_i,\mu_i)$ is $(h,d,N)$-decorated.
Assume that $D$ is an infinite graph with $V(D) = V(T_1 \cup T_2 \cup ...)$ such that for $i'>i \geq 1$, if $u \in V(T_i)$ is adjacent to $w \in V(T_{i'})$ in $D$, and $v \in V(T_{i'})$ precedes $w$ in $T_{i'}$ with respect to $(\phi_{i'},\tau_{i'},\mu_{i'})$, then $u$ is adjacent to $v$ in $D$.
If the roots of $T_1,T_2,...$ form a stable set in $D$, then there exists an infinite stable set $S$ of $D$ such that $\lvert S \cap V(T_i) \rvert \leq 1$ for each $i \geq 1$ and such that the set of heads of all edges of $T_1 \cup T_2 \cup ...$ with tails in $S$ is rich in $D$.
\end{theorem}

\begin{pf}
Let $b=\min\{\mu_i(e): i \geq 1, e \in E(T_i)\}$, and let $c=\min\{\lvert \phi_i(e) \rvert: i \geq 1, e \in E(T_i)\}$.
We shall prove this theorem by induction on $N+h-b-c$.
Note that $b \leq N$ and $c \leq h$, so $N+h-b-c \geq 0$.
When $b=N$, we know $\mu_i(e)=N$ for any $i \geq 1$ and $e \in E(T_i)$, so this theorem follows from Lemma \ref{tree lemma same mu} for any $c$.
So we may assume $b<N$ and assume that this lemma holds when $N+h-b-c$ is smaller.

For each $i \geq 1$, we define $F_i=\{e \in E(T_i): \mu_i(e)=b,\lvert \phi_i(e) \rvert=c\}$ and define $(T_i',\phi_i',\tau_i',\mu_i')$ to be the $F_i$-contraction of $(T_i,\phi_i,\tau_i,\mu_i)$.
Define $D'$ to be the $(\bigcup_{i \geq 1}F_i)$-contraction of $D$. 

For each $i \geq 1$, since $\lvert \phi_i'(e) \rvert=c$ and $\mu_i(e) = b$ for all $e \in E(T_i')$, a node $v$ of $T_i'$ precedes a node $w$ of $T_i'$ in $T_i'$ with respect to $(\phi_i',\tau_i',\mu_i')$ if and only if the root of the component of $T_i-F_i$ contracted into $v$ precedes the root of the component of $T_i-F_i$ contracted into $w$ in $T_i$ with respect to $(\phi_i,\tau_i,\mu_i)$.
Hence, if $x \in V(T_i'), y \in V(T_j')$ with $i<j$ and $x$ is adjacent to $y$ in $D'$, then $x$ is adjacent in $D'$ to all nodes of $T_{j}'$ preceding $y$.
By Lemma \ref{tree lemma same mu}, there exists an infinite stable set $S' \subseteq V(D')$ in $D'$ such that $\lvert S' \cap V(T_i') \rvert \leq 1$ for each $i \geq 1$ and such that the set of the children of the members of $S'$ is rich in $D'$.

Let $i_1<i_2<...$ be the indices $j$ such that $S' \cap V(T_j') \neq \emptyset$.
For each $j \geq 1$, define $T_j''$ to be the component of $T_{i_j}-F_{i_j}$ contracted into the node in $S' \cap V(T_{i_j}')$, $\phi_i'' = \phi_i|_{E(T_j'')}$, $\tau_i'' = \tau_i|_{E(T_j'')}$ and $\mu_i'' = \mu_i|_{E(T_j'')}$.
Define $D''=D[V(T_1'' \cup T_2'' \cup ...)]$.

Then by Lemma \ref{still rich}, it suffices to show that there exists an infinite stable set $S''$ in $D''$ such that $\lvert S'' \cap V(T_i'') \rvert \leq 1$ for all $i \geq 1$ and the set of heads of all edges of $T_1'' \cup T_2'' \cup ...$ with tails in $S''$ is rich in $D''$.

For each $i \geq 1$, we define $W_i=\{e \in E(T_i''): \lvert \phi_i(e) \rvert=c\}$ and define $(T_i^W,\phi_i^W,\tau_i^W,\mu_i^W)$ to be the $W_i$-contraction of $(T''_i,\phi''_i,\tau''_i,\mu''_i)$.
Define $D^W$ to be the $(\bigcup_{i \geq 1}W_i)$-contraction of $D''$. 

By the definition of $W_i$, for each $i \geq 1$, a node $v$ of $T_i^W$ precedes a node $w$ of $T_i^W$ in $T_i^W$ with respect to $(\phi_i^W,\tau_i^W,\mu_i^W)$ if and only if the root of the component of $T_i''-W_i$ contracted into $v$ precedes the root of the component of $T_i''-W_i$ contracted into $w$ in $T_i''$ with respect to $(\phi_i'',\tau_i'',\mu_i'')$.
So if $x \in V(T_i^W), y \in V(T_j^W)$ with $i<j$ and $x$ is adjacent to $y$ in $D^W$, then $x$ is adjacent in $D^W$ to all nodes of $T_{j}^W$ preceding $y$.
Note that for each $i \geq 1$ and each edge $e$ in a component of $T_i-F_i$, $\mu_i(e)+\lvert \phi_i(e) \rvert \geq b+c+1$. 
So for each $i \geq 1$ and $e \in W_i$, $\mu_i(e) \geq b+1$, so $\min\{\mu_i^W(e): i \geq 1, e \in E(T_i^W)\} + \min\{\lvert \phi_i^W(e) \rvert: i \geq 1, e \in E(T_i^W)\} \geq (b+1)+c$.
Hence by the induction hypothesis, there exists an infinite stable set $S^W \subseteq V(D^W)$ in $D^W$ such that $\lvert S^W \cap V(T_i^W) \rvert \leq 1$ for each $i \geq 1$ and such that the set of the children of the members of $S^W$ is rich in $D^W$.

Let $i^W_1<i^W_2<...$ be the indices $j$ such that $S^W \cap V(T_j^W) \neq \emptyset$.
For each $j \geq 1$, define $T_j^*$ to be the component of $T_{i_j}''-W_{i_j}$ contracted into the node in $S^W \cap V(T_{i_j}^W)$, $\phi_i^* = \phi_i|_{E(T_j^*)}$, $\tau_i^* = \tau_i|_{E(T_j^*)}$ and $\mu_i^* = \mu_i|_{E(T_j^*)}$.
Define $D^*=D[V(T_1^* \cup T_2^* \cup ...)]$.

Note that for each $i \geq 1$ and each edge $e$ of $T_i^*$, we know $e \not \in W_i$, so $\lvert \phi^*_i(e) \rvert \geq c+1$. 
Hence $\min\{\mu_i^*(e): i \geq 1, e \in E(T_i^*)\} + \min\{\lvert \phi_i^*(e) \rvert: i \geq 1, e \in E(T_i^*)\} \geq b+(c+1)$.
By the induction hypothesis, there exists an infinite stable set $S$ in $D^*$ such that $\lvert S \cap V(T_i^*) \rvert \leq 1$ for all $i \geq 1$ and the set of heads of all edges of $T_1^* \cup T_2^* \cup ...$ with tails in $S$ is rich in $D^*$.
Then this theorem follows from Lemma \ref{still rich}.
\end{pf}

\section{Assemblages, encodings and simulations}
\label{sec: new wqo bounded depth}

The goal of this section is to show how to reduce the well-quasi-ordering problem with respect to the topological minor relation for graphs with given tree-decompositions to the one that focuses on their bags. 
We need a number of new terminologies to achieve this goal.
We will provide the intuition about those terminologies in Section \ref{subsec:assemblages_intuition} and provide the formal description in the remaining subsections.
We first mention well-known results about well-quasi-ordering that we will use in this paper in Section \ref{subsec:easy_wqo}.

\subsection{Preliminary for well-quasi-ordering} \label{subsec:easy_wqo}

We say that $(S,\preceq)$ is a \defn{well-quasi-ordered set} if $\preceq$ is a well-quasi-ordering on $S$.
Note that if $(S_1,\preceq_1)$ and $(S_2, \preceq_2)$ are two well-quasi-ordered sets, then $S_1 \times S_2$ is well-quasi-ordered by $\preceq_3$, where $(s_1,s_2) \preceq_3 (s_1',s_2')$ if and only if $s_1 \preceq s_1'$ and $s_2 \preceq s_2'$.
We call $(S_1 \times S_2, \preceq_3)$ the \defn{well-quasi-ordered set obtained from $(S_1,\preceq_1), (S_2,\preceq_2)$ by Cartesian product}, and denote it by $(S_1 \times S_2, \preceq_1 \times \preceq_2)$.
For any two sets $A,B$, we define $A \uplus B$ to be the union of $A$ and a disjoint copy of $B$.
Then $S_1 \uplus S_2$ is well-quasi-ordered by $\preceq_4$, where $s \preceq_4 s'$ if and only if either $s,s' \in S_1$ and $s \preceq_1 s'$ , or $s,s' \in S_2$ and $s \preceq_2 s'$.
We call $(S_1 \uplus S_2, \preceq_4)$ the \defn{well-quasi-ordered set obtained from $(S_1,\preceq_1),(S_2,\preceq_2)$ by disjoint union}.

The following theorem was proved by Higman and gave another way to obtain another well-quasi-ordered set from a well-quasi-ordered set. 

\begin{theorem}[\cite{h}] \label{Higman's lemma}
Let $(S,\preceq)$ be a well-quasi-ordered set.
For every finite sequences $A=(a_1,a_2,...,a_n)$ and $B=(b_1,b_2,...,b_m)$ over $S$, we say that $A \preceq' B$ if there exist $1 \leq i_1 < i_2 < ... < i_n \leq m$ such that $a_j \preceq b_{i_j}$ for every $j \in [n]$.
Then the finite sequences over $S$ are well-quasi-ordered by $\preceq'$.
\end{theorem}

We call the new well-quasi-ordered set mentioned in Theorem \ref{Higman's lemma} the \defn{well-quasi-ordered set obtained from $(S,\preceq)$ by Higman's lemma}.

\subsection{Intuition} \label{subsec:assemblages_intuition}
To apply the minimal bad sequence argument to prove well-quasi-ordering results based on the tree-structure of the graphs given by their rooted tree-decompositions, we need two main properties for the tree-decomposition: ``linkedness property'' and ``absorption property''.
Those properties are simply conceptual and we will not include a precise definition for them.
We have obtained a tree-decomposition with the above two properties in the previous sections.
The goal for this section is to show how to use them to prove well-quasi-ordering results.

The linkedness property roughly says that given a rooted tree-decomposition of a graph $G$, ``whenever'' a tree node $a$ is an ancestor of another tree node $b$, we can find a homeomorphic embedding from the subgraph of $G$ induced by the union of all bags at the descendants of $b$ into the subgraph of $G$ induced by the union of all bags at the descendants of $a$, and this homeomorphic embedding ``preserves the roots''.
Here the ``roots'' mean the common vertices in the bags at $a$ and its parent (and common vertices in the bags at $b$ and its parent, respectively).
By a homeomorphic embedding ``preserving the roots'', we roughly mean that the roots of the first graph are mapped to the roots of the second graph.
So we consider rooted graphs, which consist of a graph and a sequence of vertices with no repeated entries.
This sequence is called the root march.
On the other hand, to use the linkedness property to prove our well-quasi-ordering result with respect to the topological minor relation, we have to make it more flexible in the sense that we allow some roots of the first graph to be mapped to disjoint paths in the second graph each containing a root of the second graph instead of being mapped to roots.
To handle it, we have to declare which roots are allowed to be flexible.
It is the motivation of the ``essential number'' associated with each vertex in the root march.
Then we consider the ``rooted extension'' of a rooted graph, which is obtained by adding a copy for each root and adding a certain number of edges between each root and its copy according to the essential number.
By doing so, we can describe the aforementioned flexible homeomorphic embedding between rooted graphs in terms of the usual homeomorphic embedding between rooted extensions.
We formally define rooted graphs and related notions in Section \ref{subsec:rooted_graphs}

The absorption property roughly says that given a rooted tree-decomposition of a graph $G$ and a tree node $t$, we can ``encode'' the subgraph of $G$ induced by the union of all bags at the descendants of $t$ into the subgraph of $G$ induced by the bag at $t$ such that by simply seeing the encodings of two given graphs at their root bags in their tree-decomposition, we can decide whether one graph is a topological minor of another.
In other words, whenever one encoding is a topological minor of another encoding, we should be able to construct a homeomorphic embedding from the entire first graph to the entire second graph.
The main technicality in this section arises from this part.

First, the encoding has to record, given a tree node $t$, what vertices in the bag at $t$ are also contained in the bags at its children.
This leads to the notion of an assemblage, which consists of a rooted graph and a multiset of marches.
Then we define homeomorphic embeddings between (labelled) assemblages, which are called ''simulations'' in later sections.
We formally define assemblages and simulations in Section \ref{subsec:assemblages}.

Second, given a tree node $t$ and its child $c$, the encoding has to record the subgraph induced by the union of all bags at the descendant of $c$.
We should treat this subgraph as a rooted graph, where the roots are the common vertices in the bags at $t$ and $c$.
We call this subgraph ``the branch at $c$''.
The essential number associated with a root of the branch, which represents the ``flexibility'' mentioned above when discussing the linkedness property, involves the flexibility of this vertex in the original tree-decomposition and involves whether it also belongs to a bag at a descendant of another child of $t$.
We address the formal definition of branches in Section \ref{subsec:branches} and prove that our definition for branches are ``well-defined'' in the sense that a branch of a branch is also a branch.

We formally define encodings in Section \ref{subsec:encoding}.
Roughly speaking, an encoding is an assemblage that records the aforementioned information about branches. 
The main technical result in this section (Lemma \ref{encoding simulation}) is that the simulation relation between encodings can recover a homeomorphic embedding between the original graphs.
A proof sketch for Lemma \ref{encoding simulation} will be provided later.

\subsection{Rooted graphs} \label{subsec:rooted_graphs}

A \defn{march} in a graph is either the empty set or a sequence of distinct vertices of the graph such that each entry is associated with a number in $\{0,1,2\}$, called the \defn{essential number}.

We say that $(G,\gamma)$ is a \defn{rooted graph} if $G$ is a graph and $\gamma$ is a march in $G$. 
In this case, we say that $\gamma$ is the \defn{root march} of $(G,\gamma)$.
We denote the set of entries of $\gamma$ by \defn{$V(\gamma)$}.
Let $\gamma=(v_1,v_2,...,v_k)$ and assume that $v_i$ is associated with essential number $a_i$ in $\gamma$ for each $i \in [k]$.
The \defn{rooted extension} of $(G,\gamma)$ is the graph $G'$ obtained from $G$ by 
	\begin{itemize}
		\item adding vertices $u_1,u_2,...,u_k$, and 
		\item for each $i \in [k]$, adding $a_i$ parallel edges between $v_i$ and $u_i$.
	\end{itemize}
We call the sequence $(u_1,u_2,...,u_k)$ the \defn{indicator} of $G'$.

Recall that if $f$ is a function and $\sigma=(x_1,x_2,...,x_n)$ is a sequence whose entries are in the domain of $f$, then we define $f(\sigma)=(f(x_1),f(x_2),...,f(x_n))$.

Let $(G_1,\gamma_1)$ and $(G_2,\gamma_2)$ be rooted graphs.
Let $G_1',G_2'$ be the rooted extensions of $(G_1,\gamma_1),(G_2,\gamma_2)$ with indicators $I_1,I_2$, respectively.
We say that $\eta$ is a \defn{homeomorphic embedding} from $(G_1,\gamma_1)$ to $(G_2,\gamma_2)$ if the following hold.
	\begin{itemize}
		\item $\gamma_1$ and $\gamma_2$ have the same length, and for every $i \in [\lvert V(\gamma_1) \rvert]$, the essential number associated with the $i$-th vertex in $\gamma_1$ equals the essential number associated with the $i$-th vertex in $\gamma_2$.
		\item $\eta$ is a homeomorphic embedding from $G_1'$ to $G_2'$ such that $\eta(I_1)=I_2$.
		\item If $v$ is the $i$-th vertex in $\gamma_2$ for some $i \in [\lvert V(\gamma_2) \rvert]$, and $v$ is a vertex in $\eta(e)-\eta(V(G_1'))$ for some edge $e$ of $G_1'$, then either 
			\begin{itemize}
				\item $e$ is an edge incident with the $i$-th vertex in $I_1$, or
				\item the essential number associated with the $i$-th vertex in $\gamma_1$ is 0, and $e$ is an edge incident with the $i$-th vertex in $\gamma_1$.
			\end{itemize}
		\item If $v$ is the $i$-th vertex in $\gamma_2$ for some $i \in [\lvert V(\gamma_2) \rvert]$, and $v=\eta(v')$ for some vertex $v'$ of $G_1'$, then $v'$ is the $i$-th vertex in $\gamma_1$.
	\end{itemize}
In this case, we say that $(G_1,\gamma_1)$ is a \defn{rooted topological minor} of $(G_2,\gamma_2)$ and write $\eta: (G_1,\gamma_1) \hookrightarrow (G_2,\gamma_2)$.

\subsection{Assemblages and the simulation relation} \label{subsec:assemblages}

We say that a tuple $(G,\gamma_0,\Gamma)$ is an \defn{assemblage} if $(G,\gamma_0)$ is a rooted graph and $\Gamma$ is a finite multiset of marches in $G$.
We also call $\gamma_0$ the \defn{root march} of the assemblage $(G,\gamma_0,\Gamma)$.

For simplicity of notations, for a quasi-order $Q$, we also denote its ground set by $Q$, and denote the relation by $\leq_Q$.

Let $Q$ be a quasi-order.
We say that $(G,\gamma_0,\Gamma,f,\phi)$ is a \defn{$Q$-assemblage} if 
	\begin{itemize}
		\item $(G,\gamma_0,\Gamma)$ is an assemblage, and 
		\item $f: \Gamma \rightarrow Q$ and $\phi: V(G) \rightarrow Q$ are functions.
	\end{itemize}
We call $(G,\gamma_0,\Gamma)$ the \defn{underlying assemblage} of the $Q$-assemblage $(G,\gamma_0,\Gamma,f,\phi)$.
We say that a $Q$-assemblage $(G',\gamma_0',\Gamma',f',\phi')$ \defn{simulates} $(G,\gamma_0,\Gamma,f,\phi)$ if there exist $\eta: (G,\gamma_0) \hookrightarrow (G',\gamma_0')$ and an injection $\iota: \Gamma \rightarrow \Gamma'$ such that
	\begin{itemize}
		\item $\phi(v) \leq_Q \phi'(\eta(v))$ for every $v \in V(G)$, and 
		\item $\eta(\sigma)=\iota(\sigma)$ and $f(\sigma) \leq_Q f'(\iota(\sigma))$ for every $\sigma \in \Gamma$.
	\end{itemize}
In this case, we write \defn{$(G,\gamma_0,\Gamma,f,\phi) \preceq (G',\gamma_0',\Gamma',f',\phi')$}.
We call $\preceq$ the \defn{simulation relation}.

\subsection{Rooted tree-decomposition of assemblages and branches} \label{subsec:branches}

We say that $(T,\X,\alpha)$ is a \defn{rooted tree-decomposition} of an assemblage $(G,\gamma_0,\Gamma)$ if $(T,\X)$ is a rooted tree-decomposition of $G$ such that $V(\gamma_0)$ is contained in the bag of the root of $T$, and $\alpha$ is a function from $\Gamma$ to $V(T)$ such that $V(\sigma) \subseteq X_{\alpha(\sigma)}$ for each $\sigma \in \Gamma$. 
A \defn{rooted tree-decomposition} of a $Q$-assemblage (for some quasi-order $Q$) is a rooted tree-decomposition of its underlying assemblage.

Let $Q$ be a quasi-order and let $(T,\X,\alpha)$ be a rooted tree-decomposition of a $Q$-assemblage $(G,\gamma_0,\Gamma,f,\phi)$.
Let $t$ be a non-root node of $T$ and let $p$ be the parent of $t$.
Assume that there exists an ordering $\pi_t$ on $X_t \cap X_p$, and assume that $\gamma_{p}$ is defined, where $\gamma_p=\gamma_0$ if $p$ is the root of $T$.
We define $\gamma_t,\Gamma_t,f_t,\phi_t,S_t$ as follows.
	\begin{itemize}
		\item $\gamma_t$ is the march such that $V(\gamma_t)=X_t \cap X_p$ with the ordering $\pi_t$, where for each vertex $v \in V(\gamma_t)$, the essential number $j$ associated with $v$ is defined as follows.
			\begin{itemize}
				\item[(BR0)] $j=0$, if the following hold.
					\begin{itemize}
						\item[(BR01)] $v \not \in V(\sigma)$ for every $\sigma \in \Gamma$ in which $\alpha(\sigma)$ is a non-descendant of $t$. 
						\item[(BR02)] There exists no edge incident with $v$ whose other end is in $V(G)-\uparrow t$.
						\item[(BR03)] Either $v \not \in V(\gamma_0)$, or $v \in V(\gamma_0)$ and $v$ is associated with essential number 0 in $\gamma_0$. 
					\end{itemize}
				\item[(BR1)] $j=1$, if the following hold.
					\begin{itemize}
						\item[(BR11)] $v \not \in V(\sigma)$ for every $\sigma \in \Gamma$ in which $\alpha(\sigma)$ is a non-descendant of $t$. 
						\item[(BR12)] either
							\begin{itemize}
								\item $v$ is incident with exactly one edge whose other end is in $V(G)-\uparrow t$, and either $v \not \in V(\gamma_0)$, or $v \in V(\gamma_0)$ and $v$ is associated with essential number 0 in $\gamma_0$, or  
								\item there exists no edge incident with $v$ whose other end is in $V(G)-\uparrow t$, $v \in V(\gamma_0)$, and $v$ is associated with essential number 1 in $\gamma_0$. 
							\end{itemize}
					\end{itemize}
				\item[(BR2)] $j=2$, otherwise.
			\end{itemize}
		\item $\Gamma_t=\{\sigma \in \Gamma: \alpha(\sigma)$ is a descendent of $t\}$. 
		\item $f_t=f|_{\Gamma_t}$. 
		\item $\phi_t = \phi|_{\uparrow t}$. 
		\item $S_t$ is the $Q$-assemblage $(G[\uparrow t],\gamma_t,\Gamma_t,f_t,\phi_t)$. 
	\end{itemize}
We call $S_t$ the \defn{$(f,\phi)$-branch of $(T,\X,\alpha)$ at $t$ (with respect to $\pi_t$)}, and we call the underlying assemblage of $S_t$ the \defn{branch of $(T,\X,\alpha)$ at $t$ (with respect to $\pi_t$)}.

\begin{lemma} \label{branch of brnach}
Let $Q$ be a quasi-order and let $(T,\X,\alpha)$ be a rooted tree-decomposition of a $Q$-assemblage $(G,\gamma_0,\Gamma,f,\phi)$.
Assume that for every node $t'$ other than the root of $T$, there exists an ordering $\pi_{t'}$ on $X_{t'} \cap X_{p'}$, where $p'$ is the parent of $t'$. 
Let $t$ be a non-root node of $T$. 
Let $(G[\uparrow t], \gamma_t, \Gamma_t,f_t,\phi_t)$ be the $(f,\phi)$-branch at $t$ with respect to $\pi_t$.
Let $T'$ be the rooted subtree of $T$ induced by the descendants of $t$ rooted at $t$.
Let $\X' = (X_t: t \in V(T'))$.
Let $\alpha' = \alpha|_{\Gamma_t}$.
Then $(T',\X',\alpha')$ is a rooted tree-decomposition of $(G[\uparrow t], \gamma_t, \Gamma_t,f_t,\phi_t)$ such that for every descendant $t'$ of $t$ with $t' \neq t$, the $(f_t,\phi_t)$-branch of $(T',\X',\alpha')$ at $t'$ with respect to $\pi_{t'}$ is the same as the $(f,\phi)$-branch of $(T,\X,\alpha)$ at $t'$ with respect to $\pi_{t'}$.
\end{lemma}

\begin{pf}
It is clear that $(T',\X',\alpha')$ is a rooted tree-decomposition of $(G[\uparrow t], \gamma_t, \Gamma_t,f_t,\phi_t)$.
Let $t'$ be a descendant of $t$ with $t' \neq t$.
It suffices to show that the $(f_t,\phi_t)$-branch of $(T',\X',\alpha')$ at $t'$ with respect to $\pi_{t'}$ is the same as the $(f,\phi)$-branch of $(T,\X,\alpha)$ at $t'$ with respect to $\pi_{t'}$.
Note that by definition, the only possible difference between these two branches are at their root marches.
Denote the $(f,\phi)$-branch of $(T,\X,\alpha)$ at $t'$ with respect to $\pi_{t'}$ by $(G[\uparrow t'], \gamma_{t'},\Gamma_{t'},f_{t'},\phi_{t'})$.
Denote the $(f_t,\phi_t)$-branch of $(T',\X',\alpha')$ at $t'$ with respect to $\pi_{t'}$ by $(G[\uparrow t'], \gamma'_{t'},\Gamma_{t'},f_{t'},\phi_{t'})$.
It suffices to show that $\gamma_{t'} = \gamma'_{t'}$.

Note that $V(\gamma'_{t'})=V(\gamma_{t'})$.
Let $v \in V(\gamma_{t'})$.
To show $\gamma_{t'} = \gamma'_{t'}$, it suffices to show that the essential number associated with $v$ in $\gamma'_{t'}$ equals the essential number associated with $v$ in $\gamma_{t'}$.
It obvious holds unless $v \in V(\gamma_t)$ by the definition of $\gamma_{t'}$ and $\gamma_{t'}'$.

So we may assume $v \in V(\gamma_t)$.

Let $\ell$ be the essential number associated with $v$ in $\gamma_{t'}$. 
Let $\ell'$ be the essential number associated with $v$ in $\gamma_{t'}'$. 
Let $\ell_t$ be the essential number associated with $v$ in $\gamma_t$.
Recall that $\gamma_t$ is the root march of $(G[\uparrow t], \gamma_t, \Gamma_t,f_t,\phi_t)$.

Suppose to the contrary that $\ell \neq \ell'$.

\medskip

\noindent{\bf Claim 1:} There exists no $\sigma \in \Gamma$ such that $v \in V(\sigma)$, and $\alpha(\sigma)$ is a non-descendant of $t'$.

\noindent{\bf Proof of Claim 1:}
Suppose that there exists $\sigma \in \Gamma$ such that $v \in V(\sigma)$, and $\alpha(\sigma)$ is a non-descendant of $t'$.
Then $\ell=2$.
If $\alpha(\sigma)$ is a descendant of $t$, then $\sigma \in \Gamma_t$, so $\ell=\ell'=2$, a contradiction.
So $\alpha(\sigma)$ is a non-descendant of $t$.
Hence $\ell_t=2$. 
This implies that (BR03) and (BR12) are violated when considering $\gamma'_{t'}$, so $\ell'=2$. 
Hence $\ell=2 = \ell'$, a contradiction.
$\Box$

\medskip

\noindent{\bf Claim 2:} There exists no $\sigma \in \Gamma_t$ such that $v \in V(\sigma)$, and $\alpha(\sigma)$ is a non-descendant of $t'$.

\noindent{\bf Proof of Claim 2:}
This claim follows from Claim 1 since $\Gamma_t \subseteq \Gamma$.
$\Box$

\medskip

By Claims 1 and 2, (BR01) and (BR11) hold when considering $\ell,\ell_t,\ell'$.

\medskip

\noindent{\bf Claim 3:} There exists an edge of $G$ incident with $v$ whose other end is in $V(G)-(T,\X)\uparrow t'$. 

\noindent{\bf Proof of Claim 3:}
Suppose that there exists no edge of $G$ incident with $v$ whose other end is in $V(G)-(T,\X)\uparrow t'$.
So there exists no edge of $G[\uparrow t]$ incident with $v$ whose other end is in $V(G[\uparrow t])-(T',\X')\uparrow t'$.
Hence (BR02) holds when considering $\ell,\ell_t,\ell'$.

Hence if $v \not \in V(\gamma_0)$, then $\ell_t=\ell=0$, and it implies that $\ell'=0=\ell$, a contradiction.
So $v \in V(\gamma_0)$.
Hence both $\ell$ and $\ell_t$ equal the essential number associated with $v$ in $\gamma_0$ by (BR03), (BR12) and (BR2). 
Similarly, $\ell'$ equal the essential number associated with $v$ in $\gamma_t$, which is $\ell_t$, by (BR03), (BR12) and (BR2). 
So $\ell'=\ell$, a contradiction.
$\Box$

\medskip

\noindent{\bf Claim 4:} There exist at least two edges of $G$ incident with $v$ whose other ends are in $V(G)-(T,\X)\uparrow t'$.

\noindent{\bf Proof of Claim 4:}
Suppose to the contrary that this claim does not hold.
By Claim 3, there exists exactly one edge $e$ of $G$ incident with $v$ whose other end is in $V(G)-(T,\X)\uparrow t'$.
So $\ell \in \{1,2\}$.

Suppose that the end of $e$ other than $v$ is in $V(G)-(T,\X)\uparrow t$.
Then $\ell_t=\ell$, and there exists no edge of $G[\uparrow t]$ incident with $v$ whose other end is in $V(G[\uparrow t])-(T',\X')\uparrow t'$.
Since $v \in V(\gamma_t)$, $\ell'=\ell_t$.
Hence $\ell=\ell'$, a contradiction.

So the end of $e$ other than $v$ is in $(T,\X)\uparrow t-(T,\X)\uparrow t'$.
Hence there exists no edge of $G$ incident with $v$ whose other end is in $V(G)-(T,\X)\uparrow t$.
If $v \not \in V(\gamma_0)$, then $\ell_t=0$ and $\ell=1$, so $\ell'=1=\ell$, a contradiction.
So $v \in V(\gamma_0)$.
If $v$ is associated with essential number 0 in $\gamma_0$, then $\ell_t=0$ and $\ell=1$, so $\ell'=1=\ell$, a contradiction.
So $v$ is associated with essential number 1 or 2 in $\gamma_0$, then $\ell_t \geq 1$ and $\ell=2$, so $\ell'=2=\ell$, a contradiction.
$\Box$

\medskip

Claim 4 implies $\ell=2$.

If there exists no edge of $G$ incident with $v$ whose other end is in $V(G)-(T,\X)\uparrow t$, then there exist at least two edges of $G$ incident with $v$ whose other ends are in $V(G[\uparrow t])-(T',\X')\uparrow t'$, so $\ell'=2=\ell$, a contradiction.
So there exists an edge of $G$ incident with $v$ whose other end is in $V(G)-(T,\X)\uparrow t$.
In particular, $\ell_t \geq 1$.
So if there exists an edge of $G[\uparrow t]$ incident with $v$ whose other end is in $V(G[\uparrow t])-(T,\X)\uparrow t'$, then $\ell'=2=\ell$, a contradiction.
Hence there exists no edge of $G[\uparrow t]$ incident with $v$ whose other end is in $V(G[\uparrow t])-(T,\X)\uparrow t'$.
This implies that there exist at least two edges of $G$ incident with $v$ whose other ends are in $V(G)-(T,\X)\uparrow t$.
So $\ell_t=2$, and hence $\ell'=2=\ell$, a contradiction.
\end{pf}

\subsection{Encoding} \label{subsec:encoding}

Let $Q$ be a quasi-order and let $(T,\X,\alpha)$ be a rooted tree-decomposition of a $Q$-assemblage $(G,\gamma_0,\Gamma,f,\phi)$.
For every non-root node $t$ of $T$, let $S_t$ be the $(f,\phi)$-branch at $t$ (with respect to an ordering $\pi_t$ of its root march), define $b_t$ to be the sequence with length $\lvert V(\gamma_t) \rvert$ such that for every integer $i \in [\lvert V(\gamma_t) \rvert]$, the $i$-th entry of $b_t$ is 
	\begin{itemize}
		\item 1 if the $i$-th vertex in $\gamma_t$ is adjacent to a vertex in $\uparrow t-V(\gamma_t)$, and 
		\item 0 otherwise.
	\end{itemize}
Now let $t$ be a node of $T$ (possibly the root of $T$).
Let $S=\{S_c: c$ is a child of $t\}$ and let it be ordered by the simulation relation.
So $S$ is a quasi-order.
Let $S'=\{b_c: c$ is a child of $t\}$ and let it be ordered by the equality relation.
So $S'$ is a quasi-order.
Let $S''$ be the quasi-order obtained by the Cartesian product of $S$ and $S'$.
Define $Q'$ to be the quasi-order obtained from $Q$ and $S''$ by disjoint union.
The \defn{encoding of $(T,\X,\alpha)$ at $t$ (with respect to $\pi_t$ and $\pi_c$ for all children $c$ of $t$)} is the $Q'$-assemblage $(H,\gamma_H,\Gamma_H,f_H,\phi_H)$ such that the following hold.
	\begin{itemize}
		\item $H=G[X_t]$.
		\item $\gamma_H=\gamma_t$ if $t$ is not the root of $T$; $\gamma_H=\gamma_0$ if $t$ is the root of $T$.
		\item $\Gamma_H=\Lambda_1 \uplus \Lambda_2$, where $\Lambda_1=\{\gamma_c: c$ is a child of $t\}$ and $\Lambda_2 = \{\sigma \in \Gamma: \alpha(\sigma)=t\}$. 
			(Note that we keep both elements in $\Gamma_H$ if some element appears in both $\Lambda_1$ and $\Lambda_2$.)
		\item $f_H(\sigma)=(S_c,b_c)$ if $\sigma \in \Lambda_1$ and $\sigma=\gamma_c$; $f_H(\sigma)=f(\sigma)$ if $\sigma \in \Lambda_2$.
		\item $\phi_H=\phi|_{V(H)}$.
	\end{itemize}
We compare encodings by the simulation relation (with respect to $Q'$).

The following lemma (Lemma \ref{encoding simulation}) is the main result of this section, which states that if the encoding of a rooted tree-decomposition of a labelled assemblage at its root simulates the encoding of another rooted tree-decomposition of another labelled assemblage at its root, then the first labelled assemblage simulates the second one.
We sketch its proof here.
The simulation between encodings gives us a homeomorphic embedding $\eta$ between the root bags, and for each child $c$ of the root of the first tree, a homeomorphic embedding $\eta_c$ from the branch at $c$ to the branch at some child $c'$ of the root of the second tree.
We shall construct a desired homeomorphic embedding from the first whole assemblage to the second one by using $\eta$ and $\eta_c$.
Clearly, vertices contained in the root bag but not in any other bag should be mapped according to $\eta$, and each vertex not contained in the root bag is contained in exactly one branch (say at branch at $c$) and should be mapped according to $\eta_c$.
However, it is unclear how to map the common vertices of the root bag and a child bag; more seriously, some vertex can be contained in multiple child bags.
The first part of the proof (Claims 1-10) shows that we can resolve this ambiguity by determining where those vertices are mapped to based on a simple rule in a way that we can further embed edges of the first assemblage into paths in the second one.
Then we formally construct a homeomorphic embedding from the first assemblage to the second one and verify the correctness (Claims 11-15).
Finally, there is a natural way to extend the homeomorphic embedding to a simulation (Claims 16 and 17) to complete the proof.

\begin{lemma} \label{encoding simulation}
Let $Q$ be a quasi-order.
Let $(G,\gamma_0,\Gamma,f,\phi)$ and $(G',\gamma_0',\Gamma',f',\phi')$ be $Q$-assemblages with rooted tree-decompositions $(T,\X,\alpha)$ and $(T',\X',\alpha')$, respectively.
Let $r$ and $r'$ be the roots of $T$ and $T'$, respectively.
Let $\gamma_r=\gamma_0$ and $\gamma_{r'}=\gamma_0'$.
Assume that for every $t \in \{r,r'\}$, there exists an ordering $\pi_t$ of $V(\gamma_t)$ which is the same as the ordering of $\gamma_t$, and for every child $c$ of $t$, there exists an ordering $\pi_c$ of $X_t \cap X_c$.
If the encoding of $(T',\X',\alpha')$ at $r'$ (with respect to $\pi_{r'}$ and $\pi_c$ for all children $c$ of $r'$) simulates the encoding of $(T,\X,\alpha)$ at $r$ (with respect to $\pi_r$ and $\pi_c$ for all children of $r$), then $(G',\gamma_0',\Gamma',f',\phi')$ simulates $(G,\gamma_0,\Gamma,f,\phi)$.
\end{lemma}

\begin{pf}
Let $(H,\gamma_0,\Gamma_H,f_H,\phi_H)$ and $(H',\gamma_0',\Gamma_H',f_H',\phi_H')$ be the encodings of $(T,\X,\alpha)$ and $(T',\X',\alpha')$ at $r$ and $r'$, respectively.
Let $\eta: (H,\gamma_0) \hookrightarrow (H',\gamma_0')$ and $\iota$ be the functions witnessing the simulation between these two encodings.
For each child $c$ of $r$ (or $r'$), we define $S_c,b_c$ to be the $(f,\phi)$-branch (or $(f',\phi')$-branch) at $c$ and the sequence, respectively, as mentioned in the definition of encodings, and define $\gamma_c$ to be the march in $\Gamma_H$ (or $\Gamma_{H'}$) mapped to $(S_c,b_c)$ by $f_H$ (or $f_{H'}$).
For each child $c$ of $r$, define $\eta_c,\iota_c$ to be the functions that witness the simulation $f_H(\gamma_c) \leq_{Q'} f'_H(\iota(\gamma_c))$, where $Q'$ is the quasi-order mentioned in the definition of encodings.
Let $\overline{G}, \overline{G'}, \overline{H}, \overline{H'}$ be the rooted extensions of $G,G',H,H'$, respectively.

\medskip

\noindent{\bf Claim 1:} If $c$ is a child of $r$ and $v \in V(\gamma_c)$ such that $\eta_c(v) \not \in V(\overline{H'})$, then the essential number associated with $v$ in $\gamma_c$ is 0 or 1.

\noindent{\bf Proof of Claim 1:}
Note that $V(\iota(\gamma_c)) \subseteq V(\overline{H'})$.
Since $\eta_c(v) \not \in V(\overline{H'})$, $\eta_c(v) \not \in V(\iota(\gamma_c))$.
If the essential number associated with $v$ in $\gamma_c$ is 2, then there are two edges between $v$ and the corresponding indicator, so $\eta_c(v) \in V(\iota(\gamma_c))$, a contradiction.
So $v$ is associated with essential number 0 or 1 in $\gamma_c$.
$\Box$

\medskip

Let $A= \bigcup_c V(\gamma_c)$, where the union is over all children $c$ of $r$.
For each vertex $v \in A$, define $L_v$ to be the set of children $c$ of $r$ such that $v \in V(\gamma_c)$ and $\eta_c(v) \not \in V(\overline{H'})$. 

\medskip

\noindent{\bf Claim 2:} If $v \in A$, then the following statements hold.
	\begin{itemize}
		\item If $c$ is a node with $c \in L_v$ such that $v$ is associated with essential number 1 in $\gamma_c$, then $v$ is adjacent to a vertex in $\uparrow c- V(\gamma_c)$, and $\lvert L_v \rvert \leq 2$.
		\item If $c$ is a node with $c \in L_v$ such that $v$ is adjacent to a vertex in $\uparrow c - V(\gamma_c)$ and $v$ is associated with essential number 0 in $\gamma_c$, then $L_v = \{c\}$.
		\item If for every $c \in L_v$, $v$ is associated with essential number 0 in $\gamma_c$ and $v$ is not adjacent to any vertex in $\uparrow c - V(\gamma_c)$, then for every $c \in L_v$, $V(\gamma_c)$ contains all neighbors of $v$ in $G$.
	\end{itemize}

\noindent{\bf Proof of Claim 2:}
For each child $c$ of $r$, let $c'$ be the child of $r'$ such that $\gamma_{c'}=\iota(\gamma_c)$.

We first assume that $c \in L_v$ and $v$ is associated with essential number 1 in $\gamma_c$ and prove the first statement of this claim.
Since $\eta_{c}(v) \not \in V(\overline{H'})$ and $v$ is associated with essential number 1 in $\gamma_c$, there exists a path in $G'[(T',\X')\uparrow c']$ from $\eta_{c}(v)$ to $\eta(v)$ disjoint from $V(\eta(\gamma_c))-\{\eta(v)\}$.
So $\eta(v)$ is adjacent to a vertex in $(T',\X')\uparrow c'-V(\gamma_{c'})$.
Since $b_{c} = b_{c'}$, $v$ is adjacent to a vertex in $(T,\X)\uparrow c - V(\gamma_{c})$.
So $v$ is associated with essential number 1 or 2 in $\gamma_d$ for every child $d$ of $r$ with $v \in V(\gamma_d)$ and $d \neq c$.
If there exists $d \in L_v-\{c\}$, then since $v$ is adjacent to some vertex in $(T,\X)\uparrow c - V(\gamma_{c})$, the essential number associated with $v$ in $\gamma_d$ is 1 by Claim 1, so a similar argument shows that $v$ is adjacent to a vertex in $\uparrow d - V(\gamma_d)$; since $v$ is adjacent to one vertex in $(T,\X)\uparrow c - V(\gamma_{c})$ and one vertex in $(T,\X)\uparrow d - V(\gamma_{d})$, $v$ is associated with essential number 2 in $\gamma_{d'}$ for every child $d'$ of $r$ other than $c$ and $d$, so $L_v = \{c,d\}$ by Claim 1.
This proves the first statement.

Now we assume that $c \in L_v$, $v$ is adjacent to a vertex in $\uparrow c - V(\gamma_c)$, and $v$ is associated with essential number 0 in $\gamma_c$.
Suppose that there exists $d \in L_v-\{c\}$.
Since $v$ is adjacent to a vertex in $\uparrow c- V(\gamma_c)$, $v$ is associated with essential number 1 in $\gamma_d$ by Claim 1.
By Statement 1 of this claim, $v$ is adjacent to a vertex in $\uparrow d-V(\gamma_d)$, contradicting that $v$ is associated with essential number 0 in $\gamma_c$ by (BR02).
This proves the second statement.

Finally, we assume that for every $c \in L_v$, $v$ is associated with essential number 0 in $\gamma_c$, and $v$ is not adjacent to any vertex in $\uparrow c-V(\gamma_c)$.
So for every $c \in L_v$, $v$ is not adjacent to any vertex in $V(G)-\uparrow c$ (by (BR02)) and is not adjacent to any vertex in $\uparrow c- V(\gamma_c)$, so all neighbors of $v$ in $G$ belong to $V(\gamma_c)$.
This proves the claim.
$\Box$

\medskip

\noindent{\bf Claim 3:} If $c_1,c_2$ are different children of $r$ such that $v \in V(\gamma_{c_1}) \cap V(\gamma_{c_2})$, $\eta_{c_1}(v) \not \in V(\overline{H'})$, $\eta_{c_2}(v) \not \in V(\overline{H'})$, and at least one of the essential numbers associated with $v$ in $\gamma_{c_1},\gamma_{c_2}$ is non-zero, then the following statements hold.
	\begin{itemize}
		\item $v$ is associated with essential number 1 in both $\gamma_{c_1}$ and $\gamma_{c_2}$, and $v$ is adjacent to a vertex in $\uparrow c_1 - V(\gamma_{c_1})$ and adjacent to a vertex in $\uparrow c_2 -V(\gamma_{c_2})$.
		\item For every neighbor $u \in V(\gamma_{c_1}) \cap V(\gamma_{c_2})$ of $v$, $\eta_{c_1}(u)$ and $\eta_{c_2}(u)$ are either both in $V(\overline{H'})$ or both not in $V(\overline{H'})$.
	\end{itemize}

\noindent{\bf Proof of Claim 3:}
For $i \in \{1,2\}$, let $c_i'$ be the child of $r'$ such that $\iota(\gamma_{c_i})=\gamma_{c_i'}$.
We know $\{c_1,c_2\} \subseteq L_v$ since $v \in V(\gamma_{c_1}) \cap V(\gamma_{c_2})$, $\eta_{c_1}(v) \not \in V(\overline{H'})$ and $\eta_{c_2}(v) \not \in V(\overline{H'})$.
Since at least one of the essential numbers associated with $v$ in $\gamma_{c_1},\gamma_{c_2}$ is non-zero, there exists $i \in [2]$ such that the essential number associated with $v$ in $\gamma_{c_i}$ is non-zero (and hence equals 1 by Claim 1).
By Statement 1 of Claim 2, $L_v=\{c_1,c_2\}$ and $v$ is adjacent to a vertex in $\uparrow c_i-V(\gamma_{c_i})$.
Hence $v$ is associated with essential number 1 in $\gamma_{c_{3-i}}$ (by Claim 1 and (BR02)) and hence in both $\gamma_{c_1}$ and $\gamma_{c_2}$.
Again by Statement 1 of Claim 2, $v$ is adjacent to a vertex in $\uparrow c_1-V(\gamma_{c_1})$ and a vertex in $\uparrow c_2-V(\gamma_{c_2})$.
So the first statement of this claim holds.

Let $u \in V(\gamma_{c_1}) \cap V(\gamma_{c_2})$ be a neighbor of $v$.
Suppose that $\eta_{c_1}(u) \not \in V(\overline{H'})$ and $\eta_{c_2}(u) \in V(\overline{H'})$.
Since $v$ is associated with essential number 1 in $\gamma_{c_2}$ (by the first statement of this claim) and $\eta_{c_2}(v) \not \in V(\overline{H'})$ and $\eta_{c_2}(u) \in V(\overline{H'})$, $\eta_{c_2}(u)$ is adjacent to a vertex in $(T',\X') \uparrow c_2'-V(\gamma_{c_2'})$.
Since $b_{c_2}=b_{c_2'}$, $u$ is adjacent to a vertex in $\uparrow c_2-V(\gamma_{c_2})$, so $u$ is not associated with essential number 0 in $\gamma_{c_1}$ by (BR02).
Since $\eta_{c_1}(u) \not \in V(\overline{H'})$, $u$ is associated with essential number 1 in $\gamma_{c_1}$ by Claim 1. 
So there exists a unique vertex $w \in \uparrow c_2 -V(\gamma_{c_2})$ adjacent to $u$ by (BR12).
Note that $\eta(u)$ is the vertex in $V(\gamma_{c_1'}) \cap V(\eta_{c_1}(e))$, where $e$ is the edge in $S_{c_1}$ between $u$ and the indicator adjacent to $u$.
So essential number associated with $\eta(u)$ in $\gamma_{c_1'}$ is 1.
Since $u$ is adjacent to a vertex in $\uparrow c_2-V(\gamma_{c_2})$, $\eta(u)$ is adjacent to a vertex in $(T',\X')\uparrow c_2'-V(\gamma_{c_2'})$.
Hence there exists a unique vertex $w' \in (T',\X')\uparrow c_2'-V(\gamma_{c_2'})$ adjacent to $\eta(u)$.
Since $\eta_{c_2}(v),\eta_{c_2}(w) \not \in V(\overline{H'})$, we know $w' \in \eta_{c_2}(uv) \cap \eta_{c_2}(wu)$, so $v=w$.
But $v \in V(\gamma_{c_2})$ and $w \not \in V(\gamma_{c_2})$, so $v \neq w$, a contradiction.
This proves that $\eta_{c_1}(u)$ and $\eta_{c_2}(u)$ are either both in $V(\overline{H'})$ or both not in $V(\overline{H'})$.
$\Box$

\medskip

\noindent{\bf Claim 4:} If $v \in V(\gamma_r) \cap A$ and $v$ is associated with essential number 2 in $\gamma_r$, then $L_v=\emptyset$.

\noindent{\bf Proof of Claim 4:}
Let $c$ be a child of $r$ such that $v \in V(\gamma_c)$.
To prove this claim, it suffices to prove that $\eta_c(v) \in V(\overline{H'})$.
Since $v$ is associated with essential number 2 in $\gamma_r$, $v$ is associated with essential 2 in $\gamma_c$ by definition by (BR03) and (BR12).
So $\eta_c(v) \in V(\overline{H'})$ by Claim 1.
$\Box$

\medskip

Claims 1 and 4 imply that if $v \in A$ and $L_v \neq \emptyset$, then for every $c \in L_v$, the essential number associated with $v$ in $\gamma_c$ is in $\{0,1\}$ and equals the number of edges incident with $v$ whose other end is in $V(\overline{G})-\uparrow c$.

\medskip

\noindent{\bf Claim 5:} If $c$ is a child of $r$ and $v_1,v_2 \in V(\gamma_c)$ are adjacent vertices with $c \in L_{v_2}-L_{v_1}$, then either $v_1$ is adjacent to a vertex in $\uparrow c-V(\gamma_c)$, or all of the following statements hold.
	\begin{itemize}
		\item $v_2$ is associated with essential number 0 in $\gamma_c$,
		\item $v_2$ is not adjacent to any vertex in $V(G)-\uparrow c$, 
		\item there exists a unique edge $e$ of $G$ between $v_1$ and $v_2$, and
		\item there exists a path in $\eta_c(e)$ from $\eta_{c}(v_2) \not \in V(\overline{H'})$ to $\eta(v_2)$ disjoint from $V(\eta(\gamma_c))-\{\eta(v_2)\}$ such that $\eta_c(e)$ is the union of this path and an edge between $\eta(v_1)$ and $\eta(v_2)$. 
	\end{itemize}

\noindent{\bf Proof of Claim 5:}
Let $c'$ be the child of $r'$ such that $\iota(\gamma_{c})=\gamma_{c'}$.
Since $c \not \in L_{v_1}$, $\eta_{c}(v_1)=\eta(v_1)$.
Since $c \in L_{v_2}$, $\eta_c(v_2) \not \in V(\eta(\gamma_c))$.
Let $v_1v_2$ be an arbitrary edge of $G$ between $v_1$ and $v_2$.
Since $v_1v_2 \in E(G)$, $\eta_c(v_1v_2)$ is a path in $G'[(T',\X') \uparrow c']$ from $\eta_{c}(v_2) \not \in V(\overline{H'})$ to $\eta_{c}(v_1)=\eta(v_1)$. 
Note that $\eta_c(v_1v_2)$ is disjoint from $V(\eta(\gamma_c))-\{\eta(v_1),\eta(v_2)\}$. 

We may assume that $v_1$ is not adjacent to a vertex in $\uparrow c-V(\gamma_c)$, for otherwise we are done.

If $\eta(v_2) \not \in V(\eta_c(v_1v_2))$, then $\eta(v_1)$ is adjacent to a vertex in $\uparrow c'-V(\gamma_{c'})$; since $b_{c}=b_{c'}$, $v_1$ is adjacent to a vertex in $\uparrow c-V(\gamma_{c})$, a contradiction.
So $\eta(v_2) \in V(\eta_c(v_1v_2))$.
Since $\eta(v_2) \in V(\overline{H'})$, $\eta_c(v_2) \neq \eta(v_2)$.
So $\eta(v_2) \in V(\gamma_{c'})$ is an internal vertex of $\eta_c(v_1v_2)$.
Hence the essential number associated with $v_2$ in $\gamma_c$ is 0, and $v_1v_2$ is the unique edge of $G$ between $v_1$ and $v_2$ by the definition of homomorphic embeddings of rooted graphs.
So there exists no edge incident with $v_2$ whose other end is in $V(G)-\uparrow c$ by (BR02).
Since $\eta(v_2)$ is an internal vertex of $\eta_c(v_1v_2)$, there exists a path $P$ in $G'[(T',\X') \uparrow c']$ from $\eta_{c}(v_2) \not \in V(\overline{H'})$ to $\eta(v_2)$ disjoint from $V(\eta(\gamma_c))-\{\eta(v_2)\}$.
Since $v_1$ is not adjacent to a vertex in $\uparrow c-V(\gamma_c)$, $\eta_c(v_1v_2)$ is the union of $P$ and an edge between $\eta(v_1)$ and $\eta(v_2)$.
$\Box$

\medskip

\noindent{\bf Claim 6:} If $v_1,v_2 \in A$ are distinct adjacent vertices, $c_1$ and $c_2$ are distinct children of $r$ with $c_1 \in L_{v_1},c_2 \in L_{v_2}$ and $\{v_1,v_2\} \subseteq V(\gamma_{c_1}) \cap V(\gamma_{c_2})$ such that $v_1$ is associated with essential number 1 in $\gamma_{c_1}$, and $v_2$ is not adjacent to any vertex in $V(G)-\uparrow c_2$, then  
	\begin{itemize}
		\item $L_{v_2} = \{c_1,c_2\}$, 
		\item $v_2$ is associated with essential number 0 in $\gamma_{c_1}$,
		\item all neighbors of $v_2$ are contained in $V(\gamma_{c_1}) \cap V(\gamma_{c_2})$, and
		\item either $L_{v_1}=\{c_1\}$ or $L_{v_1}=\{c_1,c_2\}$.
	\end{itemize}

\noindent{\bf Proof of Claim 6:}
For $i \in \{1,2\}$, let $c_i'$ be the child of $r'$ such that $\iota(\gamma_{c_i})=\gamma_{c_i'}$.

Suppose $\eta_{c_1}(v_2)=\eta(v_2)$.
Since $c_1 \in L_{v_1}$ and $v_1$ is associated with essential number 1 in $\gamma_{c_1}$, $\eta_{c_1}(v_1v_2)$ is a path in $G'[\uparrow c_1']$ from $\eta_{c_1}(v_1) \not \in V(\gamma_{c_1'})$ to $\eta_{c_1}(v_2)=\eta(v_2)$ disjoint from $V(\gamma_{c_1'})-\{\eta_{c_1}(v_2)\}$.
Since $b_{c_1}=b_{c_1'}$, $v_2$ is adjacent to a vertex in $\uparrow c_1-V(\gamma_{c_1}) \subseteq V(G)-\uparrow c_2$, a contradiction.

So $\eta_{c_1}(v_2) \neq \eta(v_2)$.
Hence $c_1 \in L_{v_2}$.
Since $v_2$ is not adjacent to any vertex in $V(G)-\uparrow c_2$, $v_2$ is associated with essential number 0 or 2 in $\gamma_{c_1}$ by Statement 1 of Claim 2.
Since $\eta_{c_1}(v_2) \neq \eta(v_2)$, $v_2$ is associated with essential number 0 in $\gamma_{c_1}$.
By (BR02), $v_2$ is not adjacent to any vertex in $V(G)-\uparrow c_1$.
Hence all neighbors of $v_2$ are contained in $\uparrow c_1 \cap \uparrow c_2 \subseteq V(\gamma_{c_1}) \cap V(\gamma_{c_2})$.

To prove this claim, it suffices to prove that $L_{v_2}=\{c_1,c_2\}$, and either $L_{v_1}=\{c_1\}$ or $L_{v_1}=\{c_1,c_2\}$.

Since $c_1 \in L_{v_2}$, $\{c_1,c_2\} \subseteq L_{v_2}$.
Since $c_1 \in L_{v_1}$ and $v_1$ is associated with essential number 1 in $\gamma_{c_1}$, Statement 1 of Claim 2 implies $\lvert L_{v_1} \rvert \leq 2$.
If $L_{v_2} \subseteq L_{v_1}$, then $L_{v_1}=L_{v_2} = \{c_1,c_2\}$ and we are done.

So we may assume that $L_{v_2}-L_{v_1} \neq \emptyset$.

Suppose that there exists $d \in L_{v_2}$ such that $v_1 \not \in V(\gamma_d)$.
Then the assumption $d \in L_{v_2}$ and the existence of an edge $v_1v_2$ implies that $v_2$ is associated with essential number 1 in $\gamma_d$ by Claim 1 and (BR02).
But it implies that $v_2$ is adjacent to a vertex in $\uparrow d-V(\gamma_d)$ by Statement 1 of Claim 2, a contradiction.

Hence for every $d \in L_{v_2}$, $v_1 \in V(\gamma_d)$.

Let $c \in L_{v_2}-L_{v_1}$.
Note that $v_1 \in V(\gamma_c)$ and $c \neq c_1$.
Since all neighbors of $v_2$ are contained in $V(\gamma_{c_1}) \cap V(\gamma_{c_2})$, by the fact $b_c=b_{c'}$, where $c'$ is the children of $r'$ such that $\iota(\gamma_c)=\gamma_{c'}$, we know $\eta(v_2)$ is not adjacent to any vertex in $\uparrow c'-V(\gamma_{c'})$.
This together with the assumption $c \in L_{v_2}-L_{v_1}$, we know that $\eta_c(v_1v_2)$ is disjoint from $\eta(v_2)$, so $\eta(v_1)$ is adjacent to a vertex in $\uparrow c'-V(\gamma_{c'})$.
Since $b_c=b_{c'}$, $v_1$ is adjacent to a vertex in $\uparrow c-V(\gamma_c)$.
Since $v_1$ is associated with essential number 1 in $\gamma_{c_1}$, $\uparrow c_1 \cup \uparrow c$ contains all neighbors of $v_1$, and $c$ is the unique element in $L_{v_2}-L_{v_1}$. 
Hence either $c_2 \in L_{v_1} \cap L_{v_2}$, or $c_2$ is the unique element in $L_{v_2}-L_{v_1}$.
If $c_2 \in L_{v_1}$, then $L_{v_1} = \{c_1,c_2\}$ by Statement 1 of Claim 2, so $c \neq c_2$ and $v_1$ is adjacent to a vertex in $\uparrow c_1-V(\gamma_{c_1})$ (by Statement 1 of Claim 2) and a vertex in $\uparrow c-V(\gamma_c)$, and hence $v_1$ is associated with essential number 2 in $c_2$ and $c_2 \not \in L_{v_1}$, a contradiction.
So $c_2$ is the unique element in $L_{v_2}-L_{v_1}$.
Since $v_1$ is adjacent to a vertex in $\uparrow c_1-V(\gamma_{c_1})$ (by Statement 1 of Claim 2) and a vertex in $\uparrow c_2-V(\gamma_{c_2})$, we know that $v_1$ is associated with essential number 2 in $\gamma_{d}$ for every child $d \not \in \{c_1,c_2\}$ of $r$ with $v_1 \in V(\gamma_{d})$, so $L_{v_1} \subseteq \{c_1,c_2\}$.
Since $c_2 \not \in L_{v_1}$, $L_{v_1}=\{c_1\}$.
And by the uniqueness of $c$, $L_{v_2}=\{c_1,c_2\}$.
$\Box$

\medskip

\noindent{\bf Claim 7:} If $v_1,v_2 \in A$ are adjacent vertices, $L_{v_1} \neq L_{v_2}$, $c_1$ and $c_2$ are distinct children of $r$ such that $\{v_1,v_2\} \subseteq V(\gamma_{c_1}) \cap V(\gamma_{c_2})$, and for every $i \in \{1,2\}$, $c_i \in L_{v_i}$ and $v_i$ is associated with essential number 0 in $\gamma_{c_i}$, then 
	\begin{itemize}
		\item there exists a unique edge $e$ of $G$ between $v_1$ and $v_2$, and 
		\item for every $i \in \{1,2\}$, $L_{v_i}=\{c_i\}$ and there exists a path in $\eta_{c_i}(e)$ from $\eta_{c_i}(v_i)$ to $\eta(v_i)$ disjoint from $V(\eta(\gamma_{c_i}))-\{\eta(v_i)\}$ .
	\end{itemize}

\noindent{\bf Proof of Claim 7:}
We first assume that there exists $c \in L_{v_1}-L_{v_2}$.
If $v_2 \not \in V(\gamma_c)$, then $c \neq c_1$ and $v_2$ is a vertex in $V(G)-\uparrow c$ adjacent to $v_1$, so $v_1$ is associated with essential number 1 in $\gamma_c$ by Claim 1 and (BR02); by Statement 1 of Claim 2, $v_1$ is adjacent to a vertex in $\uparrow c-V(\gamma_c)$, so $v_1$ is not associated with essential number 0 in $\gamma_{c_1}$, a contradiction.
So $v_2 \in V(\gamma_c)$.
Since $c_2 \in L_{v_2}$, $c \neq c_2$.
Since $v_2$ is associated with essential number 0 in $\gamma_{c_2}$, $v_2$ is not adjacent to a vertex in $\uparrow c - V(\gamma_c)$.
Since $c \in L_{v_1}-L_{v_2}$, by Claim 5, $v_1$ is not adjacent to any vertex in $V(G)-\uparrow c$.
Let $v_1v_2$ be an arbitrary edge of $G$ between $v_1$ and $v_2$, and let $c'$ be the child of $r'$ such that $\iota(\gamma_{c})=\gamma_{c'}$.
Since $c \in L_{v_1}-L_{v_2}$, $\eta_c(v_1v_2)$ is a path in $G'[(T',\X')\uparrow c']$ from $\eta_c(v_1) \not \in V(\gamma_{c'})$ to $\eta_c(v_2)=\eta(v_2)$. 
Since $v_2$ is not adjacent to a vertex in $\uparrow c - V(\gamma_c)$ and $b_c=b_{c'}$, $\eta(v_1)$ is an internal vertex of $\eta_c(v_1v_2)$.
Hence there exists a path in $\eta_{c}(v_1v_2)$ from $\eta_{c}(v_1)$ to $\eta(v_1)$ disjoint from $V(\gamma_{c})-\{\eta(v_1)\}$.
So $\eta(v_1)$ is adjacent to a vertex in $\uparrow c'-V(\gamma_{c'})$, and the edge of $G$ between $v_1$ and $v_2$ is unique.
Since $b_c=b_{c'}$, $v_1$ is adjacent to a vertex in $\uparrow c-V(\gamma_{c})$.
Since $v_1$ is associated with essential number 0 in $\gamma_{c_1}$, $c=c_1$.
By Statement 2 of Claim 2, $L_{v_1}=\{c_1\}$.
Since $c=c_1$, there exists a path in $\eta_{c_1}(v_1v_2)$ from $\eta_{c_1}(v_1)$ to $\eta(v_1)$ disjoint from $V(\gamma_{c_1})-\{\eta(v_1)\}$.

Hence if $L_{v_1} \not \subseteq L_{v_2}$, then there exists a unique edge $e$ of $G$ between $v_1$ and $v_2$, $L_{v_1}=\{c_1\}$, and there exists a path in $\eta_{c_1}(e)$ from $\eta_{c_1}(v_1)$ to $\eta(v_1)$ disjoint from $V(\gamma_{c_1})-\{\eta(v_1)\}$.
Similarly, if $L_{v_2} \not \subseteq L_{v_1}$, then there exists a unique edge $e$ of $G$ between $v_1$ and $v_2$, $L_{v_2}=\{c_2\}$, and there exists a path in $\eta_{c_2}(e)$ from $\eta_{c_2}(v_2)$ to $\eta(v_2)$ disjoint from $V(\gamma_{c_2})-\{\eta(v_2)\}$.

Since $L_{v_1} \neq L_{v_2}$, by symmetry, we may assume that $L_{v_1} \not \subseteq L_{v_2}$.
So there exists a unique edge $e$ of $G$ between $v_1$ and $v_2$, $L_{v_1}=\{c_1\}$, and there exists a path in $\eta_{c_1}(e)$ from $\eta_{c_1}(v_1)$ to $\eta(v_1)$ disjoint from $V(\gamma_{c_1})-\{\eta(v_1)\}$.

Since $c_2 \not \in \{c_1\}=L_{v_1}$, $L_{v_2} \not \subseteq L_{v_1}$.
So $L_{v_2}=\{c_2\}$, and there exists a path in $\eta_{c_2}(e)$ from $\eta_{c_2}(v_2)$ to $\eta(v_2)$ disjoint from $V(\gamma_{c_2})-\{\eta(v_2)\}$. 
$\Box$

\medskip

\noindent{\bf Claim 8:} If $v_1,v_2 \in A$ are distinct adjacent vertices, $c_1$ and $c_2$ are distinct children of $r$ such that $c_1 \in L_{v_1},c_2 \in L_{v_2}$ and $\{v_1,v_2\} \subseteq V(\gamma_{c_1}) \cap V(\gamma_{c_2})$, then either 
	\begin{itemize}
		\item $L_{v_1}=L_{v_2}$, or 
		\item there exists a unique edge $e$ of $G$ between $v_1$ and $v_2$, and either
			\begin{itemize}
				\item for every $i \in \{1,2\}$, $v_i$ is associated with essential number 0 in $\gamma_{c_i}$, $L_{v_i}=\{c_i\}$, and there exists a path in $\eta_{c_i}(e)$ from $\eta_{c_i}(v_i)$ to $\eta(v_i)$ disjoint from $V(\eta(\gamma_{c_i}))-\{\eta(v_i)\}$, or
				\item there exists $i^* \in \{1,2\}$ such that $L_{v_{i^*}}=\{c_{i^*}\}$, $L_{v_{3-i^*}}=\{c_1,c_2\}$, and all neighbors of $v_{3-i^*}$ are contained in $V(\gamma_{c_1}) \cap V(\gamma_{c_2})$.
			\end{itemize}
	\end{itemize}

\noindent{\bf Proof of Claim 8:}
For $i \in \{1,2\}$, let $c_i'$ be the child of $r'$ such that $\iota(\gamma_{c_i})=\gamma_{c_i'}$.
We may assume $L_{v_1} \neq L_{v_2}$, for otherwise we are done.

If for every $i \in \{1,2\}$, $v_i$ is associated with essential number 0 in $\gamma_{c_i}$, then by Claim 7, the second statement of this claim holds.

Hence by symmetry and Claim 1, we may assume that $v_1$ is associated with essential number 1 in $\gamma_{c_1}$.
By Statement 1 of Claim 2, $v_1$ is adjacent to a vertex in $\uparrow c_1-V(\gamma_{c_1})$ and $\lvert L_{v_1} \rvert \leq 2$.
So $v_1$ is associated with essential number 1 or 2 in $\gamma_{c_2}$ by (BR02).

We first assume that $v_2$ is not adjacent to any vertex in $V(G)-\uparrow c_2$.
Then by Claim 6, $c_1 \in L_{v_1} \subseteq L_{v_2}=\{c_1,c_2\}$ and all neighbors of $v_2$ are contained in $V(\gamma_{c_1}) \cap V(\gamma_{c_2})$.
Since $L_{v_1} \neq L_{v_2}$, $L_{v_1}=\{c_1\}$.
Since $v_1$ is associated with essential number 1 in $\gamma_{c_1}$, there exists at most one edge incident with $v_1$ whose other end is in $\uparrow c_2-V(\gamma_{c_2})$.
So $\eta(v_1)$ is associated with essential number 1 in $\gamma_{c_1'}$ and there exists at most one edge $e'$ incident with $\eta(v_1)$ whose other end is in $\uparrow c'_2-V(\gamma_{c'_2})$.
This together with the fact $c_2 \in L_{v_2}-L_{v_1}$ imply that there are at most two edges of $G$ between $v_1$ and $v_2$, and if there are two edges $e_1,e_2$ of $G$ between $v_1$ and $v_2$, one of $\eta_{c_2}(e_1)$ and $\eta_{c_2}(e_2)$ contains $e'$, and the other contains $\eta(v_2)$ and an edge between $\eta(v_1)$ and $\eta(v_2)$.
But all neighbors of $v_2$ are contained in $V(\gamma_{c_1}) \cap V(\gamma_{c_2})$.
So $v_2$ has no neighbor in $\uparrow c_2-V(\gamma_{c_2})$.
Since $b_c=b_{c'}$, $\eta(v_2)$ has no neighbor in $\uparrow c_2'-V(\gamma_{c_2'})$.
Since $c_2 \in L_{v_2}$, if both $e_1$ and $e_2$ exist, then some of $\eta_{c_2}(e_1)$ and $\eta_{c_2}(e_2)$ contains $\eta(v_2)$ and an edge between $\eta(v_2)$ and a vertex in $\uparrow c_2'-V(\gamma_{c_2'})$, a contradiction.
So there exists an unique edge of $G$ between $v_1$ and $v_2$.
Hence Statement 2 of this claim holds.

So we may assume that $v_2$ is adjacent to some vertex in $V(G)-\uparrow c_2$.
In particular, $v_2$ is associated with essential number 1 or 2 in $\gamma_{c_2}$.
Since $c_2 \in L_{v_2}$, $v_2$ is associated with essential number 1 in $\gamma_{c_2}$ by Claim 1.

Suppose that $v_1$ is not adjacent to a vertex in $\uparrow c_2-V(\gamma_{c_2})$.
Since $v_1$ is associated with essential number 1 or 2 in $\gamma_{c_2}$, by Statement 1 of Claim 2, either $c_2 \not \in L_{v_1}$, or $v_1$ is associated with essential number 2 in $\gamma_{c_2}$.
So $\eta_{c_2}(v_1)=\eta(v_1)$.
That is, $c_2 \not \in L_{v_1}$.
So $c_2 \in L_{v_2}-L_{v_1}$.
By Claim 5, since $v_1$ is not adjacent to a vertex in $\uparrow c_2-V(\gamma_{c_2})$, $v_2$ is not adjacent to any vertex in $V(G)-\uparrow c_2$, a contradiction.

Hence $v_1$ is adjacent to a vertex in $\uparrow c_2-V(\gamma_{c_2})$.
Since $v_1$ is associated with essential number 1 in $\gamma_{c_1}$, there exists a unique vertex, say $w$, in $\uparrow c_2-V(\gamma_{c_2})$ adjacent to $v_1$.
Since $b_{c_1}=b_{c_1'}$ and $\eta(v_1)$ is associated with essential number 1 in $\gamma_{c_1'}$, there exists a unique vertex $w'$ in $\uparrow c_2'-V(\gamma_{c_2'})$ adjacent to $\eta(v_1)$.
If $c_2 \not \in L_{v_1}$, then since $c_2 \in L_{v_2}$ and $v_2$ is associated with essential number 1 in $\gamma_{c_2}$, $w'$ belongs to $\eta_{c_2}(v_1w) \cap \eta_{c_2}(v_1v_2)$, a contradiction.
So $c_2 \in L_{v_1}$.
Since $\lvert L_{v_1} \rvert \leq 2$, $L_{v_1}=\{c_1,c_2\}$.

Since $v_1$ is adjacent to $v_2$, Statement 2 of Claim 3 (taking $v=v_1$ and $u =v_2$) implies that either $\{c_1,c_2\} \subseteq L_{v_2}$, or $\{c_1,c_2\} \cap L_{v_2}=\emptyset$.
Since $c_2 \in L_{v_2}$, $\{c_1,c_2\} \subseteq L_{v_2}$.
Since $v_2$ is associated with essential number 1 in $\gamma_{c_2}$ and $c_2 \in L_{v_2}$, by Statement 1 of Claim 2, $\lvert L_{v_2} \rvert \leq 2$.
So $L_{v_2}=\{c_1,c_2\} = L_{v_1}$, a contradiction.
$\Box$

\medskip

\noindent{\bf Claim 9:} If $v_1,v_2 \in A$ are adjacent vertices, $c_1$ and $c_2$ are distinct children of $r$ with $c_1 \in L_{v_1},c_2 \in L_{v_2}$ and $L_{v_1} \neq L_{v_2}$, then there exists a unique edge $e$ of $G$ between $v_1$ and $v_2$, and either  
	\begin{itemize}
		\item for every $i \in \{1,2\}$, there exists an edge $f_i$ of the rooted extension of $(G[\uparrow c_i],\gamma_{c_i})$ incident with $v_i$ and there exists a path in $\eta_{c_i}(f_i)$ from $\eta_{c_i}(v_i)$ to $\eta(v_i)$ disjoint from $V(\eta(\gamma_{c_i}))-\{\eta(v_i)\}$ such that if $v_1,v_2 \in V(\gamma_{c_i})$, then $f_i=e$, and if $v_{3-i} \not \in V(\gamma_{c_i})$, then $f_i$ is the edge between $v_i$ and a vertex in the indicator, or
		\item $\{v_1,v_2\} \subseteq V(\gamma_{c_1}) \cap V(\gamma_{c_2})$, and there exists $i^* \in \{1,2\}$ such that $L_{v_{i^*}}=\{c_{i^*}\}$, $L_{v_{3-i^*}}=\{c_1,c_2\}$ and all neighbors of $v_{3-i^*}$ are contained in $V(\gamma_{c_1}) \cap V(\gamma_{c_2})$.
	\end{itemize}

\noindent{\bf Proof of Claim 9:}
If $\{v_1,v_2\} \subseteq V(\gamma_{c_1}) \cap V(\gamma_{c_2})$, then this claim immediately follows from Claim 8.

So we may assume that at least one of $v_1,v_2$ does not belong to $V(\gamma_{c_1}) \cap V(\gamma_{c_2})$.
If $v_1 \not \in V(\gamma_{c_2})$ and $v_2 \not \in V(\gamma_{c_1})$, then for each $i \in \{1,2\}$, $v_{3-i}$ is a neighbor of $v_i$ not contained in $\uparrow c_i$, so $v_i$ is associated with essential number 1 in $\gamma_{c_i}$ (by Claim 1 and (BR02)), and hence there exists a unique edge $e$ of $G$ between $v_1$ and $v_2$, and for each $i \in \{1,2\}$, there exists a path in $\eta_{c_i}(f_i)$ from $\eta_{c_i}(v_i)$ to $\eta(v_i)$ disjoint from $V(\eta(\gamma_{c_i}))-\{\eta(v_i)\}$, where $f_i$ is the edge of the rooted extension of $(G[\uparrow c_i],\gamma_{c_i})$ between $v_i$ and the vertex in the indicator adjacent to $v_i$, and hence Statement 1 of this claim holds.

Hence we may assume that $v_1 \not \in V(\gamma_{c_2})$ and $v_2 \in V(\gamma_{c_1})$ by symmetry.
Since $v_1$ is adjacent to $v_2$ and $c_2 \in L_{v_2}$, $v_2$ is associated with essential number 1 in $\gamma_{c_2}$ by Claim 1 and (BR02).
So $v_2$ is not adjacent to any vertex in $\uparrow c_1-V(\gamma_{c_1})$, and there exists a unique edge $e$ of $G$ between $v_1$ and $v_2$.
Since $c_2 \in L_{v_2}$ and $v_2$ is associated with essential number 1 in $\gamma_{c_2}$, there exists a path in $\eta_{c_2}(f)$ from $\eta_{c_2}(v_2)$ to $\eta(v_2)$ disjoint from $V(\eta(\gamma_{c_2}))-\{\eta(v_2)\}$, where $f$ is the edge of the rooted extension of $(G[\uparrow c_2],\gamma_{c_2})$ between $v_2$ and the vertex in the indicator adjacent to $v_2$.

Since $v_2$ is associated with essential number 1 in $\gamma_{c_2}$, $v_2$ is adjacent to a vertex in $\uparrow c_2-V(\gamma_{c_2})$ by Statement 1 of Claim 2.
Since $e$ is an edge incident with $v_2$ whose other end is in $V(G)-\uparrow c_2$, $v_2$ is not adjacent to any vertex in $\uparrow c_1-V(\gamma_{c_1})$.
If $c_1 \in L_{v_1}-L_{v_2}$, then by Claim 5 (taking $c=c_1$, $v_1=v_2$ and $v_2=v_1$), there exists a path in $\eta_{c_1}(e)$ from $\eta_{c_1}(v_1) \not \in V(\overline{H'})$ to $\eta(v_1)$ disjoint from $V(\eta(\gamma_{c_1}))-\{\eta(v_1)\}$, so we are done.

Hence we may assume that $c_1 \not \in L_{v_1}-L_{v_2}$. 
Since $c_1 \in L_{v_1}$, $c_1 \in L_{v_2}$.
Since $v_2$ is adjacent to a vertex in $\uparrow c_2-V(\gamma_{c_2})$, $v_2$ is associated with essential number 1 in $\gamma_{c_1}$ by Claim 1 and (BR02). 
By Statement 1 of Claim 2, $v_2$ is adjacent to some vertex in $\uparrow c_1-V(\gamma_{c_1})$, a contradiction.
$\Box$

\medskip

\noindent{\bf Claim 10:} There exists a function $g$ that maps each vertex $v \in A$ with $L_v \neq \emptyset$ to an element in $L_v$ such that the following statements hold.
	\begin{itemize}
		\item If $x,y \in A$ with $L_x=L_y$, then $g(x)=g(y)$.
		\item If $v_1,v_2 \in A$ are adjacent vertices with $L_{v_1} \neq \emptyset \neq L_{v_2}$ and $g(v_1) \neq g(v_2)$, then there exists a unique edge $e$ of $G$ between $v_1$ and $v_2$, and either
			\begin{itemize}
				\item for every $i \in \{1,2\}$, $v_i$ is adjacent to a vertex in $\uparrow g(v_i)-V(\gamma_{g(v_i)})$, and there exist an edge $f_i$ of the rooted extension of $(G[\uparrow g(v_i)],\gamma_{g(v_i)})$ incident with $v_i$ and a path in $\eta_{g(v_i)}(f_i)$ from $\eta_{g(v_i)}(v_i)$ to $\eta(v_i)$ disjoint from $V(\eta(\gamma_{g(v_i)}))-\{\eta(v_i)\}$ such that if $v_1,v_2 \in V(\gamma_{g(v_i)})$, then $f_i=e$, and if $v_{3-i} \not \in V(\gamma_{g(v_i)})$, then $f_i$ is the edge between $v_i$ and a vertex in the indicator, or
				\item $\{v_1,v_2\} \subseteq V(\gamma_{g(v_1)}) \cap V(\gamma_{g(v_2)})$, and there exist $i^* \in \{1,2\}$ and $c^* \in L_{v_{i^*}}$ such that $L_{v_{i^*}}=\{g(v_{i^*}),c^*\}$, $v_{i^*}$ is not adjacent to a vertex in $\uparrow g(v_{i^*})-V(\gamma_{g(v_{i^*})})$, and for every neighbor $z \in A$ of $v_{i^*}$ with $L_z \neq \emptyset$ and $g(z) \neq g(v_{i^*})$, we have
					\begin{itemize}
						\item $L_z = \{c^*\}$, $g(z)=c^*$,
						\item there exists a unique edge of $G$ between $z$ and $v_{i^*}$, 
						\item $z$ is adjacent to a vertex in $\uparrow g(z)-V(\gamma_{g(z)})$,
						\item $z$ is associated with essential number 1 in $\gamma_{g(z)}$,
						\item $\eta_{g(v_{i^*})}(v_{i^*}z)$ is a path in $G$ from $\eta_{g(v_{i^*})}(v_{i^*})$ to $\eta_{g(v_{i^*})}(z)=\eta(z)$ disjoint from $V(\eta(\gamma_{g(v_{i^*})}))-\{\eta(z)\}$, and
						\item there exists a path in $\eta_{c^*}(f_z)$ from $\eta_{c^*}(z)$ to $\eta(z)$ disjoint from $V(\eta(\gamma_{c^*}))-\{\eta(z)\}$, where $f_z$ is the edge of the rooted extension of $(G[\uparrow c^*],\gamma_{c^*})$ between $z$ and a vertex in the indicator. 
					\end{itemize}
			\end{itemize}
	\end{itemize}

\noindent{\bf Proof of Claim 10:}
Clearly, there exists a function that maps each vertex $v \in A$ with $L_v \neq \emptyset$ to an element in $L_v$ such that if $x,y \in A$ with $L_x=L_y$, then $g(x)=g(y)$.
We shall prove that $g$ satisfies the second statement of this claim.

Let $A_1 = \{x \in A: L_x \neq \emptyset, x$ is adjacent to a vertex in $\uparrow g(x)-V(\gamma_{g(x)})\}$.
Let $A_0 = \{x \in A: L_x \neq \emptyset, x$ is not adjacent to a vertex in $\uparrow g(x)-V(\gamma_{g(x)})\}$.

Let $v_1,v_2 \in A$ be adjacent vertices with $g(v_1) \neq g(v_2)$.
Note that $g(v_1) \in L_{v_1}$, $g(v_2) \in L_{v_2}$, and $g(v_1) \neq g(v_2)$.
So $L_{v_1} \neq L_{v_2}$ by Statement 1.
By Claim 9, there uniquely exists an edge of $G$ between $v_1$ and $v_2$.

We first assume that $\{v_1,v_2\} \subseteq A_1$.
Since $\{v_1,v_2\} \subseteq A_1$, for every $i \in \{1,2\}$, $v_i$ is adjacent to a vertex in $\uparrow g(v_i)-V(\gamma_{g(v_i)})$, and by Claim 9, there exist an edge $f_i$ of the rooted extension of $(G[\uparrow g(v_i)],\gamma_{g(v_i)})$ incident with $v_i$ and a path in $\eta_{g(v_i)}(f_i)$ from $\eta_{g(v_i)}(v_i)$ to $\eta(v_i)$ disjoint from $V(\eta(\gamma_{g(v_i)}))-\{\eta(v_i)\}$ such that if $v_1,v_2 \in V(\gamma_{g(v_i)})$, then $f_i=e$, and if $v_{3-i} \not \in V(\gamma_{g(v_i)})$, then $f_i$ is the edge between $v_i$ and a vertex in the indicator, so this claim holds.

Hence we may assume that at least one of $v_1$ and $v_2$ is not in $A_1$.
So at least one of $v_1$ and $v_2$ is in $A_0$.
By Claim 9, there exists a unique edge $e$ of $G$ between $v_1$ and $v_2$, $\{v_1,v_2\} \subseteq V(\gamma_{g(v_1)}) \cap V(\gamma_{g(v_2)})$, and there exists $i^* \in \{1,2\}$ such that $L_{v_{3-i^*}}=\{g(v_{3-i^*})\}$, $L_{v_{i^*}}=\{g(v_1),g(v_2)\}$ and all neighbors of $v_{i^*}$ are contained in $V(\gamma_{g(v_1)}) \cap V(\gamma_{g(v_2)})$.
Since all neighbors of $v_{i^*}$ are contained in $V(\gamma_{g(v_1)}) \cap V(\gamma_{g(v_2)})$, $v_{i^*} \in A_0$.
Since $g(v_{i^*}) \in L_{v_{i^*}}-L_{v_{3-i^*}}$ and $v_{i^*} \in A_0$, $\eta_{g(v_{i^*})}(e)$ is a path in $G$ from $\eta_{g(v_{i^*})}(v_{i^*})$ to $\eta_{g(v_{i^*})}(v_{3-i^*})=\eta(v_{3-i^*})$ disjoint from $V(\gamma_{g(v_{i^*})})-\{\eta(v_{3-i^*})\}$ and shows that $v_{3-i^*}$ is adjacent to a vertex in $\uparrow g(v_{i^*})-V(\gamma_{g(v_{i^*})})$.
Since $g(v_{3-i^*}) \in L_{v_{3-i^*}}$, $v_{3-i^*}$ is associated with essential number 1 in $\gamma_{g(v_{3-i^*})}$ by Claim 1 and (BR02).
So $v_{3-i^*} \in A_1$ by Statement 1 of Claim 2, and there exists a path in $\eta_{g(v_{3-i^*})}(f_{v_{3-i^*}})$ from $\eta_{g(v_{3-i^*})}(v_{3-i^*})$ to $\eta(v_{3-i^*})$ disjoint from $V(\eta(\gamma_{g(v_{3-i^*})}))-\{\eta(v_{3-i^*})\}$, where $f_{v_{3-i^*}}$ is the edge of the rooted extension of $(G[\uparrow g(v_{3-i^*})],\gamma_{g(v_{3-i^*})})$ between $v_{3-i^*}$ and a vertex in the indicator. 

Let $c^* = g(v_{3-i^*})$.

Let $Z = \{z \in A: zv_{i^*} \in E(G), L_z \neq \emptyset, g(z) \neq g(v_{i^*})\}$.
Since $v_{i^*} \in A_0$, for every $z \in Z$, the above argument (by taking $\{v_1,v_2\}=\{v_{i^*},z\}$) shows that $L_{z}=\{g(z)\}$ and $L_{v_{i^*}}=\{g(v_{i^*}),g(z)\}$; since $L_{v_{i^*}}=\{g(v_1),g(v_2)\}$, we know $g(z)=c^*$; hence the above argument (by taking $\{v_1,v_2\}=\{v_{i^*},z\}$) shows that $\{z,v_{i^*}\} \subseteq V(\gamma_{g(z)}) \cap V(\gamma_{g(v_{i^*})})$, and $L_{z}=\{c^*\}$, $L_{v_{i^*}}=\{g(v_{i^*}),g(z)\}$, all neighbors of $v_{i^*}$ are contained in $V(\gamma_{g(v_{i^*})}) \cap V(\gamma_{g(z)})$, $z \in A_1$ and is associated with essential number 1 in $\gamma_{c^*}$, there exists a unique edge $e_z$ of $G$ between $v_{i^*}$ and $z$, $\eta_{g(v_{i^*})}(e_z)$ is a path in $G$ from $\eta_{g(v_{i^*})}(v_{i^*})$ to $\eta_{g(v_{i^*})}(z)=\eta(z)$ disjoint from $V(\gamma_{g(v_{i^*})})-\{\eta(z)\}$, and there exists a path in $\eta_{c^*}(f_z)$ from $\eta_{c^*}(z)$ to $\eta(z)$ disjoint from $V(\eta(\gamma_{c^*}))-\{\eta(z)\}$, where $f_z$ is the edge of the rooted extension of $(G[\uparrow c^*],\gamma_{c^*})$ between $z$ and a vertex in the indicator. 
So this claim holds.
$\Box$

\medskip

Let $g$ be a function that maps each vertex $v \in A$ with $L_v \neq \emptyset$ to an element in $L_v$ satisfying Claim 10. 
Note that for every vertex $v \in V(G)-V(H)$, there exists a unique child $c$ of $r$ such that $v \in \uparrow c-V(\gamma_c)$, and we define $g(v)=c$.
Hence $g$ is a function whose domain is $\{v \in A: L_v \neq \emptyset\} \cup (V(G)-V(H))$.

Note that if $v \in A$ with $L_v=\emptyset$, then $\eta_c(v)$ equals $\eta(v)$ for all children $c$ of $r$ with $v \in V(\gamma_c)$.

Now we define $\pi_V:V(\overline{G}) \rightarrow V(\overline{G'})$ as follows.
	\begin{itemize}
		\item For each vertex $v \in V(G)-V(H)$, there exists a unique child $c$ of $r$ such that $v \in \uparrow c$, and we define $\pi_V(v)=\eta_c(v)$. 
		\item For each vertex $v \in A$ with $L_v \neq \emptyset$, define $\pi_V(v)=\eta_{g(v)}(v)$.
		\item For each vertex $v \in A$ with $L_v=\emptyset$, define $\pi_V(v)=\eta(v)$. (Recall that $\eta(v)=\eta_c(v)$ for every child $c$ of $r$ with $v \in V(\gamma_c)$ in this case.)
		\item For each vertex $v \in V(H)-A$, define $\pi_V(v)=\eta(v)$.
		\item For each $j \in [\lvert \gamma_0 \rvert]$, $\pi_V$ maps the $j$-th entry of the indicator of $\overline{G}$ to the $j$-th entry of the indicator of $\overline{G'}$. 
	\end{itemize}
Note that if $v \in A$ and $\pi_V(v) \in V(\overline{H'})$, then $L_v=\emptyset$ and $\pi_V(v)=\eta(v)$.

Clearly, $\pi_V$ is an injection.

We define $\pi_E$ to be a function that maps each edge $e$ of $\overline{G}$, say with ends $u,v$, to a subgraph of $\overline{G'}$ as follows.
(Note that $u=v$ when $e$ is a loop.)
	\begin{itemize}
		\item If $\pi_V(u),\pi_V(v) \in V(\overline{H'})$, then we know $\pi_V(u)=\eta(u)$ and $\pi_V(v)=\eta(v)$, and we define $\pi_E(e)=\eta(e)$.
		\item If $\pi_V(u) \in V(\overline{H'})$ and $\pi_V(v) \not \in V(\overline{H'})$, then we know $\pi_V(u)=\eta(u)$, and $g(v)$ is the unique child $c$ of $r$ such that $\pi_V(v)=\eta_c(v)$, and we define $\pi_E(e)$ as follows.
			\begin{itemize}
				\item If $u \in \uparrow g(v)$, then we know $\eta_{g(v)}(u)=\eta(u)=\pi_V(u)$, and we define $\pi_E(e)$ to be $\eta_{g(v)}(e)$. 
				\item If $u \not \in \uparrow g(v)$, then we know that $v \in V(\gamma_{g(v)}) \subseteq V(\overline{H})$, $v$ is associated with essential number 1 in $\gamma_{g(v)}$, $u \in V(\overline{H})$ and there exists a path in the image of $\eta_{g(v)}$ from $\eta_{g(v)}(v)$ to $\eta(v)$ disjoint from $V(\eta(\gamma_{g(v)}))-\{\eta(v)\}$, and we define $\pi_E(e)$ to be the path obtained by concatenating the path $\eta(e)$ and the path in the image of $\eta_{g(v)}$ from $\eta_{g(v)}(v)$ to $\eta(v)$ just mentioned.
			\end{itemize}
		\item If $\pi_V(u),\pi_V(v) \not \in V(\overline{H'})$, then $g(u),g(v)$ are the unique children $c_u,c_v$ of $r$, respectively, such that $\pi_V(u)=\eta_{c_u}(u)$ and $\pi_V(v)=\eta_{c_v}(v)$, and we define $\pi_E(e)$ as follows.
			\begin{itemize}
				\item If $g(u)=g(v)$, then define $\pi_E(e)=\eta_{g(u)}(e)$.
				\item If $g(u) \neq g(v)$ and $e \in E(H)$, then 
					\begin{itemize}
						\item if $u$ is adjacent to a vertex in $\uparrow g(u)-V(\gamma_{g(u)})$, and $v$ is adjacent to a vertex in $\uparrow g(v)-V(\gamma_{g(v)})$, then define $\pi_E(e)$ to be the path obtained by concatenating the path in the image of $\eta_{g(u)}$ from $\eta_{g(u)}(u)$ to $\eta(u)$ disjoint from $V(\eta(\gamma_{g(u)}))-\{\eta(u)\}$ mentioned in Claim 10, the path $\eta(e)$, and the path in the image of $\eta_{g(v)}$ from $\eta(v)$ to $\eta_{g(v)}(v)$ disjoint from $V(\eta(\gamma_{g(v)}))-\{\eta(v)\}$ mentioned in Claim 10,
						\item otherwise, we know $\{u,v\} \subseteq V(\gamma_{g(u)}) \cap V(\gamma_{g(v)})$ by Claim 10, and we may assume by symmetry that $u$ is not adjacent to any vertex in $\uparrow g(u)-V(\gamma_{g(u)})$, and we define $\pi_E(e)$ to be the path obtained by concatenating the path $\eta_{g(u)}(e)$ and the path in the image of $\eta_{g(v)}$ from $\eta_{g(v)}(v)$ to $\eta(v)$ disjoint from $V(\gamma_{g(v)})-\{\eta(v)\}$ mentioned in Claim 10.
					\end{itemize}
				\item If $g(u) \neq g(v)$ and $e \not \in E(H)$, then by symmetry we may assume that $u \in V(H)$ and $v \in \uparrow g(v)-V(H)$, so $u$ is associated with essential number 1 in $\gamma_{g(u)}$ and with essential number 1 or 2 in $\gamma_{g(v)}$ (by Claim 2), and we define $\pi_E(e)$ to be the path obtained by concatenating the path in the image of $\eta_{g(u)}$ from $\eta_{g(u)}(u)$ to $\eta(u)$ disjoint from $V(\eta(\gamma_{g(u)}))-\{\eta(u)\}$, the path in the image of $\eta_{g(v)}$ from $\eta(u)$ to $\eta_{g(v)}(u)$ disjoint from $V(\eta(\gamma_{g(v)}))-\{\eta(u)\}$, and $\eta_{g(v)}(e)$.
			\end{itemize}
	\end{itemize}

Now we show that $(\pi_V,\pi_E)$ is a homoemorphic embedding.

It is straightforward to check that for every $e \in E(\overline{G})$, say with ends $u$ and $v$, the intersection of $\pi_E(e)$ and the image of $\pi_V$ is $\{\pi_V(u),\pi_V(v)\}$.
Note that for each $e \in E(G)$, $\pi_E(e)$ is contained in a union of $\eta(e)$ (if $\eta(e)$ is defined) and subpaths of $\eta_c(e')$ intersecting $V(\overline{H'})$ only at vertices in $\{\eta(u),\eta(v)\}$ (ignore $\eta(u)$ or $\eta(v)$ if it is undefined), for some children $c$ of $r$ and edges $e'$ of the rooted extension of $(G[\uparrow c], \gamma_c)$ with ends $u$ and $v$ such that either $e'=e$ or $e'$ is the edge between $\{u,v\}$ and a vertex in the indicator by Claim 10.

\medskip

\noindent{\bf Claim 11:} If $e_1,e_2 \in E(G)$ such that $\pi_E(e_1) \cap \pi_E(e_2) \not \subseteq \bigcup_{x \in V(e_1 \cap e_2)} \pi_V(x)$, then each of $e_1$ and $e_2$ is incident with a vertex in $V(H)$.

\noindent{\bf Proof of Claim 11:}
By symmetry, suppose to the contrary that $e_1$ is not incident with a vertex in $V(H)$.
So there exists a child $c$ of $r$ such that both ends of $e_1$ are contained in $\uparrow c-V(\gamma_c)$.
Let $u_1$ and $v_1$ be the ends of $e_1$.
Hence $\pi_V(u_1)=\eta_c(u_1) \not \in V(H')$, $\pi_V(v_1)=\eta_c(v_1) \not \in V(H')$, and $\pi_E(e_1)=\eta_c(e_1)$ is contained in $G[\uparrow c-V(\gamma_c)]$.
Since $\eta_c$ is a homeomorphic embedding, $\pi_E(e_1) \cap \pi_E(e_2) \subseteq \bigcup_{x \in V(e_1 \cap e_2)}\pi_V(x)$, a contradiction.
$\Box$

\medskip

\noindent{\bf Claim 12:} If $e_1,e_2 \in E(G)$ such that $\pi_E(e_1) \cap \pi_E(e_2) \not \subseteq \bigcup_{x \in V(e_1 \cap e_2)} \pi_V(x)$, then each of $e_1$ and $e_2$ is an edge of $H$.

\noindent{\bf Proof of Claim 12:}
By symmetry, suppose to the contrary that $e_1$ is not an edge of $H$.
By Claim 11, $e_1$ has exactly one end in $V(H)$.
Since $e_1$ is not an edge of $H$, $e_1$ is not a loop.
Let $u$ and $v$ be the ends of $e_1$.
By symmetry, we may assume that there exists a child $c$ of $r$ such that $v \in \uparrow c-V(\gamma_c)$ and $u \in V(H)$.
Note that $\pi_V(v) \not \in V(H')$ and $c=g(v)$.

Suppose $\pi_E(e_1)=\eta_{g(v)}(e_1)$.
Then $V(\pi_E(e_1)) \cap V(H') \subseteq \{\eta(u)\}$.
Since $\pi_E(e_1) \cap \pi_E(e_2) \not \subseteq \bigcup_{x \in V(e_1 \cap e_2)} \pi_V(x)$, $V(\pi_E(e_1)) \cap V(H') \neq \emptyset$ and $u$ is an end of $e_2$.
So $V(\pi_E(e_1)) \cap V(H') = \{\eta(u)\}$.
If $\eta(u)$ is not an internal vertex of $\pi_E(e_1)$, then since $V(\pi_E(e_1)) \cap V(H') = \{\eta(u)\}$, $\pi_V(u)=\eta(u)$, so $\pi_E(e_1) \cap \pi_E(e_2) \subseteq \bigcup_{x \in V(e_1 \cap e_2)} \pi_V(x)$, a contradiction.
Hence $\eta(u)$ is an internal vertex of $\pi_E(e_1)$.
So $\eta_{g(v)}(u) = \pi_V(u) \not \in V(H')$, $g(u)=g(v)$, and $u$ is associated essential number 0 in $\gamma_{g(u)}$.
Hence the end of $e_2$ other than $u$, denoted by $v_2$, belongs to $\uparrow g(u)$.
If $\pi_E(e_2)=\eta_{g(u)}(e_2)$, then $\pi_E(e_1) \cap \pi_E(e_2) \subseteq \bigcup_{x \in V(e_1 \cap e_2)} \pi_V(x)$, a contradiction.
So $\pi_E(e_2) \neq \eta_{g(u)}(e_2)$.
Hence $\pi_V(v_2) \not \in V(H')$, $g(v_2)$ is defined, and $g(v_2) \neq g(u)$.
Since $v_2 \in \uparrow g(u)$ and $g(v_2) \neq g(u)$, $e_2 \in E(H)$.
Note that $u$ is adjacent to a vertex in $\uparrow g(u)-V(\gamma_{g(u)})$.
If $v_2$ is adjacent to a vertex in $\uparrow g(v_2)-V(\gamma_{g(v_2)})$, then by Lemma 10, $\eta(u)$ is an internal vertex of both $\eta_{g(u)}(e_1)$ and $\eta_{g(u)}(e_1')$ for some edge $e_1'$ of the rooted extension of $(G[\uparrow g(u)],V(\gamma_{g(u)}))$ different from $e_1$, a contradiction. 
So $v_2$ is not adjacent to a vertex in $\uparrow g(v_2)-V(\gamma_{g(v_2)})$.
By Claim 10, $u$ is associated with essential number 1 in $\gamma_{g(u)}$, a contradiction.

Hence $\pi_E(e_1) \neq \eta_{g(v)}(e_1)$.
Since $u \in \uparrow g(v)$, $\pi_V(u) \not \in V(H')$.
So $g(u)$ is defined and $g(u) \neq g(v)$.
Since $v \in \uparrow g(v)-V(\gamma_{g(v)})$ and $g(u) \neq g(v)$, $u$ is associated with essential number 1 in $\gamma_{g(u)}$.
Since $g(u) \neq g(v)$ and $e_1 \not \in E(H)$, $\pi_E(e_1)$ is obtained from by concatenating the path $P$ in the image of $\eta_{g(u)}$ from $\eta_{g(u)}(u)$ to $\eta(u)$ disjoint from $V(\eta(\gamma_{g(u)}))-\{\eta(u)\}$ and a path in the image of $\eta_{g(v)}$ only intersecting $V(H')$ at $\eta(u)$.
Since $u$ is associated with essential number 1 in $\gamma_{g(u)}$, $P$ is a subpath of $\eta_{g(u)}(f_u)$, where $f_u$ is the edge of the rooted extension of $(G[\uparrow g(u)],\gamma_{g(u)})$ between $u$ and a vertex in the indicator.

So $\pi_E(e_1) \cap \pi_E(e_2) - \bigcup_{x \in V(e_1 \cap e_2)} \pi_V(x) = \{\eta(u)\}$ and $u$ is an end of $e_2$.
This together with the fact that $u$ is associated with essential number 1 in $\gamma_{g(u)}$ imply that $e_2$ is not a loop.
Let $v_2$ be the end of $e_2$ other than $u$.
Since $u$ is associated with essential number 1 in $\gamma_{g(u)}$ and $e_1$ is between $u$ and $V(G)-\uparrow g(u)$, $v_2 \in \uparrow g(u)$.
Since $P \subseteq \eta_{g(u)}(f_u)$, $\pi_E(e_2) \neq \eta_{g(u)}(e_2)$, for otherwise $\pi_E(e_1) \cap \pi_E(e_2) \subseteq \bigcup_{x \in V(e_1 \cap e_2)} \pi_V(x)$.
So $\pi_V(v_2) \not \in V(H')$.
Hence $g(v_2)$ is defined, $g(v_2) \neq g(u)$, and $v_2 \in V(\gamma_{g(u)})$.

Suppose $v_2$ is adjacent to a vertex in $\uparrow g(v_2)-V(\gamma_{g(v_2)})$.
Since $u$ is adjacent to a vertex in $\uparrow g(u)-V(\gamma_{g(u)})$ and $u,v_2 \in V(\gamma_{g(u)})$, by Claim 10, $\eta(u) \in V(P) \cap V(\eta_{g(u)}(e_2)) \subseteq V(\eta_{g(u)}(f_u)) \cap V(\eta_{g(u)}(e_2)) = \{\eta_{g(u)}(u)\}$, a contradiction.

So $v_2$ is not adjacent to any vertex in $\uparrow g(v_2)-V(\gamma_{g(v_2)})$.
By Claim 10, $\{u,v_2\} \subseteq V(\gamma_{g(u)}) \cap V(\gamma_{g(v_2)})$, and $g(v_2) \not \in L_u$. 
Since $g(v_2) \not \in L_u$ and $v_2$ is not adjacent to any vertex in $\uparrow g(v_2)-V(\gamma_{g(v_2)})$, there exists an edge $e_2'$ of $\eta_{g(v_2)}(e_2)$ incident with $\eta(u)$ whose other end is in $\uparrow g(v_2)' - V(\gamma_{g(v_2)'})$, where $g(v_2)'$ is the child of $r'$ with $\iota(\gamma_{g(v_2)})=\gamma_{g(v_2)'}$.
Since $u$ is associated with essential number 1 in $\gamma_{g(u)}$, $g(v_2)=g(v)$.
Hence $\eta(u)$ is incident with an edge $e_1'$ in $\eta_{g(v_2)'}(e_1)$. 
Since $e_1 \neq e_2$, $e_1'$ and $e_2'$ are distinct.
So $u$ is not associated with essential number 1 in $\gamma_{g(u)}$, a contradiction. 
$\Box$

\medskip

\noindent{\bf Claim 13:} There exist no distinct edges $e_1,e_2 \in E(G)$ such that $\pi_E(e_1) \cap \pi_E(e_2) \not \subseteq \bigcup_{x \in V(e_1 \cap e_2)} \pi_V(x)$.

\noindent{\bf Proof of Claim 13:}
Suppose to the contrary that there exist distinct $e_1,e_2 \in E(G)$ such that $\pi_E(e_1) \cap \pi_E(e_2) \not \subseteq \bigcup_{x \in V(e_1 \cap e_2)} \pi_V(x)$.
By Claim 12, $e_1$ and $e_2$ are edges of $H$.
For each $i \in \{1,2\}$, let $u_i,v_i$ be the ends of $e_i$.
Recall that $\pi_E(e_1) \cap \pi_E(e_2) - \bigcup_{x \in V(e_1 \cap e_2)} \pi_V(x) \subseteq \{\eta(u_1),\eta(v_1),\eta(u_2),\eta(v_2)\}$.
And for each $i \in \{1,2\}$, $V(\pi_E(e_i)) \cap V(H') \subseteq V(\eta(e_i))$.
So $\{u_1,v_1\} \cap \{u_2,v_2\} \neq \emptyset$.
By symmetry, we may assume that $u_1=u_2$.

Suppose $\pi_V(u_1) \in V(H')$.
If $\pi_V(v_1) \in V(H')$, then $\pi_V(u_1)=\eta(u_1)$ and $\pi_V(v_1)=\eta(v_1)$, so $\pi_E(e_1) \cap \pi_E(e_2) \subseteq \bigcup_{x \in V(e_1 \cap e_2)} \pi_V(x)$.
So $\pi_V(v_1) \not \in V(H')$.
Similarly, $\pi_V(v_2) \not \in V(H')$.
Hence $g(v_1)$ and $g(v_2)$ are defined.
If $u_1 \in \uparrow g(v_1) \cap \uparrow g(v_2)$, then $\pi_E(e_1) \cap \pi_E(e_2) = \eta_{g(v_1)}(e_1) \cap \eta_{g(v_2)}(e_2) \subseteq \bigcup_{x \in V(e_1 \cap e_2)} \pi_V(x)$.
So by symmetry, we may assume that $u_1 \not \in \uparrow g(v_1)$.
Hence $v_1$ is associated with essential number 1 in $\gamma_{g(v_1)}$, so $e_1$ is the unique edge between $v_1$ and $V(G)-\uparrow g(v_1)$, and hence $v_1 \neq v_2$.
So $\pi_E(e_1)$ is contained in $\eta(e_1) \cup \eta_{g(v_1)}(f_1)$, where $f_1$ is the edge of the rooted extension of $(G[\uparrow g(v_1)], \gamma_{g(v_1)})$ between $v_1$ and a vertex in the indicator.
Since $\pi_E(e_1) \cap \pi_E(e_2) \not \subseteq \bigcup_{x \in V(e_1 \cap e_2)} \pi_V(x)$ and $v_1 \neq v_2$, we have $u_1 \not \in \uparrow g(v_2)$, $g(v_1) \neq g(v_2)$, and $v_2$ is associated with essential number 1 in $\gamma_{v_2}$.
Since $\pi_V(u_1)=\eta(u_1)$, $\{\eta(v_1), \eta(v_2)\} \cap V(\pi_E(e_1)) \cap V(\pi_E(e_2)) \neq \emptyset$.
By symmetry, we may assume that $\eta(v_1) \in V(\pi_E(e_1)) \cap V(\pi_E(e_2))$.
Since $g(v_1) \neq g(v_2)$, $v_1$ is incident with an edge whose other end is in $\uparrow g(v_2)-V(\gamma_{g(v_2)})$ which does not contain $u$.
So $v_1$ is not associated with essential number 1 in $\gamma_{g(v_1)}$, a contradiction.

So $\pi_V(u_1) \not \in V(H')$.
Suppose $\pi_V(v_1) \in V(H')$ and $\pi_V(v_2) \in V(H')$.
Since $\pi_E(e_1) \cap \pi_E(e_2) \not \subseteq \bigcup_{x \in V(e_1 \cap e_2)} \pi_V(x)$, either $v_1 \not \in \uparrow g(u_1)$ or $v_2 \not \in \uparrow g(u_1)$.
So $u_1$ is associated with essential number 1 in $\gamma_{g(u_1)}$.
Hence exactly one of $v_1$ and $v_2$ is not in $\uparrow g(u_1)$.
But it implies that $\pi_E(e_1) \cap \pi_E(e_2) \subseteq \bigcup_{x \in V(e_1 \cap e_2)} \pi_V(x)$ since $u_1$ is associated with essential number 1 in $\gamma_{g(u_1)}$, a contradiction.

Hence by symmetry, we may assume that $\pi_V(u_1) \not \in V(H')$ and $\pi_V(v_1) \not \in V(H')$.

Suppose $\pi_E(e_1)=\eta_{g(u_1)}(e_1)$ and $\pi_V(v_2) \in V(H')$.
In particular, $g(u_1)=g(v_1)$.
Then $v_2 \not \in \uparrow g(u_1)$, for otherwise $\pi_E(e_1) \cap \pi_E(e_2) \subseteq \bigcup_{x \in V(e_1 \cap e_2)} \pi_V(x)$.
So $u_1$ is associated with essential number 1 in $\gamma_{g(u_1)}$, and $e_2$ is the unique edge incident with $u_1$ whose other end is in $V(G)-\uparrow g(u_1)$.
Since $g(u_1)=g(v_1)$, $v_1 \in V(\gamma_{g(u_1)})$, so $v_1 \neq v_2$.
So $\pi_E(e_1) \cap \pi_E(e_2) \subseteq \bigcup_{x \in V(e_1 \cap e_2)} \pi_V(x)$, a contradiction.

Hence either $\pi_E(e_1) \neq \eta_{g(u_1)}(e_1)$ or $\pi_V(v_2) \not \in V(H')$.
Suppose $\pi_E(e_1) = \eta_{g(u_1)}(e_1)$.
In particular, $g(u_1)=g(v_1)$ and $\pi_V(v_2) \not \in V(H')$.
Hence $g(u_1) \neq g(v_2)$, for otherwise $\pi_E(e_1) \cap \pi_E(e_2) \subseteq \bigcup_{x \in V(e_1 \cap e_2)} \pi_V(x)$.
Since $g(v_1)=g(u_1) \neq g(v_2)$, $v_1 \neq v_2$.
Since $\pi_E(e_1) = \eta_{g(u_1)}(e_1)$ and $g(u_1) \neq g(v_2)$, by Claim 10, $\pi_E(e_1) \cap \pi_E(e_2) \subseteq \bigcup_{x \in V(e_1 \cap e_2)} \pi_V(x)$, a contradiction.

Hence $\pi_E(e_1) \neq \eta_{g(u_1)}(e_1)$.
So $g(u_1) \neq g(v_1)$, and hence $u_1 \neq v_1$.

Suppose $v_2 \not \in \uparrow g(u_1)$.
So $u_1$ is associated with essential number 1 in $\gamma_{g(u_1)}$, and $e_2$ is the unique edge of $G$ incident with $u_1$ whose other end is in $V(G)-\uparrow g(u_1)$.
Hence $v_1 \in \uparrow g(u_1)$, and $\eta(u_1)$ is contained in $\eta_{g(u_1)}(f_u)$, where $f_u$ is the edge of the rooted extension of $(G[\uparrow g(u_1)],\gamma_{g(u_1)})$ between $u_1$ and a vertex in the indicator.
Since $e_1 \in E(H)$, $v_1 \in V(\gamma_{g(u_1)})$.
Since $g(u_1) \in L_{u_1}$, $u_1$ is adjacent to a vertex in $\uparrow g(u_1)-V(\gamma_{g(u_1)})$ by Claim 2.
If $v_1$ is adjacent to a vertex in $\uparrow g(v_1)-V(\gamma_{g(v_1)})$, then by Claim 10, since $v_1 \in V(\gamma_{u_1})$, $\eta(u_1) \in V(\eta_{g(u_1)}(e_1))$, contradicting $\eta(u_1) \in V(\eta_{g(u_1)}(f_u))-\{\eta_{g(u_1)}(u_1)\}$.
So $v_1$ is not adjacent to a vertex in $\uparrow g(v_1)-V(\gamma_{g(v_1)})$.
Then $\pi_E(e_1)$ contains $\eta_{g(v_1)}(e_1)$ which contains an edge incident with $\eta(u_1)$ whose other end is in $\uparrow g(v_1)'-V(\gamma_{g(v_1)'})$ by Claim 10, where $g(v_1)'$ is the child of $r'$ with $\iota(\gamma_{g(v_1)})=\gamma_{g(v_1)'}$.
Since $b_{g(v_1)}=b_{g(v_1)'}$, $u_1$ is incident with an edge whose other end is in $\uparrow g(v_1)-V(\gamma_{g(v_1)})$ which is different from $e_2$, a contradiction.

Hence $v_2 \in \uparrow g(u_1)$.
Suppose $\pi_V(v_2) \in V(H')$.
Then $\pi_E(e_2)=\eta_{g(u_1)}(e_2)$ and $v_1 \neq v_2$.
So by Claim 10, $\pi_E(e_1) \cap \pi_E(e_2) \subseteq \bigcup_{x \in V(e_1 \cap e_2)} \pi_V(x)$.

So $\pi_V(v_2) \not \in V(H')$.
Note that we showed that $\pi_V(u_1) \not \in V(H')$ and $\pi_V(v_1) \not \in V(H')$ imply $g(u_1) \neq g(v_1)$ and $v_2 \in \uparrow g(u_1)$.
As $\pi_V(u_1) \not \in V(H')$ and $\pi_V(v_2) \not \in V(H')$, a similar argument shows that $g(u_1) \neq g(v_2)$ and $v_1 \in \uparrow g(u_1)$.
By Claim 10, for each $i \in \{1,2\}$, there unique exists an edge between $u_1$ and $v_i$.
So $v_1 \neq v_2$.

Since $v_1 \neq v_2$, if $u_1$ is not adjacent to a vertex in $\uparrow g(u_1)-V(\gamma_{g(u_1)})$, then by Claim 10, $V(\pi_E(e_1)) \cap V(\pi_E(e_2)) = \{\pi_V(u_1)\}$, a contradiction.
So $u_1$ is adjacent to a vertex in $\uparrow g(u_1)-V(\gamma_{g(u_1)})$.

Suppose for every $i \in \{1,2\}$, $v_i$ is adjacent to a vertex in $\uparrow g(v_i)-V(\gamma_{g(v_i)})$.
Since $\{v_1,v_2\} \subseteq \uparrow g(u_1)$, by Claim 10, $\eta(u_1) \in V(\eta_{g(u_1)}(e_1)) \cap V(\eta_{g(u_1)}(e_2)) = \{\eta_{g(u_1)}(u_1)\}$, contradicting the definition of $g(u_1)$.

So by symmetry, we may assume that $v_1$ is not adjacent to a vertex in $\uparrow g(v_1)-V(\gamma_{g(v_1)})$.
By Claim 10, $\eta(u_1)$ is incident with an edge $e_1^*$ in $\eta_{g(v_1)}(e_1)$ whose other end is in $\uparrow g(v_1)'-V(\gamma_{g(v_1)'})$, where $g(v_1)'$ is the child of $r'$ with $\iota(\gamma_{g(v_1)}) = \gamma_{g(v_1)'}$.
So $u_1$ is associated with essential number 1 in $\gamma_{g(u_1)}$, and $\eta(u_1)$ is contained in $\eta_{g(u_1)}(f_u)$, where $f_u$ is the edge of the rooted extension of $(G[\uparrow g(u_1)],\gamma_{g(u_1)})$ between $u_1$ and a vertex in the indicator.
If $v_2$ is adjacent to a vertex in $\uparrow g(v_2)-V(\gamma_{g(v_2)})$, then by Claim 10, $\eta(u_1) \in V(\eta_{g(u_1)}(e_2)) \cap V(\eta_{g(u_1)}(f_u)) = \{\eta_{g(u_1)}(u_1)\}$, a contradiction.
So $v_2$ is not adjacent to a vertex in $\uparrow g(v_2)-V(\gamma_{g(v_2)})$.
By Claim 10, $\eta(u_1)$ is incident with an edge $e_2^*$ in $\eta_{g(v_2)}(e_2)$ whose other end is in $\uparrow g(v_2)'-V(\gamma_{g(v_2)'})$, where $g(v_2)'$ is the child of $r'$ with $\iota(\gamma_{g(v_2)}) = \gamma_{g(v_2)'}$.
Since $e_1 \neq e_2$, $e_1^* \neq e_2^*$.
Since $g(u_1) \neq g(v_1)$ and $g(u_1) \neq g(v_2)$, $u_1$ is not associated with essential number 1 in $\gamma_{g(u_1)}$, a contradiction.
This proves the claim.
$\Box$

\medskip

By Claim 13, $(\pi_V,\pi_E)$ defines a homeomorphic embedding from $\overline{G}$ to $\overline{G'}$ such that for every $i \in [\lvert V(\gamma_0) \rvert]$, $\pi_V$ maps the $i$-th entry in the indicator of $\overline{G}$ to the $i$-th entry in the indicator of $\overline{G'}$. 

\medskip

\noindent{\bf Claim 14:} If $v$ is the $i$-th vertex of $\gamma_{0}'$ for some $i \in [\lvert V(\gamma_{0}') \rvert]$ and is a vertex in $\pi_E(e)-\pi_V(V(\overline{G}))$ for some edge $e$ of $\overline{G}$, then either $e$ is an edge incident with the $i$-th vertex in the indicator of $\overline{G}$, or the essential number associated with the $i$-th vertex in $\gamma_0$ is 0 and $e$ is an edge incident with the $i$-th vertex in $\gamma_0$.

\noindent{\bf Proof of Claim 14:}
We may assume that $e$ is not an edge incident with the $i$-th vertex in the indicator of $\overline{G}$, for otherwise we are done.
Since $(\pi_V,\pi_E)$ defines a homeomorphic embedding from $\overline{G}$ to $\overline{G'}$ such that for every $j \in [\lvert V(\gamma_0) \rvert]$, $\pi_V$ maps the $j$-th entry in the indicator of $\overline{G}$ to the $j$-th entry in the indicator of $\overline{G'}$, the essential number associated with the $i$-th vertex in $\gamma_0$ is 0.

So it suffices to show that $e$ is an edge incident with the $i$-th vertex in $\gamma_0$.
Let $x$ and $y$ be the ends of $e$.

If $e \in E(\overline{H})$ and $\pi_E(e) = \eta(e)$, then since $\eta$ is a homeomorphic embedding from $(H,\gamma_0)$ to $(H',\gamma_0')$, the essential number associated with the $i$-th vertex in $\gamma_0$ is 0 and $e$ is an edge incident with the $i$-th vertex in $\gamma_0$, so we are done.
Hence we may assume that either $e \not \in E(\overline{H})$, or $e \in E(\overline{H})$ and $\pi_E(e) \neq \eta(e)$.
So either $\pi_V(x) \not \in V(\overline{H'})$, or $\pi_V(y) \not \in V(\overline{H'})$.

By symmetry, we may assume that $\pi_V(y) \not \in V(\overline{H'})$.
So $g(y)$ is defined.
Let $g(y)'$ be the child of $r'$ with $\iota(\gamma_{g(y)})=\gamma_{g(y)'}$.

We first assume $\pi_V(x) \in V(\overline{H'})$.
Since $\eta$ is a homeomorphic embedding from $(H,\gamma_0)$ to $(H',\gamma_0')$ and $\eta_{g(y)}$ is a homeomorphic embedding from $(G[\uparrow g(y)],\gamma_{g(y)})$ to $(G'[\uparrow g(y)'],\gamma_{g(y)'})$, we know that $v = \eta(u)$ for some $u \in V(H)$ and $e$ is incident with $u$.
Hence $v$ is the $i$-th vertex in $\gamma_0'$ with $v=\eta(u)$ for some $u \in V(H)$.
By the definition of a homeomorphic embedding of rooted graphs, $u$ is the $i$-th vertex in $\gamma_0$.
So $e$ is incident with the $i$-th vertex in $\gamma_0$ and we are done.

So we may assume $\pi_V(x) \not \in V(\overline{H'})$.
Hence $g(x)$ is defined.

We first assume $g(x)=g(y)$.
Since $\pi_V(x),\pi_V(y) \not \in V(\overline{H'})$, $v \in V(\pi_E(e)) = V(\eta_{g(y)}(e))$.
So $v \in V(\gamma_{g(y)'})$, and there exists $u \in V(\gamma_{g(y)})$ with $\eta(u)=v$ such that $e$ is incident with $u$.
Hence $v$ is the $i$-th vertex in $\gamma_0'$ with $v=\eta(u)$ for some $u \in V(H)$.
By the definition of a homeomorphic embedding of rooted graphs, $u$ is the $i$-th vertex in $\gamma_0$.
So $e$ is incident with the $i$-th vertex in $\gamma_0$ and we are done.

Hence we may assume $g(x) \neq g(y)$.
Since $\eta,\eta_{g(x)}$ and $\eta_{g(y)}$ are homeomorphic embeddings, $e$ is an edge incident with the $i$-th vertex in $\gamma_0$.
$\Box$

\medskip

\noindent{\bf Claim 15:} $(\pi_V,\pi_E)$ defines a homeomorphic embedding from $(G,\gamma_0)$ to $(G',\gamma_0')$. 

\noindent{\bf Proof of Claim 15:}
Since $(H',\gamma_0',\Gamma_H',f_H',\phi_H')$ simulates $(H,\gamma_0,\Gamma_H,f_H,\phi_H)$, we know that $\gamma_0$ and $\gamma_0'$ have the same length, and for every $i \in [|V(\gamma_0)|]$, the essential number associated with the $i$-th vertex in $\gamma_0$ equals the essential number associated with the $i$-th vertex in $\gamma_0'$.
By the definition of $(\pi_V,\pi_E)$ and Claim 13, $(\pi_V,\pi_E)$ gives a homeomorphic embedding from $\overline{G}$ to $\overline{G'}$ such that for every $i \in [|V(\gamma_0)|]$, $\pi_V$ maps the $i$-th entry in the indicator of $\overline{G}$ to the $i$-th entry in the indicator of $\overline{G'}$.
Moreover, if $v'$ is the $i$-th vertex in $\gamma_0'$ for some $i \in [\lvert V(\gamma_0') \rvert]$ with $\pi_V(v)=v'$ for some vertex $v$ of $\overline{G}$, then $v'=\pi_V(v)=\eta(v)$, so $v$ is the $i$-th vertex in $\gamma_0$ since $\eta: (H,\gamma_0) \hookrightarrow (H',\gamma_0')$.
The above facts together with Claim 14 imply this claim. 
$\Box$

\medskip

It is clear that $\phi(v) \leq_{Q} \phi'(\pi_V(v))$ for every $v \in V(G)$.
So to prove that $(G',\gamma_0',\Gamma',f',\phi')$ simulates $(G,\gamma_0,\Gamma,f,\phi)$, it suffices to define an injection $\iota^*: \Gamma \rightarrow \Gamma'$ such that for every $\sigma \in \Gamma$, $\pi_V(\sigma)=\iota^*(\sigma)$ and $f(\sigma) \leq_{Q} f'(\iota^*(\sigma))$.

Define $\iota^*: \Gamma \rightarrow \Gamma'$ such that for every $\sigma \in \Gamma$,
	\begin{itemize}
		\item if $\alpha(\sigma)=r$, then $\iota^*(\sigma)=\iota(\sigma)$;
		\item if $\alpha(\sigma)$ is a descendant of $c$ for some child $c$ of $r$, then $\iota^*(\sigma)=\iota_c(\sigma)$. (Note that such child $c$ is unique, so $\iota^*(\sigma)$ is well-defined.)
	\end{itemize}

\noindent{\bf Claim 16:} $\iota^*$ is an injection.

\noindent{\bf Proof of Claim 16:}
Note that for distinct children $c$ and $c'$ of $r$, we know that $\iota$, $\iota_c$ and $\iota_{c'}$ are injections with disjoint images.
So $\iota^*$ is an injection.
$\Box$

\medskip

\noindent{\bf Claim 17:} For every $\sigma \in \Gamma$, $\pi_V(\sigma)=\iota^*(\sigma)$.

\noindent{\bf Proof of Claim 17:}
Suppose to the contrary that there exists $\sigma \in \Gamma$ such that $\pi_V(\sigma) \neq \iota^*(\sigma)$.

Suppose $\alpha(\sigma)=r$.
Then $V(\sigma) \subseteq X_r$.
So for each vertex $v$ in $V(\sigma)$, either $v \not \in V(\gamma_c)$ for every child $c$ of $r$, or $v$ is associated with essential number 2 in $\gamma_c$ for every child $c$ of $r$ with $v \in V(\gamma_c)$ by (BR01), (BR11) and (BR2) since $r$ is a non-descendant of every child of $r$.
So for every $v \in V(\sigma)$, $\pi_V(v)=\eta(v)$. 
Hence $\pi_V(\sigma)=\eta(\sigma)=\iota(\sigma)=\iota^*(\sigma)$, a contradiction.

So there exists a child $c^*$ of $r$ such that $\alpha(\sigma)$ is a descendant of $c^*$.
If $\pi_V(\sigma)=\eta_{c^*}(\sigma)$, then $\pi_V(\sigma)=\eta_{c^*}(\sigma)=\iota_{c^*}(\sigma)=\iota^*(\sigma)$, a contradiction.
Hence there exists a vertex $v \in V(\sigma)$ such that $\pi_V(v) \neq \eta_{c^*}(v)$.
So $v \in V(\gamma_{c^*})$, $L_v \neq \emptyset$ and $g(v) \neq c^*$.
Hence, $v$ is in $V(\gamma_{g(v)})$ and is associated with essential number 0 or 1 in $\gamma_{g(v)}$. 
However, since $\alpha(\sigma)$ is a descendant of $c^* \neq g(v)$, and $c^*$ is a non-descendant of $g(v)$, $\alpha(\sigma)$ is a non-descendant of $g(v)$.
So $v$ is associated with essential number 2 in $\gamma_{g(v)}$, a contradiction.
This proves $\pi_V(\sigma)=\iota^*(\sigma)$.
$\Box$

\medskip

Since $\pi_V(\sigma)=\iota^*(\sigma)$ for all $\sigma \in \Gamma$, it is straightforward to verify that $f(\sigma) \leq_{Q} f'(\iota^*(\sigma))$ by the properties of $f_H, \eta$, $f_c$ and $\eta_c$ for all children $c$ of $r$.
This proves that $(G',\gamma_0',\Gamma',f',\phi')$ simulates $(G,\gamma_0,\Gamma,f,\phi)$.
\end{pf}

\section{Well-behaved assemblages} \label{sec: well-behaved}

In this section, we develop the last piece of our machinery and prove the main theorem of this paper (Theorem \ref{robertson conj bounded tree width}).
We shall demonstrate how to formally use the linkedness property and the absorption property of a tree-decomposition mentioned in earlier sections to prove well-quasi-ordering results.

In Section \ref{subsec:warm_up_well_behaved}, we define ``well-behaved'' sets of assemblages, which can be considered sets of assemblages that make the absorption property hold.
In Section \ref{subsec:using_linking}, we show how to use the linkedness property of a tree-decomposition to reduce a well-quasi-ordering problem on assemblages to the one on their encodings.
We introduce, in Sections \ref{subsec:node_realizers} and \ref{subsec:unimpeded}, the form of linkedness property that we will use in Section \ref{subsec:using_linking}.
Finally, we apply all machinery developed in this paper to prove Theorem \ref{robertson conj bounded tree width} in Section \ref{subsec:assemblages_to_graphs}.

\subsection{Warm up} \label{subsec:warm_up_well_behaved}
We say that a set $\F$ of assemblages is \defn{well-behaved} if for every well-quasi-order $Q$, for every infinite sequence of $Q$-assemblages $S_1,S_2,...$ each having underlying assemblage in $\F$, there exist $1 \leq i < i'$ such that $S_{i'}$ simulates $S_i$.

\begin{lemma} \label{well-behaved bounded size}
Let $n$ be a nonnegative integer.
If $\F^n$ is the set of assemblages $(G,\gamma_0,\Gamma)$ with $\lvert V(G) \rvert \leq n$, then $\F^n$ is well-behaved.
\end{lemma}

\begin{pf}
Let $Q$ be a well-quasi-order.
Let $((G_1,\gamma_1,\Gamma_1,f_1,\phi_1), (G_2, \gamma_2,\Gamma_2,f_2,\phi_2),...)$ be an infinite sequence of $Q$-assemblages whose underlying assemblages belong to $\F^n$.
Since $n$ is finite, we may assume that $G_1,G_2,...$ have the same number of vertices and may assume that $V(\overline{G_1})=V(\overline{G_2})= \cdots$, where each $\overline{G_i}$ is the rooted extension of $G_i$.
Since we can record the number of loops incident with the vertices of $G_i$ and the number of parallel edges between two vertices of $G_i$ as a sequence over nonnegative integers with at most $|V(G_1)|^2$ entries, Higman's lemma implies that we may further assume that for any $i<j$, and for every pair of distinct elements $v,v'$ of $V(G_1)$, the number of loops incident with $v$ in $G_i$ is at most the number of loops incident with $v$ in $G_j$, and the number of edges in $G_i$ between $v$ and $v'$ is at most the number of edges in $G_j$ between $v$ and $v'$.
Therefore, for $i<j$, $G_i$ is a subgraph of $G_j$.
Since there are only finitely many distinct marches on $V(G_1)$, we may assume that $\gamma_1=\gamma_2= \cdots$.
Furthermore, since $V(G_1)$ is finite, we may assume that for every $v \in V(G_1)$, $\phi_1(v) \leq_Q \phi_2(v) \leq_Q \cdots$.
For each $i \geq 1$ and for each march $\sigma$ with $V(\sigma) \subseteq V(G_1)$, let $\sigma_1,\sigma_2,...$ be the elements of the multiset $\Gamma_i$ with $\sigma_j=\sigma$ for all $j$, and we define $f'_i(\sigma)=(f_i(\sigma_1),f_i(\sigma_2),...)$.
Note that $\Gamma_i$ is finite, so $f_i'(\sigma)$ is a finite sequence.
Hence $f_i'$ is a function from the set of marches on $V(G_1)$ to the well-quasi-order set, denoted by $Q'$, obtained from $Q$ by Higman's Lemma.
Since the domain of each $f_i'$ is finite, there exist $1 \leq j < j'$ such that $f_j'(\sigma) \leq_{Q'} f_{j'}'(\sigma)$ for all marches $\sigma$ with $V(\sigma) \subseteq V(G_1)$.
So there exists an injection $\iota: \Gamma_j \rightarrow \Gamma_{j'}$ such that $\iota(\gamma)$ is the same march as $\gamma$, and $f_j(\gamma) \leq_Q f_{j'}(\iota(\gamma))$ for all $\gamma \in \Gamma_j$.
Define $\eta$ to be the identity homeomorphic embedding from $(G_j,\gamma_j)$ to $(G_{j'},\gamma_{j'})$.
Then $(G_{j'},\gamma_{j'},\Gamma_{j'},f_{j'},\phi_{j'})$ simulates $(G_j.\gamma_j,\Gamma_j,f_j,\phi_j)$ witnessed by $\eta,\iota$.
\end{pf}

\subsection{Node-realizers} \label{subsec:node_realizers}

Let $(T,\X,\alpha)$ be a rooted tree-decomposition of an assemblage $(G,\gamma_0,\Gamma)$.
Denote $\X$ by $(X_t: t \in V(T))$.
The \defn{node-realizer} of $(T,\X,\alpha)$ is the rooted tree-decomposition $(T',\X',\alpha)$ obtained from $(T,\X,\alpha)$ by subdividing each edge $xy$ of $T$ once, defining the bag of the corresponding new node to be $X_x \cap X_y$, and adding a new vertex, which is the root of $T'$, adjacent to the root of $T$ and defining its bag to be $V(\gamma_0)$.
Denote $\X'$ by $(X'_t: t \in V(T'))$.
For each edge $e$ of $T$, $e$ is a node of $T'$, and we let $t_e$ be the head of $e$ in $T$, and let $\gamma_{t_e}$ be the root march of the underlying assemblage of the $(f,\phi)$-branch of $(T,\X,\alpha)$ at $t_e$ (with respect to an arbitrary ordering of $X'_{e}$).
Note that the definition of $\gamma_{t_e}$ is independent with $f$ and $\phi$.
For each edge $e$ of $T$ and each $Z \subseteq V(G)$, we say that the node $e$ of $T'$ corresponds to a \defn{$\Gamma$-pseudo-edge-cut modulo $Z$} if every vertex $v$ in $V(\gamma_{t_e})-Z$ is associated with essential number 0 or 1 in $\gamma_{t_e}$. 
The \defn{$\Gamma$-elevation} of the node-realizer $(T',\X',\alpha)$ is the elevation of $(T',\X')$ but pseudo-edge-cuts in the definition for elevation, $(Z,s)$-strips and $(Z,s)$-depth are replaced by $\Gamma$-pseudo-edge-cuts; namely, the $\Gamma$-elevation of $(T',\X',\alpha)$ is the maximum positive integer $h$ such that there exist $Z \subseteq V(G)$, a positive integer $s$ and a sequence $(t_1,t_2,...,t_h)$ of nodes of $T'$ such that 
	\begin{itemize}
		\item $t_i$ is a precursor of $t_{i+1}$ for every $i \in [h-1]$,
		\item there exists a directed path in $T$ passing through $t_1,t_2,...,t_h$ in the order listed,
		\item $Z \subseteq \bigcap_{i=1}^h X'_{t_i}$,
		\item $X'_{t_1}-Z,X'_{t_2}-Z,...,X'_{t_h}-Z$ are pairwise disjoint non-empty sets with size $s$,
		\item there exist $|X_{t_1}'|$ disjoint paths in $G$ from $X_{t_1}'$ to $X_{t_h}'$, and
		\item there exists no node $t$ in $t_1T't_h$ such that $|X_t'|=|X_{t_1}'|$ and the separation given by $t$ in $(T',\X')$ is a $\Gamma$-pseudo-edge-cut modulo $Z$.
	\end{itemize}

\subsection{Unimpeded tree-decomposition} \label{subsec:unimpeded}

Let $(T,\X,\alpha)$ be a rooted tree-decomposition of an assemblage $(G,\gamma_0,\Gamma)$.
Denote $\X$ by $(X_t: t \in V(T))$.
Let $N$ be a positive integer.
We say that $(T,\X,\alpha)$ is \defn{$N$-unimpeded} if there exist $\lvert Z_1 \rvert$ disjoint paths in $G$ from $Z_1$ to $Z_2$ whenever $Z_1,Z_2,...,Z_{N+1}$ are pairwise distinct sets with the same size such that 
	\begin{itemize}
		\item for each $i \in [N+1]-[1]$, there exists an edge $s_it_i$ of $T$ such that $Z_i = X_{s_i} \cap X_{t_i}$,
		\item $s_2t_2,s_3t_3,...,s_{N+1}t_{N+1}$ are distinct edges of $T$ such that some directed path in $T$ passing through them in the order listed,
		\item either $Z_1 = V(\gamma_0)$, or $Z_1 = X_{s_1} \cap X_{t_1}$ for some edge $s_1t_1$ of $T$ with $s_1t_1 \neq s_2t_2$ such that some directed path in $T$ passing through $s_1t_1,s_2t_2,...,s_{N+1}t_{N+1}$ in the order listed, and
		\item $\lvert X_x \cap X_y \rvert \geq \lvert Z_1 \rvert$ for all edges $xy$ of $T$ contained in the path $s_1Tt_{N+1}$, where $s_1$ is the root of $T$ if $Z_1=V(\gamma_0)$, and $s_1$ is the tail of $s_1t_1$ otherwise.
	\end{itemize}

\begin{lemma} \label{unimpeded_weakly_linked}
Let $N$ be a positive integer.
Let $(T,\X,\alpha)$ be a rooted tree-decomposition of an assemblage $(G,\gamma_0,\Gamma)$.
Let $(R,\Y,\alpha)$ be the node-realizer of $(T,\X,\alpha)$.
If $(T,\X,\alpha)$ is $N$-unimpeded, then $(R,\Y)$ is weakly $N$-linked.
\end{lemma}

\begin{pf}
Denote $\X$ by $(X_t: t \in V(T))$ and denote $\Y$ by $(Y_t: t \in V(R))$.

Suppose to the contrary that $(R,\Y)$ is not weakly $N$-linked.
Then there exist nodes $t_1,t_2,...,t_{N+1}$ of $R$ such that $t_i$ is a precursor of $t_{i+1}$ for every $i \in [N]$ and the sets $Y_{t_1}, Y_{t_2},...., \allowbreak Y_{t_{N+1}}$ are distinct, but there exist no $|Y_{t_1}|$ disjoint sets in $G$ from $Y_{t_1}$ to $Y_{t_2}$.
We choose the nodes $t_1,t_2,...,t_{N+1}$ such that $\lvert \{t_i \in V(R)-V(T): i \in [N+1]\} \rvert$ is as large as possible. 
Let $P$ be the directed path in $R$ from $t_1$ to $t_{N+1}$.

Since $(R,\Y,\alpha)$ is the node-realizer of $(T,\X,\alpha)$, for every edge $xy$ of $R$, either $Y_x \subseteq Y_y$ or $Y_y \subseteq Y_x$.
So if there exists $i \in [N]$ such that $t_i$ is the parent of $t_{i+1}$, then $Y_{t_i}$ and $Y_{t_{i+1}}$ are two distinct sets with the same size such that one of them is a subset of the other, a contradiction.
Hence for every $i \in [N]$, $t_i$ is not the parent of $t_{i+1}$.

Suppose to the contrary that $t_1 \in V(T)$.
Let $c$ be the child of $t_1$ in $R$ contained in $P$.
Since $t_1 \in V(T)$, $c$ is obtained by subdividing an edge of $T$, so $Y_c \subseteq Y_{t_1}$.
Since $t_1$ is a precursor of $t_2$, $|Y_c| \geq |Y_{t_1}|$.
But $Y_c \subseteq Y_{t_1}$, so $Y_c=Y_{t_1}$.
Since $t_1$ is not the parent of $t_2$, $c \neq t_2$.
Then $c,t_2,t_3,...,t_{N+1}$ is a better choice than $t_1,t_2,...,t_{N+1}$, a contradiction.

So $t_1 \in V(R)-V(T)$.
Hence either $t_1$ is the root of $R$, or $t_1 \in E(T)$.
For the former, we know $Y_{t_1}=V(\gamma_0)$ and we let $s_1=t_1$ and $Z_1 = Y_{t_1}=V(\gamma_0)$; for the latter, let $s_1$ be the parent of $t_1$ and let $q_1$ be the node of $T$ such that $t_1=s_1q_1$.
For every $i \in [N+1]-[1]$, let $s_i$ be the parent of $t_i$ in $R$.
Note that $\{s_i: i \in [N+1]-[1]\} \cap \{t_i: i \in [N+1]\} = \emptyset$.

For every $i \in [N+1]-[1]$, if $t_i \in V(T)$, then $s_i \in E(T)$, and since $t_1$ is a precursor of $t_i$, we know $Y_{s_i} = Y_{t_i}$, so replacing $t_i$ by $s_i$ gives a better choice than $t_1,...,t_{N+1}$, a contradiction.

Hence for every $i \in [N+1]-[1]$, $t_i \in E(T)$, so there exists a node $q_i$ of $T$ such that $t_i = s_iq_i$ and $X_{s_i} \cap X_{q_i} = Y_{t_i}$ for some node $q_i$ of $T$.
So $s_2q_2,...,s_{N+1}q_{N+1}$ are distinct edges of $T$ passed through by some directed path in $T$ in the order listed. 

Moreover, for every edge $xy$ of $T$ contained in the directed path $s_1Tq_{N+1}$, we know that $xy$ is a node of $R$ contained in $P$, so $|X_x \cap X_y| = |Y_{xy}| \geq |Y_{t_1}|$.
Since $(T,\X,\alpha)$ is $N$-unimpeded, there exist $|Y_{t_1}|$ disjoint paths in $G$ from $Y_{t_1}$ to $Y_{t_2}$, a contradiction.
\end{pf}

\subsection{Using tree-decomposition} \label{subsec:using_linking}

Let $\F$ be a family of assemblages.
We say that a \defn{rooted tree-decomposition $(T,\X,\alpha)$ of an assemblage is over $\F$} if for every $t \in V(T)$ and for every ordering $\pi$ of the vertices in the bags, the underlying assemblage of the encoding of $(T,\X,\alpha)$ at $t$ (with respect to $\pi)$ belongs to $\F$.

The \defn{adhesion} of a rooted tree-decomposition $(T,\X,\alpha)$ of an assemblage is the adhesion of $(T,\X)$.

\begin{theorem} \label{well-behaved tree-decomposition}
Let $\F$ be a well-behaved family of assemblages.
Let $h,d,N$ be positive integers.
Let $\F_{h,d,N}$ be the family consisting of all assemblages $(G,\gamma_0,\Gamma)$ that has an $N$-unimpeded rooted tree-decomposition over $\F$ of adhesion at most $h$ such that the $\Gamma$-elevation of its node-realizer is at most $d$.
Then $\F_{h,d,N}$ is well-behaved.
\end{theorem}

\begin{pf}
Let $Q$ be a well-quasi-order.
For each positive integer $i$, let $(G_i,\gamma_i,\Gamma_i,f_i,\phi_i)$ be a $Q$-assemblage with $(G_i,\gamma_i,\Gamma_i) \in \F_{h,d,N}$.
By the definition of $\F_{h,d,N}$, for every $i \geq 1$, there exists an $N$-unimpeded rooted tree-decomposition $(T^i,\X^i,\alpha^i)$ of $(G_i,\gamma_i,\Gamma_i)$ over $\F$ of adhesion at most $h$ such that the node-realizer of $(T^i,\X^i,\alpha^i)$, denoted by $(R^i,\Y^i,\alpha^i)$, has $\Gamma$-elevation at most $d$.
To prove this theorem, it suffices to prove that there exist $i' > i \geq 1$ such that $(G_{i'},\gamma_{i'},\Gamma_{i'},f_{i'},\phi_{i'})$ simulates $(G_i,\gamma_i,\Gamma_i,f_i,\phi_i)$.

For every $i \geq 1$, let $\overline{G_i}$ be the rooted extension of $(G_i,\gamma_i)$, and we denote $\X^i$ by $(X^i_t: t \in V(T^i))$ and denote $\Y^i$ by $(Y^i_t: t \in V(R^i))$.

We call a node $t$ of $R^i$ a \defn{chopper} in $(R^i,\Y^i)$ if either $t$ has no precursor, or there do not exist $\lvert Y^i_t \rvert$ disjoint paths in $G$ from $Y^i_{t'}$ to $Y^i_t$, where $t'$ is the precursor of $t$ closest to $t$.

We define the \defn{level} of each node $t$ of $R^i$, denote by $\mu_i(t)$, recursively according to the breadth-first-search order of $R^i$ as follows.
	\begin{itemize}
		\item If $t$ has no precursor, then define $\mu_i(t)=0$.
		\item If $t$ has a precursor, then let $t'$ be the precursor of $t$ closest to $t$, and define the level of $t$ as follows.
			\begin{itemize}
				\item If $t$ is a chopper, then define $\mu_i(t)=\mu_i(t')+1$.
				\item If $t$ is not a chopper, then some precursor of $t$ is a chopper, and let $t''$ be such a chopper closest to $t$ and define
					\begin{itemize}
						\item $\mu_i(t)=\mu_i(t')$, if $Y^i_t \cap Y^i_{t''}=Y^i_{t'} \cap Y^i_{t''}$, and for every vertex $v \in Y^i_t \cap Y^i_{t''}$, 
							\begin{itemize}
								\item if there exists $\sigma \in \Gamma_i$ with $v \in V(\sigma)$ such that $\alpha^i(\sigma)$ is a non-descendant of $t$, then there exists $\sigma' \in \Gamma_i$ with $v \in V(\sigma')$ such that $\alpha^i(\sigma')$ is a non-descendant of $t'$, and 
								\item the number of edges incident with $v$ whose other ends are in $V(\overline{G_i})-(R^i,\Y^i)\uparrow t$ and the number of edges incident with $v$ whose other ends are in $V(\overline{G_i})-(R^i,\Y^i)\uparrow t'$ are either both at least two or both equal to a number $j$ with $j \in \{0,1\}$;
							\end{itemize}
						\item $\mu_i(t)=\mu_i(t')+1$, otherwise.
					\end{itemize}
			\end{itemize}
	\end{itemize} 

Let $N'=(3^h(h+1)^2+2)(N+1)$.

\medskip

\noindent{\bf Claim 1:} The level of every node is at most $N'$.

\noindent{\bf Proof of Claim 1:}
Suppose to the contrary that some node of $R^i$ has level at least $N'+1$.
So there exist nodes $t_0,t_1,...,t_{N'},t_{N'+1}$ of $R^i$ such that $\mu_i(t_0)=0$, and for every $j \in [0,N']$, $t_j$ is a precursor of $t_{j+1}$ with $\mu_i(t_{j+1})=\mu_i(t_j)+1$ such that $\mu_i(t)=\mu_i(t_j)$ for every precursor $t$ of $t_{j+1}$ contained in the path between $t_j$ and $t_{j+1}$.

Suppose that there exist $0 \leq j_1<j_2<...<j_{N+1} \leq N'+1$ such that $t_{j_\ell}$ is a chopper for each $\ell \in [N+1]$.
So for each $\ell \in [N]$, there exist no $\lvert Y^i_{j_1} \rvert$ disjoint paths in $G$ from $Y^i_{j_\ell}$ to $Y^i_{j_{\ell+1}}$.
In particular, $Y^i_{j_1}, ...,Y^i_{j_{N+1}}$ are pairwise distinct sets with the same size.
Since $(T^i,\X^i,\alpha^i)$ is $N$-unimpeded, $(R^i,\Y^i)$ is weakly $N$-linked by Lemma \ref{unimpeded_weakly_linked}, so there exist $\lvert Y^i_{j_1} \rvert$ disjoint paths in $G$ from $Y^i_{j_1}$ to $Y^i_{j_2}$, a contradiction.

Hence at most $N$ nodes in $\{t_j: 0 \leq j \leq N'+1\}$ are choppers.
Let $c=3^h(h+1)^2$.
Since $N' \geq (c+1)(N+1)+N$, there exists $a \in [N'+1-c]$ such that $t_a,t_{a+1},...,t_{a+c}$ are not choppers.
Let $t''$ be the chopper that is the precursor of $t_a$ closest to $t_a$.
Since $\mu_i(t)=\mu_i(t_j)=j \geq 1$ for each $j \in [N']$ and for each precursor $t$ of $t_{j+1}$ contained in the path between $t_j$ and $t_{j+1}$, there exist $\lvert Y^i_{t_1} \rvert$ disjoint paths in $G$ from $Y^i_{t_a}$ to $Y^i_{t_{a+c}}$.
So $t''$ is the chopper that is the precursor of $t_{a+j}$ closest to $t_{a+j}$ for each $j \in [0,c]$.
Since $(R^i,\Y^i)$ is a tree-decomposition of adhesion at most $h$, there are at most $h+1$ different possibilities for $Y^i_{t_{a+j}} \cap Y^i_{t''}$ for $0 \leq j \leq c$.
And for each $v \in Y^i_{t''}$ and $0 \leq j \leq c$, either $v \not \in Y^i_{t_{a+j}}$, or $v \in Y^i_{t_{a+j}} \cap Y^i_{t''}$ and $v$ is incident with at least two edges whose other ends are in $V(\overline{G_i})-(R^i,\Y^i)\uparrow t_{a+j}$, or $v \in Y^i_{t_{a+j}} \cap Y^i_{t''}$ and $v$ is incident with exactly $\ell$ edge whose other end is in $V(\overline{G_i})-(R^i,\Y^i)\uparrow t_{a+j}$ for some $\ell \in \{0,1\}$.
Since $c+1>3^h(h+1)^2$, there exist $a \leq q_1<q_2<...<q_{h+2} \leq a+c$ such that for every $1 \leq j_1<j_2 \leq h+2$, $Y^i_{t_{q_{j_1}}} \cap Y^i_{t''}=Y^i_{t_{q_{j_2}}} \cap Y^i_{t''}$ and for every vertex $v \in Y^i_{t_{q_{j_1}}} \cap Y^i_{t''}$, the number of edges incident with $v$ whose other ends are in $V(\overline{G_i})-(R^i,\Y^i)\uparrow t_{q_{j_1}}$ and the number of edges incident with $v$ whose other ends are in $V(\overline{G_i})-(R^i,\Y^i)\uparrow t_{q_{j_2}}$ are either both at least two or both equal to some number $\ell$ with $\ell \in \{0,1\}$.

By the definition of $\mu_i$, for every $\ell \in [h+1]$, since $\mu_i(t_{q_{\ell}}) < \mu_i(t_{q_{{\ell+1}}})$, there exist $v_\ell \in Y^i_{t_{q_{1}}} \cap Y^i_{t''}$ and $\sigma_\ell \in \Gamma_i$ such that $v_\ell \in V(\sigma_\ell)$ and $\alpha^i(\sigma_\ell)$ is a non-descendant of $t_{q_{{\ell+1}}}$ but there exists no $\sigma' \in \Gamma_i$ with $v \in V(\sigma')$ such that $\alpha^i(\sigma')$ is a non-descendant of $t_{q_{\ell}}$.
Note that $v_1,v_2,...,v_{h+1}$ are distinct vertices in $Y^i_{q_{1}} \cap Y^i_{t''}$ which has size at most $h$, a contradiction.
$\Box$

\medskip

For every $i \geq 1$, if a node $t$ of $R^i$ has a precursor, then since $(R^i,\Y^i)$ is a node-realizer, $Y^i_t$ is a subset of the bag of its parent. 
For each node $t$ of $R^i$, we define $\pi'_t$ to be an ordering of the vertices in $Y^i_t$ such that the following hold.
	\begin{itemize}
		\item If $t$ is the root of $R^i$, then $\pi'_t$ is the ordering same as $\gamma_i$.
		\item If $t$ is a chopper in $(R^i,\Y^i)$ but not the root of $R^i$, then $\pi_t'$ is an arbitrary ordering.
		\item If $t$ is not a chopper in $(R^i,\Y^i)$, then there exist $\lvert Y^i_t \rvert$ disjoint paths $P_1,P_2,...,P_{\lvert Y^i_t \rvert}$ in $G_i$ from $Y^i_t$ to $Y^i_{t'}$, where $t'$ is the precursor of $t$ closest to $t$, such that for each $j \in [\lvert Y^i_t \rvert]$, the ends of $P_j$ are the $j$-th vertices in $Y^i_t$ and $Y^i_{t'}$ with respect to $\pi'_t$ and $\pi'_{t'}$, respectively.
	\end{itemize}
Then for each non-root node $t$ of $T^i$ with parent $p$, we define $\pi_t=\pi'_{pt}$, where $pt$ is the node of $R^i$ obtained by subdividing the edge $pt$ of $T^i$, and we define $\gamma_t$ to be the root march of the underlying assemblage of the $(f_i,\phi_i)$-branch of $(T^i,\X^i,\alpha^i)$ at $t$ with respect to $\pi_t$.  

For each $i \geq 1$, define $\psi_i$ and $\tau_i$ to be the functions from $E(T^i)$ to $[V(G_i)]^{\leq h}$ such that for all edges $pt$ of $T^i$, where $t$ is a child of $p$, the following hold.
	\begin{itemize}
		\item $\psi_i(pt)=X^i_p \cap X^i_t$.
		\item If $t$ is a chopper in $(R^i,\Y^i)$, then $\tau_i(pt)=\psi_i(pt)$; otherwise, $\tau_i(pt)$ is defined to be the set of all vertices $v \in V(\gamma_t)$ associated with essential number 2 in $\gamma_t$. 
	\end{itemize}
Note that when $t$ is not a chopper, $\tau_i(pt)$ contains all vertices $v$ in $V(\gamma_t)=\psi_i(pt)$ satisfying that either $v$ is incident with at least two edges whose other ends are in $V(\overline{G_i})-(R^i,\Y^i)\uparrow pt$, or $v \in V(\sigma)$ for some $\sigma \in \Gamma_i$ in which $\alpha^i(\sigma)$ is a non-descendant of $t$ in $R^i$.
In addition, $\mu_i$ is defined on the nodes of $R^i$, so its domain contains $E(T^i)$.

Recall the definition of a node preceding another node stated in Section \ref{subsec:decorated_trees}.

\medskip

\noindent{\bf Claim 2:} For each $i \geq 1$, if a node $v$ of $T^i$ precedes another node $w$ of $T^i$ with respect to $(\psi_i,\tau_i,\mu_i|_{E(T^i)})$, then the following statements hold.
	\begin{itemize}
		\item There exist $\lvert V(\gamma_v) \rvert$ disjoint paths $P_1,P_2,...,P_{\lvert V(\gamma_v) \rvert}$ in $G_i$ from $V(\gamma_v) \subseteq X^i_v$ to $V(\gamma_w) \subseteq X^i_w$ such that for each $\ell \in [\lvert V(\gamma_v) \rvert]$, the ends of $P_\ell$ are the $\ell$-th vertices in $\gamma_v$ and $\gamma_w$.
		\item If $e_v,e_w$ are the edges of $T^i$ with heads $v,w$, respectively, then for every $j \in [\lvert V(\gamma_v) \rvert]$, either 
			\begin{itemize}
				\item the $j$-th entry of $\gamma_v$ is associated with essential number 2 in $\gamma_v$ and the $j$-th entry of $\gamma_w$ is associated with essential number 2 in $\gamma_w$, or 
				\item there exists $\ell \in \{0,1\}$ such that the $j$-th entry of $\gamma_v$ is associated with essential number $\ell$ in $\gamma_v$, the $j$-th entry of $\gamma_w$ is associated with essential number $\ell$ in $\gamma_w$, and $\ell_v=\ell_w=\ell$, where $\ell_v$ is the number of edges incident with the $j$-th entry of $\gamma_v$ whose other ends are in $V(\overline{G_i})-(R^i,\Y^i)\uparrow e_v$, and $\ell_w$ is the number of edges incident with the $j$-th entry of $\gamma_w$ whose other ends are in $V(\overline{G_i})-(R^i,\Y^i)\uparrow e_w$.
			\end{itemize}
	\end{itemize}

\noindent{\bf Proof of Claim 2:}
Let $e_v,e_w$ be the edges of $T^i$ with heads $v,w$, respectively.
Since $v$ precedes $w$ with respect to $(\psi_i,\tau_i,\mu_i|_{E(T^i)})$, $\lvert \psi_i(e) \rvert \geq \lvert \psi_i(e_v) \rvert = \lvert \psi_i(e_w) \rvert=\lvert V(\gamma_v) \rvert$ for every edge $e$ in $vT^iw$.
So $e_v$ is a precursor of $e_w$ in $(R^i,Y^i)$.
Let $t_1,t_2,...,t_c$ (for some integer $c$) be the nodes in $e_vR^ie_w$ with $\lvert Y^i_{t_j} \rvert = \lvert \psi_i(e_v) \rvert$ for every $j \in [c]$ such that they appear in $e_vR^ie_w$ in the order listed.
So $t_1=e_v$, $t_c=e_w$, and for each $j \in [c-1]$, $t_j$ is the closest precursor of $t_{j+1}$.
Hence $\mu_i(t_1)\leq \mu_i(t_2) \leq ... \leq \mu_i(t_c)$ by the definition of $\mu_i$.
Since $v$ precedes $w$, $\mu_i(t_1)=\mu_i(e_v)=\mu_i(e_w)=\mu_i(t_c)$, so $\mu_i(t_j)=\mu_i(t_1)$ for all $j \in [c]$.
By the definition $\mu_i$ and $\pi'_{t_j}$, we know that for each $j \in [c-1]$, there exist $\lvert V(\gamma_v) \rvert$ disjoint paths $P^j_1,P^j_2,...,P^j_{\lvert V(\gamma_v) \rvert}$ in $G_i$ from $Y^i_{t_j}$ to $Y^i_{t_{j+1}}$ such that for each $\ell \in [\lvert V(\gamma_v) \rvert]$, the ends of $P^j_\ell$ are the $\ell$-th vertices of $\gamma_{t_j}$ and $\gamma_{t_{j+1}}$.
By concatenating those paths, we obtain paths $P_1,...,P_{\lvert V(\gamma_v) \rvert}$ in $G_i$ from $Y^i_{t_1}$ to $Y^i_{t_c}$ such that for each $\ell \in [\lvert V(\gamma_v) \rvert]$, the ends of $P_\ell$ are the $\ell$-th vertices in $\gamma_v$ and $\gamma_w$.
This proves the first statement of the claim since $Y^i_{t_1}=V(\gamma_v) \subseteq X_v^i$ and $Y^i_{t_c}=V(\gamma_w) \subseteq X_w^i$.

Now we prove Statement 2 of this claim.
Suppose that $x$ and $y$ are the $j$-th vertices of $\gamma_v$ and $\gamma_w$, respectively, for some $j \in [\lvert V(\gamma_v) \rvert]$, such that $x$ and $y$ violate Statement 2 of this claim. 
Let $\ell_v$ be the number of edges incident with $x$ whose other ends are in $V(\overline{G_i})-(R^i,\Y^i)\uparrow e_v$.
Let $\ell_w$ be the number of edges incident with $y$ whose other ends are in $V(\overline{G_i})-(R^i,\Y^i)\uparrow e_w$.

We first suppose that $\tau_i(e_v) \neq \psi_i(e_v)$.
Since $v$ precedes $w$ in $T^i$ with respect to $(\psi_i,\tau_i, \allowbreak \mu_i|_{E(T^i)})$, we know $\tau_i(e_v)=\tau_i(e_w)$ and $\lvert \psi_i(e_v) \rvert = \lvert \psi_i(e_w) \rvert$.
So $\tau_i(e_w) \neq \psi_i(e_w)$.
In particular, $e_v$ and $e_w$ are not choppers in $(R^i,\Y^i)$.
Since $e_v$ and $e_w$ are not choppers in $(R^i,\Y^i)$, there exist $e_v'$ and $e_w'$ such that $e_v'$ is the precursor that is the chopper of $e_v$ closest to $e_v$, and $e_w'$ is the precursor of $e_w$ that is the chopper closest to $e_w$ in $(R^i,\Y^i)$. 
Since $\mu_i(e_v)=\mu_i(e_w)$, by Statement 1 of this claim, no internal node of the path $e_vR^ie_w$ with bag size $\lvert \psi_i(e_v) \rvert$ is a chopper in $(R^i,\Y^i)$.
So $e_v'=e_w'$.
Since $\mu_i(e_v)=\mu_i(e_w)$ and $e_v$ is a precursor of $e_w$, $Y^i_{e_v} \cap Y^i_{e_v'} = Y^i_{e_w} \cap Y^i_{e_v'}$.
Since $\tau_i(e_v)=\tau_i(e_w) \subseteq \psi_i(e_v) \cap \psi_i(e_w)$, either $x \in \tau_i(e_v)$ and $y \in \tau_i(e_w)$, or $x \not \in \tau_i(e_v)$ and $y \not \in \tau_i(e_w)$.
Since $e_v$ and $e_w$ are not choppers in $(R^i,\Y^i)$, if $x \in \tau_i(e_v)$ and $y \in \tau_i(e_w)$, then $x$ is associated with essential number 2 in $\gamma_v$ and $y$ is associated with essential number 2 in $\gamma_w$, a contradiction.
So $x \not \in \tau_i(e_v)$ and $y \not \in \tau_i(e_w)$.
Hence $\ell_v \leq 1$ and $\ell_w \leq 1$.
Since the essential number associated with vertices in $\gamma_i$ are represented by edges in $\overline{G_i}$, we know that $\ell_v \leq 1$ and $\ell_w \leq 1$ implies that $x$ is associated with essential number $\ell_v$ in $\gamma_v$, and $y$ is associated with essential number $\ell_w$ in $\gamma_w$.
If $x \in Y^i_{e_v} \cap Y^i_{e_v'}$, then since $Y^i_{e_v} \cap Y^i_{e_v'} = Y^i_{e_w} \cap Y^i_{e_v'}$, we have $x=y$, and by the definition of $\mu_i$, we have $\ell_v=\ell_w=\ell$ for some $\ell \in \{0,1\}$, a contradiction. 
So $x \not \in Y^i_{e_v} \cap Y^i_{e_v'}$.
Since $e_v \in e'_vR^ie_w$, $y \not \in Y^i_{e_w} \cap Y^i_{e_v'}$.
Since $e_v$ and $e_w$ are not choppers, $\ell_v \neq 0 \neq \ell_w$.
So $\ell_v=\ell_w=1$, a contradiction.

Hence $\tau_i(e_v)=\psi_i(e_v)$.
Since $\tau_i(e_v)=\tau_i(e_w)$ and $\lvert \psi_i(e_v) \rvert = \lvert \psi_i(e_w) \rvert$, $\tau_i(e_w)=\psi_i(e_w)$.
So $Y^i_{e_v}=Y^i_{e_w}$, $x=y$, and no internal node of the path $e_vR^ie_w$ with bag size $\lvert \psi_i(e_v) \rvert$ is a chopper.
Hence the precursor $e^*$ that is the chopper closest to $e_v$ is the precursor that is the chopper closest to $e_w$.
Since $\mu_i(e_v)=\mu_i(e_w)$, if $x=y \in Y^i_{e^*} \cap Y^i_{e_v}$, then by the definition of $\mu_i$, the number of edges incident with $x$ whose other ends are in $V(\overline{G_i})- (R^i,\Y^i)\uparrow e_v$ and the number of edges incident with $y$ whose other ends are in $V(\overline{G_i})-(R^i,\Y^i)\uparrow e_w$ are either both at least two or both equal to a number $\ell$ with $\ell \in \{0,1\}$, and either there exists $\sigma \in \Gamma_i$ with $x=y \in V(\sigma)$ in which $\alpha^i(\sigma)$ is a non-descendant of $e_v$, or there exists no $\sigma' \in \Gamma_i$ with $x=y \in V(\sigma')$ in which $\alpha^i(\sigma')$ is a non-descendant of $e_w$, so $x$ and $y$ do not violate Statement 2 of this claim, a contradiction.
So $x=y \not \in Y^i_{e^*} \cap Y^i_{e_v}$.
In particular, $e^* \neq e_v$.
So $e_v$ is not a chopper.
Since $e_v$ and $e_w$ are not choppers and $x=y \in \tau_i(e_v)=\tau_i(e_w)$, $x$ is associated with essential number 2 in $\gamma_v$ and $y$ is associated with essential number 2 in $\gamma_w$, a contradiction.
This proves the second statement of this claim.
$\Box$

\medskip

\noindent{\bf Claim 3:} For each $i \geq 1$, $(T^i,\psi_i,\tau_i,\mu_i|_{E(T^i)})$ is $(h,d+1,N')$-decorated.

\noindent{\bf Proof of Claim 3:}
By Claim 1, the image of $\mu_i$ is contained in $\{0,1,...,N'\}$.
By the definition of a tree-decomposition, it is easy to see that if $e,e',e''$ are edges of $E(T^i)$ appearing on a directed path in $T^i$ in the order listed, then $\psi_i(e) \cap \psi_i(e'') \subseteq \psi_i(e')$.

Suppose that $(T^i,\psi_i,\tau_i,\mu_i|_{E(T^i)})$ is not $(h,d+1,N')$-decorated.
Then there exist a directed path $P$ in $T^i$, distinct edges $e_1,e_2,...,e_{d+2}$ in $P$ appeared in the order listed with $\lvert \psi_i(e_1) \rvert = \lvert \psi_i(e_2) \rvert = ... =\lvert \psi_i(e_{d+2}) \rvert$, and a set $Z$ with $\psi_i(e_j) \cap \psi_i(e_\ell)=Z$ for all $1 \leq j < \ell \leq d+2$ such that for every $e \in E(P)$, $\lvert \psi_i(e) \rvert \geq \lvert \psi_i(e_1) \rvert$, and for every edge $e' \in E(P)$ with $\lvert \psi_i(e') \rvert = \lvert \psi_i(e_1) \rvert$, we have $\tau_i(e') \not \subseteq Z$ and $\mu_i(e') = \mu_i(e_1)$.

Note that $Z \neq \psi_i(e_1)$, for otherwise $\psi_i(e_1) = Z=\psi_i(e_1) \cap \psi_i(e_2)$, so $\psi_i(e_2) = \psi_i(e_1)=Z$ and $e_2$ is an edge $e'$ of $E(P)$ with $\lvert \psi_i(e') \rvert = \lvert \psi_i(e_1) \rvert$ and $\tau_i(e_2) \subseteq \psi_i(e_2)=Z$.
In addition, no node in $e_2R^ie_{d+2} \subseteq e_1R^ie_{d+2}-\{e_1\}$ with bag size $\lvert \psi_i(e_1) \rvert$ is a chopper.

So in $(R^i,\Y^i)$, $e_2,e_3,...,e_{d+2}$ are nodes of $R^i$ such that 
	\begin{itemize}
		\item $e_j$ is a precursor of $e_{j+1}$ for every $j \in [d+1]-[1]$,
		\item $Z \subseteq Y^i_{e_j}$ for all $j \in [d+2]-[1]$,
		\item $Y^i_{e_j}-Z$ are pairwise disjoint nonempty sets with the same size for all $j \in [d+2]-[1]$, and
		\item there exist $\lvert Y^i_{e_2} \rvert$ disjoint paths in $G$ from $Y^i_{e_2}$ to $Y^i_{e_{d+2}}$.
	\end{itemize}
Since the $\Gamma_i$-elevation of $(R^i,\Y^i,\alpha^i)$ is at most $d$, there exists a node $e^*$ of $R^i$ belonging to the directed path in $R^i$ from $e_2$ to $e_{d+2}$ such that $\lvert Y^i_{e^*} \rvert = \lvert Y^i_{e_2} \rvert$, and $(A_{e^*},B_{e^*})$ (i.e.\ the separation of $G$ given by the node $e^*$ in $(R^i,\Y^i)$) is a $\Gamma_i$-pseudo-edge-cut modulo $Z$.
We choose $e^*$ such that $e^*$ is an edge of $T^i$ if possible.
Since $(R^i,\Y^i)$ is the node-realizer of $(T^i,\X^i)$, and $e_2$ and $e_{d+2}$ are edges of $T^i$, we know that $e^*$ is an edge of $T^i$ by our choice of $e^*$.
So $e^*$ is an edge in $P$.
Hence, $\lvert \psi_i(e^*) \rvert = \lvert \psi_i(e_2) \rvert$. 
So $\tau_i(e^*) \not \subseteq Z$.

Hence there exists $v \in \tau_i(e^*)-Z$.
Since $\lvert Y^i_{e^*} \rvert = \lvert Y^i_{e_2} \rvert$ and $e^* \in e_2R^ie_{d+2}$, $e^*$ is not a chopper.
Since $v \in \tau_i(e^*)$, $v$ is associated with essential number 2 in $\gamma_{t_{e^*}}$, where $t_{e^*}$ is the head of $e^*$ in $T^i$.
But $v \not \in Z$, contradicting that $(A_{e^*},B_{e^*})$ is a $\Gamma_i$-pseudo-edge-cut modulo $Z$.
$\Box$

\medskip

Define $D$ to be an infinite graph with $V(D) = \bigcup_{i \geq 1} V(T^i)$ such that for $i'>i \geq 1$, if a vertex $x \in V(T^i)$ is adjacent to a vertex $y \in V(T^{i'})$, then either
	\begin{itemize}
		\item $x$ is the root of $T^i$, $y$ is the root of $T^{i'}$, and $(G_{i'},\gamma_{i'},\Gamma_{i'},f_{i'}, \phi_{i'})$ simulates $(G_i,\gamma_i,\Gamma_i,f_i,\phi_i)$, or
		\item $x$ is not the root of $T^i$, $y$ is not the root of $T^{i'}$, and the $(f_{i'},\phi_{i'})$-branch of $(T^{i'},\X^{i'},\alpha^{i'})$ at $y$ with respect to $\pi_y$ simulates the $(f_i,\phi_i)$-branch of $(T^i,\X^i,\alpha^i)$ at $x$ with respect to $\pi_x$.
	\end{itemize}

\medskip

\noindent{\bf Claim 4:} If $i' > i \geq 1$, $u \in V(T^i)$ is adjacent in $D$ to $w \in V(T^{i'})$, and $v \in V(T^{i'})$ precedes $w$ in $T^{i'}$ with respect to $(\psi_{i'},\tau_{i'},\mu_{i'}|_{E(T^{i'})})$, then $u$ is adjacent in $D$ to $v$.

\noindent{\bf Proof of Claim 4:}
We may assume $v \neq w$, for otherwise we are done.
Since $v$ precedes $w$ with respect to $(\psi_{i'},\tau_{i'},\mu_{i'}|_{E(T^{i'})})$, we know that $v$ and $w$ are not the root of $T^{i'}$.
Since $u$ is adjacent to $w$, $u$ is not the root of $T^i$.
Since the simulation relation for $Q$-assemblages is transitive, it suffices to prove that the $(f_{i'},\phi_{i'})$-branch of $(T^{i'},\X^{i'},\alpha^{i'})$ at $v$ with respect to $\pi_v$ simulates the $(f_{i'},\phi_{i'})$-branch of $(T^{i'},\X^{i'},\alpha^{i'})$ at $w$ with respect to $\pi_w$.

Since $v$ precedes $w$ with respect to $(\psi_{i'},\tau_{i'},\mu_{i'}|_{E(T^{i'})})$, by Claim 2, there exist $\lvert V(\gamma_v) \rvert$ disjoint paths $P_1,...,P_{\lvert V(\gamma_v) \rvert}$ in $G_{i'}[\uparrow v]$ from $V(\gamma_v)$ to $V(\gamma_w)$ such that for each $j \in [|V(\gamma_v)|]$, the ends of $P_j$ are the $j$-th entries of $\gamma_v$ and $\gamma_w$.
Let $e_v,e_w$ be the edges of $T^{i'}$ with heads $v,w$, respectively.
So $e_w$ is not a chopper.
Hence $\tau_{i'}(e_w)$ is the set of vertices in $\gamma_w$ associated with essential number 2 in $\gamma_w$. 
Since $\tau_{i'}(e_v)=\tau_{i'}(e_w)$, $\tau_{i'}(e_v) = \tau_{i'}(e_w) \subseteq \psi_{i'}(e_v) \cap \psi_{i'}(e_w) \subseteq V(\gamma_v) \cap V(\gamma_w)$. 
In particular, every vertex in $\gamma_w$ associated with essential number 2 is contained in $V(\gamma_w) \cap V(\gamma_v)$.
Note that by Statement 2 of Claim 2, for each $j \in [|V(\gamma_v)|]$, there exists $\ell_j \in \{0,1,2\}$ such that the vertex in $V(P_j) \cap V(\gamma_v)$ is associated with essential number $\ell_j$ in $\gamma_v$, and the vertex in $V(P_j) \cap V(\gamma_w)$ is associated with essential number $\ell_j$ in $\gamma_w$.

Let $H_v,H_w$ be the rooted extensions of $(G_{i'}[\uparrow v],\gamma_v),(G_{i'}[\uparrow w],\gamma_w)$, respectively.
Define $\pi_V: V(H_w) \rightarrow V(H_v)$ such that $\pi_V|_{\uparrow w}$ is the identity map and $\pi_V$ maps the $j$-th entry in the indicator of $H_w$ to the $j$-th entry in the indicator of $H_v$ for each $j \in [|V(\gamma_v)|]$.
Define $\pi_E$ to be the function with domain $E(H_w)$ such that for every $e \in E(H_w)$, 
	\begin{itemize}
		\item if $e \in E(G_{i'}[\uparrow w])$ or $e$ is incident with an entry in the indicator which is adjacent to a vertex associated with essential number $2$ in $\gamma_w$, then $\pi_E(e)=e$;
		\item otherwise, $e$ is an edge incident with an entry in the indicator of $H_w$ with degree one, say the $j$-th entry, then we define $\pi_E(e)$ to be the path obtained by concatenating $P_j$ with the edge incident with the $j$-th entry in the indicator of $H_v$.
	\end{itemize}
It is clear that $(\pi_V,\pi_E)$ is a homeomorphic embedding from $(G_w,\gamma_w)$ to $(G_v,\gamma_v)$, since every vertex in $\gamma_w$ associated with essential number 2 is contained in $V(\gamma_w) \cap V(\gamma_v)$, and for each $j \in [|V(\gamma_v)|]$, there exists $\ell_j \in \{0,1,2\}$ such that the vertex in $V(P_j) \cap V(\gamma_v)$ is associated with essential number $\ell_j$ in $\gamma_v$, and the vertex in $V(P_j) \cap V(\gamma_w)$ is associated with essential number $\ell_j$ in $\gamma_w$.
Define $\iota: \Gamma_w \rightarrow \Gamma_v$ to be the identity map.
Then $(\pi_V,\pi_E)$ and $\iota$ witness that the $(f_{i'},\phi_{i'})$-branch of $(T^{i'},\X^{i'},\alpha^i)$ at $v$ simulates the one at $w$.
$\Box$

\medskip

We may assume that the roots of $T^1,T^2,...$ form a stable set in $D$, for otherwise this theorem is proved.
By Claims 3, 4 and Theorem \ref{decorated tree lemma}, there exists an infinite stable set $S$ of $D$ such that $\lvert S \cap V(T^i) \rvert \leq 1$ for every $i \geq 1$, and the set, denoted by $C$, of heads of all edges of $T^1 \cup T^2 \cup \cdots$  with tails in $S$ is rich in $D$.
By removing some trees from the sequence $T^1,T^2,...$, we may assume that $\lvert S \cap V(T^i) \rvert=1$ for all $i \geq 1$.

For each $i \geq 1$, let $s_i$ be the vertex in $S \cap V(T^i)$.
Since $S$ is infinite, either infinitely many elements of $S$ are the roots of some trees in $\{T^1,T^2,...\}$, or infinitely many elements of $S$ are not the roots of some trees in $\{T^1,T^2,...\}$.
So by removing some trees from the sequence $T^1,T^2,...$, we may assume that either $s_i$ is the root of $T^i$ for every $i \geq 1$, or $s_i$ is a non-root node of $T^i$ for every $i \geq 1$.

Let $C'$ be the set of $(f_i,\phi_i)$-branches of $(T^i,\X^i,\alpha^i)$ at $c$ among all $c \in C \cap V(T^i)$ and all $i \geq 1$.
Since $C$ is rich, $C'$ is well-quasi-ordered by the simulation relation. 

For each $i \geq 1$, 
	\begin{itemize}
		\item if $s_i$ is the root of $T^i$, then let $B_i=(G_i,\gamma_i,\Gamma_i,f_i,\phi_i)$, and let $(T^*_i,\X^*_i,\alpha^*_i) = (T^i,\X^i,\alpha^i)$,
		\item otherwise, let $B_i$ be the $(f_i,\phi_i)$-branch of $(T^i,\X^i,\alpha^i)$ at $s_i$ with respect to $\pi_{s_i}$, and let $(T^*_i,\X^*_i,\alpha^*_i)$ be the rooted tree-decomposition of $B_i$ by taking the maximal subtree of $T^i$ rooted at $s_i$.
	\end{itemize}
By Lemma \ref{branch of brnach}, $C'$ is the set of $(f_i,\phi_i)$-branches of $(T^*_i,\X^*_i,\alpha^*_i)$ at $c$ among all $c \in C \cap V(T^*_i)$ and all $i \geq 1$.

For each $i \geq 1$, we define $W_i$ to be the encoding of $(T^*_i,\X^*_i,\alpha^*_i)$ at $s_i$.
Since each $(T^i,\X^i,\alpha^i)$ is over $\F$, each $(T^*_i,\X^*_i,\alpha^*_i)$ is over $\F$ by Lemma \ref{branch of brnach}.
So the underlying assemblage of $W_i$ is in $\F$ for each $i \geq 1$.
Since $Q$ is a well-quasi-order, and $C'$ is well-quasi-ordered by the simulation relation, there exists a well-quasi-order $Q^*$ such that for each $i \geq 1$, $W_i$ is a $Q^*$-assemblage.
Since $\F$ is well-behaved, there exist $i' > i \geq 1$ such that $W_{i'}$ simulates $W_i$.
Since $s_{i'}$ is the root of $T^*_{i'}$ and $s_i$ is the root of $T^*_i$, by Lemma \ref{encoding simulation}, $B_{i'}$ simulates $B_i$. 

Recall that either $s_i$ is the root of $T^i$ for every $i \geq 1$, or $s_i$ is a non-root node of $T^i$ for every $i \geq 1$.
For the former, $B_{i'}=(G_{i'},\gamma_{i'},\Gamma_{i'},f_{i'},\phi_{i'})$ simulates $B_i=(G_i,\gamma_i,\Gamma_i,f_i,\phi_i)$, so we are done.
For the latter, $B_{i'}$ simulates $B_i$, and $s_{i'}$ and $s_i$ are not roots, so $s_i$ is adjacent to $s_{i'}$ in $D$, contradicting that $S$ is stable in $D$.
This proves the theorem.
\end{pf}

\subsection{Application to graphs} \label{subsec:assemblages_to_graphs}

\begin{lemma} \label{graph to assemblage}
Let $w,k,d,N$ be positive integers.
Let $\F^w$ be the set of assemblages on at most $w+1$ vertices.
If $G$ is a graph that has an $N$-linked rooted tree-decomposition $(T,\X)$ of width at most $w$ and elevation at most $d$ such that for every edge $xy$ of $T$, either $X_x \subseteq X_y$ or $X_y \subseteq X_x$, then the assemblage $(G,\emptyset,\emptyset)$ has a $2N$-unimpeded rooted tree-decomposition over $\F^w$ of adhesion at most $w+1$ such that its node-realizer has $\emptyset$-elevation at most $d$.
\end{lemma}

\begin{pf}
Define $\alpha$ to be the function with empty domain.
Then $(T,\X,\alpha)$ is a rooted tree-decomposition of the assemblage $(G,\emptyset,\emptyset)$.
By Lemma \ref{N-linked disjoint paths}, since $(T,\X)$ is $N$-linked, $(T,\X)$ is weakly $N$-linked.
Since $(T,\X)$ is weakly $N$-linked and for each edge $xy$ of $T$, either $X_x \subseteq X_y$ or $X_y \subseteq X_x$, $(T,\X,\alpha)$ is a $2N$-unimpeded rooted tree-decomposition of $(G,\emptyset,\emptyset)$, and the $\emptyset$-elevation of the node-realizer of $(T,\X,\alpha)$ equals the elevation of $(T,\X)$. 
Furthermore, the width of $(T,\X)$ is at most $w$, so $(T,\X,\alpha)$ is over $\F^w$ and has adhesion at most $w+1$.
\end{pf}

\bigskip

Now we are ready to prove the labelled version of Robertson's conjecture for graphs with bounded tree-width.

\bigskip

\noindent{\bf Proof of Theorem \ref{robertson conj bounded tree width}:}
Let $k,w,G_1,G_2,..., Q, \phi_1,\phi_2,...$ be the ones as in the statement of Theorem \ref{robertson conj bounded tree width}.

Let $N$ and $f(k,w)$ be the numbers mentioned in Theorem \ref{bounded depth}.
Let $\F^w$ be the set of assemblages on at most $w+1$ vertices.
By Lemma \ref{well-behaved bounded size}, $\F^w$ is well-behaved.
Let $\F$ be the set consisting of all assemblages $(G,\gamma,\Gamma)$ that have a $2N$-unimpeded rooted tree-decomposition over $\F^w$ of adhesion at most $w+1$ such that its node-realizer has $\Gamma$-elevation at most $f(k,w)$.
By Theorem \ref{well-behaved tree-decomposition}, $\F$ is well-behaved.

By Theorem \ref{bounded depth}, every graph $G$ of tree-width at most $w$ not containing the Robertson chain of length $k$ as a topological minor has an $N$-linked rooted tree-decomposition $(T,\X)$ of width at most $w$ and of elevation at most $f(k,w)$ such that for every edge $xy$ of $T$, either $X_x \subseteq X_y$ or $X_y \subseteq X_x$.
By Lemma \ref{graph to assemblage}, $(G,\emptyset,\emptyset)$ belongs to $\F$.
Hence $(G_1,\emptyset,\emptyset,\emptyset,\phi_1)$, $(G_2,\emptyset,\emptyset,\emptyset,\phi_2)$, ... are $Q$-assemblages whose underlying assemblages are in $\F$.
Since $\F$ is well-behaved, there exist $i' > i \geq 1$ such that $(G_{i'},\emptyset,\emptyset,\emptyset,\phi_{i'})$ simulates $(G_i,\emptyset,\emptyset,\emptyset,\phi_i)$.
Therefore, there exists a homeomorphic embedding $\eta$ from $G_i$ to $G_{i'}$ such that $\phi_i(v) \leq_Q \phi_{i'}(\eta(v))$ for every $v \in V(G_i)$.
This completes the proof.
$\Box$

\printindex

\end{document}